%% file: EFT4.tex
\documentclass[10pt]{amsart}
\pdfoutput=1
\usepackage{ mathrsfs }
\usepackage{amsmath, amsthm, amssymb,slashed,stmaryrd}
\usepackage{mathtools}
\usepackage{ifpdf}
\usepackage[pdftex]{graphicx}
\usepackage{tikz}
\usetikzlibrary{}
\usetikzlibrary{decorations.markings,matrix,arrows,decorations.pathmorphing,backgrounds,patterns,shapes.geometric,calc}

\usetikzlibrary{matrix,arrows,calc}
\usepackage[all]{xy}
\UseRawInputEncoding

\usepackage[pdftex,plainpages=false,hypertexnames=false,pdfpagelabels]{hyperref}
\usepackage{tikzit}
\input{sample.tikzstyles}

\usepackage{amsthm}
\theoremstyle{definition}
 \newtheorem{thm}{Theorem}[section]
    
 \newtheorem{cor}[thm]{Corollary}
 \newtheorem{lem}[thm]{Lemma}

    \newtheorem{hyp}[thm]{Hypothesis}
 \newtheorem{expect}[thm]{Expectation}
 \newtheorem{prop}[thm]{Proposition}
 \newtheorem{conj}[thm]{Conjecture}
 \newtheorem{defn}[thm]{Definition}
  
 \newtheorem{notation}[thm]{Notation}
 \newtheorem{ex}[thm]{Example}
  
 \newtheorem*{thm*}{Theorem}
 \theoremstyle{remark}
 \newtheorem{rmk}[thm]{Remark}
 
\def\beq{\begin{eqnarray}}
\def\eeq{\end{eqnarray}}
 \newcommand{\bp}{\begin{proof}[Proof]}
 \newcommand{\ep}{\end{proof}}

\DeclareSymbolFont{bbold}{U}{bbold}{m}{n}
\DeclareSymbolFontAlphabet{\mathbbold}{bbold}
\def\one{\mathbbold{1}}

\def\sL{{\sf L}}
\def\sR{{\sf R}}

\def\sS{{\sf S}}

\def\dR{{ d}}
\def\cutoff{{\rm co}}
\def\rg{{\sf rg}}
\def\im{{\rm im}}
\def\u{{\sf u}}
\def\s{{\sf s}}

\def\t{{\sf t}}
\def\m{{\sf m}}

\def\sS{{\sf S}}

\def\sC{{\sf C}}

\def\sR{{\sf R}}
\def\sL{{\sf L}}

\def\c{{\sf c}}
\def\F{{\sf F}}
\def\HBR{{\mathcal{E}}}
\def\Tran{{L}}
\def\Wit{{\mathbb W}{\rm it}}
\def\Witt{{\rm Wit}}
\def\pair{{\beta}}

\def\Mell{\mathcal{M}_{\rm ell}}
\def\Sym{{\rm Sym}}
\def\p{{\rm p}}
\def\a{{\rm a}}

\def\rS{{\rm S}}
\def\rC{{\rm C}}
\def\rL{{\rm L}}
\def\rR{{\rm R}}
\def\rT{{\rm T}}

\def\fl{{\rm fl}}
\def\Fl{{\rm Fl}}

\def\M{{\mathbb{M}}}
\def\V{\mathcal{V}}
\def\SMfld{{\sf SMfld}}
\def\Mfld{{\sf Mfld}}
\def\inn{{\sf in}}

\def\Or{{\sf Or}}
\def\orr{{\sf or}}
\def\c{{\sf c}}
\def\inv{{\sf inv}}

\def\out{{\sf out}}
\def\sTr{{\sf sTr}}

\def\Tr{{\sf Tr}}
\def\deg{{\rm deg}}

\def\MString{{\rm MString}}
\def\BString{{\rm BString}}
\def\MSO{{\rm MSO}}
\def\MSpin{{\rm MSpin}}

\def\KO{{\rm KO}}
\def\MF{{\rm MF}}
\def\mmf{{\rm mf}}
\def\Ind{{\rm Ind}}
\def\K{{\rm K}}
\def\Mod{{\sf Mod}}

\def\Ch{{\rm Ch}}
\def\CS{{\rm CS}}
\def\Eu{{\rm Eu}}
\def\EE{\hbox{-{\sf Euc}}}
\def\Pf{{\rm Pf}}
\def\Pff{{\sf Pf}}
\def\Mf{{\rm Mf}}

\def\Cl{{\rm Cl}}
\def\Fer{{\rm Fer}}
\def\cCl{\C{\rm l}}

\def\act{{\sf act}}
\def\proj{{\sf pr}}
\def\twist{\mathscr{T}}
\def\Ob{{\rm Ob}}

\def\Mor{{\rm Mor}}

\def\H{{\rm H}}
\def\HH{{\mathbb H}}
\def\HHz{{\mathfrak H}}
\def\iHH{{\mathring{\HH}}}
\def\Spin{{\rm Spin}}

\def\SU{{\rm SU}}
\def\SO{{\rm SO}}
\def\BSO{{\rm BSO}}

\def\Map{{\sf Map}}
\def\RG{{\sf RG}}

\def\EFT{{\sf EFT}}
\def\Vect{{\sf Vect}}

\def\vol{{\rm vol}}

\def\Lat{{\sf Lat}}

\def\sLat{s{\sf Lat}}
\def\EFT{ \hbox{-{\sf EFT}}}
\def\eft{ \hbox{-{\sf eft}}}
\def\Bord{\hbox{-{\sf Bord}}}
\def\EBord{\hbox{-}{\sf EBord}}
\def\CEBord{\hbox{-}{\sf CEBord}}
\def\Ann{\mathscr {C}}
\def\tAnn{\widetilde{\mathscr {C}}}
\def\ebord{\hbox{-{\sf ebord}}}
\def\TMF{{\rm TMF}}

\def\pt{{\rm pt}}

\def\ev{{\rm ev}}
\def\odd{{\rm odd}}

\def\bS{{\mathbb{S}}}
\def\A{{\mathbb{A}}}
\def\B{{\mathbb{B}}}
\def\R{{\mathbb{R}}}

\def\TA{{\sf {TA}}}

\def\E{{\mathbb{E}}}

\def\N{{\mathbb{N}}}
\def\id{{{\rm id}}}
\def\C{{\mathbb{C}}}

\def\Z{{\mathbb{Z}}}

\def\GL{{\rm GL}}
\def\End{{\rm End}}
\def\SL{{\rm SL}}
\def\MP{{\rm MP}}
\def\Hom{{\sf Hom}}

\def\Spec{{\rm Spec}}

\def\aa{{\sf a}}
\def\b{{\sf b}}
\newcommand{\op}{{\sf{op}}}   
\newcommand{\sq}{/\!\!/}

\def\leftsquigarrow{\ensuremath{\rotatebox[origin=c]{-180}{$\rightsquigarrow$}}}
\def\acts{ \ \ensuremath{\rotatebox[origin=c]{90}{$\curvearrowright$}} \ }

\def\twocommute{\ensuremath{\rotatebox[origin=c]{30}{$\Rightarrow$}}}

\newcommand{\hocolim}  
{\operatornamewithlimits{{hocolim}}}

\newcommand\nc{\newcommand}
\nc\mf\mathfrak
\nc\mc\mathcal
\nc\mb\mathbb

\begin{document}

\title{How do field theories detect the torsion in topological modular forms?}
\author{Daniel Berwick-Evans}

\begin{abstract}
We construct deformation invariants of $2|1$-dimensional Euclidean field theories valued in a cohomology theory approximating topological modular forms. This implies several results anticipated by Stolz and Teichner and gives the first torsion invariants of field theories valued in $\pi_*\TMF$. The framework leads to a version of the elliptic Euler class described in terms of field-theoretic data, as well as a partial construction of the supersymmetric sigma model with target a string manifold. 
\end{abstract}


\maketitle 

\vspace{-.3in}
\section{Introduction and statement of results}

Stolz and Teichner describe categories of $2|1$-dimensional Euclidean field theories over a smooth manifold $M$ indexed by a degree $n\in \Z$. They conjecture that such field theories determine geometric cocycles for topological modular forms~\cite[Conjecture~1.17]{ST11}
\beq
&&\begin{tikzpicture}[baseline=(basepoint)];
\node (A) at (0,0) {$\left\{\begin{array}{c} {\rm fully}\hbox{-}{\rm extended} \ 2|1\hbox{-}{\rm Euclidean} \\ {\rm field\ theories\ of \ degree} \ n \ { \rm over} \ M\end{array} \right\}$};
\node (B) at (7,0) {$\TMF^{n}(M).$};
\draw[->,dashed] (A) to node [above] {\small cocycle}  (B);
\path (0,-.05) coordinate (basepoint);
\end{tikzpicture}\label{eq:conjecture}
\eeq
It is further expected that this cocycle model inherits an analytic pushforward from quantization of field theories. The ultimate goal is to prove a TMF-generalization of the Atiyah--Singer index theorem. 
In this paper we formulate and prove a 1-extended approximation to~\eqref{eq:conjecture} via categories $2|1\eft^n(M)$ that encode the values of degree~$n$ field theories over~$M$ on (Ramond) supercylinders and supertori. Our definition also imposes a Wilsonian-style ellipticity condition on objects of $2|1\eft^n(M)$, see Remark~\ref{rmk:cutoffs}. On the homotopical side,~$\KO_\MF$ is a height~$\le 1$ approximation of $\TMF$ with a map of ring spectra~$\TMF\to \KO_\MF$, see~\eqref{eq:pullbacksquare}.

\begin{thm}\label{thm1}
Partially-defined $2|1$-Euclidean field theories determine geometric cocycles
\beq\label{eq:cocycle}
2|1\eft^n(M)\xrightarrow{{\rm cocycle}} \KO^n_\MF(M).
\eeq
\end{thm}

The Ando--Hopkins--Rezk $\sigma$-orientation of TMF~\cite{AHR} provides a map of ring spectra
\beq\label{eq:ogstringorientation}
\MString\xrightarrow{\sigma} \TMF\to \KO_\MF,
\eeq
implying that $\KO_\MF$ has Euler classes for vector bundles with string structures. The following approximates Stolz and Teichner's proposed elliptic Euler class~\cite[Theorem 1.0.3]{ST04}.

\begin{thm}\label{thm2}
A real vector bundle $V\to M$ with geometric string structure determines an object $\Eu(V)\in 2|1\eft^{{\rm dim}(V)}(M)$ whose image under~\eqref{eq:cocycle} is the $\KO_\MF$-Euler class of $V$.
\end{thm} 

Hopkins--Mahowald proved that the string orientation detects all torsion in $\pi_*\TMF$ \cite[Theorem 6.25]{HopkinsICM2002}, and 
Bunke--Naumann~\cite[\S1.4]{BunkeNaumann} show that the 3-torsion in degrees $4k-1$ injects along $\pi_*\TMF\to \pi_*\KO_\MF$ and that $\Z/24\simeq \pi_3\TMF\hookrightarrow \pi_3 \KO_\MF$ injects.
 The following approximates the anticipated $\pi_*\TMF$-valued analytic index while also giving the first $\TMF$-torsion detected by structures in supersymmetric Euclidean field theories. 

\begin{thm}\label{thm3}
For an $n$-dimensional string manifold $X$, there is an object $\sigma(X)\in 2|1\eft^{-n}(\pt)$ that is a cocycle refinement (using~\eqref{eq:cocycle}) of the class in the image of 
$$
[X]\in \pi_n(\MString)\to \pi_n\TMF\to \pi_n\KO_\MF=\KO^{-n}_\MF(\pt).
$$ 
In particular, $S^3$ with string structure specified by $1\in \Z\simeq \H^3(S^3;\Z)$ determines a cocycle $\sigma(S^3)\in 2|1\eft^{-3}(\pt)$ representing the generator $\nu$ of $\Z/24\simeq \pi_3(\TMF)\hookrightarrow \pi_3(\KO_\MF)$. 
\end{thm}

\subsection{Context}
In the late 1980s, Witten proposed a surprising connection between elliptic cohomology and quantum field theory~\cite{Witten_Elliptic,Witten_Dirac}. One focal point of his work is a cobordism invariant of string manifolds valued in integral modular forms
\beq\label{eq:thewittengenus}
\left\{\!\!\begin{array}{c}{\rm String \ manifolds} \\ {\rm up\ to \ cobordism} \end{array}\!\!\right\}\to \MF^\Z
\eeq
now called the \emph{Witten genus}. In physical language, Witten argued that partition functions of quantum field theories with $\mathcal{N}=(0,1)$ supersymmetry determine integral modular forms. The invariant~\eqref{eq:thewittengenus} then arises from the partition function of the $\mathcal{N}=(0,1)$ supersymmetric sigma model with target a string manifold. Segal~\cite{Segal_Elliptic,SegalCFT} and subsequently Stolz and Teichner~\cite{ST04,ST11} provided mathematical frameworks in which aspects of Witten's physical argument could be made precise, leading to the conjecture~\eqref{eq:conjecture}.

In parallel to these developments, Hopkins and collaborators clarified the homotopy-theoretic origins of the Witten genus~\cite{HopkinsICM94,AHSI,HopkinsICM2002,Lurie_Elliptic,AHR}. This culminated in the universal elliptic cohomology theory of topological modular forms $(\TMF)$ and its string orientation, 
\beq\label{eq:stringfactorsWG}
&&\begin{tikzpicture}[baseline=(basepoint)];
\node (X) at (-4,0) {$\MString$};
\node (Y) at (-2,0) {$\TMF,$};
\node (A) at (0,0) {$\pi_*\MString$};
\node (B) at (3,0) {$\pi_*\TMF$};
\node (C) at (6,0) {$\MF^\Z$.};
\draw[->] (X) to node [above] {$\sigma$} (Y);
\draw[->] (A) to node [below] {$\pi_*\sigma$} (B);
\draw[->] (B) to node [below] {$\phi$} (C);
\draw[->,bend left=10] (A) to node [above] {\eqref{eq:thewittengenus}} (C);
\path (0,0) coordinate (basepoint);
\end{tikzpicture}
\eeq
The map $\sigma$ refines the Witten genus as indicated on the right above, using
the canonical (edge) homomorphism~$\phi\colon \pi_*\TMF\to \MF^\Z$. The map~$\phi$ is a rational isomorphism, with kernel and cokernel nontrivial 2- and~3-torsion groups. It is natural to ask if Witten's bridge between quantum field theory and modular forms extends to a link with topological modular forms. This question has proven difficult: a description of $\TMF$-torsion in the language of quantum field theory has eluded understanding for decades. 

Part of the difficultly comes from our meager physical understanding of $\mathcal{N}=(0,1)$-theories. With enough supersymmetry, a quantum field theory becomes so constrained that its path integral is effectively controlled by algebra rather than analysis. This often leads to exact computations grounded in rigorous mathematics. For example, 2-dimensional theories with $\mathcal{N}=(2,2)$ supersymmetry admit the $A$- and $B$-twists from mirror symmetry; the correlation functions in these theories can be evaluated using symplectic geometry and algebraic geometry, respectively~\cite{Wittentopological}. By comparison, $\mathcal{N}=(0,1)$ is the minimal supersymmetry in dimension~2, and so it imposes relatively few constraints on the path integral. Correspondingly few quantities have been computed by physicists. Due to this more challenging landscape, $\mathcal{N}=(0,1)$-theories have received comparatively little attention. And yet, minimally supersymmetric theories are in some sense the most fundamental. 

Recent work in physics has used~\eqref{eq:conjecture} to predict when a pair of $\mathcal{N}=(0,1)$-theories lie in the same deformation class~\cite{GJF,TopMoon,Tachikawa,TachikawaYamashita,GPPV,LinPei,TMFmoon}. One can then ask for physically meaningful invariants witnessing these predictions.
Gaiotto, Johnson-Freyd and Witten~\cite{GJFW} consider the $\mathcal{N}=(0,1)$-supersymmetric sigma model with target~$S^3$. This is a supersymmetric Wess--Zumino--Witten model for the group $S^3 \simeq \SU(2)$ depending on a level $k\in \Z\simeq \H^4(B\SU(2);\Z)$. Gaiotto--Johnson-Freyd--Witten argue that $k$ mod 24 is a deformation invariant of the theory. In light of the isomorphism $\pi_3(\TMF)\simeq \Z/24$, this matches with the prediction from conjecture~\eqref{eq:conjecture} and the $\pi_*\TMF$-valued refinement~\eqref{eq:stringfactorsWG} of the Witten genus. Subsequent work of Gaiotto and Johnson-Freyd generalizes this invariant to $\mathcal{N}=(0,1)$ supersymmetric sigma models with target a $(4k-1)$-dimensional string manifold, where the level is replaced by more general  anomaly cancelation data \cite{GJF2}. 
 Theorems~\ref{thm1} and~\ref{thm3} make the Gaiotto--Johnson-Freyd--Witten invariant precise and identify its values with $\TMF$-torsion in the image of $\pi_{4k-1}\TMF\to \pi_{4k-1}\KO_\MF$. This gives a partial answer to the title question of this paper. 


\subsection{Relation to prior work}

Theorem~\ref{thm1} implies several results anticipated by Stolz and Teichner. Specifically, \cite[Theorem 1.0.2]{ST04} describes a map (contingent on~\cite[Hypothesis 3.3.13]{ST04})
\beq\label{eq:KTateEFT}
\left\{\!\!\!\begin{array}{c} 2|1\hbox{-}{\rm Euclidean\ field\ theories} \\ {\rm of \ degree} \ n \ { \rm over} \ M\end{array} \!\!\!\right\}\longrightarrow \KO^{n}(M)(\!(q)\!)
\eeq
as well as a modularity property for the image of~\eqref{eq:KTateEFT} when $M=\pt$. A result closely related to~\eqref{eq:KTateEFT} was proved by Pokman Cheung~\cite{Cheung}, and a special case of the modularity property is stated as \cite[Theorem 1.15]{ST11} where a proof is outlined. The map of ring spectra
$$
\KO_\MF\to \KO(\!(q)\!),
$$
allows one to deduce a version of~\eqref{eq:KTateEFT} from Theorem~\ref{thm1}. Furthermore, the modularity encoded by $\KO_\MF$ reduces to Stolz and Teichner's claimed modularity for~\eqref{eq:KTateEFT} when $M=\pt$. Building on~\eqref{eq:KTateEFT}, Stolz and Teichner describe a field theory in \cite[Theorem 1.0.3]{ST04} that is expected to provide a cocycle representative of the $\TMF$-Euler class of a string vector bundle under~\eqref{eq:conjecture}. In some sense, this object is the heart of the Stolz--Teichner program, providing a description of the string orientation of $\TMF$ in terms of 2-dimensional supersymmetric field theories. Theorem~\ref{thm2} provides a partial construction and verifies that $\Eu(V)\in 2|1\eft^n(M)$ is a cocycle representative of the $\KO_\MF$-Euler class, and hence approximates the $\TMF$-Euler class. This realizes a 1-categorical version of the Stolz--Teichner program; see also Remark~\ref{rmk:comparetoSTEuler}. We caution that the precise definition of field theory differs in all three of \cite{ST04}, \cite{Cheung} and \cite{ST11}, although the spirit of the definitions are quite similar. The definitions in this paper are closely related to the ones in~\cite{ST11}, see Proposition~\ref{rmk:comparetoST}. 

Our results also shed some light on Bunke and Naumann's torsion invariants of $(4k-1)$-dimensional string manifolds~\cite{BunkeNaumann}
\beq\label{EqBNinvariant}
\pi_{4k-1}\MString\to \C(\!(q)\!)/(\Z(\!(q)\!)+\MF_{2k}).
\eeq
In Theorem~\ref{thm:string} and Corollary~\ref{cor:BN1}, we show that the composition of maps of spectra~\eqref{eq:ogstringorientation} determines the invariant~\eqref{EqBNinvariant} by passing to coefficients. 
This description implies the conjecture~\cite[\S1.5, Open problem~3]{BunkeNaumann}, see Corollary~\ref{cor:BN2}. Theorem~\ref{thm3} then identifies the Bunke--Naumann invariant with the Gaiotto--Johnson-Freyd--Witten invariant, completing the argument outlined in~\cite[\S2.5]{GJF2}; see Remark~\ref{rmk:GJFcompare}.

Theorems~\ref{thm1}-\ref{thm3} forge new connections between algebraic topology and Dirac--Ramond operators in the physics literature. 
In particular, Alvarez and Windey~\cite{AlvarezWindey} study families of Dirac--Ramond operators~\cite{Witten_Dirac} associated to a bundle of string manifolds $X\to M$, extracting modular invariants at the level of cohomology of $M$ with coefficients over~$\C$. A natural extension of their construction involves the \emph{Bismut--Ramond superconnection}, which gives both integral and torsion refinements of Alvarez and Windey's invariants; see~\S\ref{sec:BRsuperconn}. The Bismut--Ramond superconnection is the analytical input to the construction of the objects $\sigma(X)\in 2|1\eft^{-n}(\pt)$ in Theorem~\ref{thm3}. 

The geometry behind Theorems~\ref{thm1}-\ref{thm3} expand upon a verification of the conjecture~\eqref{eq:conjecture} with~$\C$ coefficients in~\cite{DBE_MQ,DBEChern}, i.e., a  height~0 confirmation of the Stolz--Teichner program. These previous results show that the value of a $2|1$-Euclidean field theory on supertori in~$M$ produces a class in $\H(M;\MF^\omega)$, ordinary cohomology valued in weak modular forms; see Proposition~\ref{prop:maintorusprop}. The present paper extends this by also including the data of a $2|1$-Euclidean field theory on supercylinders in~$M$. The trace of the values on supercylinders recovers the values on supertori, so that Theorems~\ref{thm1}-\ref{thm3} directly generalize the results from~\cite{DBE_MQ,DBEChern}.

A crucial new ingredient developed in this paper and the companion~\cite{DBEEFT} is an ellipticity-type condition for field theories, namely that energy cutoffs vary smoothly and admit a type of Chern--Simons form, see Definition~\ref{defn:admitscutoffsEFT}. This Wilsonian approach to mathematical quantum field theory is influenced by Costello's work~\cite{costello_WG1,costbook,costello_WG2}. In our framework, the existence of a cutoff is a \emph{property} of a Stolz--Teichner field theory, not additional data; see Remark~\ref{rmk:cutoffs}.

An enduring challenge in modern mathematical physics is the search for a precise and complete definition of quantum field theory. In the early days of the subject, there were no definitions at all; now mathematicians are spoiled for choice, making it unclear which is the ``correct" definition. Each mathematical framework has its own strengths and insights, e.g.,~\cite{SegalCFT,ST04,CDOII,HuKriz,Cheung,HST,costello_WG1,costello_WG2,ST11,DouglasHenriques,KitchlooII,GradyPavlov}.
Conjecture~\eqref{eq:conjecture} is an excellent proving ground, as it requires enough sophistication to identify structures in quantum field theory that encode an intricate object like~$\TMF$. With these shifting landscapes in mind, the cocycle map~\eqref{eq:cocycle} is designed to be as robust as possible to changes in the mathematical definition of field theory: an object of $2|1\eft^n(M)$ is the data attached to certain explicit moduli spaces of $2|1$-Euclidean supermanifolds, and such data is unambiguously part of the geometry of a $\mathcal{N}=(0,1)$-theory.
Restricting the data of a fully-defined theory to its values on these specific moduli spaces gives a diagram
\beq
&&\begin{tikzpicture}[baseline=(basepoint)];
\node (A) at (0,0) {$\left\{\begin{array}{c} {\rm fully}\hbox{-}{\rm extended} \\ 2|1\hbox{-}{\rm Euclidean\ field\ theories} \\ {\rm of \ degree} \ n \ { \rm over} \ M\end{array} \right\}$};
\node (B) at (7,0) {$\TMF^n(M)$};
\node (C) at (0,-1.75) {$2|1\eft^n(M)$};
\node (D) at (7,-1.75) {$\KO^n_\MF(M).$};
\draw[->,dashed] (A) to node [above] {\small conjecture~\eqref{eq:conjecture}}  (B);
\draw[->,dashed] (A) to node [left] {restriction} (C);
\draw[->] (B) to (D);
\draw[->] (C) to node [above] {\small Theorem~\ref{thm1}}  (D);
\path (0,-.8) coordinate (basepoint);
\end{tikzpicture}\label{eq:conjectureprime}
\eeq
Theorem~\ref{thm1} therefore helps tame the zoo of possible definitions of field theory in~\eqref{eq:conjecture}: extending the partially-defined field theories in $2|1\eft^n(M)$ to fully-extended ones must lift a class in $\KO_\MF$ to one in $\TMF$. Physicists have recently compiled lists of the $\mathcal{N}=(0,1)$-theories expected to represent specific $\TMF$-classes, e.g., see~\cite[Table 2]{GPPV}. This provides explicit cases where one can study lifting problems for $\KO_\MF$-classes to $\TMF$-classes and compare with extending down in a given quantum field theory.

\begin{rmk} \label{rmk:concordance}
Stolz and Teichner's~\cite[Conjecture~1.17]{ST11} states that a cocycle map~\eqref{eq:conjecture} induces a bijection between concordance classes of field theories and $\TMF$ classes. 
It is easy to see that the map~\eqref{eq:cocycle} factors through concordance classes, but we do not expect a bijection for the definitions used in this paper. The basic problem is that the dimensions of the vector bundles~$\HBR_k$ (see~\eqref{eq:sequenceofsuperconn} below) are concordance invariants in our framework, but that are not invariants of the underlying class in $\KO_\MF$. This issue is already present when considering models for $\K$-theory in terms of superconnections on vector bundles. We know of two plausible approaches that could lead to the desired bijection on concordance classes. The first fixes an infinite dimensional separable Hilbert space~$\mathcal{H}$, and considers subspaces $\HBR\subset M\times \mathcal{H}$ where $\HBR$ is no longer required to be locally trivial, but the family of projection operators determined by $\HBR$ is required to vary continuously with~$M$. This generalizes the definition of a rigged Hilbert space. A second option considers families of field valued in sheaves of topological vector spaces with no local freeness condition; see~\cite[Remark 3.16]{ST11}. A sheafy generalization of Theorem~\ref{thm1} would require one to enhance the tools from index theory to apply in this more general setting, e.g., adapt~\cite[Ch.~9-10]{BGV} to consider superconnections and index bundles supported on general sheaves of topological vector spaces. While this seems quite interesting, it is tangential to the goals of the present paper: our examples of interest (namely, the $\KO_\MF$-Euler class and the supersymmetric $\sigma$-model) are defined on vector bundles, and so for simplicity we work with this class of field theories. 
\end{rmk}

\subsection{Detailed overview of results}
In the remainder of this introduction we describe the objects in Theorems~\ref{thm1}-\ref{thm3} and outline proofs, highlighting some of the key intermediary results about $2|1$-Euclidean field theories and $\KO_\MF$. 

\subsubsection{A $\KO_\MF$-index theorem}
Let $\MF$ denote the graded ring of weakly holomorphic modular forms (see~\S\ref{Sec:MF} for our conventions) and $\H_{\MF}$ the commutative ring spectrum representing ordinary cohomology with coefficients in $\MF$. Similarly, let $\H_{\C(\!(q)\!)[u^{\pm 2}]}$ denote the spectrum representing ordinary cohomology with coefficients in $\C(\!(q)\!)[u^{\pm 2}]$ where $|u^2|=-4$. Define $\KO_\MF$ as the pullback in commutative ring spectra,
\beq
&&\begin{tikzpicture}[baseline=(basepoint)];
\node (T) at (-2,.75) {$\TMF$};
\node (A) at (0,0) {$\KO_\MF$};
\node (B) at (4,0) {$\KO(\!(q)\!)$};
\node (C) at (0,-1) {$\H_\MF$};
\node (D) at (4,-1) {$\H_{\C(\!(q)\!)[u^{\pm 2}]}$};
\node (P) at (.7,-.4) {\scalebox{1.5}{$\lrcorner$}};
\draw[->,] (A) to (B);
\draw[->] (A) to  (C);
\draw[->] (C) to node [below] {$q$-expand} (D);
\draw[->] (B) to node [right] {Ch} (D);
\draw[->,bend left=10] (T) to node [above] {Miller} (B);
\draw[->,bend right=10] (T) to (C);
\draw[->,dashed] (T) to (A);
\path (0,-.25) coordinate (basepoint);
\end{tikzpicture}\label{eq:pullbacksquare}
\eeq
where $\Ch$ is the Pontryagin character, and the lower horizontal arrow is determined by $q$-expansion of coefficients. By the universal property of the pullback, there is a map $\TMF\to \KO_\MF$ determined by smashing the Miller character with $\H_\C$~\cite{Miller}. The Miller character comes from evaluation of $\TMF$ at the Tate curve, which essentially extracts the height~1 part of $\TMF$. Similarly, the Chern--Dold character $\TMF\to \TMF\otimes \C\simeq \H_\MF$ extracts the height~0 part. In this way, $\KO_\MF$ is a height $\le 1$ approximation to~$\TMF$. 

Bismut's families index theorem~\cite{Bismutindex} together with the analytic description of the $\KO(\!(q)\!)$-valued Witten genus from~\cite{Witten_Dirac} leads to the following, proved in~\S\ref{sec:thm:KOMForientation}. 

\begin{prop}\label{thm:KOMForientation}
Given a family of Riemannian string manifolds $\pi\colon X\to M$, the analytic index of the Bismut--Ramond superconnection determines a class $\sigma(X)\in \KO^{-n}_\MF(M)$. 
\end{prop}

The string orientation~\eqref{eq:ogstringorientation} also determines a topological index: for a family $\pi\colon X\to M$ of string manifolds with fiber dimension~$n$, consider the composition
\beq\label{eq:topindex}
&&[X]\in \MString^{-n}(M)\xrightarrow{\sigma} \TMF^{-n}(M)\to \KO_\MF^{-n}(M).
\eeq
A characterization of $[\MString,\KO_\MF]$ in terms of orientations of $\KO(\!(q)\!)$ and $\H_\MF$ leads to the following index theorem for $\KO_\MF$, proved in~\S\ref{sec:proofofthm:string}.

\begin{thm}\label{thm:string}
The topological index~\eqref{eq:topindex} equals the analytic index in Proposition~\ref{thm:KOMForientation}.
\end{thm}

This index theorem recovers a previously-studied invariant of Bunke and Naumann and proves a related conjecture, as we now explain. When $M=\pt$, Proposition~\ref{thm:KOMForientation} determines analytically-constructed invariants of string manifolds,
\beq\label{eq:stringinvariants}
\left\{\!\!\begin{array}{c}{\rm String \ manifolds} \\ {\rm up\ to \ cobordism} \end{array}\!\!\right\}\simeq \pi_*\MString\to \pi_*\KO_\MF.
\eeq
The homotopy groups of $\KO_\MF$ are computed in Proposition~\ref{prop:htpygrps}, and in particular 
$$
\pi_{4k}\KO_\MF\simeq \MF_{2k}^\Z,\qquad \pi_{4k-1}\simeq \C(\!(q)\!)/(\Z(\!(q)\!)+\MF_{2k})
$$
recovering integral modular forms in degrees $4k$ and the target of the Bunke--Naumann invariants~\eqref{EqBNinvariant} in degrees $4k-1$.

\begin{cor}\label{cor:BN1}
In degree $*=4k$, \eqref{eq:stringinvariants} agrees with the Witten genus as an integral modular form.
In degrees $4k-1$, it agrees with the Bunke--Naumann torsion invariants of string manifolds~\cite{BunkeNaumann}. 
\end{cor}

The following verifies~\cite[\S1.5, Open problem~3]{BunkeNaumann}.

\begin{cor} \label{cor:BN2}
The Bunke--Naumann torsion invariant of $(4k-1)$-dimensional string manifolds is determined by the topological pushforward~\eqref{eq:topindex}. 
\end{cor}
\bp
This follows directly from Theorem~\ref{thm:string} and Corollary~\ref{cor:BN1}. 
\ep

\subsubsection{Stolz and Teichner's $2|1$-Euclidean field theories}\label{Sec:partiallydefinedFT}
The categories $2|1\eft^n(M)$ in Theorem~\ref{thm1} come from distilling the geometric information of a $2|1$-Euclidean field theory of  degree~$n$ in the framework developed by Stolz and Teichner in~\cite{ST11}.
We briefly review their setup below.

The definition of degree~$n$ field theory in conjecture~\eqref{eq:conjecture} starts with a category of $2|1$-Euclidean bordisms over a smooth manifold~$M$. A 1-extended version of this bordism category, denoted $2|1\EBord(M)$, is constructed in~\cite{ST11}. Roughly, a morphism in~$2|1\EBord(M)$ is 
\beq
&&\begin{tikzpicture}[baseline=(basepoint)];
\node (A) at (0,0) {$Y_\out^{1|1}$};
\node (B) at (5,0) {$Y_\inn^{1|1}$};
\node (C) at (2.5,0) {$\Sigma^{2|1}$};
\node (D) at (2.5,-1) {$M$};
\draw[->,right hook-latex] (A) to (C);
\draw[->,left hook-latex] (B) to (C);
\draw[->] (C) to node [left] {$\Phi$} (D);
\path (0,-.5) coordinate (basepoint);
\end{tikzpicture}\label{eq:abordism}
\eeq
for a $2|1$-Euclidean manifold $\Sigma^{2|1}$, a map $\Phi\colon \Sigma^{2|1} \to M$, and $1|1$-dimensional sub supermanifolds~$Y_{\inn},Y_{\out}\hookrightarrow \Sigma^{2|1}$ that determine the source and target of the bordism. The precise definition of $2|1\EBord(M)$ requires one consider families of bordisms and their isometries (as a symmetric monoidal stack), but we ignore this for now; see~\S\ref{sec:21EB} for details.

A \emph{degree~$n$ field theory over~$M$} is a natural transformation
\beq
&&\begin{tikzpicture}[baseline=(basepoint)];
\node (A) at (0,0) {$2|1\EBord(M)$};
\node (B) at (5,0) {$\TA$};
\node (C) at (2.5,0) {$E \Downarrow$};
\draw[->,bend left=15] (A) to node [above] {$\one$} (B);
\draw[->,bend right=15] (A) to node [below] {$\Fer_n$} (B);
\path (0,0) coordinate (basepoint);
\end{tikzpicture}\label{eq:twistedEFT}
\eeq
where $\TA$ is the Morita category of (topological) algebras, bimodules, and bimodule maps. The functor $\one$ is the constant functor to the monoidal unit of $\TA$, and
\beq\label{eq:maintwist}
\Fer_n\colon 2|1\EBord(M)\to  \TA,\qquad \Fer_n\otimes \Fer_m\simeq \Fer_{n+m}
\eeq
is a $\otimes$-invertible functor equipped with the indicated isomorphisms for $n,m\in \Z$. A complete construction of $\Fer_n$ has yet to appear, but it is expected that~\eqref{eq:maintwist} is the unique supersymmetric extension of the $n$th tensor power of the chiral free fermion as constructed in~\cite[\S5-6]{ST11}.
Hypothesis~\ref{hyp:chiral} gives a precise formulation of this extension, and we assume this hypothesis below. 

\begin{rmk}
We emphasize that Theorems~\ref{thm1}-\ref{thm3} do not depend on Hypothesis~\ref{hyp:chiral}, but rather it allows us to compare our results with~\cite{ST11}. In~\S\ref{sec:susyextunique} we prove the existence of a unique supersymmetric extension of $\Fer_n$ for the class of $2|1$-Euclidean bordisms relevant in this paper, partially verifying Hypothesis~\ref{hyp:chiral}.
\end{rmk}

It turns out to be important to demand an additional structure on degree~$n$ field theories~\eqref{eq:twistedEFT} called \emph{reflection positivity}; this structure is implicit in the ``adjunction transformations" that were part of the definition of $2|1$-dimensional field theory in \cite{ST04}. Below we follow the approach of Freed--Hopkins~\cite[\S3]{FreedHopkins}: a \emph{reflection structure} is $\Z/2$-equivariance data for a natural transformation~\eqref{eq:twistedEFT} (Definition~\ref{defn:RP0}), and the \emph{positivity property} is that a certain hermitian pairing is positive (Definition~\ref{defn:RPFT}). 

\begin{defn}[Definition~\ref{defn:STEFT}]\label{defn:21EFT}
Let $2|1\EFT^n(M)$ denote the groupoid whose objects are reflection positive natural transformations~\eqref{eq:twistedEFT}. 
\end{defn}

\begin{rmk}\label{rmk:KTate}
We note that reflection positivity is unnecessary in the $d=1$ version of conjecture~\eqref{eq:conjecture}: unoriented, real $1|1$-dimensional field theories recover real K-theory~\cite{HST}.\footnote{In~\cite{HST}, a positivity condition is imposed that is distinct from reflection positivity; it is unclear whether this condition is necessary to recover $\KO$-theory.} However, unoriented and real structures do not generalize to the~$2|1$-dimensional setting: chiral supersymmetry is not preserved under orientation reversal, nor is it preserved by complex conjugation. The Miller character $\TMF\to \KO(\!(q)\!)$ necessitates a map from $2|1$-Euclidean field theories to real K-theory if conjecture~\eqref{eq:conjecture} is to hold. One mechanism (and indeed, the only one we know) to obtain real structures on Hilbert spaces of states in dimension~2 assumes reflection positivity, e.g., see~\cite[3.2.2]{GPPV}. 
\end{rmk}

\subsubsection{The categories $2|1\eft^n(M)$ and $\KO_\MF$-cocycles}
The cocycle map in Theorem~\ref{thm1} comes from restricting~\eqref{eq:twistedEFT} to certain explicit subcategories of $2|1\EBord(M)$,
\beq
\label{eq:twistedC}\label{eq:twistedrestriction}
&&\begin{tikzpicture}[baseline=(basepoint)];
\node (A) at (-1.5,0) {$\Ann^{2|1}(M)$};
\node (B) at (4,-1) {$2|1\EBord(M)$};
\node (C) at (8,-1) {$\TA.$};
\node (D) at (6.2,-1) {$E\Downarrow$};
\node (E) at (-1.5,-1.5) {$\null$};
\draw[->,dashed] (A) to node [left] {trace} (E);
\draw[->] (A) to node [above=5pt] {supercylinders}  (B);
\draw[->,bend left=12] (B) to node [above] {$\one$} (C);
\draw[->,bend right=12] (B) to node [below] {$\Fer_n$}  (C);
\path (0,-.05) coordinate (basepoint);
\end{tikzpicture}
\eeq
\vspace{-.55in}
\beq
\label{eq:restricttotori}&&\begin{tikzpicture}[baseline=(basepoint)];
\node (A) at (-1.5,0) {$\mathcal{L}_0^{2|1}(M)$};
\node (B) at (2.9,.8) {$\null$};
\node (C) at (8,0) {$\null$};
\node (D) at (6.2,0) {$\null$};
\draw[->] (A) to node [below=3pt] {supertori} (B);
\path (0,-.05) coordinate (basepoint);
\end{tikzpicture}
\eeq
The restriction of $E$ along~\eqref{eq:twistedrestriction} affords a twisted representation of the category $\Ann^{2|1}(M)$ of nearly constant supercylinders in~$M$, viewed as bordisms between supercircles in~$M$. 
The restriction of $E$ along~\eqref{eq:restricttotori} determines a line bundle with section over the stack $\mathcal{L}^{2|1}_0(M)$ of nearly constant supertori in~$M$, viewed as bordisms from the empty set to itself. The images of $\Ann^{2|1}(M)$ and $\mathcal{L}^{2|1}_0(M)$ in $2|1\EBord(M)$ are compatible via a (categorical) trace, indicated by the dashed arrow. The moduli spaces~$\mathcal{L}^{2|1}_0(M)$ and~$\Ann^{2|1}(M)$ are finite-dimensional, and so the restrictions~\eqref{eq:twistedrestriction} and~\eqref{eq:restricttotori} extract data from a field theory that can be analyzed using standard techniques in super geometry. 

\begin{rmk}
The Euclidean supercylinders and supertori in subcategories $\Ann^{2|1}(M)$ and  $\mathcal{L}^{2|1}_0(M)$ have nonbounding (i.e., odd) spin structures. In physical language, the restrictions~\eqref{eq:twistedrestriction} and~\eqref{eq:restricttotori} therefore only probe the Ramond sector of a $2|1$-Euclidean field theory, ignoring the Neveu--Schwarz sector. The Ramond sector is the part of the $\mathcal{N}=(0,1)$-supersymmetric sigma model relevant to Witten's construction of the Witten genus~\cite{Witten_Dirac}. 
\end{rmk}

Restriction along~\eqref{eq:twistedrestriction} extracts familiar objects from index theory. 

\begin{prop}[Propositions~\ref{prop:superconn}, \ref{prop:selfadjoint},~\ref{prop:realA}] \label{prop:mainprop}
The restriction~\eqref{eq:twistedrestriction} determines:
\begin{enumerate}
\item[i.] a Fr\'echet vector bundle $\HBR\to M$ with fiberwise $S^1\times \cCl_n$-action, where the eigenvalues of the $S^1$-action are bounded below; 
\item[ii.] a smooth, 1-parameter family of superconnections $\A_\ell$, $\ell\in \R_{>0}$ on $\HBR$ that are $S^1$-invariant and $\cCl_n$-linear.
\end{enumerate}
The restriction along~\eqref{eq:twistedrestriction}  of a reflection positive degree~$n$ field theory further determines
\begin{enumerate}
\item[iii.] a hermitian metric on $\HBR$ for which the $S^1$-action is unitary, the $\cCl_n$-action is self-adjoint and the superconnection $\A$ is self-adjoint;
\item[iv.] a real structure on $\HBR$ compatible with the $S^1$-action and superconnection for which the Clifford action is real relative to the standard real structure on $\cCl_n$. 
\end{enumerate}
\end{prop}

The $S^1$-eigenspaces extract a sequence of vector bundles with superconnections (Lemma~\ref{lem:weightspaces})
\beq\label{eq:sequenceofsuperconn}
\HBR\simeq \bigoplus_{k\in \Z} \HBR_k,\qquad \A_\ell=\bigoplus_{k\in \Z} \A_{\ell,k}
\eeq
where $\HBR_k=\{0\}$ for $k\ll 0$. We emphasize that the $\HBR_k$ need not be finite rank. The geometric objects in Proposition~\ref{prop:mainprop} generalize standard data in $\mathcal{N}=(0,1)$ supersymmetric quantum field theory (reviewed in~\S\ref{sec:motivateQFT}) via the $M$-family of operators
\beq\label{eq:Loperators}
&&L_0:=\bigoplus_k (k +\frac{\ell}{2\pi} \A_{\ell,k}^2),\qquad \bar{L}_0:=\bigoplus_k \frac{\ell}{2\pi}\A_{\ell,k}^2,\qquad \bar{G}_0:=\bigoplus_k \frac{\sqrt{\ell}}{\sqrt{2\pi}} \A_{\ell,k}.
\eeq
We observe that these operators satisfy the condition that $L_0-\bar L_0$ has integer eigenvalues (and so generates an $S^1$-action), $[L_0,\bar L_0]=0$, and~$\bar{G}_0^2=\bar L_0$. Assuming existence and uniqueness of solutions to differential equations (see Remark~\ref{rmk:ODEfail}) the degree~$n$ representation of $\Ann^{2|1}(M)$ is given by
\beq\label{eq:potentialsemigroup}
&&q^{L_0}{\bar q}^{\bar L_0+\theta \sqrt{2\pi} \bar G_0},\qquad q=e^{2\pi i z/\ell}, \quad \ell\in \R_{>0}, \ z/\ell\in \HHz\subset \C
\eeq
recovering the usual super semigroup formula when $M=\pt$, e.g.,~\cite[\S4]{SegalCFT},~\cite[\S3.3]{ST04}, and \cite[Eq.~2.1.17]{Cheung}. Geometrically, $\ell$ denotes the circumference of a supercircle, and $\tau=z/\ell$ is the conformal modulus of a supercylinder; see~\S\ref{sec:supercyl}.

When attempting to connect $2|1$-Euclidean field theories over $M$ with the topology of~$M$, Proposition~\ref{prop:mainprop} suggests that we extract a sequence of $\KO$-classes
\beq\label{Eq:eftindex}
{\rm Index}\colon 2|1\EFT^n(M)\dashrightarrow \KO^n(M)(\!(q)\!),\qquad E\mapsto \sum q^k\Ind(\A_{\ell,k})
\eeq
by applying the families index construction (e.g., see~\cite[Ch.~10]{BGV}) to each $\A_{\ell,k}$ in~\eqref{eq:sequenceofsuperconn} extracted from a $2|1$-Euclidean field theory~$E$. The validity of the families index construction is a \emph{property} of a superconnection, and hence the index~\eqref{Eq:eftindex} is defined on the following subcategory of~$2|1\EFT^n(M)$. 
 
 \begin{defn} \label{defn:admitscutoffsEFT}
A degree~$n$ representation of $\Ann^{2|1}(M)$ \emph{admits cutoffs} if it is infinitesimally generated (see Definition~\ref{defn:infgen}) and the superconnections~\eqref{eq:sequenceofsuperconn} admit smooth index bundles in the sense of Definition~\ref{defn:smoothindexbundleEFT}. A field theory $E\in 2|1\EFT^n(M)$ \emph{admits smooth cutoffs} if the the associated representation of $\Ann^{2|1}(M)$ does. Let $2|1\EFT^n_{\cutoff}(M)\subset 2|1\EFT^n(M)$ denote the full subcategory whose objects admit cutoffs. 
 \end{defn}

\begin{rmk}\label{rmk:cutoffs}
The above terminology comes from the fact that the existence of the families index is the same as the existence of a smooth energy cutoff theory in the sense of Wilsonian effective field theory, e.g., see~\cite[\S1.3]{costbook} for a mathematical exposition where the definition of field theory also involves an ellipticity condition~\cite[Definition~9.3.1]{costbook}.
In our framework, the existence of cutoffs is an analytic condition on field theories analogous to the ellipticity property for differential operators. This analogy is tight for $1|1$-dimensional field theories where ellipticity for a family of Dirac operators implies the cutoff property for the associated object in $1|1\eft^n(M)$, see~\cite{DBEEFT}. It bears noting that from the very beginning of this subject, Witten suggested mathematicians use cutoffs to analyze 2-dimensional theories with $\mathcal{N}=(0,1)$ supersymmetry~\cite[\S1]{Witten_Dirac}. 
\end{rmk}


\begin{rmk}\label{rmk:finitetype}
If the $S^1$-eigenspaces $\HBR_k$ in~\eqref{eq:sequenceofsuperconn} are finite-dimensional vector bundles, the cutoff property in Definition~\ref{defn:admitscutoffsEFT} is automatically satisfied using the existence and uniqueness of solutions to ordinary differential equations, along with the existence of the families index bundles in finite-dimensions, e.g.,~\cite[\S9.1]{BGV}. In particular, $M$-families of vertex algebras with connection produce the operators~\eqref{eq:Loperators} and satisfy the cutoff condition. There are many proposed constructions relating (families of) vertex algebras with elliptic cohomology classes, e.g., see~\cite{LiuMF,CDOII,HuKriz,costello_WG1,PokmanCDO,KitchlooII,Andreexamples,LinPei}. The representation~\eqref{eq:potentialsemigroup} of Euclidean supercylinders offers a point of contact between these vertex algebraic constructions and the Stolz--Teichner program. 
\end{rmk}

\begin{rmk}
When $M=\pt$, the cutoff property simplifies considerably, essentially imposing a growth constraint on the eigenvalues of $\A_\ell^2$, see Remark~\ref{rmk:indexvsdiff}. 
\end{rmk}

Now we shift attention to moduli spaces of supertori and the restriction of a degree~$n$ field theory along~\eqref{eq:restricttotori}. From general considerations this restriction produces a line bundle with section over $\mathcal{L}_0^{2|1}(M)$~\cite[\S5.3]{ST11}. The following generalizes~\cite[Theorem~1.1]{DBEChern} from a statement in degree~$n=0$ to arbitrary degree. The vector spaces $\mmf^{k,l}$ are $C^\infty$-versions of weak modular forms, viewed as a function on based, oriented lattices in $\C\simeq \R^2$ that transform with weight $(-k/2,-l/2)$, see~\eqref{Eq:mmf}.

\begin{prop}\label{prop:maintorusprop}
Sections of the $n^{\rm th}$ tensor power of the Pfaffian line are given by data 
\beq\label{eq:deriveddata0}
&&\Gamma(\mathcal{L}^{2|1}_0(M);\Pff^{\otimes -n})\simeq\left\{\begin{array}{l} \displaystyle Z\in \bigoplus_{j+k=n} \Omega^j(M;\mmf^{k,0}))\\ 
\displaystyle Z_v\in \bigoplus_{j+k=n-1} \Omega^j(M;\mmf^{(k+2),2}) \\ 
\displaystyle Z_{\bar \tau}\in \bigoplus_{j+k=n-1} \Omega^k(M;\mmf^{k,-2})\end{array} \right\}
\eeq
satisfying
\beq\label{eq:deriveddata}
 \dR Z=0\quad \partial_{v}Z=\dR Z_v\quad \partial_{\bar \tau}Z=\dR Z_{\bar \tau},
\eeq
where $d$ is the exterior derivative, $\partial_v$ is the derivative with respect to the volume of a lattice, and $\partial_{\bar \tau}$ is the derivative with respect to the antiholomorphic part of the conformal modulus of the lattice. In particular, a degree~$n$ field theory determines a (de~Rham) cohomology class $[Z]$ with the properties
\beq
&&[\partial_{v}Z]=[\dR Z_v]=0,\quad [\partial_{\bar \tau}Z]=[\dR Z_{\bar \tau}]=0, 
\implies [Z]\in \H^n(M;\MF^\omega)\label{eq:underlyingclass}
\eeq
implying that $[Z]$ is a class valued in the ring $\MF^\omega$ of weak modular forms. 
\end{prop}

\begin{rmk}
The datum $Z_{\bar \tau}$ with property~\eqref{eq:underlyingclass} recovers the holomorphic anomaly equation~\cite[pages~9 and~14]{GJF2}. 
\end{rmk}

To summarize the discussion so far, restricting a degree~$n$ field theory to supertori and supercylinders in $M$ gives the indicated arrows
\beq
&&\begin{tikzpicture}[baseline=(basepoint)];
\node (T) at (-2.5,1) {$2|1\EFT_{\cutoff}^n(M)$};
\node (A) at (0,0) {$\KO_\MF^n(M)$};
\node (B) at (5,0) {$\KO^n(M)(\!(q)\!)$};
\node (CC) at (0,-1) {$\H^n(M;\MF)$};
\node (C) at (-2.75,-1) {$\H^n(M;\MF^\omega)$};
\node (D) at (5,-1) {$\H^n(M;\C(\!(q)\!)[u^{\pm 2}]).$};
\draw[->] (A) to (B);
\draw[->] (A) to  (CC);
\draw[->] (CC) to (C);
\draw[->] (CC) to node [below] {$q$-expand} (D);
\draw[->] (B) to node [right] {Ch} (D);
\draw[->,bend left=10] (T) to node [above=3pt] {Restrict to cylinders + Index } (B);
\draw[->,bend right=10] (T) to node [left=3pt] {Restrict to tori} (C);
\draw[->,dotted] (T) to (A);
\draw[->,dashed,bend right=10] (T) to (CC);
\path (0,-.12) coordinate (basepoint);
\end{tikzpicture}\label{eq:pullbacksquareEFT}
\eeq
These restrictions satisfy a compatibility property: the value of a field theory on a supertorus is the trace of the value on a supercylinder~\cite[\S2]{STTraces}, 
\beq\label{eq:itsaproperty}
Z=\sTr_{n}(q^{L_0}{\bar q}^{\bar L_0})\in \Omega^\bullet(M;C^\infty(\HHz\times \R_{>0}))
\eeq
in the notation of~\eqref{eq:potentialsemigroup} and \eqref{eq:deriveddata0}. Using that the $S^1$-eigenspaces in Propositon~\ref{prop:mainprop} are bounded below, the equality~\eqref{eq:itsaproperty} implies that $[Z]\in \H^n(M;\MF)\subset \H^n(M;\MF^\omega)$ determines a de~Rham cohomology class with values in weakly holomorphic modular forms i.e., $[Z]$ is meromorphic at $i\infty$ (see Definition~\ref{defn:modularform}). This produces the long-dashed arrow in~\eqref{eq:pullbacksquareEFT}. 

\begin{rmk}
The supertrace on the right in~\eqref{eq:itsaproperty} depends on the degree of the theory; it is a normalized version of the Clifford supertrace, see Remark~\ref{rmk:CliffvsFertrace}. \end{rmk}

The dotted arrow in~\eqref{eq:pullbacksquareEFT} is more difficult: the compatibility property~\eqref{eq:itsaproperty} is not enough information to determine a map to~$\KO^n_\MF$. Instead, one requires homotopy compatibility data relating de~Rham cocycles in $\Omega^n(M;\C(\!(q)\!)[u^{\pm 2}])$ that are extracted from the two sides of~\eqref{eq:itsaproperty}. On the left side, choosing a $\bar\partial_{\bar \tau}$-preimage $T$ of $Z_{\bar \tau}$ in~\eqref{eq:deriveddata} gives
$$
Z_{\bar \tau}=\partial_{\bar \tau} T\implies \partial_{\bar \tau}(Z-dT)=0.
$$
A choice of $T$ always exists because the half plane $\HHz\subset \C$ is Stein. By restricting to lattices of volume~1 and $q$-expanding, we obtain a differential form
\beq\label{eq;holomorphicpartition}
(Z-dT)|_{v=1}\in \Omega^\bullet(M;\C(\!(q)\!)[u^{\pm 2}])\qquad  {\rm (depends\ on \ choice\ of} \ T ).
\eeq 
 On the right hand side of~\eqref{eq:itsaproperty}, the Pontryagin character of a choice of index bundle constructs a cocycle in the Cech--de~Rham complex with coefficients in formal power series,
\beq
&&\Ch(\Ind(\bar G_0))\in \Omega^\bullet(\mathfrak{U};\C(\!(q)\!)[u^{\pm 2}]) \qquad {\rm (Chern \ character \  of\ the\ index \ bundle)}.\label{Eq:Chernofindex}
\eeq
Definition~\ref{defn:admitscutoffsEFT} guarantees that an index bundle exists for each $k\in \Z$, but its construction depends on choices of cutoffs, see~\S\ref{sec:cutoffconstrution}.
Using~\eqref{eq:itsaproperty}, the cocycles \eqref{eq;holomorphicpartition} and~\eqref{Eq:Chernofindex} represent the same cohomology class, but need not be equal as de~Rham cocycles. Bismut's local index theorem~\cite{Bismutindex} specifies a coboundary measuring their difference. This coboundary is the homotopy compatibility required for the dotted arrow in~\eqref{eq:pullbacksquareEFT}.

In~\S\ref{sec:datafromEFT}, we construct the groupoids $2|1\eft^n(M)$ that are purpose-built to capture the pieces of a $2|1$-Euclidean field theory that determine a $\KO_\MF$-cocycle via~\eqref{eq:cocycle}.

 \begin{defn}\label{defn:21eft}
The groupoid $2|1\eft^n(M)$ has objects degree~$n$ representations of~$\Ann^{2|1}(M)$ (Definition~\ref{defn:degreenAnn}) that are reflection positive (Definition~\ref{defn:rssP}), trace class (Definition~\ref{defn:tracestru}) and admit cutoffs (Definition~\ref{defn:admitscutoffsEFT}). Morphisms are isomorphisms between degree~$n$ representations of~$\Ann^{2|1}(M)$ that are compatible with the reflection structure and trace. 
\end{defn}

\begin{prop}[Lemma~\ref{prop:EFT2} and Propositions~\ref{prop:EFT3} and \ref{prop:trace}] \label{rmk:comparetoST}
There is a span
$$
2|1\EFT^n(M) \hookleftarrow 2|1\EFT_{\cutoff}^n(M)\xrightarrow{\rm restrict} 2|1\eft^n(M),
$$
i.e., the groupoid $2|1\eft^n(M)$ admits a restriction map from the subcategory of Stolz and Teichner's field theories that admit cutoffs. 
\end{prop}

By definition, the dashed arrow in~\eqref{eq:pullbacksquareEFT} factors through $2|1\eft^n(M)$. This completes the sketch of the cocycle map~\eqref{eq:cocycle} in Theorem~\ref{thm1}; the complete proof is in~\S\ref{sec:proofofthm1}.

\begin{rmk}\label{rmk:GJFcompare}
When $M=\pt$, the construction of the dotted arrow in~\eqref{eq:pullbacksquareEFT} reduces to the torsion invariant in the physical argument from~\cite[pages~14-15]{GJF2}. In short, there are two possible ways to extract a holomorphic quantity from a $\mathcal{N}=(0,1)$-theory, namely~\eqref{eq;holomorphicpartition} and~\eqref{Eq:Chernofindex}. Each way depends on a choice, and the difference between the choices defines an invariant. The choices in~\cite[pages~14-15]{GJF2} are identical to~\eqref{eq;holomorphicpartition} and~\eqref{Eq:Chernofindex}. 
\end{rmk}

An object of $2|1\eft^n(M)$ is determined by the geometric data in Propositions~\ref{prop:mainprop} and~\ref{prop:maintorusprop} satisfying various properties. Theorems~\ref{thm2} and \ref{thm3} follow from explicit construction of such data and verification of the properties. Theorem~\ref{thm3} (proved in~\S\ref{sec:partialsigmamodel}) uses the Bismut--Ramond superconnection from Proposition~\ref{thm:KOMForientation} to construct a degree~$-n$ representation of supercylinders; Bismut's local index theorem plays a vital role for verifying the cutoff property for this representation. Theorem~\ref{thm2} (proved in~\S\ref{sec:EulercocycleKOMF}) follows by first realizing the $\KO_\MF$-Euler class of a string vector bundle $V\to M$ in terms of an infinite sequence of finite rank Clifford module bundles built out of~$V$, see~\S\ref{sec:KOMFEuler}. A choice of connection on $V$ then supplies the data~\eqref{eq:sequenceofsuperconn} which (after a choice of string structure) determines an object $\Eu(V)\in 2|1\eft^{{\rm dim}(V)}(M)$. In this case, the cutoff property follows from Remark~\ref{rmk:finitetype}. 

\subsection{Acknowledgements} 
This project has been in gestation for several years and has benefited from conversations with many colleagues, including Matt Ando, Mark Behrens, David Ben-Zvi, Emily Cliff, Nora Ganter, Andr\'e Henriques, Tom Nevins, Stephan Stolz and Arnav Tripathy. I would like to express particular thanks to Theo Johnson-Freyd and Peter Teichner for many enlightening conversations about supersymmetric field theories, and to Lennart Meier and Charles Rezk for their help in analyzing the string orientation of~$\KO_\MF$. This work was supported by the National Science Foundation under grant number~DMS-2205835.

\newpage

\setcounter{tocdepth}{1} 
\tableofcontents

\section{Motivation from supersymmetric quantum field theory}\label{sec:motivate}
In this section we review Witten's physical construction of the Witten genus via 2-dimensional quantum field theory with $\mathcal{N}=(0,1)$ supersymmetry~\cite{Witten_Dirac}. His argument shows that the partition function of this type of quantum field theory is an integral modular by exploiting a compatibility property between values on cylinders and values on tori. For families of field theories this compatibility property becomes compatibility data, which is the key observation for extracting $\KO_\MF$-valued invariants; compare diagrams~\eqref{eq:pullbacksquareEFT} and~\eqref{eq:strictpullbacksquare}.

\subsection{Supersymmetric quantum mechanics and the index theorem}
Let $X$ be an $n$-dimensional Riemannian spin manifold with principal spin bundle $P\to X$. The spinor bundle $\bS:=P\times_{\Spin(n)}\cCl_n$ supports the $\cCl_{-n}$-linear Dirac operator $\slashed{D}$, reviewed in~\S\ref{sec:Cliffordandspin} below. McKean and Singer's argument~\cite{McKeanSinger} shows that the Clifford supertrace $\sTr_{\Cl_{-n}}(e^{-t\slashed{D}^2})$ is independent of $t$, and the limit
\beq\label{eq:rationalizedindex}
\lim_{t\to \infty} \sTr_{\Cl_{-n}}(e^{-t\slashed{D}^2})={\rm sdim}_{\Cl_{-n}}(\ker(\slashed{D})),
\eeq
computes the Clifford super dimension of the kernel of $\slashed{D}^2$. Under the Atiyah--Bott--Shapiro construction~\cite{ABS}, the Clifford module $\ker(\slashed{D})$ is the $\KO^{-n}(\pt)$-valued index of the Dirac operator. The Clifford superdimension~\eqref{eq:rationalizedindex} is the rationalization of this class, which equals the usual index of $\slashed{D}$ when the dimension $n$ is divisible by~4, and is zero otherwise.
On the other hand, the $t\to 0$ limit can be identified with the integral of the $\hat{A}$-form, e.g., see~\cite[Theorem~4.21]{BGV}. In summary, the $t\to 0$ and $t\to \infty$ limits give the two sides of the Atiyah--Singer index theorem
\beq\label{eq:partition function2}
&&\sTr_{\Cl_{-n}}(e^{-t\slashed{D}^2})=\left\{ \begin{array}{ll} \displaystyle{\rm sdim}_{\Cl_{-n}}(\ker(\slashed{D})) & t\to \infty \\ \displaystyle \int_X\hat{A}(X) & t\to 0.\end{array}\right.
\eeq
 Below we review a description of~\eqref{eq:partition function2} in the language of supersymmetric quantum mechanics due to Witten~\cite{susymorse} and Alvarez-Gaum\'e~\cite{Alvarez}. 

The basic data of a quantum mechanical system is a $\Z/2$-graded space of states $\mathcal{H}$ with a hermitian inner product and a self-adjoint even operator $H$ (the Hamiltonian) generating time evolution. In the Wick-rotated theory, time evolution is a trace class semigroup representation
\beq\label{eq:semigroup}
&&\R_{> 0}\to \End(\mathcal{H}),\qquad t\mapsto e^{-tH}
\qquad 
\begin{tikzpicture}[baseline=(basepoint)];
\node (D) at (-.5,0) {$\bullet$};
\node (D) at (-.5,1) {$\bullet$};
\draw[thick] (-.5,0) to (-.5,1);
\draw[thick,<->] (-1,0) to node [left] {$t$} (-1,1);
\node (G) at (.25,.5) {$I_t.$};
\path (0,.5) coordinate (basepoint);
\end{tikzpicture}
\eeq
As indicated on the right of~\eqref{eq:semigroup}, one can reformulate this semigroup representation geometrically: the parameter $t$ represents a length of an interval $I_t$. Gluing intervals adds the lengths, $I_t\circ I_s\simeq I_{t+s}$, corresponding to multiplication in the semigroup. 

A superalgebra $A$ acting on $\mathcal{H}$ is a \emph{symmetry algebra} if it commutes with $H$. A trace on $A$ (see~\S\ref{sec:supertrace}) defines the \emph{partition function}
\beq
Z:=\sTr_A(e^{-tH})\in C^\infty(\R_{>0};A/[A,A])\label{eq:partitionfunctioneasy}
\eeq
valued in the abelianization, i.e., the zeroth Hochschild homology of $A$. In the geometric picture, the trace~\eqref{eq:partitionfunctioneasy} corresponds to gluing together the ends of the length $t$ interval in~\eqref{eq:semigroup} to obtain a circle of circumference~$t$. 

\emph{Supersymmetric} quantum mechanics adds the data of odd self-adjoint operators $Q_1,\dots Q_\mathcal{N}$ satisfying~\cite{susymorse}
$$
Q_i^2=\frac{1}{2}[Q_i,Q_i]=H \ {\rm for} \ i=1,\dots ,\mathcal{N} \quad {\rm and} \quad [Q_i,Q_j]=0 \ {\rm for} \ i\ne j,
$$
where the bracket $[-,-]$ is the supercommutator. By the graded Jacobi identity, the $Q_i$ generate a symmetry algebra. 
A Riemannian spin manifold $X$ determines a quantum mechanical system with $\mathcal{N}=1$ supersymmetry via~\cite[Remark 3.2.21]{ST04}
\beq\label{eq:SQMsystem}
\mathcal{H}=\Gamma(X;\bS), \quad H=\slashed{D}^2,\quad Q=\slashed{D},
\eeq
and the Clifford algebra $\Cl_{-n}$ is a symmetry algebra commuting with $Q$. Supersymmetry promotes the semigroup representation~\eqref{eq:semigroup} to a supersemigroup representation,
\beq\label{eq:supersemigroup}
\R^{1|1}_{>0}\to \End(\mathcal{H}),\qquad (t,\theta)\mapsto e^{-tH+\theta D}.
\eeq
The geometric description in terms of intervals of length $t$ in~\eqref{eq:semigroup} generalizes to the gluing of \emph{superintervals} of length $(t,\theta)$. This can be made precise in Stolz and Teichner's formalism for $1|1$-Euclidean field theories, and we refer to the companion paper~\cite{DBEEFT} for details.

A consequence of supersymmetry is that the partition function~\eqref{eq:partitionfunctioneasy} is independent of~$t$. For the system~\eqref{eq:SQMsystem}, the $t\to \infty$ limit computes (e.g., see~\cite[pg.~5-6]{GJF2})
\beq\label{eq:Wittenindex}
&&Z=\sTr_{\cCl_{-n}}(e^{-tQ^2})=\sTr(\Gamma_{-n}\circ e^{-tQ^2})\stackrel{t\to\infty}{=}\begin{array}{c}{\rm 1\hbox{-}point\ function } \\ {\rm of \ \Gamma_{-n}}\end{array} 
\eeq
where $\Gamma_{-n}$ is a particular element of the Clifford action determining the Clifford supertrace, see~\eqref{eq:Gamma}. When $n=0$ (or after a Morita equivalence, when $n=0\ {\rm mod} \ 8$) \eqref{eq:Wittenindex} is a signed count of ground states, called the \emph{Witten index} of the quantum mechanical system. On the other hand the $t\to 0$ limit can be computed using path integral techniques,
\beq\label{eq:pathintegralAhat}
Z=\sTr_{\cCl_{-n}}(e^{-tQ^2})=\int_{LX} e^{-S(\gamma)} d\gamma\stackrel{t\to 0}{=} \int_X \hat{A}(X)
\eeq
as an integral over the (super) loop space. In the $t\to 0$ limit, this integral (formally) localizes onto the constant loops~\cite{AtiyahCircular}, ~\cite[\S4]{susymorse}.
One can then \emph{define} the path integral as a formal application of the fixed point formula in $S^1$-equivariant cohomology, yielding the integral of the $\hat{A}$-form over~$X$, e.g., see~\cite{JonesPetrack}. To summarize, the limits of~\eqref{eq:partition function2} have the interpretation in supersymmetric quantum mechanics
\beq\label{eq:ASindexthm}
&&\left(\begin{array}{c} t\to \infty: \\ {\rm count \ of } \\ {\rm ground\ states} \end{array}\right) \rightsquigarrow \quad{\rm sdim}_{\Cl_{-n}}(\ker(\slashed{D}))=\int_X \hat{A}(X) \quad \leftsquigarrow \left(\begin{array}{c} t\to 0: \\ {\rm localized}\\ {\rm path\ integral} {\rm }\end{array}\right).
\eeq
Physically, the equality~\eqref{eq:ASindexthm} is the agreement between the Hamiltonian and Lagrangian (i.e., path integral) approaches to quantum mechanics, giving a physics argument for the Atiyah--Singer index theorem. The Witten genus and its modularity properties are deduced from a generalization of~\eqref{eq:ASindexthm} in 2-dimensional supersymmetric quantum field theory. 

\subsection{2-dimensional Euclidean quantum field theory}
 The basic data in a 2-dimensional Euclidean quantum field theory is a $\Z/2$-graded state space~$\mathcal{H}$ with a hermitian inner product, a self-adjoint Hamiltonian $H\in \End(\mathcal{H})$ generating infinitesimal time evolution, and a self-adjoint momentum operator $P\in \End(\mathcal{H})$ generating infinitesimal spatial translation. For 2-dimensional quantum field theories compactified on a circle of radius~$r$, the operator~$P$ has spectrum in~$\frac{1}{r}\cdot \Z\subset \R$. Define the operators~\cite[Eq.~106]{WittenICM} 
\beq\label{eq:L0ops}
&&L_0=(H+P)/2,\qquad \bar L_0=(H-P)/2.
\eeq
Consider the semigroup $\HHz/2\pi r\Z$ determined by addition in $\HHz\subset \C$, the upper half plane. Define a representation of this semigroup via
\vspace{-.2in}
\beq\label{eq:QFTsemigroup}
&&\begin{array}{c}\displaystyle q^{L_0}\bar q^{\bar L_0}=e^{-t H+i xP}\colon \HHz/2\pi r\Z\to \End(\mathcal{H}),\\ \null \\ \tau=\frac{x+it}{2\pi}, \ q=e^{2\pi i \tau}.\end{array}\qquad 
\begin{tikzpicture}[baseline=(basepoint)];
\node (D) at (-.5,0) {$\bullet$};
\node (G) at (.3,.8) {$\bullet$};
\draw[thick] (.5,0) to (.5,1);
\draw[thick] (-.5,0) to (-.5,1);
\draw[thick,<->] (-1,0) to node [left] {$t$} (-1,1);
\draw[thick,->, bend left=30] (-.5,1.35) to node [above] {$x$} (.4,1.36);
\node [draw, thick, ellipse, minimum width=1cm,minimum height=.45cm] (b) at (0,0){};
\node [draw, thick, ellipse, minimum width=1cm,minimum height=.45cm] (c) at (0,1){};
\node (G) at (1.5,.5) {$\rC_{r,\tau}$};
\path (0,.5) coordinate (basepoint);
\end{tikzpicture}
\eeq
This generalizes the time-evolution semigroup~\eqref{eq:semigroup} in quantum mechanics. As indicated on the right of~\eqref{eq:QFTsemigroup}, one can reformulate this semigroup representation geometrically: the parameters $(r,\tau)$ determine a cylinder $\rC_{r,\tau}$ viewed as a Euclidean (i.e., flat Riemannian) bordism between circles of radius~$r$, where the complex parameter $\tau$ encodes the height~$t$ of the cylinder and the rotation $x \ {\rm mod} \ 2\pi r$ that measures the difference in parameterizations of the boundary circles~\cite[Figure~2]{WittenICM}. Composition in the semigroup then corresponds to the gluing of cylinders, 
\beq\label{eq:Euclideanpicture}\label{eq:gluecylinderpic}
&&\begin{tikzpicture}[baseline=(basepoint)];
	\begin{pgfonlayer}{nodelayer}
		\node [style=none] (0) at (1.5, 1.75) {};
		\node [style=none] (1) at (1, 1.75) {};
		\node [style=none] (2) at (1.5, 0.25) {};
		\node [style=none] (3) at (1, 0.25) {};
		\node [style=none] (4) at (1.5, -1.25) {};
		\node [style=none] (5) at (1, -1.25) {};
		\node [style=none] (6) at (-2, 2) {};
		\node [style=none] (7) at (-2.5, 2) {};
		\node [style=none] (8) at (-2, 0.5) {};
		\node [style=none] (9) at (-2.5, 0.5) {};
		\node [style=none] (10) at (-2, 0) {};
		\node [style=none] (11) at (-2.5, 0) {};
		\node [style=none] (12) at (-2, -1.5) {};
		\node [style=none] (13) at (-2.5, -1.5) {};
		\node [style=none] (14) at (-2.5, 0.5) {};
		\node [style=none] (15) at (-2.5, 0.5) {};
				\node [style=none] (21) at (-.5, .25) {$\xrightarrow{\rm glue=semigroup}$};
		\node [style=none] (20) at (.6,-1.7) {$\rC_{r,\tau'}\circ \rC_{r,\tau}=\rC_{r,\tau+\tau'}$};
\node [style=none] (19) at (-3,-.5) {$\rC_{r,\tau}$};
\node [style=none] (18) at (-3,1.5) {$\rC_{r,\tau'}$};
	\end{pgfonlayer}
	\begin{pgfonlayer}{edgelayer}
		\draw [bend left=270, looseness=0.75,thick] (0.center) to (1.center);
		\draw [bend left=270,thick,dashed] (2.center) to (3.center);
		\draw [bend left=270,thick,dashed] (4.center) to (5.center);
		\draw [thick] (0.center) to (4.center);
		\draw [thick] (1.center) to (5.center);
		\draw [bend left=90,thick] (0.center) to (1.center);
		\draw [bend left=90, looseness=0.75,thick] (2.center) to (3.center);
		\draw [bend left=90, looseness=0.75,thick] (4.center) to (5.center);
		\draw [bend left=270, looseness=0.75,thick] (6.center) to (7.center);
		\draw [bend left=270,thick,dashed] (8.center) to (9.center);
		\draw [bend left=90, looseness=0.75,thick] (6.center) to (7.center);
		\draw [bend left=90, looseness=0.75,thick] (8.center) to (9.center);
		\draw [thick] (8.center) to (6.center);
		\draw [thick] (9.center) to (7.center);
		\draw [bend left=270, looseness=0.75,thick] (10.center) to (11.center);
		\draw [bend left=270,thick,dashed] (12.center) to (13.center);
		\draw [bend left=90, looseness=0.75,thick] (10.center) to (11.center);
		\draw [bend left=90, looseness=0.75,thick] (12.center) to (13.center);
		\draw [thick] (12.center) to (10.center);
		\draw [thick] (13.center) to (11.center);
	\end{pgfonlayer}
	\path (0,0) coordinate (basepoint);
\end{tikzpicture}
\qquad\qquad\qquad\begin{tikzpicture}[baseline=(basepoint)];
\node [style=none] (19) at (-2,-1.57) {$\rC_{r,\tau}$};
\node [style=none] (19) at (1.25,-1.25) {$\rT_{r,\tau}$};
	\begin{pgfonlayer}{nodelayer}
		\node [style=none] (1) at (2.25, 0) {};
		\node [style=none] (2) at (1.75, 0) {};
		\node [style=none] (4) at (0.5, 0) {};
		\node [style=none] (5) at (0, 0) {};
		\node [style=none] (8) at (-2.5, 1.5) {};
		\node [style=none] (9) at (-2, 1.5) {};
		\node [style=none] (12) at (-2.5, -1) {};
		\node [style=none] (13) at (-2, -1) {};
						\node [style=none] (21) at (-1, .25) {$\xrightarrow{\rm glue=trace}$};
	\end{pgfonlayer}
	\begin{pgfonlayer}{edgelayer}
		\draw [bend left=270, looseness=1.25,thick] (1.center) to (5.center);
		\draw [bend right=90,thick] (2.center) to (4.center);
		\draw [bend left=90,thick,dashed] (4.center) to (5.center);
		\draw [bend right=90,thick] (4.center) to (5.center);
		\draw [bend right=90, looseness=1.25,thick] (5.center) to (1.center);
		\draw [bend right=90,thick] (4.center) to (2.center);
		\draw [bend left=90,thick] (8.center) to (9.center);
		\draw [bend right=90,thick] (8.center) to (9.center);
		\draw [bend right=90, looseness=1.25,thick] (12.center) to (13.center);
		\draw [bend left=90,thick,dashed] (12.center) to (13.center);
		\draw [thick] (9.center) to (13.center);
		\draw [thick] (8.center) to (12.center);
	\end{pgfonlayer}
	\path (0,0) coordinate (basepoint);
\end{tikzpicture}
\eeq
As explained in~\cite[\S3]{WittenICM}, gluing the ends of a cylinder together corresponds to the trace,
$$
Z_r(\tau):=\Tr(q^{L_0}\bar q^{\bar L_0})\in C^\infty(\HHz/2\pi r\Z)
$$
and the \emph{partition function} $Z_r(\tau)$ is interpreted as the value of the field theory on a torus $\rT_{r,\tau}=\R^2/(2\pi r\Z\oplus 2\pi\tau\Z)$ associated with the based lattice $2\pi r\Z\oplus 2\pi\tau\Z\subset \C\simeq \R^2$. Path integral arguments lead one to expect that the partition function only depends on the underlying torus, not the based lattice. Hence, given operators~\eqref{eq:L0ops} for each $r\in \R_{>0}$, we expect an $\SL_2(\Z)$-invariant function~$Z$, 
\beq\label{eq:SL2invariance}
Z(r,\tau):=\Tr(q^{L_0}\bar q^{\bar L_0})\in C^\infty(\R_{>0}\times \HHz)^{\SL_2(\Z)},
\eeq
$$
\left[\begin{array}{cc} a & b \\ c& d\end{array}\right]\cdot (r,\tau)=\left(r|c\tau+d|,\frac{a\tau+b}{c\tau+d}\right),\quad \left[\begin{array}{cc} a & b \\ c& d\end{array}\right]\in \SL_2(\Z)
$$
where in~\eqref{eq:SL2invariance} $L_0$ and $\bar L_0$ are now a family of operators depending on $r\in \R_{>0}$.

An enhancement of the picture~\eqref{eq:Euclideanpicture} endows the surfaces with spin structures, resulting in a \emph{spin Euclidean field theory}. In this case, there are two possible circles of radius $r$ corresponding to periodic spin structure (called the \emph{Ramond sector}) and the antiperiodic spin structure (called the \emph{Neveu--Schwarz sector}). Hence, a spin Euclidean field theory has a pair of Hilbert spaces of states and a pair of semigroup representations~\eqref{eq:QFTsemigroup}, one for each spin structure on the circle. The partition function of a spin theory comes in 4 flavors~\eqref{eq:2dQFTpartition} corresponding to the 4 spin structures on a Euclidean torus $\C/2\pi r \Z\oplus 2 \pi \tau\Z$; see~\S\ref{sec:dumbEuc}. The relationship with the semigroups of operators comes from taking the trace or the supertrace in the Ramond or Neveu-Schwarz sector. The $\SL_2(\Z)$-invariance property~\eqref{eq:SL2invariance} is refined to an $\MP_2(\Z)$-invariance property for the metaplectic double cover, where the action permutes three of the spin structures; see~\S\ref{sec:back2}. Our focus in this paper is on the Ramond sector, and the partition function associated with the supertrace, i.e., the value of the theory on nonbounding tori. 

\subsection{Symmetries, free fermions, and degree~$n$ field theories}

The generalization of a symmetry algebra $A$ in quantum mechanics is a \emph{current algebra} in 2-dimensional field theory. A degree~$n$ Euclidean field theory carries an action by the current algebra of the~$n$ (chiral) free fermions, which is essentially determined by the Clifford algebra of a loop space, see~\eqref{eq:Fernoncircle}. The discussion below is informal, but there are many mathematically rigorous approaches; see \cite[\S8]{SegalCFT} and \cite[\S2.3]{ST04}, and also \cite{DouglasHenriques,Tener,LudewigFer}. 

For a spin circle with radius $r$, consider the (infinite-dimensional) Clifford algebra
\beq\label{eq:Fernoncircle}
\Fer_n(S^1_r):=\Cl\Big(\Gamma(S^1_r;\bS\otimes \C^n), \int_{S^1_r}\langle-,-\rangle\Big)
\eeq
where $\Gamma(S^1_r;\bS\otimes \C^n)$ is the space of sections of the $\C^n$-twisted spinor bundle, $\langle-,-\rangle$ is the $\C$-bilinear extension of the standard inner product on $\R^n$, and the integral is defined using the isomorphism $\bS\otimes \bS\simeq \Omega^1_{S^1_r}$. For the nonbounding spin structure, $\bS$ is the trivial bundle; this allows the identification $\Fer_n(S^1_r)\simeq \cCl(L_r\R^n)$ with the complex Clifford algebra of the loop space of $\R^n$ for loops of radius~$r$ with inner product $\int_{S^1_r}\langle-,-\rangle dt$. 

Below we will focus on two pieces of data coming from a degree~$n$ Euclidean field theory:
\begin{enumerate}
\item[(i.)] a $\Fer_n(S^1_r)$-module structure on~$\mathcal{H}$ compatible with a semigroup representation~\eqref{eq:QFTsemigroup};
\item[(ii.)] a trace valued in a line bundle over the moduli of spin tori. 
\end{enumerate}
Each of these aspects requires a bit of setup. 

We begin with (i). It is useful to view the pair $(\Gamma(S^1_r;\bS\otimes \C^n), \int_{S^1_r}\langle-,-\rangle)$ in~\eqref{eq:Fernoncircle} as an odd symplectic vector space, and $\Fer_n(S^1_r)$ the associated Weyl algebra. From this perspective, a spin cylinder $\rC_{r,\tau}$ determines a Lagrangian subspace
\beq\label{eq:Lagsubspace}
&&\Gamma_{\rm harmonic}(\rC_{r,\tau};\bS\otimes \C^n)\xhookrightarrow{{\rm restrict}} \Gamma(S^1_r;\bS\otimes \C^n)\oplus \Gamma(S^1_r;\bS\otimes \C^n)
\eeq
by restriction of harmonic spinors (i.e., the kernel of the Dirac operator) to the boundaries of $C_{r,\tau}$. A general construction associates to the Lagrangian subspace~\eqref{eq:Lagsubspace} a Fock module $\Fer_n(\rC_{r,\tau})$ that is an irreducible $\Fer_n(S^1_r)$-bimodule~\cite[Definition 2.2.4]{ST04}. This Fock module is quite explicit, simply being $\Fer_n(S^1_r)$ as a bimodule over itself with left action modified by an algebra isomorphism~\cite[page~59]{ST11}
\beq\label{eq:algebraiso}
\varphi_\tau\colon \Fer_n(S^1_r)\to \Fer_n(S^1_r)
\eeq 
for each $\tau \in \HHz$.
For each sector in a spin Euclidean field theory, a degree~$n$ structure determines the data of a $\Fer_n(S^1_r)$-action on $\mathcal{H}$ that is compatible with the semigroup representation~\eqref{eq:QFTsemigroup} via~\eqref{eq:algebraiso}. 

\begin{rmk}
Decomposing into Fourier modes allows one to describe~\eqref{eq:Fernoncircle} and~\eqref{eq:algebraiso} in terms of an infinite tensor product of finite-dimensional Clifford algebras, see~\eqref{eq:Ferdefn},~\eqref{eq:antiperFer} and~\eqref{eq:explicitcomposition}.
\end{rmk}

Now we turn to the trace data (ii) above. In a degree~$n$ Euclidean field theory, the trace of a semigroup representation in $\Fer_n(S^1_r)$-modules takes values in the family of vector spaces~\cite[\S2]{STTraces}
\beq\label{eq:Pfaffianfibers}
&&\sTr_{\Fer_n}(q^{L_0}\bar q^{\bar L_0})\in \Fer_n(\rC_{r,0})\otimes_{\Fer_n(S^1_{r})\otimes \Fer_n(S^1_{r})^\op} \Fer_n(\rC_{r,\tau}) \quad ({\rm for} \ (r,\tau) \ {\rm fixed}),
\eeq
i.e., the zeroth Hochschild homology of $\Fer_n(S^1_{r}$ valued in the bimodule $\Fer_n(\rC_{r,\tau})$, analogously to the trace~\eqref{eq:partitionfunctioneasy}. For each $(r,\tau)$, the vector space~\eqref{eq:Pfaffianfibers} is 1-dimensional. Letting $(r,\tau)\in \R_{>0}\times \HHz$ vary, the resulting line bundle is canonically isomorphic to the restriction of the $n$th power of the Pfaffian line along the covering map~\cite[Definition 2.3.12]{ST04}
$$
\R_{>0}\times \HHz\to \R_{>0}\times \HHz\sq \MP_2(\Z),\quad r\in \R_{>0}, \ \tau\in \HHz.
$$
The $\SL_2(\Z)$-invariance condition~\eqref{eq:SL2invariance} generalizes to $\MP_2(\Z)$-equivariance data for the Pfaffian line and an $\MP_2(\Z)$-invariance property for a section. For the supetrace in the Ramond sector, this reads
\beq\label{eq:2dQFTpartition}
Z:=\sTr_{\Fer_n}(q^{L_0}\bar q^{\bar L_0})\in \Gamma(\R_{>0}\times \HHz;\Pf^{\otimes n})^{\MP_2(\Z)}.
\eeq
These $\MP_2(\Z)$-invariant sections are a $C^\infty$-version of modular forms; see~\S\ref{sec:back2}. 

\begin{rmk}
There are versions of the partition function~\eqref{eq:2dQFTpartition} for the other 3 spin structures on the torus. A choice of sector (either Ramond or Neveu--Schwarz sector) together with a choice of ordinary trace or the supertrace corresponds to the 4 spin structures, where~\eqref{eq:2dQFTpartition} is the supertrace in the Ramond sector. For the remaining 3 choices, the $\MP_2(\Z)$-action permutes the spin structures, and so the invariance condition analogous to~\eqref{eq:2dQFTpartition} involves a compatibility condition among~3 functions, e.g., see~\cite[page~60]{ST11}. 
\end{rmk}

\begin{rmk}\label{rmk:anomaly}
Symmetries of a quantum system are inextricably linked with anomalies. In the case at hand, Fock modules for the algebra~\eqref{eq:Fernoncircle} are closely related positive energy representations of the loop group $L\Spin(n)$. This leads to a connection between degree~$n$ field theories and classical Chern--Simons theory~\cite[\S5.2]{ST04}; roughly a degree~$n$ theory determines a (non-topological) boundary condition for classical Chern--Simons theory of the group $\Spin(n)$ at level~1, which is an invertible topological field theory. Boundary theories (also known as \emph{relative} or \emph{twisted} theories) can be interpreted as anomalous where the bulk theory encodes the anomaly, e.g., see~\cite[\S2.3]{Freedshortrange} and~\cite{FreedTelemanrelative}.
\end{rmk}

\subsection{2-dimensional supersymmetric quantum field theory}\label{sec:motivateQFT}

A further enhancement of a degree~$n$ spin Euclidean field theory is a \emph{supersymmetric Euclidean field theory}, which requires the data of odd square roots of the operators~\eqref{eq:L0ops}. The data of $k$ odd square roots of $L_0$ and $l$ odd square roots of $\bar L_0$ is called \emph{$\mathcal{N}=(k,l)$ supersymmetry}. We will be most interested in $\mathcal{N}=(0,1)$ supersymmetry, i.e., a single odd operator $\bar G_0$ with
$$
\bar G_0^2=\frac{1}{2}[\bar G_0,\bar G_0]=\bar L_0. 
$$
Generalizing~\eqref{eq:supersemigroup} in supersymmetric quantum mechanics, the operator $\bar G_0$ promotes the semigroup representation~\eqref{eq:QFTsemigroup} to a \emph{super semigroup} representation 
\beq\label{Eq:ogsupersemi}
q^{L_0}\bar{q}^{\bar L_0+\sqrt{2\pi i}\theta Q},\qquad \HHz^{2|1}\to \End(\mathcal{H})
\eeq
where $\HHz^{2|1}\subset \R^{2|1}$ is the upper half-space of the super Lie group $\R^{2|1}$ with multiplication
$$
(\tau,\bar\tau,\theta)\cdot (\tau',\bar\tau',\theta')=(\tau+\tau',\bar \tau+\bar \tau+\theta'\theta',\theta+\theta'), \qquad (\tau,\bar\tau,\theta), (\tau',\bar\tau',\theta')\in \R^{2|1}.
$$
In this case, the McKean--Singer argument that results in $t$-independence in~\eqref{eq:partition function2} results in $\bar q$-independence in~\eqref{eq:2dQFTpartition}~\cite[page~46]{ST04}, 
\beq
Z(r,\tau)&=&\sTr_A(q^{L_0}\bar q^{\bar L_0})=\sTr_A(q^{L_0}|_{{\rm ker}(\bar L_0)})=\sum_{k\in \Z} \sTr_A(q^{L_0}|_{{\rm ker}(\bar L_0)\bigcap \mathcal{H}_k})\nonumber\\
&=&\sum_{k\in \Z} q^k{\rm sdim}_{\Fer_n}({\rm ker}(\bar L_0)\bigcap \mathcal{H}_k)\nonumber 
\eeq
where $\mathcal{H}_k$ is the $k$th eigenspace of $L_0-\bar L_0$. The above computation has three consequences. First, $Z=Z(r,\tau)$ is independent of $\bar q$, and hence depends holomorphically on $\HHz$. Second, the $q$-expansion of $Z$ has integral coefficients, implying that the $\R_{>0}$~dependence is trivial as there is no nontrivial $\R_{>0}$-deformation of an integral power series. Third and finally, since a trace class operator is compact and the semigroup representation is valued in trace class operators, the powers of $q$ must be bounded below. This brings us to the punchline.

\begin{expect}[{\cite[Theorem 3.3.14]{ST04}}] The partition function~\eqref{eq:2dQFTpartition} of a 2-dimensional field theory with $\mathcal{N}=(0,1)$ supersymmetry and degree~$n$ is an integral modular form of weight $-n/2$.\end{expect}

\begin{rmk}
Another source of modular forms in 2-dimensional field theory are \emph{chiral conformal} field theories. 
A \emph{conformal} field theory includes the data of operators $L_k$ and $\bar L_k$ on~$\mathcal{H}$ for $k\in \Z$  generating a representation of the Virasoro algebra on~$\mathcal{H}$ that extends~\eqref{eq:2dQFTpartition}. A \emph{chiral} field theory has $\bar L_k=0$ the zero operator for all $k\in \Z$. Hence, a chiral theory is canonically endowed with a generator of $\mathcal{N}=(0,1)$ supersymmetry by taking $\bar G_0=0$. It is expected that chiral conformal theories have unique extensions to Euclidean field theories with $\mathcal{N}=(0,1)$ supersymmetry; see Hypothesis~\ref{hyp:chiral} and Remark~\ref{Eq:chiraltosusy} below. 
\end{rmk}

The geometric interpretation of the super semigroup~\eqref{Eq:ogsupersemi} involves the gluing of Euclidean \emph{supercylinders}. These can be visualized using the same pictures as \eqref{eq:Euclideanpicture}, but where the moduli parameters are specified by $(r,\tau,\bar\tau,\theta)\in \R_{>0}\times \HHz^{2|1}$, and gluing of Euclidean supercylinders gives
$$
\rC_{r,\tau,\bar\tau,\theta}\circ \rC_{r,\tau',\bar\tau',\theta'}\simeq \rC_{r,\tau+\tau',\bar\tau+\bar\tau'+\theta\theta',\theta+\theta'}.
$$
The picture in~\eqref{Eq:ogsupersemi} corresponding to the trace similarly involves a moduli space of supertori, determined by a lattice $(\ell_1,\bar\ell_1,\lambda_1)\Z\oplus (\ell_2,\bar\ell_2,\lambda_2)\Z\subset \R^{2|1}$. These ideas can be made precise in Stolz and Teichner's framework. 
The moduli spaces of these objects are remarkably rich, and subtle features of their geometry are crucial in the proof of Theorem~\ref{thm1}.

\subsection{The Dirac--Ramond operator and the Witten genus} \label{sec:WGmotivate}
Next we sketch Witten's partial construction of the $\mathcal{N}=(0,1)$-supersymmetric sigma model with target a string manifold. To start, for a real vector bundle $V\to X$, define the formal power series of vector bundles
\beq\label{eq:wittenpolynomial}
W(V)&=&\bigotimes_{k\ge 0} \Sym_{q^k} (V\otimes  \C)=\bigoplus_{k\ge 0} q^k W_k(V)
\eeq
where $\Sym^l(V)$ denotes the $l^{\rm th}$ symmetric power of $V$ and 
\beq\label{eq:Symqdefn}
&&\Sym_{q^k}(V\otimes  \C)=\bigoplus_{l\ge 0} q^{kl} \Sym^l(V\otimes \C)=\underline{\C}+q^k(V\otimes  \C)+q^{2k} \Sym^2(V\otimes\C)+\dots 
\eeq
If $V$ is equipped with a connection~$\nabla^V$, then the vector bundles $W_k(V)$ inherit a connection~$W_k(\nabla^V)$, and we write $W(\nabla^V)$ for the formal power series of these connections. 

\begin{defn}\label{defn:DiracRamond}
The \emph{Dirac--Ramond} operator of a spin manifold $X$ is the formal power series of $\cCl_{-n}$-linear twisted Dirac operators (see~\cite[Eq.~17]{Witten_Dirac})
$$
\bar G_0:=\slashed{D}\otimes W(\nabla^{TX})=\sum_{k \ge 0} q^k\slashed{D}\otimes W_k(\nabla^{TX}), \quad \bar G_0\colon \Gamma(\bS\otimes W(TX))\to \Gamma(\bS\otimes W(TX))
$$
where $\slashed{D}$ is the $\Cl_{-n}$-linear Dirac operator on $X$ and $\nabla^{TX}$ is the Levi-Civita connection on~$TX$. Define the \emph{index} of the Dirac--Ramond operator as the formal power series,
\beq\label{eq:IndofDiracR}
\Ind(\bar G_0):=\phi(q)^n \sum_{k \ge 0} q^k\Ind_{\Cl_{-n}}(\slashed{D}\otimes W_k(TX))\in \KO^{-n}(\pt)(\!(q)\!)
\eeq
where $\Ind_{\Cl_{-n}}$ is the (real) Clifford index and $\phi(q)$ is defined as
\beq\label{eq:phiproductdefn}
\phi(q)=\prod_{m>0}(1-q^m)\in \C(\!(q)\!).
\eeq
\end{defn}

\begin{rmk} \label{rmk:CliffvsFertrace}
The normalization by $\phi(q)$ in~\eqref{eq:IndofDiracR} comes the trace for $n$~copies of the chiral free fermion. The algebra $\Fer_n(S^1_r)$ is Morita equivalent to the Clifford algebra~$\Cl_n$. Consequently, the trace for the fermions is compatible with the Clifford supertrace after normalizing by the trace of the Morita bimodule~\eqref{eq:Morita}, which is precisely~$\phi(q)^n$. This normalization also agrees with~\cite[Equations~17-18]{Witten_Dirac}. 
\end{rmk}

The Ramond sector of the supersymmetric $\sigma$-model with target $X$ has as its space of states the sections $\mathcal{H}=\Gamma(X;\bS\otimes W(TX))$. For $r\in \R_{>0}$ define the operators
$$
\bar G_0=\slashed{D}\otimes W(\nabla^{TX}),\qquad P|_{W_k(TX))}=k/r,\quad \implies\bar L_0=\bar G_0^2,\quad L_0=\frac{1}{2}(H+P)=\bar G_0^2+P.
$$
The algebra $\Cl_{-n}$ acts on $\mathcal{H}$ commuting with the above operators. Using the normalized trace outlined in Remark~\ref{rmk:CliffvsFertrace}, we have the partition function
\beq\label{eq:itsWittenspartition}
Z:=\phi(q)^n\sTr_{\Cl_{-n}}(q^{L_0}\bar q^{\bar L_0})=\Ind(\bar G_0)\in \KO^{-n}(\pt)(\!(q)\!).
\eeq
Hence, $Z$ encodes a signed count of ground states of $H$ for each momentum number $k$.


One may also compute the partition function~\eqref{eq:itsWittenspartition} in the path integral formalism. To describe the result of this computation we recall the following. 

\begin{defn}\label{defn:wittengenus}
Let $X$ be a Riemannian manifold with Levi-Civita connection $\nabla^X$. The \emph{Witten form} of $X$ is
\beq\nonumber
\Witt(\nabla^X)=\exp\left(\sum_{k\ge 1} \frac{\Tr(F^{2k})E_{2k}(q)}{(2\pi i)^{2k} (2k)}\right)\in \Omega^\bullet(X;\C(\!(q)\!))
\eeq
where $F=(\nabla^X)^2$ is the curvature of the Levi-Civita connection.  The \emph{Witten genus} of $X$ is 
\beq\label{Eq:Wittengenus}
\langle \Witt(\nabla^X),[X]\rangle=\int_X \Witt(\nabla^X)\in \C(\!(q)\!).
\eeq
If $X$ is even dimensional with a string structure (so that $[\Tr(F^2)]=0$) then the Witten genus is a modular form of weight $n/2$. Note that~\eqref{Eq:Wittengenus} is automatically zero if~$n\ne 0 \ {\rm mod} \ 4$.
\end{defn}

Analogously to~\eqref{eq:pathintegralAhat}, the path integral expression for the partition function can be rigorously defined by localization to the $T^2$-fixed points~\cite[\S2.6]{WitteninStrings}, \cite{BElocalization} 
\beq\label{eq:pathint}
Z(\tau)=\int_{\Map(T_\tau,X)} e^{-S(\phi)}d\phi\stackrel{{\rm vol}\to 0}{=} \int_X \Witt(\nabla^X)
\eeq
and a formal application of the fixed point theorem results in an integral of the Witten form. However, there is some ambiguity in this computation coming from the fact that the 2nd Eisenstein series is defined by a conditionally convergent series. This is an example of an \emph{anomaly}: a priori, different choices result in different values of the path integral~\eqref{eq:pathint}. However, when $[p_1(\nabla^{X})]=0$ all choices give the same result, and the anomaly vanishes.

In summary, the Ramond sector partition function of the supersymmetric $\sigma$-model with target $X$ can be computed two ways: the first is a signed count of ground states, whereas the second evaluates the path integral using localization:
\beq
&&\begin{tikzpicture}[baseline=(basepoint)];
\node (T) at (-3,.5) {$\begin{array}{c} X\ {\rm a \ } 4k\hbox{-}{\rm dimensional} \\ {\rm string\ manifold} \end{array}$};
\node (A) at (0,0) {$\MF^\Z$};
\node (B) at (4,0) {$\Z(\!(q)\!)$};
\node (C) at (0,-1.25) {$\MF$};
\node (D) at (4,-1.25) {$\C(\!(q)\!).$};
\node (P) at (.7,-.4) {\scalebox{1.5}{$\lrcorner$}};
\draw[->,] (A) to (B);
\draw[->] (A) to  (C);
\draw[->] (C) to node [below] {$q$-expand} (D);
\draw[->,right hook-latex] (B) to (D);
\draw[|->,bend left=12] (T) to node [right=40pt] {$\begin{array}{c} {\rm signed\ count\ of \ ground} \\ {\rm  states} \ {\rm (supercylinders)} \end{array}$} (B);
\draw[|->,bend right=40] (T) to node [left=8pt] {$\begin{array}{c}{\rm path\ integral}\\{\rm (supertori)}\end{array}$} (C);
\draw[|->,dashed] (T) to (A);
\path (0,-.65) coordinate (basepoint);
\end{tikzpicture}\label{eq:strictpullbacksquare}
\eeq
The signed count of ground states is extracted from the value of the field theory on supercylinders, whereas the path integral computes the value of the field theory on supertori. 
Using the Atiyah--Singer index theorem and the identity (proved in \cite{Zagiermodular})
\beq\label{eq:Weierstrassidentity}
&&\sigma(\tau,z)=\frac{e^{z/2}-e^{-z/2}}{z/2}\prod_{k>0}\frac{(1-q^kz)(1-q^kz^{-1})}{(1-q^k)^2}=\exp\left(-\sum_{k\ge 1} \frac{E_{2k}(q)}{2k(2\pi i)^{2k}}z^{2k}\right).
\eeq
 these two computations of the partition function are compatible when $X$ has a string structure. Hence we obtain an equality of elements of~$\C(\!(q)\!)$
\beq\label{eq:LagHam}
\Z(\!(q)\!) \ni \Ind(\bar G_0)=\int \Witt(\nabla^X)\in \MF
\eeq
showing that the Witten genus is an integral modular form.

In analogy to the physics argument for the index theorem~\eqref{eq:ASindexthm}, the equality~\eqref{eq:LagHam} is an instance where the Lagrangian (i.e., path integral) agrees with the Hamiltonian formalism. When passing to families of field theories there is a surprise: the analogous computations no longer agree on the nose. Instead, the results are homotopy compatible where the field theory itself specifies the data of the homotopy. This is the basic mechanism behind Theorem~\ref{thm1} that allows one to generalize Witten's argument for invariants of field theories valued in integral modular forms as the (strict) pullback in rings~\eqref{eq:strictpullbacksquare} to invariants valued in $\KO_\MF$ as the (homotopy) pullback in commutative ring spectra~\eqref{eq:pullbacksquare}.

\section{Approximating $\TMF$ by $\KO_\MF$}\label{sec:KOMF}

\subsection{The cohomology theory $\KO_\MF$}
The cohomology theory~$\TMF$ is constructed as the derived global sections of the sheaf of elliptic cohomology theories (or rather, their representing $E_\infty$-rings) on the moduli stack of elliptic curves~$\Mell$, e.g., see~\cite[\S1.3]{Lurie_Elliptic} for an overview. Restriction of sections along the map $\Spec(\Z(\!(q)\!)) \to \Mell$ classifying the punctured formal neighborhood of the Tate curve constructs the \emph{Miller character}
\beq
\TMF\to \KO(\!(q)\!).\label{eq:Millerchar}
\eeq
Indeed, the formal group law associated with the Tate curve is isomorphic to the multiplicative group, and so restriction initially yields a map of commutative ring spectra $\TMF\to \K(\!(q)\!)$. However, the Tate curve has a $\Z/2$-action by automorphisms, implying that restriction factors through the $\Z/2$-homotopy fixed points $\K(\!(q)\!)^{h\Z/2}\simeq \KO(\!(q)\!)$. 

\begin{rmk}
Lurie's approach to elliptic cohomology constructs~\eqref{eq:Millerchar} as a map of $E_\infty$-rings, see~\cite[\S4.3]{Lurie_Elliptic} for a sketch. A map~\eqref{eq:Millerchar} of commutative ring spectra is constructed explicitly in~\cite[Appendix~A]{HillLawson}, where $\KO(\!(q)\!)$ is defined as 
\beq\label{eq:KTateascolim}
\KO(\!(q)\!)=q^{-1} \hocolim_k \KO\wedge \{1,q,q^2,\dots,q^{k-1}\}_+,
\eeq
and the map~\eqref{eq:Millerchar} is shown to be unique up to homotopy. Below, it will suffice to work with the Miller character as a map of commutative ring spectra. 
\end{rmk}

Smashing the map~\eqref{eq:Millerchar} with the spectrum $\H_\C$ gives the outer homotopy commutative square in~\eqref{eq:pullbacksquare} where $\MF$ is the ring of weakly holomorphic modular forms, and we use
$$
\pi_*(\TMF)\otimes \C\simeq \MF,\qquad \pi_*(\KO(\!(q)\!))\otimes \C\simeq \C(\!(q)\!)[u^{\pm 2}],\qquad |u^2|=-4. 
$$
\begin{defn} Define $\KO_\MF$ in~\eqref{eq:pullbacksquare} as the homotopy pullback in commutative ring spectra. Let $\K_\MF$ be the similarly defined commutative ring spectrum gotten by replacing $\KO(\!(q)\!)$ by $\K(\!(q)\!)$ in~\eqref{eq:pullbacksquare}.
\end{defn} 

The spectrum $\KO_\MF$ inherits a map from $\TMF$ by the universal property~\eqref{eq:pullbacksquare}, and in this way is an approximation to $\TMF$. 

\begin{prop}\label{prop:htpygrps}
The coefficients of $\K_\MF$ and $\KO_\MF$ are
\beq
\pi_*\K_\MF&\simeq& \left\{ \begin{array}{cl} 
\MF^\Z_k & *=2k \\ 
\C(\!(q)\!)/(\Z(\!(q)\!)+\MF_{2k}) & k=2k-1
\end{array}\right.\label{eq:KMFcoeff}\\
\pi_*\KO_\MF&\simeq& \left\{ \begin{array}{cl}  
\MF_{4k}^\Z & *=8k\\ 
\Z/2(\!(q)\!) & *=8k+1\\
\Z/2(\!(q)\!) & *=8k+2\\  
\C(\!(q)\!)/(2\Z(\!(q)\!)+\MF_{4k+2})& *=8k+3 \\
\MF^\Z_{4k+2} & *=8k+4\\
0 & *=8k+5 \\ 
0 & *=8k+6 \\ 
\C(\!(q)\!)/(\Z(\!(q)\!)+\MF_{4k})& *=8k-1 \\
\end{array}\right.\label{eq:KOMFcoeff}
\eeq
where $\C(\!(q)\!)/(\Z(\!(q)\!)+\MF_{2k})$ is the quotient of Laurent series with $\C$ coefficients by the subgroup generated by Laurent series with $\Z$ coefficients and the image of $q$-expansion.
\end{prop}
\bp
We compute using the presentations 
\beq\label{eq:Kcoeffs}
\pi_*\K\simeq \Z[u^{\pm 1}],\qquad \pi_*\KO=\Z[\eta,\alpha,\beta^{\pm 1}]/\langle 2\eta,\eta^3,\eta\alpha, \alpha^2-4\beta\rangle,
\eeq
where generators have the degrees $|u|=-2$, $|\alpha|=-4$, $|\eta|=-1$ and $|\beta|=-8$. The ring homomorphism $\pi_*\KO\to \pi_*\K$ is determined by~$\beta\mapsto u^4$, $\alpha\mapsto 2u^2$, $\eta\mapsto 0$. The vertical arrow on the right in~\eqref{eq:pullbacksquare} is determined by the Pontryagin character in real K-theory, defined as the composition of complexification and the usual Chern character,
$$
\KO \to \K\xrightarrow{\Ch} \H_{\C[u^{\pm 1}]},\qquad \pi_*(\K)\otimes \C\simeq \C[u^{\pm 1}],\qquad |u|=-2.
$$
Hence, this composition has image in the 4-periodic subring $\C[u^{\pm 2}]\subset \C[u^{\pm 1}]$. The lower horizontal arrow in~\eqref{eq:pullbacksquare} can be identified with $q$-expansion~\eqref{eq:qexpand} of modular forms~\cite{LauresPhD}. 

With the above facts in place, the result follows from the Mayer--Vietoris sequence for the homotopy pullback square~\eqref{eq:pullbacksquare}. For $\K_\MF$, this sequence reads
$$
0\to \pi_{2k}\K_\MF\to u^{2k}\cdot \Z(\!(q)\!)\oplus \MF_{2k}\to u^{2k}\cdot \C(\!(q)\!)\to \pi_{2k-1}\K_\MF\to 0
$$
using that $\pi_*\K$ is concentrated in even degrees, and so~\eqref{eq:KMFcoeff} follows immediately. The real case is argued similarly, using exactness of the sequence
\beq
\pi_{*+1}\KO(\!(q)\!)\oplus \pi_{*+1}\H_{\MF}&\to& \pi_{*+1} \H_{\C(\!(q)\!)[u^{\pm 2}]}\to \pi_{*}\KO_\MF\nonumber\\
&\to& \pi_{*}\KO(\!(q)\!)\oplus \pi_{*} \H_\MF \to \pi_{*}(\H_{\C(\!(q)\!)[u^{\pm 2}]})\nonumber
\eeq
and the 8-fold periodicity of $\KO$ to establish the 8 cases above. For example, with $*=8k$ we have 
$$
\pi_{8k+1} \H_{\C(\!(q)\!)[u^{\pm 2}]}\simeq \{0\} , \quad \pi_{8k}\KO(\!(q)\!)\oplus \pi_{8k} \H_\MF\simeq \beta^k \cdot \Z(\!(q)\!) \oplus \MF_{4k} ,$$$$ \pi_{8k}(\H_{\C(\!(q)\!)[u^{\pm 2}]})\simeq u^{4k}\cdot \C(\!(q)\!)
$$
confirming $\pi_{8k}\KO_\MF\simeq \MF_{4k}^\Z$. 
\ep

\begin{rmk} 
The ring structure on $\pi_*\KO_\MF$ takes a bit more work to spell out. However, the ring homomorphisms $\pi_*\MString\to \pi_*\TMF\to \pi_*\KO_\MF$ allow us to compute most of the products of interest using the known formulas in string cobordism or $\pi_*\TMF$. For example, the generator $\nu\in \pi_3(\TMF)\simeq \Z/24$ is sent to a (nontrivial) class $\nu\in \pi_3\KO_\MF$ with order~24~\cite[\S3]{BunkeNaumann}, and this class necessarily satisfies $\eta^3=12\nu$ for $\eta\in \pi_1(\KO_\MF)$. 
\end{rmk}

\subsection{Families of string manifolds and the families Witten form}

The following enhances a family of Riemannian spin manifolds (e.g., see \cite[Geometric Data 1.1]{Freed_Det}) to include a geometric string structure in the sense of~\cite{Waldorfstring}. 

\begin{defn} \label{defn:stringstructure}
A \emph{family of Riemannian string manifolds} is
\begin{enumerate}
\item a smooth fibration of manifolds $\pi\colon X\to M$ where the vertical tangent bundle $T(X/M)\to X$ has a chosen metric $g^{(X/M)}$ and spin structure;
\item a projection $P:TX\to T(X/M)$; and
\item a trivialization of the Chern--Simons 2-gerbe with connection associated with the spin structure on $T(X/M)$ (see \cite[Definition 1.1.5.]{Waldorfstring}). 
\end{enumerate} 
\end{defn}

\begin{rmk}\label{rmk:CS2gerbe}
A \emph{family of Riemannian spin manifolds} is the data (1) and (2); these uniquely specify a (Levi-Civita) connection $\nabla^{X/M}$ on $T(X/M)$~\cite[Lemma 1.3]{Freed_Det} which determines the Chern--Simons 2-gerbe with connection in (3) as we now describe. For each $n\in \N$, the \emph{Chern--Simons 2-gerbe} is a 2-gerbe with connection over the stack classifying $n$-dimensional spin bundles with connection, denoted $*\sq \Spin(n)^\nabla$. Given a map of stacks $X\to *\sq \Spin(n)^\nabla$ classifying a metrized vector bundle $V$ with compatible connection $\nabla^V$ and choice of spin structure, the Chern--Simons 2-gerbe with connection can be pulled back to~$X$. The curvature of the resulting 2-gerbe on $X$ is precisely~$\frac{p_1}{2}(\nabla^V)$, the fractional first Pontryagin form of $V$. A trivialization of the Chern--Simons 2-gerbe associated with $V$ is in particular a choice of coboundary for~$\frac{p_1}{2}(\nabla^V)$. In the constructions below,  the only information we will use from part (3) of Definition~\ref{defn:stringstructure} is the 3-form $H\in \Omega^3(X)$ with $dH=p_1(\nabla^{X/M})\in \Omega^4(M;\C)$ determined by the image of the integral trivialization of the degree~4 cocycle $\frac{p_1}{2}(\nabla^{X/M})\in Z^4(X;\Z)$. 
\end{rmk}

\begin{defn} 
The \emph{Witten form} of a family of Riemannian spin manifolds is
\beq\label{Eq:familiesWF}
&&\Witt(\nabla^{X/M})=\exp\left(\sum_{k\ge 1} \frac{u^{2k}\Tr(F^{2k})E_{2k}(q)}{2k(2\pi i)^{2k}}\right)\in \Omega^\bullet(X;\C(\!(q)\!)[u^{\pm 2}])
\eeq
where $F$ is the curvature of $\nabla^{X/M}$. 
\end{defn}

When $M=\pt$, the above recovers Definition~\ref{defn:wittengenus}. We emphasize that $E_2(q)$ is not a modular form, see~\S\ref{sec:etageometry}. 
Hence, for a family $\pi\colon X\to M$ that admits a string structure, the integral $[\int_{X/M} \Witt(\nabla^{X/M})]\in \H(M;\MF)$ determines a cohomology class with coefficients in modular forms. However, the differential form $\int_{X/M} \Witt(\nabla^{X/M})$ need not be valued in modular forms. Instead we have the following.

\begin{prop}\label{prop:firsthomotopy}
A choice of Riemannian string structure on the family $\pi\colon X\to M$ uniquely specifies an odd differential form $\eta_H \in \Omega^\bullet(X;\C(\!(q)\!)[u^{\pm 2}])$ with the property that
\beq\label{Eq:homotopytomodular}
&&\Wit(\nabla^{X/M}):=\Witt(\nabla^{X/M})+d\eta_H\in \Omega^\bullet(X;\MF)\subset \Omega^\bullet(X;\C(\!(q)\!)[u^{\pm 2}])
\eeq
is valued in modular forms. Integrating along the fibers $\pi\colon X\to M$,
$$
\int_{X/M} \Wit(\nabla^{X/M})=\int_{X/M} \Witt(\nabla^{X/M})+d\int_{X/M} \eta_H\in \Omega^\bullet(M;\MF),
$$
therefore also yields a differential form valued in modular forms. 
\end{prop}

\bp
We recall that $p_1(\nabla^{X/M})= -\Tr(F^2)/8\pi^2$, and a geometric string structure specifies a 3-form $H$ with $dH=p_1(\nabla^{X/M})\in \Omega^4(M;\C)$, see Remark~\ref{rmk:CS2gerbe}. Then define
\beq\label{eq:thehomotopy}
\eta_H=\left(H \cdot \Witt(\nabla^{X/M})\cdot\frac{e^{E_2p_1}-1}{p_1}\right)
\eeq
where above $p_1=p_1(\nabla^{X/M})$ and we use divisibility of $e^{E_2p_1}-1$ by~$p_1$. Then $\eta_H$ satisfies
$$
d\eta_H=\left(e^{E_2p_1}\Witt(\nabla^{X/M})-\Witt(\nabla^{X/M})\right).
$$
Since
\beq\label{eq:modularWF}
&&\Wit(\nabla^{X/M}):=e^{E_2p_1}\Witt(\nabla^{X/M})=\exp\left(\sum_{k\ge 2} \frac{u^{2k}\Tr(F^{2k})E_{2k}(q)}{2k(2\pi i)^{2k}}\right)\in \Omega^\bullet(X;\MF)
\eeq
is a differential form valued in modular forms, the result follows. \ep

For future use, we also record a version of the families Witten form using the  non-holomorphic (but modular)  2nd Eisenstein series $E_2^*(\tau,\bar\tau)$, see~\eqref{Eq:2ndEisen}. Define
\beq\label{eq:nonholoWit}
&&\Witt^*(\nabla^{X/M})=\exp\left(\frac{u^{2}\Tr(F^{2})E^*_{2}(\tau,\bar\tau)}{2(2\pi i)^{2k}}+\sum_{k\ge 2} \frac{u^{2k}\Tr(F^{2k})E_{2k}(\tau)}{2k(2\pi i)^{2k}}\right).
\eeq
We have $\Witt^*(\nabla^{X/M})\in \Omega^\bullet(X;C^\infty(\HHz)[u^{\pm 2}])^{\SL_2(\Z)}$, where $\SL_2(\Z)$ acts on $C^\infty(\HHz)[u^{\pm 1}]$ as 
$$
f(\tau,\bar\tau)\mapsto f\left(\frac{a\tau+b}{c\tau+d},\frac{a\bar\tau+b}{c\bar\tau+d}\right),\qquad u\mapsto \frac{u}{c\tau+d}.
$$
Hence, $\Witt^*(\nabla^{X/M})$ is a differential form valued in a $C^\infty$-version of modular forms, see~\eqref{eq:MFassubspace}. It compares with the families Witten form~\eqref{Eq:familiesWF} as 
\beq\label{eq:nonholoWitdelta}
&&\Witt(\nabla^{X/M})=\exp\left(\frac{2\pi i }{\tau-\bar\tau} p_1(\nabla^{X/M}) \right)\Witt^*(\nabla^{X/M}),\quad p_1(\nabla^{X/M})=-\frac{u^{2}\Tr(F^2)}{8\pi^2},
\eeq
and $\Witt^*(\nabla^{X/M})$ satisfies a formula analogous to~\eqref{Eq:homotopytomodular} where
\beq\label{eq:thehomotopy2}
\eta_H=\left(H \cdot \Witt^*(\nabla^{X/M})\cdot\frac{e^{E_2^*p_1}-1}{p_1}\right),
\eeq
where above $p_1=p_1(\nabla^{X/M})$. 

\begin{rmk}
An explicit formula for $\partial_{\bar \tau}\Witt^*(\nabla^{X/M})$ is determined by~\eqref{eq:E2isntholo}. One can use this to extend $\Witt^*(\nabla^{X/M})$ to an $\SL_2(\Z)$-invariant differential form valued in the differential graded algebra $(\Omega^{0,*}(\HHz)[u^{\pm 1}],\bar\partial)$, i.e., $\Witt^*(\nabla^{X/M})$ is part of the data of a derived global section over the stack $\HHz\sq \SL_2(\Z)$. 
\end{rmk}
\subsection{The Bismut--Ramond superconnection}\label{sec:BRsuperconn}

We recall some of the ingredients in Bismut's families index theorem~\cite{Bismutindex}, see also~\S\ref{sec:KOindex} below. For $\pi\colon X\to M$ a family of Riemannian spin manifolds with principal spin bundle $P\to X$, define the spinor bundle over $X$ as the associated bundle,
$$
\bS=P\times_{\Spin(n)} \Cl_n
$$
which carries a fiberwise left action by the algebra bundle $\Cl(T(X/M))$ and a (commuting) fiberwise right action by~$\Cl_n$. This leads to an $M$-family of Dirac operators $\slashed{D}$ acting on sections of $\bS$. Since $\slashed{D}$ graded commutes with the right action of $\Cl_n$, it commutes with the associated left action of $\Cl_{-n}=\Cl_n^\op$. For $V\to X$ a metrized real vector bundle with compatible connection~$\nabla^V$, define the Fr\'echet vector bundle $\pi_*(\bS\otimes V)\to X$ whose fiber at $x\in M$ is the Fr\'echet space $\Gamma((\bS\otimes V)|_{X_x})$ where $X_x=\pi^{-1}(\{x\})$. The Dirac operator $\slashed{D}\otimes \nabla^V$ defines an odd, Clifford-linear, self-adjoint endomorphism of $\pi_*(\bS\otimes V)\to M$. Let~$\B^V$ denote the $\Cl_{-n}$-linear Bismut superconnection adapted to~$\slashed{D}\otimes \nabla^V$, 
\beq\label{eq:Bismutsuperconn}
\B^V=\slashed{D}+\nabla^{\pi_*(\bS\otimes V)}-\frac{c(T)}{4}
\eeq
where $\nabla^{\pi_*(\bS\otimes V)}$ is a Clifford linear unitary connection on $\pi_*(\bS\otimes V)$, and $c(T)$ is Clifford multiplication by the curvature 2-form $T$ of the fiber bundle $\pi\colon X\to M$. Clifford multiplication by $T$ uses the left action of $\Cl(T(X/M))$ on the spinor bundle, which commutes with the right $\Cl_n$-action. We conclude that~\eqref{eq:Bismutsuperconn} commutes with the left $\Cl_{-n}$-action. A standard reference for the Bismut superconnection is \cite[Ch.~10]{BGV}; the Clifford linear generalization is in~\cite[\S7]{DBEIndex}. Finally, we recall the $\R_{>0}$-rescaling action 
\beq\label{eq:rescale}
\B^V(t)=t^{1/2}\slashed{D}+\nabla^{\pi_*(\bS\otimes V)}-t^{-1/2}\frac{c(T)}{4},\qquad t\in \R_{>0},
\eeq
that yields a 1-parameter family of superconnections whose behavior as $t\to 0$ and $t\to \infty$ controls Bismut's families index theorem~\cite{Bismutindex}, \cite[Ch.~10]{BGV} generalizing~\eqref{eq:partition function2}.

\begin{defn} 
Given a family $\pi \colon M\to X$ of Riemannian spin manifolds, define the \emph{Bismut--Ramond superconnection} as the direct sum of $\Cl_{-n}$-linear superconnections
\beq\label{Eq:BRconn}
\A:=\B^{W(T(X/M))}=\bigoplus_{k\ge 0} q^k \B^{W_k(T(X/M))}
\eeq
defined on the Fr\'echet vector bundle over $X$
\beq\label{eq:BRV}
\pi_*(\bS\otimes W(T(X/M)))=\bigoplus_{k\ge 0} q^k \pi_*(\bS\otimes W_k(T(X/M))).
\eeq
\end{defn}

We use the notation 
$$\HBR:=\pi_*(\bS\otimes W(T(X/M)))\qquad \HBR_k:=\pi_*(\bS\otimes W_k(T(X/M)))
$$ 
so that $\HBR=\bigoplus q^k\HBR_k$. Similarly, let $\A_k:=\A|_{\HBR_k}=\B^{W_k(T(X/M))}$ denote the summand of the Bismut--Ramond superconnection on $\HBR_k$ so $\A=\bigoplus_k \A_k$. Using the $\R_{>0}$-rescaling action on Bismut superconnections~\eqref{eq:rescale}
we obtain the rescaled Bismut--Ramond superconnection,
\beq\label{eq:rescaledBR}
&&\A(t)=\bigoplus_{k\ge 0} q^k \A_{k}(t),\qquad \A_{k}(t)=t^{1/2}\slashed{D}\otimes W_k(\nabla^{X/M})+\nabla^{\HBR_k}-t^{-1/2}\frac{c(T_k)}{4}.
\eeq

\begin{defn}
The \emph{Chern character} of the Bismut--Ramond superconnection is
$$
\Ch(\A)=\phi(q)^n\sum_{k\ge 0} q^k \sTr_{\cCl_{-n}}(e^{-u\A_{k}(u)^2})\in \Omega^\bullet(M;\C(\!(q)\!)[u^{\pm 2}])
$$
where $\phi(q)$ is defined in~\eqref{eq:phiproductdefn}. We note $\Ch(\A)$ has total degree~$-n$ (see Definition~\ref{defn:chernform}).
\end{defn}

\begin{prop}\label{prop:Bismuts}
There is a uniquely specified differential form $\eta_W\in \Omega^\bullet(X;\C(\!(q)\!)[u^{\pm 2}])$,
\beq\label{eq:BismutstheoremforWG}
\int_{X/M}\Witt(\nabla^{X/M})&=&\Ch(\A)+d\eta_W  
\eeq
that mediates between the Chern character of the Bismut--Ramond superconnection and the integral of the Witten form of the family $\pi\colon X\to M$.
\end{prop}
\bp
The identity~\eqref{eq:Weierstrassidentity} gives the equality of differential forms, 
\beq
&&\phi(q)^n\sum_{k \ge 0} q^k\int_{X/M} \hat{A}(\nabla^{X/M})\Ch(\nabla^{W_k(T(X/M))})=\int_{X/M} \exp\left(\sum_{k\ge 0} \frac{\Tr(F^{2k})E_{2k}(q)}{2k(2\pi i)^{2k}}\right).
\eeq
We further have the transgression formula (see \cite[Theorem~10.32]{BGV} when $n$ is even and~\cite[Theorem~2.10]{BismutFreed2} when~$n$ is odd)
$$
\phi(q)^n\sTr_{\cCl_{-n}}(e^{-\A_{k}(s)^2})=\phi(q)^n\int_{X/M} \hat{A}(\nabla^{X/M})\Ch(\nabla^{W_k(T(X/M))})-d\int_0^s\alpha_k(t)dt,
$$ 
that measures the change in the Chern character of the Bismut--Ramond superconnection under the rescaling~\eqref{eq:rescaledBR}, where
$$
\alpha_k(s):=\phi(q)^n\sTr_{\cCl_{-n}}\left(\frac{d\A_{k}(s)}{ds} e^{-(\A_{k}(s))^2}\right).
$$
Hence if we define 
\beq\label{eq:Whomotopy}
\eta_W=\sum_{k\ge 0} q^k \int_0^1\alpha_k(t)dt
\eeq
the result follows.
\ep

\subsection{The families index of Dirac--Ramond operators with vector bundle kernels}

\begin{prop}\label{prop:index1}
For each $k\ge 0$, suppose that the fiberwise kernels 
$$
\HBR_k^0:=\ker(\slashed{D}\otimes W_k(\nabla^{X/M})|_{M_x})\subset (\HBR_k)_x
$$ 
have locally constant rank for $x\in M$. Let $p_k\colon \HBR_k\to \HBR_k^0$ denote the fiberwise orthogonal projection to the kernel.  Then there is a uniquely determined differential form $\eta_0$ satisfying
\beq\label{eq:anothertrans}
\Ch(\A)-\phi(q)^n\sum_{k\ge 0} q^k\Ch(p_k\nabla^{\HBR_k}p_k)=d \eta_0 \in \Omega^\bullet(X;\C(\!(q)\!))
\eeq
mediating between the Chern character of the superconnection and the power series of Chern characters of the kernel bundles. 
\end{prop}

\bp 
The argument follows from standard facts about Bismut superconnections; we review the key points. The spectral projection operators for a family of Dirac operators are smoothing~\cite[Proposition~9.10]{BGV}, and in particular the projections $p_k\colon \HBR_k\to \HBR_k^0$ are smoothing for each $k$. This implies that 
\beq\label{eq:thekernelbundles}
\HBR_k^0 \to X,\qquad p_k\nabla^{\HBR_k}p_k,
\eeq
is a smooth vector bundle with smooth connection. Then using the 1-parameter family of rescaled superconnections~\eqref{eq:rescaledBR}, define
$$
\alpha_k(s):=\phi(q)^n\sTr_{\cCl_{-n}}\left(\frac{d\A_{k}(s)}{ds} e^{-(\A_{k}(s))^2}\right),\qquad \eta_0=\sum_{k\ge 0} q^k\int_1^\infty \alpha_k(t)dt.
$$
In the limit $s\to \infty$, the Chern form of $\A_{k}(s)$ approaches the Chern form of~\eqref{eq:thekernelbundles}~\cite[Theorem~9.19]{BGV}. Then~\eqref{eq:anothertrans} follows from general argument for Chern--Simons forms provided that the infinite integral converges; convergence follows form large~$s$ bounds on the $C^l$-norm of $\int_1^s \alpha(t)dt$~\cite[Theorem~9.23]{BGV}. 
\ep

\begin{prop} \label{prop:KOMFindex1}
Under the assumption of Proposition~\ref{prop:index1}, a family of Riemannian string manifolds $\pi\colon X\to M$ uniquely determines a class $\sigma(X)\in \KO^{-n}_\MF(M)$. 
\end{prop}

\bp
Since the kernels $\HBR_k^0$ of the operators $\slashed{D}\otimes W_k(\nabla^{X/M})$ are bundles of $\Cl_{-n}$-modules over $M$, they determine a class
\beq
\phi(q)^n\sum_{k\ge 0} q^k [\HBR_k^0]\in \KO^{-n}(M)(\!(q)\!)\label{eq:firstKOclass}
\eeq
by the Atiyah--Bott--Shapiro construction~\cite{ABS}. In fact, for each $k$ the kernel bundle determines a map to the point-set model for the $\KO$-spectrum constructed in~\cite[Theorem 7.1]{HST}, see~\cite[\S4]{DBEIndex}. A differential form representative for the Chern character of~\eqref{eq:firstKOclass} is given by $\phi(q)^n\sum_{k\ge 0} q^k\Ch(p_k\nabla^{\HBR_k}p_k)$. On the other hand, the integral of the (modular) Witten form gives a cocycle 
\beq\label{eq:firstMFclass}
\int_{X/M}\Wit(\nabla^{X/M})\in \Omega^\bullet(M;\MF).
\eeq
Finally, we obtain a homotopy between the Chern form of~\eqref{eq:firstKOclass} and the $q$-expansion of~\eqref{eq:firstMFclass},
\beq\label{eq:firsthomotopy}
&&\int_{X/M}\Wit(\nabla^{X/M})=\phi(q)^n\sum_{k\ge 0} q^k\Ch(p_k\nabla^{\HBR_k}p_k)+d(\eta_0+\eta_H+\eta_W)
\eeq
using Propositions \ref{prop:firsthomotopy}, \ref{prop:Bismuts} and~\ref{prop:index1}. From the definition~\eqref{eq:pullbacksquare} of $\KO_\MF$ as a homotopy pullback, the data~\eqref{eq:firstKOclass}, \eqref{eq:firstMFclass}, and \eqref{eq:firsthomotopy} determine a class $[\sigma(X)]\in \KO^{-n}_\MF(M)$. 
\ep

\subsection{The proof of Proposition~\ref{thm:KOMForientation}}\label{sec:thm:KOMForientation}

To remove the vector bundle kernel assumption in Proposition~\ref{prop:index1}, we apply the families Clifford index construction to the each of the superconnections~$\A_k$; we refer to \S\ref{sec:cutoffconstrution} for a review. By standard properties of the Bimsut superconnection, the $\A_k$ admit a smooth index bundle in the sense of Definition~\ref{defn:cutoffs} (see Proposition~\ref{prop:Bismutcutoffs}). Consequently, for each $k$ we obtain an open cover 
$$
U_k^\lambda:=\{x\in X\mid \lambda\notin{\rm Spec}((\slashed{D}\otimes W_k(\nabla^{X/M}))^2_x)\}\subset X
$$
and for each $(\lambda,k)$ there is a smooth, finite-rank vector bundle with connection over $U_k^\lambda$ 
\beq\label{eq:initialcutoff}
\HBR_k^{<\lambda}:=\bigoplus_{\mu<\lambda} \HBR_k^\mu,\qquad \nabla^{<\lambda}:=p^\lambda_k\circ \A^{[1]}|_{U_\lambda}\circ p^\lambda_k, \quad p^\lambda_k\colon \HBR_k\to \HBR_k^{<\lambda}
\eeq
with fiberwise $\Cl_{-n}$-action. Above, $\HBR_k^{<\lambda}$ is the (finite) direct sum~\eqref{eq:bundleofClnmodules} of $\mu$-eigenspaces of $(\slashed{D}\otimes W_k(\nabla^{X/M})^2$ for $\mu<\lambda$, and $p^\lambda_k$ is the projection~\eqref{eq:fiberwiseincludeproject}.
Over each intersection $U_k^\lambda\bigcap U_k^{\lambda'}$ for $\lambda<\lambda'$, there is a $\Cl_{-n}$-equivariant inclusion
$$
\HBR_k^{<\lambda}\hookrightarrow \HBR_k^{<\lambda'}
$$
and the restriction of $\slashed{D}\otimes W_k(\nabla^{X/M})$ to the orthogonal complement $(\HBR_k^{<\lambda})^\perp\subset \HBR_k^{<\lambda'}$ determines a $\Cl_{-1}$-action~\eqref{eq:compatibilitydata2} generated by an odd operator $e_k^{\lambda\lambda'}$ that commutes with the $\Cl_{-n}$-action. If we choose a locally finite subcover $\{U_k^\lambda\}_{\lambda\in \Lambda_k}$ for a discrete subset $\Lambda_k\subset \R_{>0}$, the above determines a map from $X$ to a model for the $-n$th space in the $\KO$-spectrum~\eqref{eq:KOindexspec},
\beq\label{eq:KOcocycledata}
[\{U_k^\lambda,\HBR_k^{<\lambda},e_k^{\lambda\lambda'}\}_{\lambda\in \Lambda_k}]\colon M\to \KO^{-n}.
\eeq
Different choices in this construction lead to a homotopic map (with specific homotopy). The class underlying the map~\eqref{eq:KOcocycledata} is the analytic index of the $M$-family of Dirac operators $\slashed{D}\otimes W_k(\nabla^{X/M})$.

\begin{lem} \label{lem:KOqmap}
The data~\eqref{eq:KOcocycledata} determine a map $\sum_{k\ge 0} [\Ind_{\Lambda_k}(\A_k)]\colon M\to  (\KO(\!(q)\!))^{-n}$ from~$M$ to a model for the $-n$th space in the spectrum $\KO(\!(q)\!)$. Different choices in this construction lead to a homotopic map with specific homotopy.
\end{lem}
\bp
For each $N$, define an open cover of $X$ indexed by $\vec{\lambda}=(\lambda_1,\dots,\lambda_N)\in \R_{>0}^N$
\beq\label{eq:veclambda}
U^{\vec{\lambda}}&:=&\{x\in X\mid \lambda_k\notin{\rm Spec}((\slashed{D}\otimes W_k(\nabla^{X/M})^2_x)\}\subset X \ {\rm for} \ k\le N\}\\\
&=&U_k^{\lambda_1}\bigcap \dots \bigcap U_k^{\lambda_N}\nonumber
\eeq
i.e., $U^{\vec\lambda}$ is a mutual refinement of the covers $\{U_k^\lambda\}$ for $k\le N$. For each $N$, the construction~\eqref{eq:KOcocycledata} gives a map 
\beq\label{eq:maptoholim}
\left[\{U^{\vec\lambda},\bigoplus_{k=0}^N \HBR_k^{<\lambda_k},\bigoplus_{k=0}^Nf_k^{\lambda\lambda'}\}_{\vec\lambda}\right]\colon X\to (\KO\wedge \{1,q,\dots q^N\}_+ )^{-n}
\eeq
i.e., a map to the $-n$th space in the $N$th spectrum in the homotopy colimit diagram~\eqref{eq:KTateascolim}. A refinement of an open cover determines a continuous functor between \v{C}ech categories (see Definition~\ref{ex:opencover0}); so by \cite[Proposition~2.1]{Segalclassifying} the maps constructed above come equipped with canonical homotopy compatibilities for different values of~$N$. Therefore, the maps~\eqref{eq:maptoholim} are part of the data of  a map into  the homotopy colimit diagram~\eqref{eq:KTateascolim} and hence define a map into the spectrum $\KO(\!(q)\!)$. For different choices of locally finite cover $U^\lambda_k$, we obtain homotopies between maps~\eqref{eq:KOindexspec}  by Proposition~\ref{prop:appenindexbundle}, which in turn give homotopies between maps to the homotopy colimit~\eqref{eq:KTateascolim}.
\ep

\begin{lem}\label{lem:Cqqmap}
There is a uniquely determined coboundary between the Pontryagin character of the cocycle constructed in Lemma~\ref{lem:KOqmap} and the differential form $\sum q^k \Ch(\A_k(t))\in \Omega^\bullet(M;C^\infty(\R_{>0})(\!(q)\!))$, viewing $t\in C^\infty(\R_{>0})$ as the standard coordinate function. 
\end{lem}
\bp
By Lemma~\ref{lem:anotherdamnhomotopy}, for each $(\lambda,k)$ there are uniquely determined coboundaries,
\beq\label{eq:somecoboundaries}
d\eta_k^\lambda=\Ch(\A_k(t))|_{U_\lambda}-\Ch(\nabla^{<\lambda_k}).
\eeq
The Pontryagin character of the cocycle constructed in Lemma~\ref{lem:KOqmap} is a \v{C}ech--de~Rham cocycle for each $N$ 
\beq\label{eq:limitofCherns}
&&\sum_{k \ge 0}^N q^k [U^{\vec\lambda},\Ch(\nabla^{<\lambda_k}),\CS(\nabla^{<\lambda_k},\nabla^{<\lambda'_k})]\in C^*(\{U^{\vec \lambda}\} ,\Omega^\bullet(-;{\C[q]/q^{N+1}})).
\eeq
with compatibilities between different $N$ coming from refinements of the cover. For each~$N$, restricting the coboundaries~\eqref{eq:somecoboundaries} along a refinement of covers gives a coboundary in the \v{C}ech--de~Rham complex between~\eqref{eq:limitofCherns} and $\sum_{k \ge0}^N q^k\Ch(\A_k(t))$. 
\ep

\begin{proof}[Proof of Proposition~\ref{thm:KOMForientation}]
Lemmas~\ref{lem:KOqmap} and \ref{lem:Cqqmap} give a cocycle representative of the families index $\sum_{k \ge 0} [\Ind(\slashed{D}\otimes W_k(\nabla^{X/M}))] \in \KO^{-n}(M)(\!(q)\!)$ and a coboundary mediating between a cocycle representative of the Pontryagin character of the index and the Pontryagin character $\Ch(\A)\in \Omega^\bullet(M;\C(\!(q)\!)[u^{\pm 2}])$ of the Bismut--Ramond superconnection. Then Proposition~\ref{prop:firsthomotopy} and~\ref{prop:Bismuts} give a homotopy between $\Ch(\A)$ and $\int \Wit(\nabla^{X/M})$, where $\int \Wit(\nabla^{X/M})\in \Omega^\bullet(M;\MF)$. Together this determines a continuous map $\sigma(X)\colon M\to \KO_\MF^{-n}$ with underlying class $[\sigma(X)]\in \KO_\MF^{-n}(M)$. A priori this class depends on the choice of $\Lambda_k\subset \R_{>0}$ for each $k$, but any two choices are homotopic with specified homotopy by Proposition~\ref{prop:appenindexbundle} and standard properties of Chern--Simons forms.
\ep

\begin{proof}[Proof of Corollary~\ref{cor:BN1}]
In dimension $4k$, identifying the integral modular form $[\sigma(X)]\in \KO^{-4k}_\MF(\pt)\simeq \MF_{2k}^\Z$ with the Witten genus follows from the classical argument in~\S\ref{sec:WGmotivate}. 

In dimensions $4k-1$, we compare the value of the analytic index in Proposition~\ref{thm:KOMForientation} with the invariant from~\cite[Definition~3.1]{BunkeNaumann}. Let $X$ be a $(4k-1)$-dimensional string manifold. For degree reasons, we note that $\int_X \Wit(\nabla^X)=0$ and Chern character of the Dirac--Ramond operator and the index bundle is zero,
$$
\Ch(\A)=\sum \Ch(\slashed{D}\otimes W_j(TX))=0\in \Omega^{-4k+1}(\pt;\C(\!(q)\!)[u^{\pm 2}])\simeq \{0\}.
$$
Hence the class in $[\sigma(X)]\in \KO_\MF^{-(4k-1)}(\pt)$ is determined by: (i) the index of the sequence of $\Cl_{-4k+1}$-linear Dirac operators $\Ind_{\Lambda_j}(\slashed{D}\otimes W_j(TX))\colon \pt\to \KO^{-4k+1}$, and (ii) a homotopy from the 0 differential form to itself determined by
$$
\eta_H+\eta_W\in \Omega^{-4k}(\pt;\C(\!(q)\!)[u^{\pm 2}])\simeq \C(\!(q)\!)
$$
in the notation of~\eqref{eq:thehomotopy} and~\eqref{eq:Whomotopy}. 

We perform a homotopy that trivializes the data (i) and modifies the data (ii). The index ${\rm Ker}(\slashed{D}\otimes W_j(TX))$ is a $\Cl_{-4k+1}$-module, for which there necessarily exists an extension to $\Cl_{-4k+1}\otimes \Cl_1$-module structure~\cite{ABS}. The generator of the $\Cl_1$-action provides a deformation $(\slashed{D}\otimes W_j(TX))(t):=\slashed{D}\otimes W_j(TX)+\delta_j(t)$ for a family of smoothing operators $\delta_j(t)$ where 
\beq\label{eq:thedeformation}
\slashed{D}\otimes W_j(TX)+\delta_j(t)=\left\{\begin{array}{ll} \slashed{D}\otimes W_j(TX) & t<1/2 \\ {\rm invertible} & t>1 \end{array} \right.
\eeq 
This 1-parameter family of superconnections has an associated Chern--Simons form $\eta_j$ satisfying ~\cite[Theorem~9.26]{BGV} 
$$
\Ch(\slashed{D}\otimes W_j(TX))=d\eta_j. 
$$
In Bunke's framework~\cite{bunkeindex}, the deformation~\eqref{eq:thedeformation} is the data of a \emph{taming}, and the form~$\eta_j$ has an explicit description as an integral~\cite[\S1.2.7]{bunkeindex} and \cite[3.2]{BunkeNaumann}.
The class $[\sigma(X)]\in \KO_\MF^{-4k+1}(\pt)$ is then determined by the image under the quotient map
$$
[\eta_H+\eta_W+\phi(q)^n\sum q^j\eta_j]\in \KO_\MF^{-4k+1}(\pt),\quad \C(\!(q)\!)\to \C(\!(q)\!)/(\varepsilon\Z(\!(q)\!)+\MF_{2k})\simeq \KO_\MF^{-4k+1}(\pt).
$$
This recovers Bunke and Naumann's formula~\cite[Eq.~22]{BunkeNaumann} for their analytic invariant of string manifolds.
\ep

\subsection{The $\KO_\MF$ index theorem}\label{sec:proofofthm:string}
The index theorem looks to compare dashed arrows 
\beq
&&\begin{tikzpicture}[baseline=(basepoint)];
\node (T) at (-2,.75) {$\MString$};
\node (A) at (0,0) {$\KO_\MF$};
\node (B) at (4,0) {$\KO(\!(q)\!)$};
\node (C) at (0,-1) {$\H_\MF$};
\node (D) at (4,-1) {$\H_{\C(\!(q)\!)[u^{\pm 2}]}.$};
\node (P) at (.7,-.4) {\scalebox{1.5}{$\lrcorner$}};
\draw[->,] (A) to (B);
\draw[->] (A) to  (C);
\draw[->] (C) to (D);
\draw[->] (B) to (D);
\draw[->,bend left=10] (T) to (B);
\draw[->,bend right=10] (T) to (C);
\draw[->,dashed] (T) to (A);
\path (0,-.25) coordinate (basepoint);
\end{tikzpicture}\label{eq:pullbackorientation}
\eeq
Such maps are characterized by string of $\KO(\!(q)\!)$ and $\H_\MF$ together with a homotopy witnessing their compatibility. 
The following shows that the homotopy compatibility data is unique up to homotopy when it exists. 

\begin{lem}\label{lem:indedthmlemma}
There is an exact sequence 
$$
0\to [\MString,\KO_\MF]\hookrightarrow [\MString,\KO(\!(q)\!)]\times [\MString,\H_\MF]\to [\MString,\H_{\C(\!(q)\!)}]
$$
of homotopy classes of maps of commutative ring spectra. 
\end{lem}
\bp
General properties of a homotopy pullback provide the Mayer--Vietoris sequence,
\beq
\cdots&\to& [\MString,\Sigma^{-1}\H_{\C(\!(q)\!)[u^{\pm 2}]}] \to [\MString,\KO_\MF]\nonumber\\
&&\hookrightarrow [\MString,\KO(\!(q)\!)]\times [\MString,\H_\MF]\twoheadrightarrow [\MString,\H_{\C(\!(q)\!)[u^{\pm 2}]}]\nonumber
\eeq
But $\MString$ is rationally even, so $[\MString,\Sigma^{-1}\H_{\C(\!(q)\!)[u^{\pm 2}]}] =0$ and the lemma follows. 
\ep

\begin{proof}[Proof of Theorem~\ref{thm:string}]
By Lemma~\ref{lem:indedthmlemma}, maps $\MString\to \KO_\MF$ are homotopic if and only if their composite maps $\MString \to \KO(\!(q)\!)$ and $\MString\to \H_\MF$ are homotopic.
For the topological index~\eqref{eq:ogstringorientation}, the composite map $\MString\to \KO(\!(q)\!)$ sits in a homotopy commutative diagram of commutative ring spectra~\cite{AHR}
\beq
&&\begin{tikzpicture}[baseline=(basepoint)];
\node (A) at (0,0) {$\MString$};
\node (B) at (4,0) {$\TMF$};
\node (C) at (0,-1.25) {$\MSpin$};
\node (D) at (4,-1.25) {$\KO(\!(q)\!)$};
\draw[->,] (A) to node [above] {$\sigma$} (B);
\draw[->] (A) to node [left] {forget}  (C);
\draw[->] (C) to node [below] {A-B-S} (D);
\draw[->] (B) to node [right] {Miller} (D);
\path (0,-.65) coordinate (basepoint);
\end{tikzpicture}\nonumber
\eeq
for the map $\MSpin\to \KO(\!(q)\!)$ coming from the extension of the Atiyah--Bott--Shapiro orientation whose map on coefficients is the Witten genus of a spin manifold. By the Atiyah--Singer index theorem, this is homotopic to the map corresponding to the analytic index in Proposition~\ref{thm:KOMForientation}. On the other hand, the map 
\beq
&&\MString \to \TMF\to \H_{\MF}\to \H_{\C(\!(q)\!)[u^{\pm 2}]}.\label{eq:secondstringorientation}
\eeq
is an orientation of a rational cohomology theory, and hence differs from the standard orientation $\MString\to \MSO \to \H_{\C(\!(q)\!)[u^{\pm 2}]}$ by a Riemann--Roch factor gotten from the ratio of Thom classes, e.g., see~\cite[\S3.3]{AHR}. This difference class is the invertible element
\beq\label{eq:differenceclass}
\exp\left(\sum_{k\ge 2} \frac{u^{2k}p_k E_{2k}(q)}{2k}\right)\in \H^0(\BString;\C(\!(q)\!)[u^{\pm 2}])
\eeq
where ${\rm p}_{k}$ are the pullbacks of the universal Pontryagin classes along the forgetful map $\BString\to \BSO$. This is the same class~\eqref{eq:firstMFclass} associated with the analytic index, and so the composite maps $\MString\to \H_\MF$ are homotopic. This completes the proof. 
\ep

\subsection{The $\KO_\MF$-Euler class}\label{sec:KOMFEuler}
By general considerations, the string orientation~\eqref{eq:ogstringorientation} of~$\KO_\MF$ determines Euler classes for vector bundles $V\to M$ with string structures,
\beq\label{eq:itstheKOMFEuler}
\Eu(V)\in \KO_\MF^{{\rm dim}(V)}(M).
\eeq
Below we construct a representative of~$\Eu(V)$ built from Clifford modules and differential forms, using known descriptions of Euler classes of string vector bundles in $\KO(\!(q)\!)$ and in~$\H_\MF$. The following is analogous to Definition~\ref{defn:stringstructure}.

\begin{defn} \label{defn:stringstructureV}
Given a real vector bundle $V\to M$ with metric compatible connection~$\nabla^V$ a \emph{geometric string structure on $(V,\nabla^V)$} is
\begin{enumerate}
\item a spin structure on $V$;
\item a trivialization of the Chern--Simons 2-gerbe with connection associated with the spin structure on $V$ (see \cite[Definition 1.1.5.]{Waldorfstring}). 
\end{enumerate} 
\end{defn} 

\begin{rmk} 
Analogously to Remark~\ref{rmk:CS2gerbe}, the only geometric data we require below from the trivialization of the Chern--Simons 2-gerbe is the differential form $H\in \Omega^3(M)$ that is the image of the integral trivialization of the cocycle $p_1(V)\in Z^4(M;\Z)$. 
\end{rmk}

A spin structure on an $n$-dimensional vector bundle $V\to M$ determines a bundle of left $\Cl_n=\Cl^\op_{-n}$-modules,
\beq\label{eq:defnKOEuler}
\bS_V:=P\times_{\Spin_n}\Cl_{-n},\qquad [\bS_V]\in \KO^n(M)
\eeq
that represents the $\KO$-Euler class of $V$, e.g., see \cite[Remark 3.2.22]{ST04},~\cite[\S3.6]{FHTI} or \cite[1.92]{Freedalg}. 

\begin{lem}\label{lem:EulervbKO}
For a vector bundle $V\to M$ with geometric string structure, the $\KO(\!(q)\!)$-Euler class associated with the string orientation is
\beq\label{eq:Eulervb}
\phi(q)^{-{\rm dim}(V)}[\bS_V\otimes \bigotimes_{k\ge 0} \Lambda_{q^k} V]\in \KO^n(M)(\!(q)\!)
\eeq
where
\beq\label{eq:exteriordefn}
\Lambda_{q^k} V=\C\oplus q^k \Lambda^1 V\oplus q^{2k} \Lambda^2 V\oplus \cdots
\eeq
\end{lem}
\bp
Since the $\MString$-orientation of $\KO(\!(q)\!)$ factors through $\MSpin$, the $\KO(\!(q)\!)$-Euler class for the string orientation differs from the standard Euler class~\eqref{eq:defnKOEuler} by a Riemann--Roch factor, i.e., a difference class in $\KO({\rm BSpin})(\!(q)\!)$. This difference class was computed by Miller~\cite{Miller}, yielding~\eqref{eq:Eulervb}. 
\ep

\begin{rmk}
In contrast to~\eqref{eq:Symqdefn}, we note that the sum~\eqref{eq:exteriordefn} is finite when $V$ is finite-dimensional. It is related to Witten's spinor bundle on loop space \cite[Equation~24]{Witten_Dirac} in the same manner that the usual spinor bundle~\eqref{eq:spinorsdefn} is related to the Euler class in $\KO$-theory~\eqref{eq:defnKOEuler}. \end{rmk}

Similarly, the Euler class associated to the string orientation of $\H_\MF$ differs from the standard Euler class by a Riemann--Roch factor.

\begin{lem}
For a vector bundle $V\to M$ with geometric string structure, the $\H_\MF$-Euler class associated with the string orientation is
\beq\label{eq:EulervbHMF}
[\Pf(F)\Wit^{-1}(\nabla^V)]\in \H^{{\rm dim}(V)}(M;\MF)
\eeq
where $\Pf(F)$ the cocycle representative of the ordinary Euler class of~$V$ for the standard $\MSO$-orientation of $\H_\MF$, i.e., $\Pf(F)$ is the Pfaffian of the curvature of the connection on~$V$. 
\end{lem}
\bp
The proof is an easier version of the proof of Lemma~\ref{lem:EulervbKO}, e.g., see~\cite[\S3.3]{DBE_MQ}. The difference class is the modular Witten form, $\Wit^{-1}(\nabla^V)$, see~\eqref{eq:modularWF}. 
\ep

\begin{prop} 
The $\KO_\MF$-Euler class associated with the string orientation~\eqref{eq:ogstringorientation} is determined by~\eqref{eq:Eulervb} and~\eqref{eq:EulervbHMF} together with a homotopy (constructed in the proof) that is unique up to higher homotopy. 
\end{prop}
\bp
We extend the data~\eqref{eq:Eulervb} and~\eqref{eq:EulervbHMF} to a class in $\KO_\MF$ by constructing a homotopy between differential forms in $\Omega^\bullet(M;\C(\!(q)\!)[u^{\pm 2}])$. The cocycle representing the class~\eqref{eq:Eulervb} has a Chern form determined by the geometric string structure on~$V$: the vector bundle $\bS_V\otimes \bigotimes_{k\ge 0} \Lambda_{q^k} V$ inherits a connection from the connection on~$V$. Using~\eqref{eq:Weierstrassidentity} and the splitting principle, this Chern form is
\beq\label{eq:ChernofEuler}
\phi(\tau)^{-{\rm dim}(V)}\Ch(\nabla^{\bS_V\otimes \bigotimes_{k\ge 0} \Lambda_{q^k} V})&=& \Pf(F)\exp\left(-\sum_{k\ge 1} \frac{\Tr(F^{2k})E_{2k}(q)}{(2\pi i)^{2k} (2k)}\right)\\
&=&\Pf(F)\Witt(\nabla^V)^{-1}\in \Omega^\bullet(M;\C(\!(q)\!)[u^{\pm 2}]).\nonumber
\eeq
The geometric string structure on $V$ further specifies a coboundary mediating between the Chern form~\eqref{eq:ChernofEuler} and the differential form representative of the $q$-expansion of~\eqref{eq:EulervbHMF}, 
\beq\label{eq:thehomotopyV}
&&\Pf(F)(\Wit^{-1}(\nabla^V)-\Witt^{-1}(\nabla^V))=\Pf(F)d\eta_H,\nonumber\\
&& \eta_H=\left(H \cdot \Witt(\nabla^V)^{-1}\cdot\frac{-e^{E_2p_1}-1}{p_1}\right).
\eeq
Hence, we find that the data~\eqref{eq:Eulervb},~\eqref{eq:EulervbHMF}, and~\eqref{eq:thehomotopyV} determine a cocycle a class in $\KO^{{\rm dim}(V)}_\MF(M)$. By Lemma~\ref{lem:indedthmlemma}, the homotopy~\eqref{eq:thehomotopyV} is unique up to higher homotopy, and so this class is necessarily $\KO_\MF$-Euler class~\eqref{eq:itstheKOMFEuler} associated to the string orientation~\eqref{eq:ogstringorientation}.
\ep

\section{Super Euclidean geometry in dimension~2}
\subsection{Super Euclidean geometries}\label{sec:EBord}
Before specializing to the dimensions of interest, we collect some general facts about super Euclidean geometry in arbitrary super dimension. The following generalizes Thurston's notion of \emph{model geometry} \cite[Definition~3.8.1]{Thurston}. 

\begin{defn}[{\cite[\S2.5]{ST11}}] A \emph{rigid geometry} is the data of a triple $(G,\M,\M^c)$ consisting of a supermanifold~$\M$ with a left action by a super Lie group~$G$ and a codimension~1 subsupermanifold~$\M^c\subset \M$. 
\end{defn}

\begin{defn}[{\cite[Definitions 2.26, 2.28, 2.33, and 4.4]{ST11}}]
An \emph{$S$-family of $(G,\M)$-supermanifolds} is a locally trivial fiber bundle $Y\to S$ and an open cover~$\{U_i\}$ of~$Y$ with isomorphisms to open sub supermanifolds $f_i\colon U_i\xrightarrow{\sim}  V_i\subset S\times \M$ and transition data $g_{ij}\colon V_i\bigcap V_j\to G$ compatible with the $f_i$ via the $G$-action on~$\M$ and satisfying a cocycle condition. 

An \emph{isometry} between $(G,\M)$-supermanifolds is a map $\varphi\colon Y\to Y'$ over $S$ that is locally given by the $G$-action on~$\M$ relative to the open covers~$\{U_i\}$ of~$Y$ and~$\{U_i'\}$ of~$Y'$, and satisfying a compatibility condition with the transition data $g_{ij}$ and $g_{ij}'$.

A \emph{$(G,\M)$-pair} is a pair of bundles of supermanifolds $(Y^c,Y)$ over $S$ with an inclusion $Y^c\subset Y$ over $S$, where $Y$ has a $(G,\M)$-structure with respect to a chart $\{U_i\}$, and $f_i({U_i\bigcap Y^c})\subset S\times \M^c\subset S\times \M$.  \emph{Isometries} between $(G,\M)$-pairs are isometries $\varphi\colon Y\to Y'$ such that $\varphi(Y^c)=Y'^c\subset Y$. 
\end{defn}

Pulling back bundles along base changes $S'\to S$, the category of $(G,\M)$-pairs determines a stack on the site of supermanifolds; see~\cite[Definition 2.33]{ST11} and also~\cite[Remark A.23]{Powerops}. Any $S$-family of $(G,\M)$-supermanifolds determines a $(G,\M)$-pair with $Y^c=\emptyset$, so that $(G,\M)$-supermanifolds form a full substack of the stack of $(G,\M)$-pairs.

The rigid geometries relevant to the Stolz--Teichner program are defined as follows; see also~\cite[Lecture~3]{Freed5} for the (non-Wick rotated) super Poincar\'e analog.

\begin{defn}[{\cite[\S4.2]{ST11}}]\label{eq:ordEuclidean}
A \emph{super Euclidean geometry} is 
$$
(G,\M,\M^c)=(\E^{d|\delta}\rtimes \Spin(d),\R^{d|\delta},\R^{(d-1)|\delta}),\qquad d,\delta\in \N
$$
where $\M^c=\R^{(d-1)|\delta}\subset \R^{d|\delta}=\M$ is the standard inclusion and the super Lie group $\E^{d|\delta}\rtimes \Spin(d)$ is constructed out of the data: 
\begin{enumerate}
\item[i.] $\Delta$ a complex spinor representation of $\Spin(d)$ with dimension ${\rm dim}_\C(\Delta)=\delta$, and 
\item[ii.] $\Gamma\colon \Delta\otimes_\C \Delta\to \R^d\otimes_\R \C\simeq \C^d$ a $\Spin(d)$-equivariant, nodegenerate symmetric pairing for the action on $\R^d$ through the double cover $\Spin(d)\to \SO(d)$. 
\end{enumerate}
From these data, define a super Lie group $\E^{d|\delta}$ with underlying supermanifold $\R^d\times \Pi \Delta\simeq \R^{d|\delta}$ and multiplication 
\beq\label{eq:Euctransgrp}
(v,\theta)\cdot (w,\eta)=(v+w+\Gamma(\theta,\eta),\theta+\eta),\qquad v,w\in \R^d(S), \ \theta, \eta\in \Pi \Delta(S),
\eeq
using the functor of points description of $\R^{d|\delta}$, see~\eqref{eq:RnmSpts}. 
The semidirect product $\E^{d|\delta}\rtimes \Spin(d)$ acts on~$\R^{d|\delta}$ (on the left) thereby determining a rigid geometry. 
\end{defn}

\begin{rmk}
Throughout, $\E^{d|\delta}$ will denote the super Euclidean translation group with multiplication~\eqref{eq:Euctransgrp} to distinguish from the supermanifold $\R^{d|\delta}$ on which it acts. We caution that the notation~$\E^{d|\delta}$ conceals the (possibly obstructed) data~$(\Delta,\Gamma)$. Indeed, given a super dimension~$d|\delta$, there need not exist an associated super Euclidean geometry. Notably, $d|1$-Euclidean geometries exist precisely in dimensions $d=0,1,2$. Furthermore, if a super Euclidean geometry in a particular super dimension exists, it need not be unique, e.g., there are two possible $2|1$-Euclidean geometries, see Definition~\ref{Defn:Eucgrp} below. 
\end{rmk} 

\begin{ex}[Euclidean geometry without supersymmetry]\label{ex:2Euc}
Taking $\delta=0$ in Definition~\ref{eq:ordEuclidean}, the pairing~$\Gamma$ is no additional data. This leads to a Euclidean geometry $(G,\M,\M^c)\simeq (\E^d\rtimes \Spin(d),\R^d,\R^{d-1})$ for which Euclidean manifolds are $d$-manifolds with a flat metric and spin structure. This enhances Thurston's definition of Euclidean manifold by adding a spin structure~\cite[Example~3.3.2 and \S4.2]{Thurston}. Just as in Thurston's case, examples of Euclidean manifolds come from quotients of $\R^d$ by a discrete subgroup of $\E^d\rtimes \Spin(d)$ such that the associated quotient map is a covering space, or equivalently, the discrete subgroup of $\E^d\rtimes \Spin(d)$ acts properly discontinuously. A Euclidean pair $Y^c\subset Y$ is a flat $d$-manifold~$Y$ together with a totally geodesic codimension~1 submanifold~$Y^c$. We provide explicit examples of Euclidean pairs for $d=2$ in \S\ref{sec:dumbEuc} below. 
\end{ex}

\begin{defn}
For a fixed super Euclidean geometry and a smooth manifold $M$, define $d|\delta\EE(M)$ as the stack on the site of supermanifolds whose objects are $S$-families of super Euclidean pairs $(Y^c\subset Y)$ together with a map $\Phi\colon Y\to M$. Define morphisms in $d|\delta\EE(M)$ as isometries $\varphi$ of super Euclidean pairs with the property that the triangle commutes:
\beq
&&\begin{tikzpicture}[baseline=(basepoint)];
\node (A) at (0,0) {$Y^c\subset Y$};
\node (B) at (5,0) {$Y'\supset Y'^c$};
\node (C) at (2.5,-1) {$M.$};
\draw[->] (A) to node [above] {$\varphi$} (B);
\draw[->] (A) to node [below] {$\Phi$} (C);
\draw[->] (B) to node [below] {$\Phi'$} (C);
\path (0,-.5) coordinate (basepoint);
\end{tikzpicture}\label{eq:superEucpairstack}
\eeq
We set $d|\delta\EE:=d|\delta\EE(\pt)$, and note that a smooth map $M\to N$ determines a map of stacks $d|\delta\EE(M) \to d|\delta\EE(N)$ over $d|\delta\EE$. 
\end{defn}

\begin{lem}\label{lem:superfication}
Suppose a super Euclidean geometry on $\R^{d|\delta}$ has been specified. Then there is a faithful functor
\beq\label{eq:superfication}
\mathcal{S}\colon d|0\EE(M) \to d|\delta\EE(M). 
\eeq 
In particular, a spin manifold with flat metric determines a super Euclidean manifold. 
\end{lem}

\bp
This is an elaboration on \cite[Equation~4.14]{ST11}, where the functor~\eqref{eq:superfication} is called \emph{superfication}. From Definition~\ref{eq:ordEuclidean}, the inclusion of the reduced manifold $\R^d\subset \R^{d|\delta}$ is equivariant with respect to the inclusion of the reduced Lie group
\beq\label{Eq:reducedEucgrp}
\E^d\rtimes \Spin(d)\subset \E^{d|\delta}\rtimes \Spin(d). 
\eeq 
Given a family $Y^c\subset Y\to S$ of Euclidean pairs, let $\Spin(Y)\to S$ be the principal $\Spin(d)$ bundle of the vertical tangent bundle of $Y\to S$. Then define 
$$
\mathcal{S}(Y^c\subset Y):=(\Pi(\Spin(Y)|_{Y^c} \times_{\Spin(d)} \Delta)\subset \Pi(\Spin(Y)\times_{\Spin(d)} \Delta)),
$$
where $\Pi$ is the parity reversal functor~\eqref{Eq:Pifunctor} that sends a complex vector bundle $V\to N$ over an ordinary manifold to the supermanifold $\Pi V$ that regards the fibers of $V$ as odd. The cover $\{U_i\}$ defining the super Euclidean structure on $Y$ gives a cover $\{\Pi(\Spin(Y)|_{U_i}\times_{\Spin(d)}\Delta) \}$ of $\mathcal{S}(Y)$ with gluing data from $\E^d\rtimes \Spin(d)\subset \E^{d|\delta}\rtimes \Spin(d)$. Hence, one obtains a canonical super Euclidean structure on $\mathcal{S}(Y^c\subset Y)$. Furthermore, maps between Euclidean pairs $(Y^c,Y)\to (Y'^c,Y')$ determine maps of spin bundles, leading to maps between the resulting supermanifolds in the image of $\mathcal{S}$ (using that $\Pi$ is a functor). Working in a chart defining the Euclidean structure, one finds that these maps are indeed super Euclidean isometries. Finally, we take the map $\mathcal{S}(Y)\to M$ given by the composition 
$$
\mathcal{S}(Y)=\Pi(\Spin(Y)\times_{\Spin(d)} \Delta))\to Y\xrightarrow{\Phi} M
$$
where the first arrow is induced by the projection and the second arrow is part of the data of an object in $d|0\EE(M)$. This map is compatible with super Euclidean isometries in the image of~\eqref{Eq:reducedEucgrp}, and so 
this completes the construction of $\mathcal{S}$. This functor is faithful precisely because~\eqref{Eq:reducedEucgrp} is an injection. 
\ep

\begin{lem} \label{lem:reducedEuc} Given a $d|\delta$-Euclidean pair $Y^c\subset Y\to S$, the reduced family $Y^c_{\rm red} \subset Y_{\rm red}\to S_{\rm red}$ has the structure of a $d$-dimensional Euclidean pair from Example~\ref{eq:ordEuclidean}, i.e., the fibers of $Y_{\rm red}\to S_{\rm red}$ are flat Riemannian spin manifolds and $Y^c_{\rm red}$ is  a (fiberwise) totally geodesic codimension~1 submanifold. This determines a map of stacks
$$
d|\delta\EE(M)\to d|0\EE(M)
$$
covering the reduction functor from supermanifolds to manifolds.
\end{lem}
\bp
The cover $\{U_i\}$ defining the super Euclidean structure on $Y$ gives a cover $\{(U_i)_{\rm red}\}$ of $Y_{\rm red}$ with gluing data from $\E^d\rtimes \Spin(d)$. This defines an $S_{\rm red}$-family of Euclidean manifolds in the sense of Example~\ref{ex:2Euc}. Similarly, the reduced map of~\eqref{eq:codim1} is the standard inclusion $\R^{d-1}\subset \R^d$. Hence, the reduced family of a super Euclidean pair is a family of Euclidean spin manifolds and a (fiberwise) totally geodesic codimension~1 submanifold.
\ep

\begin{lem} \label{lem:flip1} 
The stack $d|\delta\EE(M)$ has an involution that to a $d|\delta$-dimensional super Euclidean pair associates an isometry 
\beq\label{eq:flip}
&&\fl\colon (Y^c,Y)\to (Y^c,Y)
\eeq
that restricts to the spin flip automorphism on the reduced Euclidean family. \end{lem}

\bp
This is an elaboration on \cite[\S2.6]{ST11}. Locally on $Y\to S$, the involution corresponds to the automorphism $-1\in \Spin(d)$ acting on $\R^{d|\delta}$ as
\beq\label{eq:fliponmodel}
\fl(v,\theta)=(v,-\theta),\qquad (v,\theta)\in \R^{d|\delta}(S). 
\eeq
The action of $-1$ commutes with the action of $\E^{d|\delta}\rtimes \Spin(d)$ on $\R^{d|\delta}$. This implies that these local involutions glue to determine a global involution on $Y\to S$. Since the inclusion $\R^{(d-1)|\delta}\hookrightarrow \R^{d|\delta}$ is preserved by the action of $-1\in \Spin(d)$, we get an involution~\eqref{eq:flip} for every $d|\delta$-dimensional Euclidean pair. Finally, we observe that $-1\in \Spin(d)<\E^d\rtimes \Spin(d)$ acts by the spin flip on $\R^d$ (this is the trivial action on $\R^d$, but flips the sheets of the spin double cover of the frame bundle), and hence acts by the spin flip on the reduced family. 
\ep

Next we observe that the super Euclidean isometry group carries an $\R_{>0}$-action by homomorphisms determined by
\beq\label{eq:rgfunctor1}
\mu\cdot (v,\theta)=(\mu^2v,\mu\theta),\qquad \mu\in \R_{>0}(S), \ (v,\theta)\in \E^{d|\delta}(S).
\eeq
From~\eqref{eq:Euctransgrp}, this gives a homomorphism $\rg_\mu\colon \E^{d|\delta}\to \E^{d|\delta}$ that commutes with the $\Spin(d)$-action. We further note that~\eqref{eq:rgfunctor1} determines a map $\rg_\mu\colon \R^{d|\delta}\to \R^{d|\delta}$ that is equivariant relative to the homomorphism $\E^{d|\delta}\rtimes \Spin(\delta) \to \E^{d|\delta}\rtimes \Spin(\delta)$ determined by~\eqref{eq:rgfunctor1}. 

\begin{rmk}
The notation in~\eqref{eq:rgfunctor1} comes from the fact that the dilation action induces the renormalization group flow on field theories.
\end{rmk}

\begin{lem} \label{lem:RGgeneral} 
The stack $d|\delta\EE(M)$ has an $\R_{>0}$-family of endofunctors
\beq\label{eq:rgfunctor2}
\RG_\mu\colon d|\delta\EE(M)\to d|\delta\EE(M),\qquad \RG_\mu\circ \RG_{\mu'}\simeq \RG_{\mu\mu'}, \quad \mu,\mu'\in \R_{>0}
\eeq
determined by applying the dilation~\eqref{eq:rgfunctor1} to the charts $U_i\subset S\times \R^{d|\delta}$ of a super Euclidean pair and the dilation action~\eqref{eq:rgfunctor1} to a super Euclidean isometry.
\end{lem}
\bp
The proof is entirely analogous to that of Lemma~\ref{lem:flip1}.
\ep

\begin{ex} For the Euclidean geometry without supersymmetry from Example~\ref{ex:2Euc}, the functor $\RG_\mu$ from Lemma~\ref{lem:flip1} scales the metric on a flat Riemannian manifold by a global factor of~$\mu^2$, where the square root $\mu$ specifies a lift of this action to the spinor bundle. \end{ex}

\subsection{2-dimensional Euclidean geometry (no supersymmetry)} \label{sec:dumbEuc}

To gather intuition, we expand on the~$d=2$ case of Example~\ref{ex:2Euc}. This geometry is modeled on the flat Riemannian spin manifold $\R^2$ with isometry group $\R^2\rtimes \Spin(2)$ acting by translations and rotations. Euclidean pairs are therefore flat 2-manifolds with spin structure and an embedded geodesic. 

Euclidean spin tori give examples of compact 2-dimensional Euclidean manifolds. Ignoring spin structures for the moment, every Euclidean torus has $\R^2$ as its universal cover and so is a quotient of $\R^2$ by the action of a subgroup $\Z^2\hookrightarrow \E^2$ of the Euclidean translation group. Let $\Lat\subset \E^2\times \E^2$ denote the submanifold of positively-oriented generators of such $\Z^2$-subgroups, i.e., the open submanifold
\beq\label{eq:Latdefn}
&&\Lat:=\{\ell_1,\ell_2\in \E^2\simeq \R^2\mid \ell_1,\ell_2\ {\rm are\ a \ positively \ oriented\ basis\ of \ } \R^2\}.
\eeq
The vectors $\ell_1,\ell_2$ generate a (based, oriented) lattice $\{n\ell_1+m\ell_2\in \R^2 \mid (n,m)\in \Z^2\}$, and we obtain a smooth fibration over $\Lat$ whose fibers are the tori $T_{\ell_1,\ell_2}=\R^2/(\ell_1\Z\oplus\ell_2\Z)$. To equip these fibers with spin structures we require the data of a lift 
\beq
&&\begin{tikzpicture}[baseline=(basepoint)];
\node (A) at (0,0) {$\Z^2$};
\node (B) at (4,0) {$\E^2$};
\node (C) at (4,1) {$\E^2\rtimes \Spin(2)$};
\draw[->] (A) to node [below] {$\langle \ell_1,\ell_2\rangle$} (B);
\draw[->,dashed] (A) to (C);
\draw[->] (C) to (B);
\path (0,-.5) coordinate (basepoint);
\end{tikzpicture}\nonumber
\eeq
to an action of $\Z^2$ on $\R^2$ that preserves the spin structure. A short computation (compare Lemma~\ref{lem:4component} below) reveals that there are 4 such possible lifts corresponding to choices of $\pm1\in \Spin(2)$ for each generator of the fundamental group
\beq\label{eq:4choices}
\langle (\ell_1,\pm 1),(\ell_2,\pm 1)\rangle \colon \Z^2 \to \E^2\rtimes \Spin(2). 
\eeq
These 4 choices are commonly referred to as periodic-periodic for (+1,+1), periodic-antiperiodic for $(+1,-1)$, antiperiodic-periodic for $(-1,+1)$, and antiperiodic-antiperiodic for $(-1,-1)$. Consequently, the moduli stack of pointed Euclidean spin tori is 
\beq\label{eq:pointedEuctori}
(\Lat\coprod \Lat\coprod \Lat\coprod \Lat)\sq \Spin(2)\times \SL_2(\Z),
\eeq
where $\SL_2(\Z)$ changes the basis of the lattice, and $\Spin(2)$ acts by rotation of lattices through the homomorphism $\Spin(2)\to \SO(2)$. The periodic-periodic spin structure is odd (i.e., nonbounding), corresponding to the trivial double cover of the frame bundle; this component is preserved by the $\SL_2(\Z)$-action. The remaining three choices are even (or bounding) spin structur; these are permuted by the $\SL_2(\Z)$-action. Hence, the stack~\eqref{eq:pointedEuctori} has two connected components. Removing the basepoint in~\eqref{eq:pointedEuctori} amounts to including the action of a torus on itself by translation, enlarging the automorphism groups in~\eqref{eq:pointedEuctori} by $\E^2/(\ell_1\Z\oplus\ell_2\Z)$. We refer to~\S\ref{sec:back2} for an alternate description of the moduli stack of spin tori~\eqref{eq:pointedEuctori} involving the metaplectic double cover of $\SL_2(\Z)$. 

Next we enhance the stack~\eqref{eq:pointedEuctori} to include maps from Euclidean tori to a smooth manifold $M$. First we observe that for any lattice there is a diffeomorphism,
\beq\label{eq:isostandardnormy}
\R^2/\Z^2\xrightarrow{\sim} \R^2/(\ell_1\Z\oplus \ell_2\Z),\qquad (x,y)\mapsto(\ell_1 x,\ell_2y)
\eeq
with the standard torus (i.e., for the unit square lattice $\Z^2\subset \E^2$). Hence, for the periodic-periodic spin structure we get a map of stacks
$$
\mathcal{L}^2(M)\sq (\E^2\rtimes \Spin(2)\times \SL_2(\Z)) \to 2\EE(M),\quad 
(\phi,\ell_1,\ell_2)\mapsto \R^2/(\ell_1\Z\oplus \ell_2\Z)\simeq \R^2/\Z^2\xrightarrow{\phi}M,
$$
where the mapping space
\beq\label{Eq:defnofsuperdoubleloops}
\mathcal{L}^2(M):=\Map(\R^2/\Z^2,M)\times \Lat
\eeq
is regarded as a sheaf (of sets). For the even spin structures, there are analogous maps where one either restricts to the subgroup of $\SL_2(\Z)$ preserving a chosen spin structure, or considers the coproduct over the periodic-antiperiodic, antiperiodic-periodic, and antiperiodic-antiperiodic spin structures.

\begin{rmk} \label{rmk:chiralD}
The isomorphism $\R^2\simeq \C$ allows one to view a Euclidean 2-torus as a complex manifold. The orientation condition for the pair $\ell_1,\ell_2\in \R^2$ in~\eqref{eq:Latdefn} translates to the ratio of generators having positive imaginary part
\beq\label{Eq:latcondition}
\Lat\simeq \{z_1,z_2\in \C^\times \mid \im(z_2/z_1)>0\}. 
\eeq
The isomorphism $\Spin(2)\times \SL_2(\Z)\simeq U(1)\times \SL_2(\Z)$ further identifies isometries of pointed Euclidean tori with isomorphisms between complex manifolds. In the complex setting, a spin structure is the same data as a square root of the canonical bundle \cite{AtiyahSpin}. This canonical bundle is trivializable, and so describing a square root in terms of monodromy in $\{\pm 1\}$ recovers the 4 choices of spin structure~\eqref{eq:4choices}.
\end{rmk}

Next we describe $2$-dimensional Euclidean pairs given by an infinite cylinder with an embedded (closed) geodesic. These examples arise by taking quotients by $\Z$-actions on $\R^2$ compatible with the standard inclusion $\R\subset \R^2$. The stabilizer of the subspace~$\R\subset \R^2$ is 
\beq\label{eq:stabilizer}
\E^1\rtimes \Z/4\subset \E^2\rtimes \Spin(2),
\eeq 
and the further subgroup $\E^1\times \Z/2\subset \E^1\rtimes \Z/4$ preserves the orientation on~$\R\subset \R^2$. For~$(\ell,\pm 1)\in \R_{>0}\times \{\pm 1\}\subset \E^1\times \{\pm 1\}$, let 
\beq\label{eq:dumbEucpair1}
&&\rS_\ell^\p =\R/(\ell,+ 1)\Z,\quad \rS_\ell^\a =\R/(\ell,-1)\Z,\quad \rC_\ell^\p=\R^2/(\ell,+ 1)\Z,\quad \rC_\ell^\a=\R^2/(\ell,-1)\Z,\quad 
\eeq
denote the circles of circumference $\ell$ and the infinite cylinders of circumference $\ell$, with periodic or antiperiodic spin structure as indicated. For $\tau\in \E^2$, we obtain inclusions
\beq\label{eq:dumbEucpair2}
i_\tau\colon \rS_\ell^{\p/\a}=\R/(\ell,\pm1)\Z\xhookrightarrow{i_0} \R^2/(\ell\Z,\pm 1) \xrightarrow{\Tran_\tau} \R^2/(\ell\Z,\pm1)=\rC_\ell^{\p/\a}
\eeq
where $i_0$ is determined by the standard inclusion $\R\subset \R^2$ and the map $\Tran_\tau$ descends from the (left) translation action of $\tau\in \E^2$ on $\R^2$. We note that $i_\tau=i_{\tau'}$ if $\tau$ and $\tau'$ have the same image in the quotient~$\E^2/\ell\Z$. Hence, the data~\eqref{eq:dumbEucpair2} determine a family of super Euclidean pairs parameterized by $(\pm1,\ell,\tau)\in \Z/2\times (\R_{>0}\times \E^2)/\Z$ for the quotient by the $\Z$-action
\beq\label{eq:normyZqupt}
n\cdot (\ell,\tau)=(\ell,n\ell+\tau),\qquad n\in \Z, \ \ell\in \R_{>0}, \ \tau\in \E^2. 
\eeq

Analogously to~\eqref{eq:isostandardnormy} there are diffeomorphisms (that are not Euclidean isometries)
\beq\label{eq:standardEucnormy}
f_\ell\colon \rS^{\p/\a}_\ell\xrightarrow{\sim} \R/\Z,\qquad F_\ell\colon \rC^{\p/\a}_\ell\xrightarrow{\sim} \R^2/\Z
\eeq
between the circle and cylinder of circumference $\ell$ with the (standard) circle and cylinder of circumference~1. Hence for each spin structure on the circle we get a map of stacks
$$
\Map(\R^2/\Z,M)\times (\R_{>0}\times \E^2)/\ell\Z\to 2\EE(M),\qquad (\phi,\ell,\tau)\mapsto (\rS^{\p/\a}_\ell\xhookrightarrow{i_\tau} \rC^{\p/\a}_\ell\simeq \R^2/\Z\xrightarrow{\phi}M),
$$
sending an $S$-point of the source to an $S$-family of Euclidean pairs over $M$. The subgroup~\eqref{eq:stabilizer} acts through isometries on the resulting family of Euclidean pairs. 

\subsection{$2|1$-Euclidean geometry} \label{Sec:Eucgeodef}
Next we turn attention to super Euclidean geometry in dimension~$2|1$. Unpacking the data from Definition~\ref{eq:ordEuclidean} in this case, $\Delta=\C$ and the pairing~$\Gamma$ is a $\Spin(2)\simeq U(1)$-equivariant map
\beq\label{eq:GammaEuc}
\Gamma\colon \Delta\otimes_\C \Delta\simeq \C\hookrightarrow \R^2\otimes_\R \C\simeq \C_{2}\oplus \C_{-2}.
\eeq
In the target above, we have decomposed $\R^2\otimes \C$ into irreducible representations of $U(1)$ of weight~$2$ and~$-2$. Therefore, we have two possible choices of $U(1)$-action on $\Delta$ that makes~\eqref{eq:GammaEuc} an equivariant map corresponding to $\Delta=\C_{\pm 1}$, the weight $+1$ and $-1$ representations. With this choice fixed, $\Gamma$ is the inclusion of $\C_{\pm 1}\otimes \C_{\pm 1}\simeq \C_{\pm 2}$ as a summand. 

\begin{defn}[$2|1$-dimensional Euclidean geometries]\label{Defn:Eucgrp}\label{defn:sEuc}
Let~$\E^{2|(0,1)}$ denote the super Lie group associated to the pairing~\eqref{eq:GammaEuc} with the choice $\Delta=\C_{-1}$. Then $\E^{2|(0,1)}$ has underlying supermanifold $\R^{2|1}$ and multiplication (using the functor of points description~\eqref{eq:r212} of~$\R^{2|1}$)
\beq
&&(z,\bar z,\theta)\cdot (z',\bar z',\theta')=(z+z',\bar z+\bar z'+\theta\theta',\theta+\theta'), \quad (z,\bar z,\theta),(z',\bar z',\theta')\in \R^{2|1}(S).\label{eq:E21mult}
\eeq
Similarly, let $\E^{2|(1,0)}$ denote the super Lie group associated to the pairing~\eqref{eq:GammaEuc} where we choose $\Delta=\C_1$. This results in the multiplication
\beq
&&(z,\bar z,\theta)\cdot (z',\bar z',\theta')=(z+z'+\theta\theta',\bar z+\bar z',\theta+\theta'), \quad (z,\bar z,\theta),(z',\bar z',\theta')\in \R^{2|1}(S).\label{eq:E21mult2}
\eeq
Consider the semidirect products $\E^{2|(0,1)}\rtimes \Spin(2)$ and $\E^{2|(1,0)}\rtimes \Spin(2)$ defined by the action (using the functor of description~\eqref{Eq:spin2} of~$\Spin(2)$)
$$
(u ,\bar u)\cdot (z,\bar z,\theta)=(u^2 z,\bar u^2 \bar z,\bar u \theta),\qquad (u,\bar u)\in \Spin(2)(S).
$$
The above defines a pair of $2|1$-dimensional super Euclidean geometries 
$$
(G,\M,\M^c)=\left\{\begin{array}{l} (\E^{2|(0,1)}\rtimes \Spin(2),\R^{2|1},\R^{1|1})\\ (\E^{2|(1,0)}\rtimes \Spin(2),\R^{2|1},\R^{1|1})\end{array}\right.
$$ 
for the codimension~1 sub supermanifold given by the map
\beq
&&\R^{1|1}\hookrightarrow \R^{2|1},\quad (t,\theta)\mapsto (t,t,\theta)=(z,\bar z,\theta), \quad (t,\theta)\in \R^{1|1}(S), \  (z,\bar z,\theta)\in \R^{2|1}(S). \label{eq:codim1}
\eeq
\end{defn} 

The two choices in $2|1$-Euclidean geometry correspond to the two possible chiralities for supersymmetry in dimension~2, often referred to as $\mathcal{N}=(1,0)$ or $\mathcal{N}=(0,1)$ supersymmetry. Our primary focus below is on the $2|1$-Euclidean geometry determined by~\eqref{eq:E21mult}, and (following~\cite[\S4.3]{ST04}) eventually we will adopt the notation $\E^{2|1}:=\E^{2|(0,1)}$ when this choice is clear. However, for the moment we explain some interconnections between the pair of super Euclidean geometries in Definition~\ref{defn:sEuc} related to complex conjugation and orientation reversal of spin manifolds; ultimately, this is the structure required to define reflection positivity for $2|1$-Euclidean field theories in the sense of~\S\ref{sec:backFT}. 

\begin{lem} \label{lem:or}
The map of supermanifolds
\beq\label{eq:ordef}
\orr\colon \R^{2|1}\to \R^{2|1},\qquad \orr(z,\bar z,\theta)=(\bar z, z,\theta)
\eeq
is equivariant for the action of super Euclidean groups relative to the isomorphism
\beq\label{eq:orstructure}
\orr\times \inv \colon \E^{2|(0,1)}\rtimes \Spin(2) \xrightarrow{\sim}  \E^{2|(1,0)}\rtimes \Spin(2)\label{eq:chiralityiso}
\eeq
determined by the product of~\eqref{eq:ordef} and the inversion isomorphism (using $\Spin(2)$ is abelian) 
\beq\label{Eq:inv}
\inv\colon \Spin(2)\to \Spin(2),\qquad (u,\bar u)\mapsto (\bar u,u). 
\eeq
Furthermore,~\eqref{eq:ordef} reverses the orientation on the reduced manifold $(\R^{2|1})_{\rm red}=\R^2$, and restricts to the identity on the subspace $\R^{1|1}\subset \R^{2|1}$. 
\end{lem}
\bp
Directly from the definitions we see that $\orr$ determines an isomorphism of super Lie groups $\E^{2|(0,1)}\simeq \E^{2|(1,0)}$. Since precomposition with inversion exchanges the representations $\C_{+2}$ and $\C_{-2}$ of $\Spin(2)\simeq U(1)$, one obtains the claimed extension to a map of semidirect products~\eqref{eq:orstructure}. In terms of real coordinates~\eqref{eq:r211} on $\R^2\subset \R^{2|1}$, the map~\eqref{eq:ordef} is $(x,y)\mapsto (x,-y)$, and so reverses the orientation on~$\R^2$. Finally, the definition~\eqref{eq:codim1} of the inclusion $\R^{1|1}\hookrightarrow \R^{2|1}$ demonstrates that~$\orr|_{\R^{1|1}}=\id$. 
\ep

We refer to Definition~\ref{defn:realsmfld} for conjugation and real structures on supermanifolds. 

\begin{lem}\label{lem:chiralityswaps}
There is an isomorphism of supermanifolds 
\beq\label{eq:real21}
\c\colon \R^{2|1}\xrightarrow{\sim}  \overline{\R}^{2|1}
\eeq
that is equivariant for actions relative to an isomorphism
\beq\label{eq:conjiso}
\E^{2|(0,1)}\rtimes \Spin(2) \xrightarrow{\sim}  \overline{\E^{2|(1,0)}\rtimes \Spin(2)},
\eeq
where both~\eqref{eq:real21} and~\eqref{eq:conjiso} are determined by the linear map $(\overline{\phantom{A}})\colon \C_{-2}\to \overline{\C}_{+2}$ gotten from complex conjugation. Furthermore,~\eqref{eq:real21} restricts to the identity on the reduced manifold~$\R^2\subset \R^{2|1}$, and restricts to a map $\R^{1|1}\to \overline{\R}^{1|1}\subset \overline{\R}^{2|1}$ between subspaces. 
\end{lem}
\bp
Complex conjugation determines an isomorphism between the $\Spin(2)$-equivariant vector bundles $\R^2\times \C_{-2}\to \R^2\times \overline{\C}_{+2}$ over $\R^2$. Applying the functor $\Pi$ from~\eqref{Eq:Pifunctor}, one obtains the map~\eqref{eq:real21} of supermanifolds whose induced map on reduced manifolds is the identity. From the definition of the super Euclidean group (using that the options in~\eqref{eq:GammaEuc} are conjguate), we obtain a homomorphism $\E^{2|(0,1)}\xrightarrow{\sim}  \overline{\E^{2|(1,0)}}$. Finally, $\Spin(2)$-equivariance of the map of vector bundles gives the extension~\eqref{eq:conjiso} to semidirect products, where we use the canonical real structure on $\Spin(2)$. 
\ep

\begin{rmk} A \emph{real structure} for a Euclidean geometry is an isomorphism $\R^{d|\delta}\to \overline{\R}{}^{d|\delta}$ that is equivariant for an isomorphism of super Lie groups $\E^{d|\delta}\rtimes \Spin(\delta)\to \overline{\E^{d|\delta}\rtimes \Spin(\delta)}$.  Such a real structure exists for the $1|1$-Euclidean geometry, see \cite[Example 6.20]{HST}. Lemma~\ref{lem:chiralityswaps} is not a real structure for the $2|1$-Euclidean geometry, as the isomorphism of super Lie groups~\eqref{eq:conjiso} changes the chirality of the supersymmetry. Because of this chiral supersymmetry, the $2|1$-Euclidean geometry does not admit a real structure. 
\end{rmk}

\begin{prop}\label{prop:Or} The maps in Lemma~\ref{lem:or} determine a morphism of stacks
\beq\label{eq:bigor}
\Or\colon 2|(0,1)\EE(M)\to 2|(1,0)\EE(M),
\eeq
that on reduced families of 2-dimensional Euclidean manifolds reverses the fiberwise orientations. 
The maps in Lemma~\ref{lem:chiralityswaps} determine a morphism of stacks
\beq\label{eq:bigc}
{\sf C}\colon 2|(1,0)\EE(M)\to 2|(0,1)\EE(M),
\eeq
covering the complex conjugation functor on supermanifolds, i.e., an $S$-family is sent to an $\overline{S}$-family. 
\end{prop}

\bp
Applying $\orr$ to the charts of a family $Y\to S$ with $(\E^{2|(0,1)\rtimes \Spin(2)},\R^{2|1})$-structure, Lemma~\ref{lem:or} constructs a $(\E^{2|(0,1)\rtimes \Spin(2)},\R^{2|1})$-structure. Furthermore, since $\R^{1|1}\subset \R^{2|1}$ is preserved by $\orr$, this construction sends a $(\E^{2|(0,1)\rtimes \Spin(2)},\R^{2|1})$-pair $(Y^c\subset Y)$ over $S$ to a $(\E^{2|(1,0)\rtimes \Spin(2)},\R^{2|1})$-pair over $S$. These constructions change the geometric structures but not the underlying family of supermanifolds; in particular, a map $\Phi\colon Y\to M$ determines a map for either geometric structure, and hence we have constructed the functor~\eqref{eq:bigor} on objects. The values on morphisms are determined by the homomorphism~\eqref{eq:chiralityiso}, using that~\eqref{eq:ordef} is equviariant with respect to this homomorphism. 

For the second statement, we start with a general observation: for a $(G,\M)$-family $Y\to S$, the conjugate family $\overline{Y}\to \overline{S}$ has a canonical $(\overline{G},\overline{\M})$ structure. Hence, for a family $Y\to S$ with $(\E^{2|(1,0)}\rtimes \Spin(2),\R^{2|1})$-structure, the family $\overline{Y}\to \overline{S}$ automatically has a $(\overline{\E^{2|(1,0)}\rtimes \Spin(2)},\overline{\R^{2|1}})$-structure. Using Lemma~\ref{lem:chiralityswaps} we endow $\overline{Y}\to \overline{S}$ with a $(\E^{2|(0,1)}\rtimes \Spin(2),\R^{2|1})$-structure. As in the previous case, since $\c$ sends the subspace $\R^{1|1}$ to the subspace $\overline{\R}^{1|1}$, this extends to a map on pairs, i.e., $(\overline{Y}^c\subset \overline{Y})$ has the canonical structure of a super Euclidean pair. Finally, we use the canonical real structure on $M$ (as an ordinary manifold) to obtain a map $\overline{Y}\to M$ as
$\overline{Y}\stackrel{\overline{\Phi}}{\to} \overline{M}\simeq M$. 
This defines the functor~\eqref{eq:bigc} on objects. The values on morphisms are determined by the homomorphism~\eqref{eq:conjiso}, using that~\eqref{eq:real21} is equviariant. 
\ep

\begin{defn} \label{defn:dagger}
The \emph{reflection} functor on $2|1$-Euclidean supermanifolds is the composition
\beq\label{eq:daggerdefn}
2|(0,1)\EE(M)\xrightarrow{\Or} 2|(1,0)\EE(M)\xrightarrow{{\sf C}} 2|(0,1)\EE(M)
\eeq
which is a map of stacks covering complex conjugation on the category of supermanifolds. 
\end{defn}

\begin{lem}
The reflection functor generates a $\Z/2$-action on the stack $2|(0,1)\EE(M)$. 
\end{lem}
\bp
From~\eqref{eq:ordef} and~\eqref{eq:real21}, we observe that $\overline{(\c\circ \orr)}\circ (\c\circ \orr)=\id_{\R^{2|1}}$. Hence, the identity map on families of super Euclidean pairs supplies the necessary 2-morphism to define a $\Z/2$-action on the stack, and the claim follows. 
\ep

Finally, we observe some connections between $2|1$-dimensional Euclidean geometry and $1|1$-dimensional Euclidean geometry (e.g., as studied in~\cite[\S6]{HST}). Define $\E^{1|1}$ as the supermanifold $\R^{1|1}$ endowed with the multiplication
$$
(t,\theta)\cdot (s,\eta)=(t+s+i\theta\eta,\theta+\eta),\qquad (t,\theta),(s,\eta)\in \R^{1|1}(S). 
$$

\begin{lem}
There are exact sequences of super Lie groups
\beq\label{eq:itsexact}
\begin{array}{ccccccccc}
1&\to &\E&\to& \E^{2|(0,1)}&\to& \E^{1|1}&\to& 1,\\
&&t&\mapsto & (t,t,0), \ (z,\bar z,\theta) &\mapsto & (2\im(z),\theta)
\end{array}
\eeq
\beq\nonumber
\begin{array}{ccccccccc}
1&\to &\E&\to& \E^{2|(1,0)}&\to& \E^{1|1}&\to& 1,\\
&&t&\mapsto & (t,t,0), \ (z,\bar z,\theta) &\mapsto & (-2\im(z),\theta)
\end{array}
\eeq
where the indicated formulas are understood using the functor of points, and we have used the (slightly abusive) notation
$$
\im(z):=\frac{z-\bar z}{2i}\in C^\infty(S),\qquad (z,\bar z)\in \R^2(S).
$$ 
\end{lem}
\bp
First observe that the maps on $S$-points satisfy the required reality conditions to give the claimed maps between supermanifolds. To see that these are furthermore maps of super Lie groups comes down to the equality
$$
\left(\frac{z-\bar z}{i}+\frac{w-\bar w}{i}+i\theta\eta,\theta+\eta\right)=\left(\frac{z+w-(\bar z+\bar w+\theta\eta)}{i},\theta+\eta\right)
$$
for the sequence involving $\E^{2|(0,1)}$, with a similar equality for the case involving $\E^{2|(1,0)}$. Exactness of the sequences is clear: $\E$ is a normal subgroup (in fact, a central subgroup) with quotient $\E^{1|1}$.
\ep

\begin{lem}\label{lem:sqrtflanddagger}
The inclusion~\eqref{eq:codim1} is $\Z/4$-equivariant for the action by the 4th roots of unity $\Z/4\subset U(1)\simeq \Spin(2)$. Under the reflection functor, the generator $i\in U(1)$ of this $\Z/4$-action is sent to $-i\in \overline{U(1)}$. 
\end{lem}
\bp
This follows from the fact that the map~\eqref{eq:orstructure} is inversion on $\Spin(2)\simeq U(1)$. 
\ep

\begin{rmk}\label{rmk:pin} The $\Z/4$-action on $\R^{1|1}$ in the previous lemma recovers the ``pin generator" from \cite[Example 6.16]{HST}, which is an important structure in unoriented $1|1$-Euclidean field theories.\end{rmk}

\subsection{Supertori as examples of $2|1$-Euclidean supermanifolds}\label{sec:toriexamples}

\begin{defn}
Let $\Hom(\Z^2,\E^{2|1}\rtimes \Spin(2))$ denote the sheaf of sets on the site of supermanifolds whose $S$-points are families of homomorphisms
$$
\Lambda\colon S\times \Z^2\to (\E^{2|1}\rtimes \Spin(2)).
$$
Homomorphisms pull back along base changes, so that $\Hom(\Z^2,\E^{2|1}\rtimes \Spin(2))$ is a presheaf, and families of maps satisfy the gluing property equivalent to the sheaf property. Define a subsheaf 
\beq\label{eq:subspaceLat}
\Hom(\Z^2,\E^{2|1}\rtimes \Spin(2))_{\Lat}\subset \Hom(\Z^2,\E^{2|1}\rtimes \Spin(2))
\eeq
whose $S$-points are families of injective homomorphisms that on reduced manifolds determine an oriented basis for $\R^2$ pointwise in $S_{\rm red}$, see~\eqref{eq:Latdefn}. The natural $\SL_2(\Z)$-action on $\Hom(\Z^2,\E^{2|1}\rtimes \Spin(2))$ by precomposition restricts to an action on the subsheaf~\eqref{eq:subspaceLat}.
\end{defn} 

\begin{rmk}
The pointwise basis condition is equivalent to the complex-valued functions
$$
z_1\colon S_{\rm red}\times (1,0)\to \R^2\rtimes \Spin(2) \to \R^2\simeq \C \quad {\rm and}\quad  z_2\colon S_{\rm red}\times (1,0)\to \R^2\rtimes \Spin(2) \to \R^2\simeq \C
$$
having the property~\eqref{Eq:latcondition}.
\end{rmk}

The condition defining the subsheaf~\eqref{eq:subspaceLat} guarantees that the $\Z^2$-action on $S\times \R^{2|1}$ determined by an $S$-point $\Lambda$ of~\eqref{eq:subspaceLat} is free, and hence the quotient is a supermanifold,
\beq\label{eq:bigfamtori}
T_\Lambda:=(S\times \R^{2|1})/\Lambda.
\eeq
Furthermore, the cover $S\times \R^{2|1}\to T_\Lambda$ determines a $2|1$-Euclidean structure on $T_\Lambda$. Hence,~\eqref{eq:bigfamtori} determines a family of $2|1$-Euclidean supermanifolds parameterized by the subsheaf~\eqref{eq:subspaceLat}.

\begin{lem}\label{lem:4component}
The sheaf~\eqref{eq:subspaceLat} is a coproduct of four components
\beq\label{eq:fourcomp}
\Hom(\Z^2,\E^{2|1}\rtimes \Spin(2))_{\Lat} \simeq \sLat \coprod \sLat_{\p\a}\coprod \sLat_{\a\p}\coprod \sLat_{\a\a}
\eeq
where each component is characterized by
$$
C^\infty(\sLat)\simeq C^\infty(\Lat)[\lambda_1,\lambda_2]/(\lambda_1\lambda_2),\quad C^\infty(\sLat_{\p\a})\simeq C^\infty(\Lat)[\lambda_2], $$$$ C^\infty(\sLat_{\p\a})\simeq C^\infty(\Lat)[\lambda_1],\quad C^\infty(\sLat)\simeq C^\infty(\Lat_{\a\a})[\lambda_1,\lambda_2]/(\lambda_1-\lambda_2),
$$
for odd functions $\lambda_i$ to be described explicitly in the proof. 
\end{lem}
\bp
To start, we establish notation that repackages the information of $\Lambda$ an $S$-point of the subsheaf~\eqref{eq:subspaceLat}. Let
\beq	
&&(\ell_1,\bar\ell_1,\lambda_1,u_1,\bar u_1)\colon S\times (1,0)\to \E^{2|1}\rtimes U(1),\nonumber\\
&&(\ell_2,\bar\ell_2,\lambda_2,u_2,\bar u_2)\colon S\times (0,1)\to \E^{2|1}\rtimes U(1)\nonumber
\eeq
denote the restrictions of $\Lambda$ to the generators of $\Z^2$. For these $S$-points to determine an $S$-family of homomorphisms, one has the commutativity condition
\beq\label{eq:commutativity}
&&(\ell_1,\bar\ell_1,\lambda_1,u_1,\bar u_1)\cdot (\ell_2,\bar\ell_2,\lambda_2,u_2,\bar u_2)=(\ell_2,\bar\ell_2,\lambda_2,u_2,\bar u_2)\cdot(\ell_1,\bar\ell_1,\lambda_1,u_1,\bar u_1). 
\eeq
Together with the condition~\eqref{Eq:latcondition}, one finds that $u_1^2=u_2^2=1$. Hence, we find four components corresponding to the four possible solutions
\beq
(u_1,u_2)=(1,1): \ \sLat, && (u_1,u_2)=(1,-1): \ \sLat_{\p\a}, \nonumber\\ 
(u_1,u_2)=(-1,1): \ \sLat_{\a\p}, && (u_1,u_2)=(-1,-1): \ \sLat_{\a\a}, \nonumber
\eeq
proving the first statement in the lemma. For the second statement, one views each component above as determining a map 
$$
(\ell_1,\bar\ell_1,\lambda_1,\ell_2,\bar\ell_2,\lambda_2)\colon S\to \R^{2|1}\times \R^{2|1}. 
$$ 
The condition~\eqref{Eq:latcondition} requires that on reduced manifolds the above determines an $S_{\rm red}$-point of $\Lat\subset \C\times \C\simeq \R^2\times \R^2\subset \R^{2|1}\times \R^{2|1}$. The commutativity condition~\eqref{eq:commutativity} with the specified values of $u_1$ and $u_2$ imposes the further conditions
\beq
(u_1,u_2)=(1,1): \lambda_1\lambda_2=\lambda_2\lambda_1 && (u_1,u_2)=(1,-1): \ \lambda_1=0, \nonumber\\ 
(u_1,u_2)=(-1,1): \ \lambda_2=0, && (u_1,u_2)=(-1,-1): \ \lambda_1=\lambda_2. \nonumber
\eeq
Together, this gives the claimed description of functions. 
\ep

\begin{rmk}
Lemma~\ref{lem:4component} shows that $\sLat_{\p\a}, \sLat_{\a\p},$ and $\sLat_{\a\a}$ are representable subsheaves of $\R^{2|1}\times \R^{2|1}$ with reduced manifolds $\Lat$. On the other hand, $\sLat$ is not representable as its functions are not the functions on any supermanifold. This makes the geometry of the periodic-periodic supertori more subtle. 
\end{rmk}

Next we enhance the above to produce super Euclidean manifolds over a smooth manifold~$M$. We will focus on the periodic-periodic supertori over $M$, though the other cases are analogous (and ultimately a bit easier). First we observe that for any periodic-periodic lattice there is an isomorphism of supermanifolds over $S$~\cite[Lemma~3.9]{DBEChern},
$$
S\times (\R^{2|1}/\Z^2)\xrightarrow{\sim} T_\Lambda,\quad (x,y,\theta)\mapsto (\ell_1x+\ell_2y,x(\bar\ell_1+\lambda_1\theta)+y(\bar\ell_2+\lambda_2\theta),\theta+x\lambda_1+y\lambda_2)
$$
with the constant $S$-family with fiber the standard supertorus~$\R^{2|1}/\Z^2$, determined by the canonical inclusions $\Z^2\subset \R^2\subset \R^{2|1}$. This isomorphism between supertori over $S$ is determined by the above $\Z^2$-equivariant map $S\times \R^{2|1}\to S\times \R^{2|1}$, and is typically not an isomorphism of Euclidean supertori over~$S$. 

\begin{defn} \label{defn:doubleloop}
Define the \emph{super Euclidean double loop space} of $M$ as the sheaf on supermanifolds
\beq\label{eq:doubleloop}
\mathcal{L}^{2|1}(M):=\Map(\R^{2|1}/\Z^2,M)\times s\Lat.
\eeq
 Identifying an $S$-point $(\Lambda,\phi)\in (\mathcal{L}^{2|1}(M))(S)$ with 
$$
\Phi\colon T_\Lambda\simeq S\times \R^{2|1}/\Z^2\xrightarrow{\phi} M
$$
gives a functor from $\mathcal{L}^{2|1}(M)\to 2|1\EE(M)$ to the stack of $2|1$-Euclidean supermanifolds over~$M$. 
\end{defn}

There is an action of $\E^{2|1}\rtimes \Spin(2)\times \SL_2(\Z)$ on $\mathcal{L}^{2|1}(M)$ where $\E^{2|1}\rtimes \Spin(2)$ acts by isometries on supertori, and $\SL_2(\Z)$ acts on $s\Lat$ by changing the basis of a lattice. This action is computed explicitly in~\cite[\S3.3]{DBEChern}. By construction, there is a map of stacks~\cite[Lemma~3.15]{DBEChern}
$$
\mathcal{L}^{2|1}(M)\sq (\E^{2|1}\rtimes \Spin(2)\times \SL_2(\Z))\to 2|1\EE(M),\quad (\phi,\Lambda)\mapsto T_\Lambda \simeq S\times \R^{2|1}/\Z^2\xrightarrow{\phi}M.
$$

\subsection{Examples of $2|1$-Euclidean pairs}
The most basic example of a $2|1$-Euclidean pair over $S$ comes from the standard inclusion~\eqref{eq:codim1}. Further examples arise from $\Z$-quotients of this basic one. We start by showing that the moduli space of $\Z$-actions involved is the same as in the case~\eqref{eq:dumbEucpair2} without supersymmetry.

\begin{lem} \label{lem:stabilizersub}
Given an $S$-point $(\ell,\bar\ell,\lambda,u,\bar u)\in (\E^{2|1}\rtimes \Spin(2))(S)$, consider the left $\Z$-action generated by
\beq\label{Eq:action}
&&S\times \R^{2|1}\xrightarrow{\Tran_{\ell,\bar\ell,\lambda,u,\bar u}} S\times \R^{2|1},\qquad (z,\bar z,\theta)\mapsto (u^2(\ell+z),\bar u^2(\bar \ell+\bar z+\lambda\theta),\bar u(\lambda+\theta)).
\eeq
This $2|1$-Euclidean isometry preserves the subspace $S\times \R^{1|1}\subset S\times \R^{2|1}$ and is free if and only if the generator lies in the subgroup 
\beq\label{eq:Zaction}
\E^\times \times \Z/2 \subset \E^{2|1}\rtimes \Spin(2).
\eeq
Finally, the $\Z$-action generated by $(\ell,\pm 1)\in (\R^\times \times \Z/2)(S)$ is the same as the action generated by~$(-\ell,\pm 1)$. 
\end{lem}

\bp
The statement first requires the map~\eqref{Eq:action} to restrict along $S\times\R^{1|1}\subset S\times \R^{2|1}$, i.e., a dashed arrow making the diagram commute
\beq
&&\begin{tikzpicture}[baseline=(basepoint)];
\node (A) at (0,0) {$S\times \R^{1|1}$};
\node (B) at (4,0) {$S\times \R^{2|1}$};
\node (C) at (0,-1.5) {$S\times \R^{1|1}$};
\node (D) at (4,-1.5) {$S\times \R^{2|1}.$};
\draw[->,right hook-latex] (A) to (B);
\draw[->,dashed] (A) to (C);
\draw[->,right hook-latex] (C) to (D);
\draw[->] (B) to node [right] {$\Tran_{\ell,\bar\ell,\lambda,u,\bar u}$} (D);
\path (0,-.75) coordinate (basepoint);
\end{tikzpicture}\nonumber
\eeq
Using the formula~\eqref{eq:codim1} the horizontal inclusions, we require $(\ell,\bar\ell,\lambda,u,\bar u)$ to lie in the subgroup $\E^1\times \Z/4\subset \E^{2|1}\rtimes \Spin(2)$ where $\Z/4\subset U(1)\simeq \Spin(2)$ includes as the 4th roots of unity. Since the action is through the ordinary Lie group $\E\rtimes \Z/4$, freeness of the action can be checked on reduced manifolds where the statement is clear. 
\ep

Given $(\tau,\bar\tau,\eta)\in \E^{2|1}(S)$, define the inclusion $i_{\tau,\bar\tau,\eta}$ as the composition
\beq\label{Eq:tauaction}
i_{\tau,\bar\tau,\eta} \colon S\times \R^{1|1} \xhookrightarrow{i_0}  S\times \R^{2|1} \xrightarrow{\Tran_{\tau,\bar\tau,\eta}} S\times \R^{2|1},
\eeq
where $i_0\colon S\times \R^{1|1}\hookrightarrow S\times \R^{2|1}$ is the  standard inclusion, and $\Tran_{\tau,\bar\tau,\eta}$ is the left translation action by $(\tau,\bar\tau,\eta)\in \E^{2|1}(S)$.

\begin{lem} \label{lem:preserve}
The map~\eqref{Eq:tauaction} descends to the $\Z$-quotient generated by $(\ell,+1)\in (\R_{>0}\times \Z/2)(S)\subset (\E^{2|1}\rtimes \Spin(2))(S)$ for all values of $(\tau,\bar\tau,\eta)$. The map~\eqref{Eq:tauaction} descends to the $\Z$-quotient generated by $(\ell,-1)\in (\R_{>0}\times \Z/2)(S)$ if and only if~$\eta=0$. 
\end{lem}

\bp
The statement amounts to checking commutativity of the diagram
\beq
&&\begin{tikzpicture}[baseline=(basepoint)];
\node (A) at (0,0) {$S\times \R^{1|1}$};
\node (B) at (4,0) {$S\times \R^{2|1}$};
\node (C) at (0,-1.5) {$S\times \R^{1|1}$};
\node (D) at (4,-1.5) {$S\times \R^{2|1}.$};
\draw[->,right hook-latex] (A) to node [above] {$i_{\tau,\bar\tau,\eta}$} (B);
\draw[->] (A) to node [left] {$(\ell,\pm 1)\cdot$} (C);
\draw[->,right hook-latex] (C) to node [below] {$i_{\tau,\bar\tau,\eta}$}  (D);
\draw[->] (B) to node [right] {$(\ell,\pm 1)\cdot$} (D);
\path (0,-.75) coordinate (basepoint);
\end{tikzpicture}\nonumber
\eeq
For an $S$-point $(t,\theta)\in \R^{1|1}(S)$, using~\eqref{eq:codim1} this is equivalent to
$$
(t+\tau+\ell,t+\bar\tau+\eta\theta+\ell,\pm\eta\pm\theta)=(t+\ell+\tau,t+\ell+\bar\tau\pm\eta\theta,\eta\pm\theta). 
$$
The above requires $\pm\eta\pm\theta=\eta\pm\theta$, and the result follows. 
\ep

Following~\eqref{eq:dumbEucpair1}, for $(\ell,\pm 1)\in (\R_{>0}\times \Z/2)(S)$, define families of \emph{supercircles} and (infinite) \emph{supercylinders} of circumference $\ell$
\beq\label{eq:Eucpair1}
&&\rS_\ell^\p=(S\times \R^{1|1})/(\ell,+ 1)\Z,\quad \rS_\ell^\a=(S\times \R^{1|1})/(\ell,-1)\Z,\\ 
&&\rC_\ell^\p=(S\times \R^{2|1}/(\ell,+ 1)\Z,\quad \rC_\ell^\a=(S\times \R^{2|1})/(\ell,-1)\Z, \label{eq:Eucpair3}
\eeq
 with periodic or antiperiodic spin structure corresponding to the $\p$ and $\a$, respectively. For $(\tau,\bar\tau,\eta)\in \E^{2|1}(S)$, Lemma~\ref{lem:preserve} yields inclusions denoted $i_{\tau,\bar\tau,\eta}$ (compare~\eqref{eq:dumbEucpair2})
\beq\label{eq:Eucpair2}\label{Eq:defnofpairs}
&&\resizebox{.9\textwidth}{!}{$i_{\tau,\bar\tau,\eta}\colon \rS^{\p}_{\ell}=(S\times \R^{1|1})/(\ell,+ 1)\Z\xhookrightarrow{i_0} (S\times \R^{2|1})/(\ell\Z,+1) \xrightarrow{\Tran_{\tau,\bar\tau,\eta}} (S\times \R^{2|1})/(\ell\Z,+)=\rC^{\p}_{\ell}$}\\
&&\resizebox{.9\textwidth}{!}{$i_{\tau,\bar\tau,0}\colon \rS^{\a}_{\ell}=(S\times \R^{1|1})/(\ell,-1)\Z\xhookrightarrow{i_0} (S\times \R^{2|1})/(\ell\Z,-1) \xrightarrow{\Tran_{\tau,\bar\tau,0}} (S\times \R^{2|1})/(\ell\Z,-)=\rC^{\a}_{\ell}$}\label{Eq:defnofpairs2}
\eeq
where $i_0$ is determined by the standard inclusion~\eqref{eq:codim1} and $\Tran_{\tau,\bar\tau,\eta}$ comes from the translation action of $(\tau,\bar\tau,\eta)\in \E^{2|1}(S)$ on $\R^{2|1}$. We note that $i_{\tau,\bar\tau,\eta}=i_{\tau',\bar\tau',\eta'}$ if $(\tau,\bar\tau,\eta)$ and $(\tau',\bar\tau',\eta')$ have the same image in the quotient~$(S\times \E^{2|1})/\ell\Z$. Hence, the data~\eqref{eq:Eucpair2} determine a family of super Euclidean pairs parameterized by $(\ell,\tau)\in (\R_{>0}\times \E^{2|1})/\Z$, and~\eqref{Eq:defnofpairs2} determines a family of super Euclidean pairs parameterized by $(\ell,\tau,\bar\tau)\in (\R_{>0}\times \E^{2})/\Z$ where in both cases the $\Z$-action defining the quotient is determined by
\beq\label{Eq:HHZquot}
&&\Z\times \R_{>0}\times \E^{2|1}\to \R_{>0}\times \E^{2|1},\qquad n\cdot (\ell,\tau,\bar\tau,\eta)=(\ell,n\ell+\tau,n\ell+\bar\tau,\eta).
\eeq

We can enhance the above to produce Euclidean manifolds over a smooth manifold~$M$; we will describe the case for the periodic spin structure, though the antiperiodic construction is entirely analogous. Define standard versions of supercircles and supercylinders
$$
\R^{1|1}/\Z\subset \R^{2|1}/\Z 
$$
for the quotient by the action of the standard subgroup $\Z\subset \E^1\subset \E^{1|1}\subset \E^{2|1}$.

\begin{lem}\label{lem:standard}
For $\ell\in \R_{>0}(S)$, there are isomorphisms of supermanifolds over~$S$
$$
(\rS^\p_\ell\subset \rC^\p_\ell)\simeq (S\times \R^{1|1}/\Z\subset S\times \R^{2|1}/\Z),
$$
with standard families. These isomorphisms are typically not isomorphisms of super Euclidean pairs. 
\end{lem}
\bp
We begin by establishing a little notation. Viewing $\ell\in \R_{>0}(S)\subset C^\infty(S)^\ev$ as the pullback of the standard coordinate on $\R_{>0}$ along a map $S\to \R_{>0}$, there is a function $\sqrt{\ell}\in C^\infty(S)^\ev$ gotten from pulling back the (positive) square root of this coordinate. With this notation in place, define isomorphisms from the families~\eqref{Eq:defnofpairs} to constant families with fiber a standard supercylinder,
\beq\label{eq:stdiso2}
&&\begin{array}{c}
F_\ell \colon \rC^\p_\ell \xrightarrow{\sim}  S\times \R^{2|1}/\Z \qquad {\rm via} \qquad 
 (z,\bar z,\theta)\mapsto (\ell^{-1} z,\ell^{-1}\bar z,\ell^{-1/2} \theta). \label{Eq:reparame}
\end{array}
\eeq 
It is straightforward to check that the displayed map $S\times \R^{2|1}\to S\times \R^{2|1}$ descends to the claimed isomorphisms between $\Z$-quotients in~\eqref{eq:stdiso2}, and furthermore sends subspaces $\rS^\p_\ell\subset \rC^\p_\ell$ to the corresponding standard supercircle. For the maps~\eqref{eq:stdiso2} to be $2|1$-Euclidean isometries, one must have $\ell=1$ and $F_\ell=\id$. 
\ep
Using the above, define a functor
\beq
\Map(\R^{2|1}/\Z,M)\times (\R_{>0}\times \E^{2|1})/\Z&\to& 2|1\EE(M),\nonumber\\ 
(\phi,\ell,\tau,\bar\tau,\eta)&\mapsto& (\rS^\p_\ell\xhookrightarrow{i_{\tau,\bar\tau,\eta}} \rC_\ell^\p\xrightarrow{F_\ell} S\times \R^{2|1}/\Z\xrightarrow{\phi} M)
\eeq
that sends an $S$-point of the source to the $2|1$-Euclidean pair in~\eqref{Eq:defnofpairs} with the map to $M$ given by the above composition, using the isomorphism~\eqref{Eq:reparame}.

\section{Supercylinders in a smooth manifold}\label{sec:supercyl}

\subsection{Warm-up: Euclidean cylinders in a manifold (without supersymmetry)}\label{sec2dwarmup}

To start we will construct moduli spaces of Euclidean cylinders in a manifold~$M$ using the 2-dimensional Euclidean geometry in~\S\ref{sec:dumbEuc}. The basic building block is the Euclidean pair~\eqref{eq:dumbEucpair2}, given by a circle embedded in a Euclidean cylinder with spin structure. Our focus will be on the periodic spin structure, so we will drop the $\p/\a$ decoration, setting $\rS_\ell:=\rS^\p_{\ell}$ and $\rC_\ell:=\rC^\p_{\ell}$ below. Taking the standard identification $\R^2\simeq \C$, define the closed upper and lower half-spaces
\beq\label{Eq:definehalfspaces}
&&\HH_+:=\{\tau\in \R^2\simeq \C \mid \im(\tau)\ge 0\},\qquad \HH_-:=\{\tau\in \R^2\simeq \C \mid \im(\tau)\le 0\}
\eeq
and the open upper and lower half-spaces,
\beq\label{Eq:definehalfspaces2}
&&\iHH_+:=\{\tau\in \R^2\simeq \C \mid \im(\tau)> 0\},\qquad \iHH_-:=\{\tau\in \R^2\simeq \C \mid \im(\tau)< 0\}.
\eeq
Note we have the redundant notation $\HHz=\iHH_+$ for the open upper half plane; we use $\iHH_+$ when emphasizing its role as a subspace of~$\HH_+$.

\begin{defn} \label{defn:EucCyl}
A \emph{positively-oriented Euclidean cylinder in $M$} is given by data $(\Phi,\ell,\tau_\inn,\tau_\out)$ for $\ell\in \R_{>0}$, $\tau_\inn,\tau_\out\in \R^2/\ell\Z$ with $\tau_\out-\tau_\inn\in \HH_+/\ell\Z$, and $\Phi\colon C_{\ell} \to M$. A \emph{negatively-oriented Euclidean cylinder in $M$} is given by the same data but with $\tau_\out-\tau_\inn\in \HH_-/\ell\Z$ in the negative half-space.  A \emph{Euclidean cylinder} in $M$ is a positively or negatively oriented Euclidean cylinder in $M$. The circumference $\ell$ loop~$\Phi\circ i_{\tau_\inn}\colon \rS_\ell \to M$ is the \emph{incoming} boundary and~$\Phi\circ i_{\tau_\out}\colon \rS_\ell\to M$ is the \emph{outgoing} boundary of a Euclidean cylinder in~$M$ for $i_\tau$ defined in~\eqref{eq:dumbEucpair2}. 
An \emph{isometry} between a pair of Euclidean cylinders in $M$ is an isometry of Euclidean manifolds $\Tran_{z,u}\colon C_\ell\to C_\ell$ that fits into the commutative diagram 
\beq\label{diag:atriplemorphnorm}
\begin{tikzpicture}[baseline=(basepoint)];
\node (BB) at (2.5,1.5) {$\rC_\ell$};
\node (A) at (-1,.75) {$\rS_\ell$};
\node (B) at (2.5,0) {$\rC_\ell$};
\node (C) at (6,.75) {$\rS_\ell$};
\node (D) at (2.5,-1.25) {$M$};
\draw[->,left hook-latex] (A) to node [below] {$i_{\tau_\inn}$} (B);
\draw[->,right hook-latex] (C) to node [below] {$i_{\tau_\out}$} (B);
\draw[->,right hook-latex] (A) to node [above=2pt] {$i_{\tau_\inn'}$} (BB);
\draw[->,left hook-latex] (C) to node [above=2pt] {$i_{\tau_\out',}$} (BB);
\draw[->] (B) to node [right] {$\Phi$} (D);
\draw[->] (B) to node [left] {$\Tran_{z,u}$} (BB);
\path (0,0) coordinate (basepoint);
\end{tikzpicture}
\eeq
where $\Tran_{z,u}$ is determined by the left action of $(z,u)\in (\R^{2}\rtimes \Spin(2))(S)$ on $\R^2$, and the source and target of the isometry are $(\ell,\tau_\inn,\tau_\out,\Phi)$ and $(\ell,\tau_\inn',\tau_\out',\Phi')$ respectively. 
\end{defn}

\begin{rmk}
Commutativity of~\eqref{diag:atriplemorphnorm} is equivalent to $\Phi'=\Phi\circ (\Tran_{z,u})^{-1}$ and
$$
\tau_\inn'=(z,u)\cdot \tau_\inn,\qquad \tau_\out'=(z,u)\cdot \tau_\out\in \R^2.
$$
\end{rmk}
By taking the isometry $(z,u)=(-\tau_\inn,1)$, a Euclidean cylinder in $M$ is isometric to one with $\tau_\inn=0$. Hence, the space (or more precisely, sheaf) 
\beq\label{eq:bigEpath}
\Map(\R^2/\Z,M)\times (\R_{>0}\times \HH_+)/\Z\coprod \Map(\R^2/\Z,M)\times (\R_{>0}\times \HH_-)/\Z
\eeq
admits a surjective map to the category in Definition~\ref{defn:EucCyl}. Explicitly, a point of~\eqref{eq:bigEpath} determines data in Definition~\ref{defn:EucCyl} via
$$
\Map(\R^2/\Z,M)\times(\R_{>0}\times \HH_\pm)/\Z \ni (\phi,\ell,\tau) \mapsto (\ell,0,\tau,\phi\circ F_\ell)=(\ell,\tau_\inn,\tau_\out,\Phi)
$$
for $F_\ell$ the diffeomorphism in~\eqref{eq:standardEucnormy}. In the notation here and below, $\phi\colon \R^{2}/\Z\to M$ is a map from a standard cylinder, whereas $\Phi\colon \rC_\ell\to M$ is a map from a Euclidean cylinder with circumference~$\ell$, and so $\Phi=\phi\circ F_\ell$. 

The sheaf~\eqref{eq:bigEpath} inherits an action by isometries of Euclidean cylinders; such isometries are determined by the action of the stabilizer subgroup~\eqref{eq:stabilizer}. Consider the action by the generator of $\Z/4\hookrightarrow \E^1\rtimes \Z/4\subset \E^2\rtimes \Spin(2)$, i.e., the isometry descending from $\sqrt{-1}\in U(1)\simeq \Spin(2)$. This gives an order 4 automorphism of~\eqref{eq:bigEpath} that we denote by $\sqrt{\Fl}$. We note that $\sqrt{\Fl}$ exchanges positively- and negatively-oriented Euclidean cylinders in $M$.
The notation comes from the fact that $\sqrt{\Fl}\circ\sqrt{\Fl}=\Fl$ is the automorphism of~\eqref{eq:bigEpath} determined by the spin flip involution from Lemma~\ref{lem:flip1}.

The 2-dimensional Euclidean geometry affords some additional structures on Euclidean cylinders. Consider the orientation-reversing map
\beq\label{Eq:easyorr}
\R^2\to \R^2, \quad (x,y)\mapsto (x,-y)
\eeq
that preserves the inclusion $\R\subset \R^2$. The \emph{orientation reversal} of a Euclidean cylinder applies~\eqref{Eq:easyorr}, giving the assignment
$$
(\Phi,\ell,\tau_\inn,\tau_\out,\phi)\mapsto (\Phi\circ\orr^{-1},\ell,\orr(\tau_\inn),\orr(\tau_\out)).
$$
Another structure comes from the diffeomorphisms
\beq\label{Eq:easyrg}
\R^2\xrightarrow{\rg_\mu} \R^2,\qquad \rg_\mu(x,y)=(\mu^2x,\mu^2y),\qquad \mu\in \R_{>0}
\eeq
that dilate the metric on $\R^2$. The action through the square of the dilation factor is required to specify an isometry of spin manifolds. The \emph{$\mu$-dilation} of a Euclidean cylinder $(\Phi,\ell,\tau_\inn,\tau_\out)$ applies the diffeomorphism~\eqref{Eq:easyrg}, giving the assignment
$$
(\Phi,\ell,\tau_\inn,\tau_\out)\mapsto  (\Phi\circ\rg_\mu^{-1},\mu^2\ell,\rg_\mu(\tau_\inn),\rg_\mu(\tau_\out))
$$
Identifying $\orr(\tau)=\overline{\tau}\in \HH_\pm\subset \C\simeq \R^2$, and using the description~\eqref{eq:bigEpath} of Euclidean cylinders, orientation reversal and $\mu$-dilation are determined by the maps
\beq
\Map(\R^2/\Z,M)\times (\R_{>0}\times \HH_\pm)/\Z&\to& (\R_{>0}\times \HH_\mp)/\Z\times \Map(\R^2/\Z,M)\nonumber\\
(\Phi,\ell,\tau)&\mapsto&(\Phi\circ \orr^{-1},\ell,\overline{\tau})\label{eq:Estructures1}\\
\Map(\R^2/\Z,M)\times (\R_{>0}\times \HH_\pm)/\Z&\to& (\R_{>0}\times \HH_\pm)/\Z\times \Map(\R^2/\Z,M)\nonumber\\
(\Phi,\ell,\tau)&\mapsto&(\Phi\circ \rg_\mu^{-1},\mu^2\ell,\mu^2\tau).\label{eq:Estructures2}
\eeq
In particular, the orientation reversal of a positively oriented Euclidean cylinder is a negatively oriented Euclidean cylinder. The incoming and outgoing boundary are preserved by orientation reversal, and dilated by $\mu$-dilation. We emphasize that~\eqref{Eq:easyorr} and~\eqref{Eq:easyrg} do not come from isometries in the 2-dimensional Euclidean geometry.

One might consider the extent to which Euclidean cylinders in $M$ can be made into the morphisms of a category, with source and target maps determined by
$$
\inn,\out\colon \Map(\R^2/\Z,M)\times (\R_{>0}\times \HH_\pm)/\Z\rightrightarrows \Map(\R/\Z,M)\times \R_{>0},
$$
sending a Euclidean cylinder in $M$ to its incoming and outgoing boundaries. One expects composition to come from concatenation of Euclidean cylinders. The obstacle is that smooth concatenation is only partially defined: general gluings need not produce smooth maps to~$M$. 
There are various ways to resolve this problem, e.g., by introducing collars or working with a partially-defined composition. 
However, for our purposes it will suffice to restrict to a subspace of~\eqref{eq:bigEpath} where the difficulties of concatenation disappear: the \emph{constant} Euclidean cylinders in~$M$. We define this category in two steps, starting with the following Lie category, i.e., a category internal to smooth manifolds. 


\begin{defn}\label{defn:tPath} Define the Lie category of \emph{constant Euclidean cylinders in $M$}, denoted~$\tAnn(M)$, as having objects and morphisms the smooth manifolds
$$
\Ob(\tAnn(M))=\R_{>0} \times M\coprod \R_{>0} \times M,
$$
$$ 
\Mor(\tAnn(M)) =(\R_{>0}\times \HH_+)/\Z\times M\coprod (\R_{>0}\times \HH_-)/\Z\times M.
$$
The source and target maps in $\tAnn(M)$ are given by the projections.
 The unit map is determined by inclusion at $0\in \HH_\pm$, and composition is determined by addition in $\HH_\pm$. 
\end{defn}

Under the inclusion of the constant maps
\beq\label{eq:normyconstantmaps}
(\R_{>0}\times \HH_\pm)/\Z\times M \hookrightarrow (\R_{>0}\times \HH_\pm)/\Z\times \Map(\R^2/\Z,M)
\eeq
the (partially-defined) operation of concatenation corresponds to composition in the category $\tAnn(M)$. 

The order 4 automorphism~$\sqrt{\Fl}$ determines an automorphism of $\Mor(\tAnn(M))$ (given by restricting~\eqref{eq:Estructures1}), and restricts along the unit map $\u\colon \Ob(\tAnn(M))\to \Mor(\tAnn(M))$ to an order~4 automorphism of $\Ob(\tAnn(M))$ whose generator permutes the components. With respect to these $\Z/4$-actions on $\Mor(\tAnn(M))$ and $\Ob(\tAnn(M))$, the source, target, unit and composition maps are all $\Z/4$-equivariant. This permits the following; we refer to~\S\ref{sec:internalcategories} for an overview of categories internal to stacks from~\cite{ST11}. 

\begin{defn}\label{defn:Path} Define the category of \emph{constant Euclidean cylinders in~$M$}, denoted~$\Ann(M)$, as the category internal to stacks with objects and morphisms
$$
\Ob(\Ann(M))=(\R_{>0} \times M\coprod \R_{>0} \times M)\sq \Z/4,
$$
$$ 
\Mor(\Ann(M)) =\Big(((\R_{>0}\times \HH_+)/\Z\times M\coprod (\R_{>0}\times \HH_-)/\Z\times M)\Big)\sq \Z/4
$$
where the generator of $\Z/4$ acts on objects and morphisms by $\sqrt{\Fl}$. The source, target, unit and composition maps are the $\Z/4$-equivariant extensions of the ones in $\tAnn(M)$. 
\end{defn}

\begin{rmk}\label{rmk:Liegropd}
We observe that $\Ann(M)$ comes from a category strictly internal to Lie groupoids: the source, target, unit and composition functors are all smooth functors between Lie groupoids. Composition is strictly associative, and similarly the left and right unitors are trivial. The functor from Lie groupoids to stacks recovers $\Ann(M)$ as defined above. The Lie groupoid perspective can be a bit more concrete. For example, the evident $\Z/4$-cover 
\beq\label{eq:Z4descent}
\tAnn(M)\to \Ann(M).
\eeq 
is a groupoid quotient on objects and morphisms for the $\Z/4$-action. This allows one to construct a functor out of $\Ann(M)$ in terms of $\Z/4$-equivariant maps out of the objects and morphisms of~$\tAnn(M)$. 
\end{rmk}

Following the previous remark, the assignments~\eqref{eq:Estructures1} and~\eqref{eq:Estructures2} extend to functors
\beq\label{eq:lookfunctors}
\Or\colon \Ann(M)\to \Ann(M),\quad \RG_\mu \colon \Ann(M)\to \Ann(M),
\eeq
where $\RG_\mu$ for $\mu\in \R_{>0}$ determine a (strict) $\R_{>0}$-action and $\Or$ determines a strict involution of~$\Ann(M)$, i.e., $\Or\circ \Or=\id_{\Ann(M)}$.

Next we explain the relationship between Euclidean cylinders and Euclidean tori in~$M$. Our focus will be on the periodic-periodic tori (see~\S\ref{sec:dumbEuc}), though there is an analogous discussion for the other spin structures. There is a $\Z$-action on~\eqref{eq:bigEpath} given by $(\Phi,\ell,\tau)\mapsto (n_\tau\cdot \Phi,\ell,\tau)$ where $n_\tau\cdot \Phi\colon \R^2/\ell\Z\to M$ is the map
\beq\label{eq:fixedsubspacecondition}
(n_\tau\cdot \Phi)(z)=\Phi(z+n\tau),\qquad n\in \Z. 
\eeq
For $\tau\in \iHH_+\subset\HH_+$ in the open upper half plane, consider the $\Z$-fixed subspace 
$$
((\R_{> 0}\times \iHH_+)/\Z \times \Map(\R^2/\Z,M))^\Z \subset (\R_{> 0}\times \HH_+)/\Z\times \Map(\R^2/\Z,M)
$$ 
whose points $(\Phi,\ell,\tau)$ satisfy~\eqref{eq:fixedsubspacecondition}. Such points may be identified with Euclidean cylinders in $M$ whose tori have conformal modulus $\tau/\ell\in \iHH_+/\Z$ and volume $\im(\tau)\ell\in \R_{>0}$. Restricting to constant Euclidean cylinders of positive volume, the condition~\eqref{eq:fixedsubspacecondition} is trivially satisfied, giving a span
\beq\label{eq:loopspace1}
&&\Mor(\Ann(M))\hookleftarrow (\R_{> 0}\times \iHH_\pm)/\Z\times M \hookrightarrow \mathcal{L}^2(M)\sq (\E^2\rtimes \Spin(2)\times \SL_2(\Z))
\eeq
where the quotient stack on the right is defined in~\S\ref{sec:dumbEuc}.

\subsection{Super Euclidean cylinders in $M$}

Analogously to the case of (nonsuper) Euclidean cylinders, the basic ingredient below is the super Euclidean pair~\eqref{eq:Eucpair2}. We will focus on the periodic spin structure, and so drop the $\p/\a$ decoration
$$
\rS_\ell:=S^\p_{\ell}=(S\times \R^{1|1})/\ell\Z \qquad \rC_\ell:=C^\p_{\ell}=(S\times \R^{2|1})/\ell\Z
$$ 
for $S$-families of supermanifolds determined by $\ell\in \R_{>0}(S)$. Define the supermanifolds with boundary via their $S$-points
$$
\HH^{2|1}_\pm(S):=\{(\tau,\bar\tau,\eta)\in \R^{2|1}(S)\mid (\tau,\bar\tau)\in \HH_\pm(S)\subset \R^2(S)\}
$$
or equivalently, restrict the structure sheaf of $\R^{2|1}$ to the half-spaces~\eqref{Eq:definehalfspaces}. Define the supermanifolds $\iHH^{2|1}_\pm$ similarly via the open half-spaces~\eqref{Eq:definehalfspaces2}. In the following, the quotients $(\R_{>0}\times \E^{2|1})/\Z$ and $(\R_{>0}\times \HH_\pm^{2|1})/\Z$ are determined by the $\Z$-action~\eqref{Eq:HHZquot} and its restrictions to $\HH_\pm^{2|1}\subset \E^{2|1}$.

\begin{defn} \label{defn:sEucCyl}
A \emph{positively-oriented $S$-family of $2|1$-Euclidean cylinders in $M$} is given by
\beq\label{eq:cyliderdatamain}
&&(\ell,\tau_\inn,\bar\tau_\inn,\eta_\inn),(\ell,\tau_\out,\bar\tau_\out,\eta_\out)\in ((\R_{>0}\times \E^{2|1})/\Z)(S), \ \  \Phi\colon \rC_\ell=(S\times \R^{2|1})/\ell\Z \to M
\eeq
subject to the condition 
$$
(\tau_\out,\bar\tau_\out,\eta_\out)\cdot(\tau_\inn,\bar\tau_\inn,\eta_\inn)^{-1}\in ((\R_{>0}\times \HH_+^{2|1})/\Z)(S)\subset ((\R_{>0}\times \E^{2|1})/\Z)(S).
$$ 
A \emph{negatively-oriented $S$-family of $2|1$-Euclidean cylinders in $M$} is given by the same data, but with the negative half-plane condition
$$
(\tau_\out,\bar\tau_\out,\eta_\out)\cdot(\tau_\inn,\bar\tau_\inn,\eta_\inn)^{-1}\in ((\R_{>0}\times \HH_-^{2|1})/\Z)(S)\subset ((\R_{>0}\times \E^{2|1})/\Z)(S).
$$ 
 An \emph{$S$-family of $2|1$-Euclidean cylinders} in $M$ is data~\eqref{eq:cyliderdatamain} that over each component of $S$ determines either a positively or negatively oriented super Euclidean cylinder in $M$. 
The $S$-family of superloops~$\Phi\circ i_{\tau_\inn,\bar\tau_\inn,\eta_\inn}\colon \rS_\ell \to M$ is the \emph{incoming} boundary and~$\Phi\circ i_{\tau_\out,\bar\tau_\out,\eta_\out}\colon \rS_\ell\to M$ is the \emph{outgoing} boundary of a super Euclidean cylinder in~$M$ for $i_{\tau,\bar\tau,\eta}$ defined in~\eqref{eq:Eucpair2}. Hence, an $S$-family of $2|1$-Euclidean cylinders in $M$ gives a diagram
\beq\label{diag:atriple}
\begin{tikzpicture}[baseline=(basepoint)];
\node (A) at (0,0) {$\rS_\ell$};
\node (B) at (2.5,0) {$\rC_\ell$};
\node (C) at (5,0) {$\rS_\ell$};
\node (D) at (2.5,-1.25) {$M.$};
\draw[->,right hook-latex] (A) to node [above] {$i_{\tau_\inn,\bar\tau_\inn,\eta_\inn}$} (B);
\draw[->,left hook-latex] (C) to node [above] {$i_{\tau_\out,\bar\tau_\out,\eta_\out}$} (B);
\draw[->] (B) to node [right] {$\Phi$} (D);
\path (0,-.75) coordinate (basepoint);
\end{tikzpicture}
\eeq
An \emph{isometry} between a pair of Euclidean supercylinders in $M$ is an isometry of $S$-families of super Euclidean manifolds $\Tran_{z,\bar z,\theta,u,\bar u}\colon C_\ell\to C_\ell$ that fits into the commutative diagram 
\beq\label{diag:atriplemorph}
\begin{tikzpicture}[baseline=(basepoint)];
\node (BB) at (2.5,1.5) {$\rC_\ell$};
\node (A) at (-1,.75) {$\rS_\ell$};
\node (B) at (2.5,0) {$\rC_\ell$};
\node (C) at (6,.75) {$\rS_\ell$};
\node (D) at (2.5,-1.25) {$M$};
\draw[->,left hook-latex] (A) to node [below] {$i_{\tau_\inn,\bar\tau_\inn,\eta_\inn}$} (B);
\draw[->,right hook-latex] (C) to node [below] {$i_{\tau_\out,\bar\tau_\out,\eta_\out}$} (B);
\draw[->,right hook-latex] (A) to node [above=2pt] {$i_{\tau_\inn',\bar\tau_\inn',\eta_\inn'}$} (BB);
\draw[->,left hook-latex] (C) to node [above=2pt] {$i_{\tau_\out',\bar\tau_\out',\eta_\out'}$} (BB);
\draw[->] (B) to node [right] {$\Phi$} (D);
\draw[->] (B) to node [left] {$\Tran_{z,\bar z,\theta,u,\bar u}$} (BB);
\path (0,0) coordinate (basepoint);
\end{tikzpicture}
\eeq
where $\Tran_{z,\bar z,\theta,u,\bar u}$ is (locally) determined by the left translation for $(z,\bar z,\theta,u,\bar u)\in (\R^{2|1}\rtimes \Spin(2))(S)$ and the source and target supercylinders of~\eqref{diag:atriplemorph} are
$$
(\Phi\colon \rC_\ell\to M,(\tau_\inn,\bar\tau_\inn,\eta_\inn),(\tau_\out,\bar\tau_\out,\eta_\out)),\quad (\Phi\circ \Tran_{z,\bar z,\theta}^{-1}\colon \rC_\ell\to M,(\tau_\inn',\bar\tau_\inn',\eta_\inn'),(\tau_\out',\bar\tau_\out',\eta_\out')),
$$
respectively, where
\beq\label{eq:targetofcyliso}
&&\resizebox{.9\textwidth}{!}{$(\tau_\inn',\bar\tau_\inn',\eta_\inn')=(z,\bar z,\theta,u,\bar u)\cdot (\tau_\inn,\bar\tau_\inn,\eta_\inn), \  (\tau_\out',\bar\tau_\out',\eta_\out')=(z,\bar z,\theta,u,\bar u)\cdot(\tau_\out,\bar\tau_\out,\eta_\out).$}
\eeq

\end{defn}

\begin{rmk}\label{rmk:Beyond}
There is an analogous version of Definition~\ref{defn:sEucCyl} for $2|1$-Euclidean cylinders whose incoming and outgoing boundaries are superloops with anti-periodic spin structures, i.e., for the Neveu--Schwarz sector. By Lemma~\ref{lem:preserve}, the moduli spaces of antiperiodic cylinders are the reduced subspaces $(\R_{>0}\times\HH_\pm)/\Z \subset (\R_{>0}\times \HH^{2|1}_\pm)/\Z$, and so do not contain information about the odd direction in the super Lie group $\E^{2|1}$. For this reason, the categories of periodic supercylinders enjoy a richer supergeometry than the categories of antiperiodic supercylinders. In particular, the odd operators $\bar G_0$ only arise for supersemigroup representations~\eqref{eq:potentialsemigroup} associated with the odd spin circle. This matches the physics literature: when studying $\mathcal{N}=(0,1)$ supersymmetry on a Riemann surface, a (globally-defined) supersymmetry operator requires the existence of a nonvanishing harmonic spinor, which only exists if the spinor bundle is trivializable e.g., see~\cite[2.2]{Beyond}. In the case of cylinders, this supersymmetry operator $\bar G_0$ therefore demands the odd spin structure.
\end{rmk}

\begin{figure}
\begin{center}
\begin{tikzpicture}[scale=1]

\node [draw, cylinder, shape aspect=4, minimum height=8cm, minimum 
width=2cm,dashed, rotate=90] (a) {};
\node (D) at (-1,-1) {$\bullet$};
\node (E) at (-6,-1.5) {$\bullet$};
\node (F) at (-6,2.5) {$\bullet$};
\node (G) at (.5,1.6) {$\bullet$};

\draw[thick] (1,-1) to (1,2);
\draw[thick] (-1,-1) to (-1,2);
\node [draw, thick, ellipse, minimum width=2cm,minimum height=.9cm] (b) at (0,-1){};
\node [draw, thick, ellipse, minimum width=2cm,minimum height=.9cm] (c) at (0,2){};
\node [draw, thick, ellipse, minimum width=2cm,minimum height=.9cm] (d) at (-5,-1.5){};
\node [draw, thick, ellipse, minimum width=2cm,minimum height=.9cm] (e) at (-5,2.5){};
\draw[->,right hook-latex] (-3.8,2.5) to node [above=.1in] {$\begin{array}{c}{\rm translated}\\ {\rm inclusion}\\ {\rm by} \ (\tau,\bar\tau,\eta) \end{array}$} (-1.2, 2);
\draw[->,right hook-latex] (-3.8,-1.5) to node [below=.1in] {$\begin{array}{c}{\rm standard}\\ {\rm inclusion} \end{array}$} (-1.2, -1.1);
\node (A) at (-5,-2.5) {$S_\ell=\R^{1|1}/\ell\Z$};
\node (B) at (-5,3.5) {$S_\ell=\R^{1|1}/\ell\Z$};
\node (C) at (-2.3,.3) {$C_\ell=\R^{2|1}/\ell\Z$};
\draw[<->,thick] (3,-1) to (7,-1);
\draw[<->,dashed] (5,-2.5) to (5,4);
\node (H) at (5.4,3.37) {$\im$};
\draw[<->,thick] (3,2) to (7,2);
\node (H) at (5,-1) {$\bullet$};
\node (H) at (6.5,2) {$\bullet$};
\node (H) at (6,1.75) {$\R^{1|1}$};
\node (H) at (3.5,-1.5) {$\R^{2|1}$};
\node (H) at (6,-1.25) {$\R^{1|1}$};
\draw[->,bend right=30]  (3.5,0) to node [above] {$/\ell\Z$} (1.3,0);
\end{tikzpicture}
\end{center}
\caption{A picture of a $2|1$-Euclidean cylinder. On the right, $\R^{2|1}$ has two embedded copies of~$\R^{1|1}$: the lower copy is the standard embedding, and the upper copy is the translate by $(\tau,\bar\tau,\eta)\in \HH^{2|1}_+$. The marked points correspond to the image of $0\in \R^{1|1}$. 
The quotient by $\ell\Z$ maps to the supercylinder with a pair of embedded supercircles. The marked points indicate the parameterizations of these embeddings; compare~\eqref{eq:QFTsemigroup}. }
\label{fig1}
\end{figure}
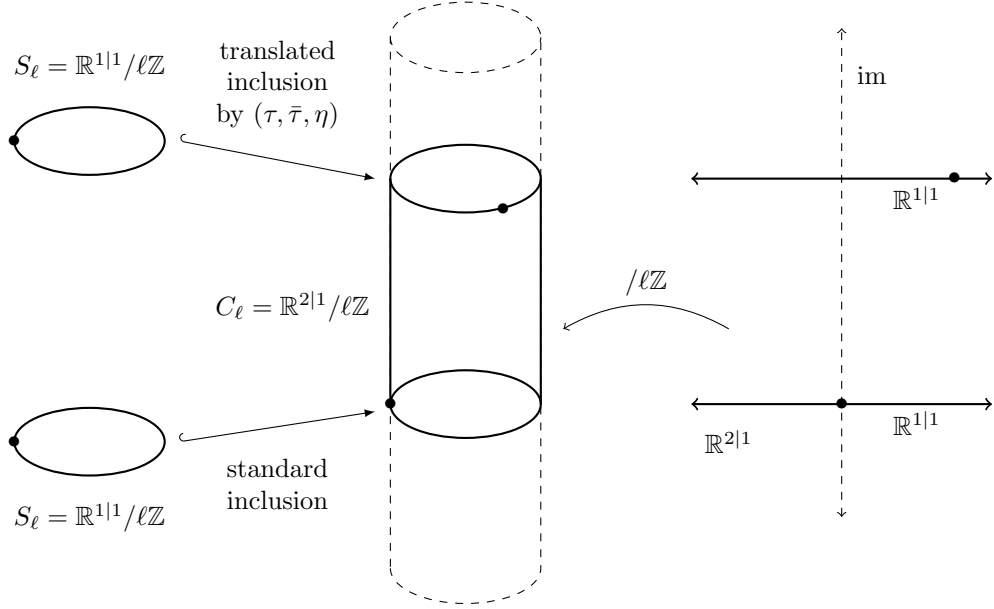



\begin{lem} The category over supermanifolds fibered in groupoids with objects~\eqref{diag:atriple} and morphisms~\eqref{diag:atriplemorph} receives an essentially surjective functor from the presheaf of sets 
\beq\label{eq:sheafofcylinders}
\Map(\R^{2|1}/\Z,M)\times (\R_{>0}\times \left(\HH^{2|1}_+\coprod \HH^{2|1}_-\right))/\Z.
\eeq
\end{lem}
\bp
An $S$-point of $(\phi,\ell,\tau,\bar\tau,\eta)$ of~\eqref{eq:sheafofcylinders} determines a family of $2|1$-Euclidean cylinders via
\beq\nonumber
\begin{tikzpicture}[baseline=(basepoint)];
\node (A) at (0,0) {$\rS_\ell$};
\node (B) at (2.5,0) {$\rC_\ell$};
\node (C) at (5,0) {$\rS_\ell$};
\node (D) at (2.5,-1.25) {$M,$};
\draw[->,right hook-latex] (A) to node [above] {$i_{0}$} (B);
\draw[->,left hook-latex] (C) to node [above] {$i_{\tau,\bar\tau,\eta}$} (B);
\draw[->] (B) to node [right] {$\Phi$} (D);
\path (0,-.75) coordinate (basepoint);
\end{tikzpicture}\qquad \Phi\colon \rC_\ell\xrightarrow{F_\ell} S\times \R^{2|1}/\Z\xrightarrow{\phi} M.
\eeq
This assignment is natural in~$S$, giving a functor between fibered categories. The functor is essentially surjective since any family of $2|1$-Euclidean cylinders is isomorphic to a family with $(\tau_\inn,\bar\tau_\inn,\eta_\inn)=0$ via the translation by $(z,\bar z,\theta)=(-\tau_\inn,-\bar\tau_\inn,-\eta_\inn)$. 
\ep

The presheaf~\eqref{eq:sheafofcylinders} inherits an action by super Euclidean isometries; by Lemma~\ref{lem:stabilizersub}, these isometries are determined by $S$-points of $\E^1\rtimes \Z/4\subset \E^{2|1}\rtimes \Spin(2)$. We will be particularly interested in the generator of $\Z/4$, corresponding to the isometry $\sqrt{-1}\in U(1)\simeq \Spin(2)<\E^{2|1}\rtimes \Spin(2)$ of the super Euclidean group. This isometry descends to an isometry between $2|1$-Euclidean cylinders denoted~$\sqrt{\fl}$,
\beq\label{sec:sqrtfldef}
\begin{tikzpicture}[baseline=(basepoint)];
\node (A) at (0,0) {$S\times \R^{2|1}$};
\node (B) at (4,0) {$S\times \R^{2|1}$};
\node (C) at (0,-1.25) {$\rC_\ell$};
\node (D) at (4,-1.25) {$\rC_\ell$};
\draw[->] (A) to node [above] {$(\sqrt{-1})\cdot $} (B);
\draw[->] (A) to (C);
\draw[->] (B) to (D);
\draw[->] (C) to node [below] {$\sqrt{\fl}$} (D);
\path (0,-.75) coordinate (basepoint);
\end{tikzpicture}\qquad \ell\in \R_{>0}(S)
\eeq
where the vertical arrows are the $\ell\Z$-quotient maps in~\eqref{eq:Eucpair3}, and the notation comes from the fact that $\sqrt{\fl}\circ\sqrt{\fl}=\fl$ is the spin flip involution of $\rC_\ell$ from Lemma~\ref{lem:flip1}. 

\begin{lem}\label{lem:sqrtflaut}
The isometry~\eqref{sec:sqrtfldef} of $\rC_\ell$ determines an order 4 automorphism of the category from Definition~\ref{defn:sEucCyl} that when restricted to the sheaf~\eqref{eq:sheafofcylinders} is by the formula
\beq\label{Eq:fliponSpots}
&& (\Phi,\ell,(\tau,\bar\tau,\eta))\mapsto (\Phi\circ \sqrt{\fl}{}^{-1},\ell,(-\tau,-\bar\tau,\sqrt{-1}\eta)).
\eeq
This automorphism exchanges positively- and negative-oriented $2|1$-Euclidean cylinders.
\end{lem}

\bp
Consider the specialization of~\eqref{diag:atriplemorph}
\beq\label{diag:sqrtflonobject}
\begin{tikzpicture}[baseline=(basepoint)];
\node (BB) at (2.5,1.5) {$\rC_\ell$};
\node (A) at (-1,.75) {$\rS_\ell$};
\node (B) at (2.5,0) {$\rC_\ell$};
\node (C) at (6,.75) {$\rS_\ell$};
\node (D) at (2.5,-1.25) {$M$};
\draw[->,left hook-latex] (A) to node [below] {$i_{\tau_\inn,\bar\tau_\inn,\eta_\inn}$} (B);
\draw[->,right hook-latex] (C) to node [below] {$i_{\tau_\out,\bar\tau_\out,\eta_\out}$} (B);
\draw[->,right hook-latex] (A) to node [above=4pt] {$i_{-\tau_\inn,-\bar\tau_\inn,\sqrt{-1}\eta_\inn}$} (BB);
\draw[->,left hook-latex] (C) to node [above=4pt] {$i_{-\tau_\out,-\bar\tau_\out,\sqrt{-1}\eta_\out}$} (BB);
\draw[->] (B) to node [right] {$\Phi$} (D);
\draw[->] (B) to node [left] {$\sqrt{\fl}$} (BB);
\path (0,0) coordinate (basepoint);
\end{tikzpicture}
\eeq
where isometry between supercylinders is the map $\sqrt{\fl}$ from~\eqref{sec:sqrtfldef}. The inclusions are the unique arrows making the triangles commute. The diagram~\eqref{diag:sqrtflonobject} determines an assignment on the objects~\eqref{diag:atriple} that uniquely extends to a map of categories fibered in groupoids by sending the isometry~\eqref{diag:atriplemorph} determined by $(z,\bar z,\theta)$ to the isometry determined by $(-z,-\bar z,\sqrt{-1}\theta)$. Specializing to an $S$-family of $2|1$-Euclidean cylinders over $M$ gotten from an $S$-point of~\eqref{eq:sheafofcylinders}, we obtain the formula~\eqref{Eq:fliponSpots}. 
\ep

\subsection{Additional structures on $2|1$-Euclidean cylinders in $M$}
The $2|1$-Euclidean geometry has some additional structures that enrich the geometry of $2|1$-Euclidean cylinders. 

\begin{lem}\label{lem:RGaut}
The automorphism~\eqref{eq:rgfunctor2} of super Euclidean manifolds determines an $\R_{>0}$-action on~\eqref{eq:sheafofcylinders} that (for $\mu\in \R_{>0}(S)$) is given by the formula
\beq\label{eq:RGonSpts}
&& \mu\cdot (\Phi,\ell,(\tau,\bar\tau,\eta))= (\Phi\circ \rg_\mu^{-1},\mu^2\ell,(\mu^2\tau,\mu^2\bar\tau,\mu\eta)).
\eeq
\end{lem}
\bp The argument is analogous to that for Lemma~\ref{lem:sqrtflaut}, but using the diagram 
\beq\label{diag:rgonobject}
\begin{tikzpicture}[baseline=(basepoint)];
\node (BB) at (2.5,1.5) {$C_{\mu^2\ell}$};
\node (A) at (-1,0) {$\rS_\ell$};
\node (AA) at (-1,1.5) {$S_{\mu^2\ell}$};
\node (B) at (2.5,0) {$\rC_\ell$};
\node (C) at (6,0) {$\rS_\ell$};
\node (CC) at (6,1.5) {$S_{\mu^2\ell}$};
\node (D) at (2.5,-1.25) {$M$};
\draw[->,right hook-latex] (A) to node [below] {$i_{\tau_\inn,\bar\tau_\inn,\eta_\inn}$} (B);
\draw[->,left hook-latex] (C) to node [below] {$i_{\tau_\out,\bar\tau_\out,\eta_\out}$} (B);
\draw[->,right hook-latex] (AA) to node [above=2pt] {$i_{\mu^2\tau_\inn,\mu^2\bar\tau_\inn,\mu\eta_\inn}$} (BB);
\draw[->,left hook-latex] (CC) to node [above=2pt] {$i_{\mu^2\tau_\out,\mu^2\bar\tau_\out,\mu\eta_\out}$} (BB);
\draw[->] (B) to node [right] {$\Phi$} (D);
\draw[->] (B) to node [left] {$\rg_\mu$} (BB);
\draw[->] (A) to node [left] {$\rg_\mu$} (AA);
\draw[->] (C) to node [left] {$\rg_\mu$} (CC);
\path (0,0) coordinate (basepoint);
\end{tikzpicture}
\eeq
in which the vertical arrows are determined by the map
\beq\label{Eq:21RGdefn}
\rg_\mu\colon \R^{2|1}\to \R^{2|1},\qquad (z,\bar z ,\theta)\mapsto (\mu^2z,\mu^2\bar z,\mu\theta)
\eeq
gotten from specializing~\eqref{eq:rgfunctor1}. 
\ep

Next we analyze the effect of orientation reversal~\eqref{eq:bigor} and complex conjugation~\eqref{eq:bigc}. For an $S$-point $f\colon S\to N$ of a supermanifold $N$, we use the notation $\overline{f}\colon \overline{S}\to \overline{N}$ to denote the corresponding $\overline{S}$-point of $\overline{N}$ from conjugation of supermanifolds. Hence, the assignment $N\mapsto \overline{N}$ determines a map of sheaves that covers the complex conjugation functor on supermanifolds; we emphasize that it does not give a map between $N$ and $\overline{N}$ (this would require $N$ to have a real structure). The notation gets a little clunky for $S$-points of $\R^{2|1}$ but we know of no great alternative; e.g., we will write
$$
(\tau,\bar\tau,\eta)\in \R^{2|1}(S)\mapsto \overline{(\tau,\bar\tau,\eta)}\in \overline{\R}^{2|1}(\overline{S}). 
$$
Finally, for a sheaf $\mathcal{F}$ on supermanifolds, define the conjugate sheaf by
\beq\label{eq:conjsheaf}
\overline{\mathcal{F}}(S):=\mathcal{F}(\overline{S}).
\eeq 

\begin{rmk}
The functor on presheaves determined by~\eqref{eq:conjsheaf} is the Kan extension of the conjugation functor on supermanifolds composed with the Yoneda embedding. In particular,~\eqref{eq:conjsheaf} is compatible with the standard notation for a representable sheaf. 
\end{rmk}

\begin{lem}\label{lem:daggeraut}
Restricted to $2|1$-Euclidean cylinders over $M$, the reflection functor~\eqref{eq:daggerdefn} is determined by a map from the sheaf~\eqref{eq:sheafofcylinders} to its conjugate 
\beq\label{eq:daggerSpts}
 (\Phi,\ell,(\tau,\bar\tau,\eta))\mapsto (\overline{\Phi}\circ\overline{\orr}^{-1},\overline{\ell},\overline{(\bar\tau,\tau,\eta)}).
\eeq
\end{lem}
\bp
The argument again runs in parallel to those above, where we consider the diagram
\beq\label{diag:orRonobject}
\begin{tikzpicture}[baseline=(basepoint)];
\node (BB) at (2.5,1.5) {$\rC_{\overline{\ell}}$};
\node (A) at (-1,0) {$\overline{\rS_{\ell}}$};
\node (AA) at (-1,1.5) {$\rS_{\overline{\ell}}$};
\node (B) at (2.5,0) {$\overline{\rC_{\ell}}$};
\node (C) at (6,0) {$\overline{\rS_{\ell}}$};
\node (CC) at (6,1.5) {$\rS_{\overline{\ell}}$};
\node (D) at (2.5,-1.25) {$\overline{M}$};
\node (DD) at (3.2,-1.28) {$\simeq M$};
\draw[->,right hook-latex] (A) to node [below] {$\overline{i_{\tau_\inn,\bar\tau_\inn,\eta_\inn}}$} (B);
\draw[->,left hook-latex] (C) to node [below] {$\overline{i_{\tau_\out,\bar\tau_\out,\eta_\out}}$} (B);
\draw[->,right hook-latex] (AA) to node [above=2pt] {$i_{\overline{{\bar\tau_\inn,\tau_\inn,\eta_\inn}}}$} (BB);
\draw[->,left hook-latex] (CC) to node [above=2pt] {$i_{\overline{{\bar\tau_\out,\tau_\out,\eta_\out}}}$} (BB);
\draw[->] (B) to node [right] {$\overline{\Phi}$} (D);
\draw[->] (B) to node [left] {$\overline{\orr}$} (BB);
\draw[->] (A) to node [left] {$\overline{\orr}$} (AA);
\draw[->] (C) to node [left] {$\overline{\orr}$} (CC);
\path (0,0) coordinate (basepoint);
\end{tikzpicture}
\eeq
where $S_{\bar \ell}=(\overline{S}\times \R^{1|1})/\bar\ell \Z$ and $C_{\bar \ell}=(\overline{S}\times \R^{2|1})/\bar\ell \Z$. The middle row is the value of the conjugation functor on a bordism over $M$, and the isomorphisms with the upper row use the orientation reversing map. Using~\eqref{eq:bigc} we obtain a $2|(1,0)$-Euclidean structure on the $\overline{S}$-family $\overline{\rC_\ell}\simeq (\overline{S}\times \R^{2|1})/\overline{\ell}\Z$ in the middle row, and then~\eqref{eq:bigor} equips the upper row with a $2|(0,1)$-Euclidean structure. 
The isomorphism $\overline{M}\simeq M$ uses the canonical real structure on an ordinary manifold.
\ep

\begin{lem}\label{lem:commuteuptoflip} The reflection functor~\eqref{eq:daggerdefn} sends the isometry $\sqrt{\fl}\colon C_\ell\to C_\ell$ defined in~\eqref{sec:sqrtfldef} to the isometry $\sqrt{\fl}{}^{-1}\colon \rC_{\bar \ell} \to C_{\bar \ell}$. Equivalently, the maps~\eqref{Eq:fliponSpots} and~\eqref{eq:daggerSpts} commute up to the spin flip automorphism of the family $\rC_\ell$, as $\sqrt{\fl}{}^{-1}=\fl\circ \sqrt{\fl}$. 
\end{lem}
\bp
The statement essentially follows from Lemma~\ref{lem:sqrtflanddagger}. To see this more explicitly, apply the conjugation functor to the diagram~\eqref{diag:sqrtflonobject} and then using the diagram~\eqref{diag:orRonobject}, we obtain an isometry of $\overline{S}$-families of $2|1$-Euclidean cylinders 
$$
\rC_{\bar \ell}\xrightarrow{\overline{\sqrt{\fl}}}\rC_{\bar \ell}\xrightarrow{\overline{\orr}} \rC_{\bar \ell},\qquad (z,\bar z,\theta)\mapsto (-\bar z,-z,-i\theta),
$$
whereas if we postcompose~\eqref{diag:orRonobject} with $\sqrt{\fl}$ we obtain an isometry 
$$
\rC_{\bar \ell}\xrightarrow{\overline{\orr}} \rC_{\bar \ell}\xrightarrow{\sqrt{\fl}} \rC_{\bar \ell},\qquad (z,\bar z,\theta)\mapsto (-\bar z,-z,i\theta),
$$
These maps differ by the action of $-1\in \Spin(2)$, i.e., the spin flip automorphism of~$\rC_\ell$. 
\ep


\subsection{The category of nearly constant supercylinders}
Restricting a super Euclidean cylinder to its incoming and outgoing boundaries determines maps of sheaves
\beq\label{eq:biginout}
\inn,\out\colon \Map(\R^{2|1}/\Z,M)\times (\R_{>0}\times \HH^{2|1}_\pm)/\Z\rightrightarrows \Map(\R^{1|1}/\Z,M)\times\R_{>0}.
\eeq
Attempting to construct a category with source and target data determined by the maps~\eqref{eq:biginout} runs into the same problems as in the case without supersymmetry: composition of $2|1$-Euclidean cylinders is only partially-defined, as concatenations need not be smooth. However, there is a (finite-dimensional) subsheaf of~\eqref{eq:sheafofcylinders} on which smooth concatenations exist, defined as follows.

\begin{defn} 
A \emph{nearly constant} $S$-family of $2|1$-Euclidean cylinders in $M$ is an $S$-point of~\eqref{eq:sheafofcylinders} lying in the subsheaf 
\beq\label{eq:nearlyconstcyl}
\Map(\R^{0|1},M)\times (\R_{>0}\times \HH^{2|1}_\pm)/\Z&\subset& \Map(\R^{2|1}/\Z,M)\times (\R_{>0}\times \HH^{2|1}_\pm)/\Z\\
(\Phi_0,\ell,\tau,\bar\tau,\eta)&\mapsto& (\Phi_0\circ p,\ell,\tau,\bar\tau,\eta)\nonumber
\eeq
for the map $p$ 
$$
p\colon \rC_\ell=(S\times \R^{2|1})/\ell\Z\to S\times\R^{0|1}
$$
determined by the standard projection $\R^{2|1}\to \R^{0|1}$, $(z,\bar z,\theta)\mapsto \theta$.
\end{defn}

The restriction of~\eqref{eq:biginout} to nearly constant $2|1$-Euclidean cylinders gives maps 
\beq\label{eq:littleinout}
\inn,\out\colon \Map(\R^{0|1},M)\times (\R_{>0}\times \HH^{2|1}_\pm)/\Z\rightrightarrows \Map(\R^{0|1},M)\times\R_{>0},
\eeq
whose values are 
\beq\label{eq:inout}
\inn(\Phi_0,\ell,\tau,\bar\tau,\eta)=(\Phi_0,\ell),\qquad \out(\Phi_0,\ell,\tau,\bar\tau,\eta)=(\Phi_0\circ L^{-1}_{\eta},\ell).
\eeq
where $\Tran_\eta\colon S\times \R^{0|1}\to S\times \R^{0|1}$ is left translation by $\eta\in \E^{0|1}(S)$. 

\begin{lem}\label{lem:concat} Suppose that $\Sigma=(\Phi,\ell,\tau,\bar\tau,\eta)$ and $\Sigma'=(\Phi',\ell,\tau',\bar\tau',\eta')$ are $S$-families of positively (respectively, negatively) oriented nearly constant $2|1$-Euclidean cylinders in~$M$ and that $\out(\Sigma)=\inn(\Sigma')$. Then there is a well-defined concatenation $\Sigma'*\Sigma$ that uniquely determines an $S$-point of~\eqref{eq:sheafofcylinders} that has image in the subsheaf of nearly constant $2|1$-Euclidean cylinders in $M$.
\end{lem}
\bp 
In the universal case (i.e., $S$ is the representable subsheaf in~\eqref{eq:nearlyconstcyl}), concatenation of positively or negatively oriented supercylinders is given by
\beq
&& (\Map(\R^{0|1},M)\times (\R_{>0}\times \HH^{2|1}_\pm)/\Z)\times_{\Map(\R^{0|1},M)\times\R_{>0}}  (\Map(\R^{0|1},M)\times (\R_{>0}\times \HH^{2|1}_\pm)/\Z) \nonumber\\
&&\simeq \Map(\R^{0|1},M)\times (\R_{>0}\times \HH^{2|1}_\pm\times \HH^{2|1}_\pm)/\Z  \stackrel{\id\times m}{\longrightarrow}\Map(\R^{0|1},M)\times (\R_{>0}\times \HH^{2|1}_\pm)/\Z \label{eq:compose}
\eeq
where $m$ is determined by the restriction of multiplication on $\E^{2|1}$ to $\HH^{2|1}_\pm$. The result follows. 
\ep

\begin{defn}\label{defn:tsC} The super Lie category of $\tAnn^{2|1}(M)$, has objects and morphisms,
$$
\Ob(\tAnn^{2|1}(M))=\Map(\R^{0|1},M)\times\R_{>0}\coprod \Map(\R^{0|1},M)\times\R_{>0},$$
$$
 \Mor(\tAnn^{2|1}(M)) =\Map(\R^{0|1},M)\times (\R_{>0}\times \HH^{2|1}_+)/\Z\coprod \Map(\R^{0|1},M)\times (\R_{>0}\times \HH^{2|1}_-)/\Z
$$
where source and target are given by the incoming and outgoing boundaries~\eqref{eq:littleinout}, while unit and compositions are determined by the super semigroup structure on $\HH^{2|1}_\pm$ (using Lemma~\ref{lem:concat}).
\end{defn}

Lemma~\ref{lem:sqrtflaut} gives a $\Z/4$-action on $\Mor(\tAnn^{2|1}(M))$ generated by $\sqrt{\Fl}$ and given by the restriction of the formula~\eqref{Eq:fliponSpots} to the subsheaf~\eqref{eq:nearlyconstcyl}. This restricts along the unit map in $\tAnn^{2|1}(M)$ to give a $\Z/4$-action on $\Ob(\tAnn^{2|1}(M))$ corresponding to the pin flip on supercircles, see Remark~\ref{rmk:pin}. The source, target, unit, and composition maps in $\tAnn^{2|1}(M)$ are all $\Z/4$-equivariant. This allows one to form the following category internal to stacks. 

\begin{defn}\label{defn:sC} The internal category $\Ann^{2|1}(M)$ of \emph{nearly constant Euclidean supercylinders in $M$}, has objects and morphisms,
$$
\Ob(\Ann^{2|1}(M))=(\Map(\R^{0|1},M)\times\R_{>0}\coprod \Map(\R^{0|1},M)\times\R_{>0})\sq \Z/4,$$
$$
 \Mor(\Ann^{2|1}(M)) =(\Map(\R^{0|1},M)\times (\R_{>0}\times \HH^{2|1}_+)/\Z\coprod \Map(\R^{0|1},M)\times (\R_{>0}\times \HH^{2|1}_-)/\Z)\sq \Z/4
$$
where source, target, unit, and composition are determined by the $\Z/4$-equivariant structure maps in~$\tAnn^{2|1}(M)$.
\end{defn}

\begin{rmk} Following Remark~\ref{rmk:Liegropd}, $\Ann^{2|1}(M)$ comes from a category strictly internal to super Lie groupoids for which composition is strictly associative. The functor from Lie groupoids to stacks recovers $\Ann^{2|1}(M)$ as defined above. The $\Z/4$-cover of internal categories
\beq\label{eq:sZ4descent}
\tAnn^{2|1}(M)\to \Ann^{2|1}(M)
\eeq 
is a groupoid quotient on objects and morphisms. This allows one to construct a functor out of $\Ann^{2|1}(M)$ in terms of $\Z/4$-equivariant maps out of the objects and morphisms of~$\tAnn^{2|1}(M)$. 
\end{rmk}

\subsection{Characterizing structure maps of $\Ann^{2|1}(M)$ via differential forms on $M$}

The isomorphism of supermanifolds $\Map(\R^{0|1},M)\simeq \Pi TM_\C$ with the odd tangent bundle gives descriptions of functions on the objects and morphisms of $\tAnn^{2|1}(M)$ in terms of differential forms:
\beq
C^\infty(\Map(\R^{0|1},M)\times\R_{>0})&\simeq& C^\infty(\Map(\R^{0|1},M))\otimes C^\infty(\R_{>0})\nonumber\\
&\simeq& \Omega^\bullet(M;C^\infty(\R_{>0})),\label{eq:objfunction}\\ 
C^\infty(\Map(\R^{0|1},M)\times (\R_{>0}\times \HH^{2|1}_\pm)/\Z)&\simeq&C^\infty(\Map(\R^{0|1},M))\otimes C^\infty((\R_{>0}\times \HH^{2|1}_\pm)/\Z) \nonumber\\
&\simeq &\Omega^\bullet(M;C^\infty((\R_{>0}\times \HH^{2|1}_\pm)/\Z)).\label{eq:morfunction}
\eeq
where $\otimes$ is the projective tensor product (e.g., see~\cite[Example~49]{HST}). 
We further observe that Taylor expansion in the odd coordinate $\eta\in C^\infty(\HH^{2|1}_\pm)$ provides an algebra isomorphism
\beq\label{Eq:halfTaylor}
C^\infty((\R_{>0}\times \HH^{2|1}_\pm)/\Z)\simeq C^\infty((\R_{>0}\times \HH_\pm)/\Z)[\eta], 
\eeq
and that there is an injective algebra homomorphism
\beq\label{eq:projfun}
C^\infty(\R_{>0})\hookrightarrow C^\infty((\R_{>0}\times \HH^{2|1}_\pm)/\Z)
\eeq
induced by the projection $(\R_{>0}\times \HH^{2|1}_\pm)/\Z\to \R_{>0}$, i.e., regarding a function on $\R_{>0}$ as a function on $(\R_{>0}\times \HH^{2|1}_\pm)/\Z$ that is constant in $\HH^{2|1}_\pm$. Consider the projection and action maps
$$
\proj\colon ((\R_{>0}\times \HH^{2|1}_\pm)/\Z)\times \Map(\R^{0|1},M) \to \R_{>0}\times \Map(\R^{0|1},M)
$$
$$
\act\colon ((\R_{>0}\times \HH^{2|1}_\pm)/\Z)\times \Map(\R^{0|1},M)\xrightarrow{p}  \R_{>0}\times \R^{0|1}\times \Map(\R^{0|1},M)\xrightarrow{\id_{\R_{>0}}\times a} \R_{>0}\times \Map(\R^{0|1},M)
$$
where $p$ is an intermediate projection, and $a$ is the left action of $\R^{0|1}$ on $\Map(\R^{0|1},M)$ by precomposition. 

\begin{lem}\label{lem:stmap}
Using~\eqref{eq:objfunction} and~\eqref{eq:morfunction}, the source and target maps $\s,\t\colon \Mor(\tAnn^{2|1}(M))\rightrightarrows \Ob(\tAnn^{2|1}(M))$ are determined by the superalgebra maps
$$
\s^*,\t^*\colon \Omega^\bullet(M;C^\infty(\R_{>0}))\to \Omega^\bullet(M;C^\infty((\R_{>0}\times \HH^{2|1}_\pm)/\Z))
$$
\beq\label{eq:projact}
&&\s^*(\alpha)=\proj^*(\alpha)=\alpha,\qquad \t^*(\alpha)=\act^*(\alpha)=(\alpha-\eta d\alpha),
\eeq
where we use the injective map $C^\infty(\R_{>0})\hookrightarrow C^\infty((\R_{>0}\times \HH^{2|1}_\pm)/\Z$ induced by the projection $\R_{>0}\times \HH^{2|1}_\pm \to \R_{>0}$ to regard a $C^\infty(\R_{>0})$-valued form as a $C^\infty((\R_{>0}\times \HH^{2|1}_\pm)/\Z$-valued form.
\end{lem}
\bp 
The $\R^{0|1}$-action on $\Map(\R^{0|1},M)$ can be characterized in terms of a map on differential forms~\cite[Lemma~3.8]{HKST}
\beq
\Omega^\bullet(M)\simeq C^\infty(\Map(\R^{0|1},M))&\to& C^\infty(\R^{0|1}\times \Map(\R^{0|1},M))\simeq \Omega^\bullet(M;C^\infty(\R^{0|1}))\nonumber\\
\alpha&\mapsto& \alpha-\eta  d\alpha,\quad \alpha\in \Omega^\bullet(M), \ C^\infty(\R^{0|1})\simeq C^\infty(\R)[\eta]\label{eq:actdeRham}
\eeq
where $d$ is the de~Rham differential on $\Omega^\bullet(M)$. Then from the description of source and target maps in~\eqref{eq:inout}, the result follows. 
\ep

\begin{rmk}
The formula~\eqref{eq:actdeRham} differs from~\cite[Lemma~3.8]{HKST} by a sign; their formula is for the right action of $\R^{0|1}$ on $\Pi TM_\C$, whereas~\eqref{eq:actdeRham} is for the left action.
\end{rmk}

\begin{lem} The $\Z/4$-action on $\Ob(\tAnn^{2|1}(M))$ and $\Mor(\tAnn^{2|1}(M))$ given by~\eqref{Eq:fliponSpots} is determined by
\beq
&&\sqrt{\Fl} \colon \Mor(\tAnn^{2|1}(M)) \to \Mor(\tAnn^{2|1}(M)),\quad (\Phi,\ell,\tau,\bar\tau,\eta) \mapsto (\Phi\circ \sqrt{\fl},\ell,-\tau,-\bar\tau,\sqrt{-1}\eta).  \label{eq:formFl}
\eeq
In terms of a map of superalgebras,
\beq
\sqrt{\Fl}^*\colon \Omega^\bullet(M;C^\infty((\R_{>0}\times \HH^{2|1}_\pm)/\Z)) &\to& \Omega^\bullet(M;C^\infty((\R_{>0}\times \HH^{2|1}_\mp)/\Z))\nonumber\\
\sqrt{\Fl}^*(\alpha)&=&(\sqrt{-1})^{-|\alpha|}\alpha(\ell,-\tau,-\bar\tau,\sqrt{-1}\eta)\label{Eq:sqflipfun}
\eeq

\end{lem}
\bp 
This follows from Lemma~\ref{lem:sqrtflaut} and the characterization of the $\R^{0|1}\rtimes \C^\times$-action on $\Map(\R^{0|1},M)$ in terms of an algebra map~\cite[Lemma~3.8]{HKST}
\beq\label{eq:HKSTaction}
&&\resizebox{.9\textwidth}{!}{$\Omega^\bullet(M)\simeq C^\infty(\Map(\R^{0|1},M))\to C^\infty(\R^{0|1}\rtimes \C^\times\times \Map(\R^{0|1},M))\simeq \Omega^\bullet(M;C^\infty(\C^\times)[\eta])$}
\eeq
where~\eqref{Eq:sqflipfun} restricts to $\Z/4\subset \C^\times$.
\ep

\begin{lem}\label{lem:structureonsP}
The maps \eqref{eq:RGonSpts} and \eqref{eq:daggerSpts} determine smooth functors 
\beq
\RG_\mu \colon \Ann^{2|1}(M) \to \Ann^{2|1}(M), && (\Phi,\mu,\ell,\tau,\bar\tau,\eta,\sqrt{\Fl})\mapsto (\Phi\circ \rg_\mu^{-1},\mu^2\ell,\mu^2\tau,\mu^2\bar\tau,\mu\eta,\sqrt{\Fl}) \nonumber\label{eq:formRG}\\ 
\dagger \colon \Ann^{2|1}(M) \to \Ann^{2|1}(M),&& (\Phi,\ell,\tau,\bar\tau,\eta,\sqrt{\Fl})\mapsto \overline{(\Phi\circ\orr^{-1},\ell,\bar\tau,\tau,\eta,\sqrt{\Fl}{}^{-1})}\label{eq:formdag} \nonumber
\eeq
given by the indicated maps on $S$-points of objects of $\Mor(\Ann^{2|1}(M))$. Furthermore, we have the equalities of smooth functors 
\beq\label{eq:compositionsofstructure}
\RG_\mu\circ \RG_{\mu'}=\RG_{\mu\mu'},\qquad \overline{\dagger}\circ \dagger=\id,
\eeq
i.e., the functors $\RG_\mu$ determine an $\R_{>0}$ action and $\dagger$ determines an involution covering conjugation of supermanifolds. The functor $\dagger$ determines a reflection structure on $\Ann^{2|1}(M)$ in the sense of Definition~\ref{defn:interreflection}. 
%
%
\end{lem}
\bp 
Functors on the $\Z/4$-cover $\tAnn^{2|1}(M)$ are determined by Lemmas \ref{lem:RGaut} and \ref{lem:daggeraut}. Since $\sqrt{\Fl}$ commutes with the dilation action, the $\Z/4$-equivariance property required for the functor $\RG_\mu \colon \Ann^{2|1}(M) \to \Ann^{2|1}(M)$ is immediate. Using Lemma~\ref{lem:commuteuptoflip}, the functor $\dagger$ arises from a $\Z/4$-equivariance property relative to the inversion homomorphism on~$\Z/4$, 
\beq\label{eq:compositionsofstructure2}
\overline{\dagger}\circ \dagger=\id,\quad \overline{\sqrt{\fl}}\circ \dagger =\dagger \circ \sqrt{\fl}{}^{-1}
\eeq
By definition, this determines a reflection structure on $\Ann^{2|1}(M)$.
\ep

\begin{cor} The functors $\RG_\mu$ and $\dagger$ are determined by the superalgebra maps
\beq
\RG_\mu^*\colon \Omega^\bullet(M;C^\infty((\R_{>0}\times \HH^{2|1}_\pm)/\Z)) &\to& \Omega^\bullet(M;C^\infty((\R_{>0}\times \HH^{2|1}_\pm)/\Z))\nonumber\\
\RG_\mu^*(\alpha)&=&\mu^{-|\alpha|}\alpha\label{Eq:RGfun}\\
\dagger^*\colon \Omega^\bullet(M;C^\infty((\R_{>0}\times \HH^{2|1}_\pm)/\Z)) &\to& \overline{\Omega^\bullet(M;C^\infty((\R_{>0}\times \HH^{2|1}_\mp)/\Z))}\nonumber\\
\dagger^*(\alpha)&=&\overline{\alpha}(\ell,\bar\tau,\tau,\eta)\label{eq:daggerfun}
\eeq
where $|\alpha|$ denotes the $\Z$-grading of $\alpha$ (as a differential form), $\overline{\alpha}$ is the usual complex conjugate of $\alpha$ (as a differential form defined over $\C$), and $\ell\in C^\infty(\R_{>0})$ and $(\tau,\bar\tau,\eta)\in C^\infty(\HH^{2|1}_\pm)$ are the standard coordinates. 
\end{cor}
\bp
We obtain~\eqref{Eq:RGfun} from restricting~\eqref{eq:HKSTaction}  to $\R_{>0}\subset \C^\times$, and~\eqref{eq:daggerfun} uses the fact that the conjugate supermanifold to $\Map(\R^{0|1},M)\simeq \Pi TM_\C$ corresponds to the opposite complex structure on $\C$-valued differential forms. 
\ep

From Definition~\ref{defn:doubleloop}, we recall that $\mathcal{L}^{2|1}(M):=\Map(\R^{2|1}/\Z^2,M)\times s\Lat$. 

\begin{lem} There is a span of inclusions 
\beq\label{eq:supercyltosuperloop}
&&\mathcal{L}^{2|1}(M) \hookleftarrow \Map(\R^{0|1},M)\times (\R_{>0}\times \iHH_+)/\Z\hookrightarrow \Mor(\Ann^{2|1}(M))
\eeq
identifying a subspace of morphisms in $\Ann^{2|1}(M)$ with $S$-points of the super double loop space of $M$. 
\end{lem}
\bp
The inclusion on the right is evident from the definition. The inclusion on the left is determined by
\beq\nonumber
(\Phi_0,\ell,\tau,\bar\tau,0)\mapsto (\Phi_0\circ p,(\ell,\ell,0),(\tau,\bar\tau,0))=(\Phi,\Lambda) \in (\Map(\R^{2|1}/\Z^2,M)\times s\Lat)(S).
\eeq
We check that the above is well-defined. The restriction $(\tau,\bar\tau,\theta)\in \iHH_+(S)\subset \HH_+(S)$ to the open half plane guarantees that this indeed gives a lattice, and the upper half plane condition implies that the lattice is oriented. Hence, the formula above gives an $S$-point $\Lambda\in s\Lat(S)$. The inclusion on mapping spaces comes from restricting~\eqref{eq:nearlyconstcyl}: in view of~\eqref{eq:inout} (taking $\eta=0$) the resulting map $\Phi=\Phi_0\circ p\colon (S\times \R^{2|1})/\ell\Z\to M$ descends to a map out of the $\Z^2$-quotient by $\Lambda$. This completes the proof. 
\ep

\section{Supercylinders and supertori as $2|1$-Euclidean bordisms}\label{sec:iota}

The main goal of this section is to construct an internal functor
\beq\label{eq:maininclude}
\Ann^{2|1}(M)\to 2|1\EBord(M)
\eeq
from the category of nearly constant supercylinders in~$M$ (see Definition~\ref{defn:sC}) to Stolz and Teichner's $2|1$-Euclidean bordism category over~$M$. Roughly, \eqref{eq:maininclude} categorifies the construction~\cite[Lemma 3.15]{DBEChern}
\beq\label{eq:superloopsinclude}
\mathcal{L}_0^{2|1}(M)\to 2|1\EBord(M),
\eeq
that sends a nearly constant super Euclidean torus in $M$ (see Definition~\ref{defn:doubleloop}) to a bordism whose source and target are both the empty $2|1$-Euclidean manifold in~$M$.

To explain how~\eqref{eq:maininclude} and~\eqref{eq:superloopsinclude} are related, we construct variations on~\eqref{eq:maininclude} that regard supercylinders as bordisms with different choices of source and target data, e.g., as bordisms from the empty set to pairs of supercircles in $M$. Compositions of these bordisms corresponds to concatenation of supercylinders or gluing to form a supertorus, and therefore generalizes the basic operations in 2-dimensional Euclidean field theory depicted in~\eqref{eq:Euclideanpicture}. The upshot is that the geometry $2|1$-Euclidean bordisms affords a families generalization of Witten's argument reviewed in~\S\ref{sec:motivate}. 

\begin{rmk}
The constructions below are all phrased within the framework of~\cite{ST11}.
However,~\eqref{eq:maininclude} and its variations are totally geometric constructions that we expect to hold in any reasonable definition of the $2|1$-Euclidean bordism category. 
\end{rmk}

\begin{rmk}\label{rmk:ultimate}
Ultimately, Theorems~\ref{thm1},~\ref{thm2} and~\ref{thm3} only involve the supercylinders whose map to $M$ is nearly constant, see~\S\ref{sec:smallbordisms}. However, for completeness we construct families of supercylinders with an arbitrary map to~$M$. Under disjoint union and composition, these families are expected to generate the full Ramond sector of~$2|1\EBord(M)$, i.e., the subcategory of bordisms whose objects are disjoint unions of superloops in $M$ with the periodic spin structure. In this way, the results of this section give an explicit geometric description of ``half" of Stolz and Teichner's $2|1$-Euclidean bordism category over~$M$~\cite[\S4.3]{ST11}. The other half  (i.e., the Neveu--Schwarz sector) can also be described using the techniques below. Following Remark~\ref{rmk:Beyond}, the Neveu--Schwarz sector has a somewhat less interesting structure.
\end{rmk}

\subsection{Families of bordisms in $2|1\EBord(M)$}\label{sec:superpathbordism}


We refer to~\S\ref{sec:21EB} for a brief review of the definition and notation for  Stolz and Teichner's $2|1$-Euclidean bordism category over~$M$. 

\begin{defn} \label{defn:21EBspt}
Define maps of stacks
\beq\label{eq:sSM}
\sS^+_M,\sS^-_M\colon \Map(\R^{2|1}/\Z,M)\times \R_{>0}\to \Ob(2|1\EBord(M))
\eeq
that to an $S$-point $(\phi,\ell)$ of the source assigns the object in the groupoid $\Ob(2|1\EBord(M))(S)$ 
\beq\label{eq:collardataI}\label{eq:stdinclude}
\begin{array}{c}
\rS_\ell=(S\times \R^{1|1})/\ell\Z=Y^c \subset Y= S\times \R^{2|1}/\ell\Z=\rC_\ell \xrightarrow{\Phi} M\\
Y\setminus Y^c \simeq (S\times \iHH^{2|1}_+)/\Z\coprod (S\times \iHH^{2|1}_-)/\Z =\left\{\begin{array}{ll} Y^+\coprod Y^- & (\rS_\ell^+,\Phi) \\ Y^-\coprod Y^+ & (\rS_\ell^-,\Phi)\end{array}\right.
\end{array}
\eeq
where $F_\ell$ is the isomorphism~\eqref{eq:stdiso2} of supermanifolds and $\Phi=\phi\circ F_\ell$. We adopt the abbreviated notation $(\rS_\ell^\pm,\Phi)$ for objects in~\eqref{eq:stdinclude}. As indicated in the second line, the difference between $(\rS_\ell^+,\Phi)$ and $(\rS_\ell^-,\Phi)$ is the partition of the collar. By pulling back families, the above assignment is natural in~$S$ and hence leads to the claimed maps of stacks~\eqref{eq:sSM}. 
\end{defn}

\begin{rmk}\label{rmk:phiandPhi}
We continue to use the lowercase notation $\phi\colon S\times \R^{2|1}/\Z\to M$ for a map from the $S$-family of standard cylinders, and $\Phi=\phi\circ F_\ell\colon (S\times \R^{2|1})/\ell\Z\to M$ for the map out of the Euclidean family determined by $\ell$. In particular, $S$-points of $\Map(\R^{2|1}/\Z,M)\times \R_{>0}$ are specified either as $(\ell,\phi)$ or $(\ell,\Phi)$, depending on which (equivalent) parameterization of $\Map(\R^{2|1}/\Z,M)\times \R_{>0}$ is more convenient. 
\end{rmk}

\begin{rmk} The families $\sS^+_M$ and $\sS^-_M$ are isomorphic in $2|1\EBord(M)$ via the isometry~$\sqrt{\fl}$ defined in~\eqref{sec:sqrtfldef} that exchanges the partitioning of the collar; see Lemma~\ref{lem:orientation11restrict}. For our purposes, it turns out to be more convenient to remember this isomorphism as additional data: the isometry~$\sqrt{\fl}$ will ultimately endow state spaces with a real structure, essentially by the argument in~\cite[3.2.2]{GPPV}. 
\end{rmk}


\begin{rmk} \label{rmk:germs1}
A pair of $S$-points of $\R_{>0}\times \Map(\R^{2|1}/\Z,M)$ determine isomorphic objects under $\sS^\pm_M$ if the maps $S\times \R^{2|1}/\Z\to M$ agree on the 1-sided neighborhood of the inclusion $S\times \R^{1|1}/\Z\hookrightarrow 
S\times \R^{2|1}/\Z$ corresponding to $Y^+$. Hence, the families in Definition~\ref{defn:21EBspt} pull back from a quotient of $\R_{>0} \times \Map(\R^{2|1}/\Z,M)$ (in sheaves) that parameterizes supercircles with the 1-sided germ of a supercylinder in $M$. 
\end{rmk}


\begin{defn}\label{defn:21EBscyl}
Define maps of stacks
\beq\label{eq:CM}
\begin{array}{c} \sC^+_M\colon \Map(\R^{2|1}/\Z,M)\times (\R_{>0}\times \HH^{2|1}_+)/\Z \to \Mor(2|1\EBord(M))\\ 
\sC^-_M\colon \Map(\R^{2|1}/\Z,M)\times (\R_{>0}\times \HH^{2|1}_-)/\Z \to \Mor(2|1\EBord(M))\end{array}
\eeq
that to an $S$-point $(\phi,\ell,\tau,\bar\tau,\eta)$ of the source assigns the $S$-families of supermanifolds over~$M$,
\beq\label{eq:morphism21}
&&\rS_\ell\coprod \rS_\ell=(S\times \R^{1|1})/\ell\Z \coprod (S\times \R^{1|1})/\ell\Z \xrightarrow{i_0 \coprod i_{\tau,\bar\tau,\eta}} \rC_\ell=(S\times \R^{2|1})/\ell\Z\xrightarrow{\Phi} M
\eeq
where $i_{\tau,\bar\tau,\eta}$ is the inclusion from~\eqref{Eq:defnofpairs} and  $\Phi:=\phi\circ F_\ell$ for $F_\ell$ defined in~\eqref{Eq:reparame}.
To complete the construction of an object in the groupoid $\Mor(2|1\EBord(M))(S)$ from the data~\eqref{eq:morphism21}, we require a choice of source and target. Take
\beq\label{eq:morphisms21EB}
&&\resizebox{.9\textwidth}{!}{$
(Y_{\inn}^c, Y_{\out}^c)=\left\{\begin{array}{lll} \Big((\rS_\ell^+,\Phi), (\rS_\ell^+,\Phi\circ \Tran_{\tau,\bar\tau,\eta}^{-1})\Big) &\implies (\rC_{\ell,\tau,\bar\tau,\eta}^+,\Phi), & (\ell,\tau,\bar\tau,\eta)\in ((\R_{>0}\times \HH^{2|1}_+)/\Z)(S)\\ 
\Big((\rS_\ell^-,\Phi), (\rS_\ell^-,\Phi\circ \Tran_{\tau,\bar\tau,\eta}^{-1})\Big) &\implies (\rC_{\ell,\tau,\bar\tau,\eta}^-,\Phi), & (\ell,\tau,\bar\tau,\eta)\in ((\R_{>0}\times \HH^{2|1}_-)/\Z)(S) \end{array}\right.$}
\eeq
i.e., the ordering in the coproduct in~\eqref{eq:morphism21} lists the source followed by the target. We adopt the abbreviated notation $(\rC_{\ell,\tau,\bar\tau,\eta}^\pm,\Phi)$ for objects in~\eqref{eq:morphisms21EB}. 
\end{defn}

\begin{rmk}\label{rmk:germs2} In parallel to Remark~\ref{rmk:germs1}, $S$-points of $\Map(\R^{2|1}/\Z,M)\times (\R_{>0}\times \HH^{2|1}_\pm)/\Z$ determine isomorphic objects in $\Mor(2|1\EBord(M))$ under $\sC^\pm_M$ if the associated maps $S\times \R^{2|1}/\Z\to M$ agree on a 1-sided neighborhood of the core of the bordism, and so the family pulls back from the associated quotient in sheaves.
\end{rmk}

Changing the source and target data~\eqref{eq:morphisms21EB} give variations on the families~\eqref{eq:CM}.

\begin{defn}\label{defn:21EBscyl2}
Define maps of stacks
\beq\label{eq:LM}
\begin{array}{c} \sL^\pm_M\colon \Map(\R^{2|1}/\Z,M)\times (\R_{>0}\times \HH^{2|1}_\pm)/\Z \to \Mor(2|1\EBord(M))\\ 
\sR^\pm_M\colon \Map(\R^{2|1}/\Z,M)\times (\R_{>0}\times \iHH^{2|1}_\pm)/\Z\to \Mor(2|1\EBord(M))
\end{array}
\eeq
as follows. Given an $S$-point $(\Phi,\ell,\tau,\bar\tau,\eta)$ of the source, take~\eqref{eq:morphism21} as above, but choose source and target data
\beq\label{eq:spath11EB}
&&
\resizebox{.93\textwidth}{!}{$(Y_{\inn}^c, Y_{\out}^c)=\left\{\begin{array}{lll} \Big(\big((\rS_\ell^\pm,\Phi) \coprod (\rS_\ell^\mp,\Phi\circ \Tran_{\tau,\bar\tau,\eta}^{-1})\big),\emptyset\Big) &\implies (\rL_{\ell,\tau,\bar\tau,\eta}^\pm,\Phi) & (\ell,\tau,\bar\tau,\eta)\in ((\R_{>0}\times \HH^{2|1}_\pm)/\Z)(S) \\
\Big(\emptyset,\big((\rS_\ell^\mp,\Phi)\coprod (\rS_\ell^\pm,\Phi\circ \Tran_{\tau,\bar\tau,\eta}^{-1})\big)\Big)& \implies(\rR_{\ell,\tau,\bar\tau,\eta}^\pm,\Phi) & (\ell,\tau,\bar\tau,\eta)\in ((\R_{>0}\times \iHH^{2|1}_\pm)/\Z)(S), \\
\end{array}\right.$}
\eeq
i.e., in the first line the target is empty, and in the second line the source is empty. 
\end{defn}

\begin{rmk}\label{rmk:collarexists}
There are implicit choices of additional collar data in Definitions~\ref{defn:21EBscyl} and~\ref{defn:21EBscyl2}; these choices are required to satisfy the technical \cite[Condition (+), page 30]{ST11} in the definition of the stack $\Mor(2|1\EBord(M))$. When choices of collar data exist, any pair of choices lead to isomorphic objects in $\Mor(2|1\EBord(M))$ essentially by intersection of collars. In the families above, the requisite choices exist and so the functors~\eqref{eq:CM} and~\eqref{eq:LM} are unique up to canonical isomorphism. We note that the restriction on $(\tau,\bar\tau,\eta)$ to the open half plane in the case of $(\rR_{\ell,\tau,\bar\tau,\eta}^\pm,\Phi)$ is necessary for the existence of collar data:
 the families $(\rR_{\ell,\tau,\bar\tau,\eta}^\pm,\Phi)$ in~\eqref{eq:LM} do not have a limit in $\Mor(2|1\EBord(M))$ as $(\tau,\bar\tau,\eta)$ approaches the boundary~$\R^{1|1}\subset \HH^{2|1}_\pm$; compare the moduli of ``right cylinders" in~\cite[\S3.2]{ST11}.
\end{rmk}

\subsection{Source, target, unit and composition of bordisms}\label{sec:compositionEucbord}

The source and target of the families of bordisms classified by $\sC^\pm_M$, $\sL^\pm_M$ and $\sR^\pm_M$ are specified in their definitions. This source and target data can equivalently be phrased in terms of 2-commuting diagrams of stacks where the 2-commutativity data is the identity. For example, we have
\beq
&&\resizebox{.9\textwidth}{!}{$
\begin{tikzpicture}[baseline=(basepoint)];
\node (A) at (0,0) {$\Map(\R^{2|1}/\Z,M)\times (\R_{>0}\times \HH^{2|1}_\pm)/\Z$};
\node (B) at (9,0) {$\Mor(2|1\EBord(M))$};
\node (C) at (0,-1.5) {$\Map(\R^{2|1}/\Z,M)\times \R_{>0}\times \Map(\R^{2|1}/\Z,M)\times \R_{>0}$};
\node (D) at (9,-1.5) {$\Ob(2|1\EBord(M))\times \Ob(2|1\EBord(M))$};
\node (E) at (5.25,-.75) {$\twocommute$};
\draw[->] (A) to node [above] {$\sC^\pm_M$} (B);
\draw[->] (A) to node [left] {$(\proj_\pm, \act_\pm)$} (C);
\draw[->] (C) to node [below] {$\sS^\pm_M\times \sS^\pm_M$} (D);
\draw[->] (B) to node [right] {$\s\times \t$} (D);
\path (0,-.75) coordinate (basepoint);
\end{tikzpicture}$} \label{eq:cyl2comm}
\eeq
where $\proj$ and $\act$ are the projection and action maps (compare~\eqref{eq:projact})
\beq
\proj,\ \act\colon \E^{2|1} \times \Map(\R^{2|1}/\Z,M)&\to& \Map(\R^{2|1}/\Z,M)\nonumber\\ 
\proj(\Phi,\ell,\tau,\bar\tau)&=&(\Phi,\ell)\label{eq:actproj}\\ 
\act(\Phi,\ell,\tau,\bar\tau,\eta)&=&(\Phi\circ \Tran_{\tau,\bar\tau,\eta}^{-1},\ell),\nonumber
\eeq
for $(\ell,\tau,\bar\tau,\eta)\in ((\R_{>0}\times \E^{2|1}_\pm)/\Z)(S)$ and $\Phi\colon (S\times \R^{2|1})/\ell\Z \to M$, and where $\proj_+,\act_+$ (respectively, $\proj_-,\act_-$) are the restrictions of $\proj,\act$ to $\HH^{2|1}_\pm\subset \E^{2|1}$, respectively. The other structure maps in $2|1\EBord(M)$ lead to similar 2-commuting diagrams; the sequence of propositions below spell these diagrams out in the cases of interest. 

\begin{prop}\label{prop:stI}
Let $\u$ be the unit functor in $2|1\EBord(M)$. There is a 2-commutative diagrams of stacks
\beq
&&\begin{tikzpicture}[baseline=(basepoint)];
\node (A) at (0,0) {$\Map(\R^{2|1}/\Z,M)\times \R_{>0}$};
\node (B) at (0,-1.5) {$\Map(\R^{2|1}/\Z,M)\times (\R_{>0}\times \HH^{2|1}_\pm)/\Z$};
\node (C) at (7,0) {$\Ob(2|1\EBord(M))$};
\node (D) at (7,-1.5) {$\Mor(2|1\EBord(M))$};
\node (E) at (3.25,-.75) {$\twocommute$};
\draw[->] (A) to node [above] {$\sS^\pm_M$} (C);
\draw[->,right hook-latex] (A) to node [left] {$\id_{\Map(\R^{2|1}/\Z,M)}\times i_0$} (B);
\draw[->] (C) to node [right] {$\u$} (D);
\draw[->] (B) to node [above] {$\sC^\pm_M$} (D);
\path (0,-.75) coordinate (basepoint);
\end{tikzpicture} \label{eq:unitcyl}
\eeq
where $i_0\colon \R_{>0}\hookrightarrow (\R_{>0}\times \HH^{2|1}_\pm)/\Z$ is inclusion along $0\in \HH^{2|1}_\pm$. Let $\c$ be the composition functor in $2|1\EBord(M)$. There is a 2-commutative diagrams of stacks
\beq
\begin{tikzpicture}[baseline=(basepoint)];
\node (A) at (0,0) {$\Map(\R^{2|1}/\Z,M)\times (\R_{>0}\times \HH^{2|1}_\pm\times \HH^{2|1}_\pm)/\Z$};
\node (B) at (8,0) {$\Mor(2|1\EBord(M))^{[2]}$};
\node (C) at (0,-1.5) {$\Map(\R^{2|1}/\Z,M)\times (\R_{>0}\times \HH^{2|1}_\pm)/\Z$};
\node (D) at (8,-1.5) {$\Mor(2|1\EBord(M))$};
\node (E) at (4,-.75) {$\twocommute$};
\draw[->] (A) to node [above] {$\sC^\pm_M \circ p_2 \times  \sC^\pm_M \circ p_1$} (B);
\draw[->] (A) to node [left] {$\m$} (C);
\draw[->] (C) to node [below] {$\sC^\pm_M$} (D);
\draw[->] (B) to node [right] {$\c$} (D);
\path (0,-.75) coordinate (basepoint);
\end{tikzpicture} \label{eq:comp2iso}
\eeq
where the top horizontal arrow is determined by the family $\sC^\pm_M$ and the maps
\beq\label{eq:composablesuperpath}
&&
\begin{array}{l} p_i\colon \Map(\R^{2|1}/\Z,M)\times (\R_{>0}\times \HH^{2|1}_\pm\times \HH^{2|1}_\pm)/\Z \to \Map(\R^{2|1}/\Z,M)\times (\R_{>0}\times \HH^{2|1}_\pm)/\Z,\\ 
p_1\big(\Phi,\ell,(\tau',\bar\tau',\eta'),(\tau,\bar\tau,\eta)\big)=(\Phi,\ell,\tau,\bar\tau,\eta), \\
p_2\big(\Phi,\ell,(\tau,\bar\tau,\eta),(\tau',\bar\tau',\eta')\big)=(\Phi\circ \Tran_{\tau,\bar\tau,\eta}^{-1},\ell,\tau',\bar\tau',\eta')\end{array}
\eeq
and $\m$ is defined as
\beq\label{eq:sfm}
\m\big(\Phi,\ell,(\tau',\bar\tau',\eta'),(\tau,\bar\tau,\eta)\big)=(\Phi,\ell,\tau'+\tau,\bar\tau'+\bar\tau+\eta'\eta,\eta'+\eta).
\eeq
\end{prop}

\bp
The commutativity of the diagram~\eqref{eq:unitcyl} is almost immediate; the only subtlety comes the choice of collars in the construction of the family $\sC^\pm_M$. Indeed, rather than a strictly commutative diagram this choice leads to 2-commutativity data gotten from the isomorphism between objects of $\Mor(2|1\EBord(M))$ from shrinking the collar in~\eqref{eq:collardataI}; see Remark~\ref{rmk:collarexists}. 

Turning to the diagram~\eqref{eq:comp2iso}, the definition of the $S$-family $(\rC^\pm_{\ell,\tau,\bar\tau,\eta},\Phi)$ and the maps~\eqref{eq:composablesuperpath} gives a map to the $S$-point of the fibered product (the strict and weak fibered product coincide, see Remark~\ref{rmk:fibration})
$$
\big(\ell,(\tau',\bar\tau',\eta'),(\tau,\bar\tau,\eta),\Phi\big)\mapsto \big((\rC^\pm_{\ell,\tau',\bar\tau',\eta'},\Phi\circ \Tran_{\tau,\bar\tau,\eta}^{-1}),(\rC^\pm_{\ell,\tau,\bar\tau,\eta},\Phi)\big). 
$$
 The construction of 2-commutativity data comes from unpacking the definition of the functor $\c$ from~\cite[page~20]{ST11}, i.e., composition of bordisms. Consider the diagram 
\beq\label{eq:triplepath}
&&
\rS_\ell \coprod \rS_\ell\coprod \rS_\ell \xrightarrow{i_0 \coprod i_{\tau,\bar\tau,\eta}\coprod i_{(\tau',\bar\tau',\eta')\cdot (\tau,\bar\tau,\eta)}} (S\times \R^{2|1})/\Z \xrightarrow{\Phi} M
\eeq
where the inclusions are as in~\eqref{Eq:defnofpairs}. 
By forgetting any one of these three inclusions, we obtain three possible $S$-points of the source of $\sC^\pm_M$ via the input datum~\eqref{eq:morphism21}: forgetting $i_{(\tau',\bar\tau',\eta')\cdot (\tau,\bar\tau,\eta)}$ gives the $S$-family of bordisms $(\rC_{\ell,\tau,\bar\tau,\eta}^\pm,\Phi)$, forgetting $i_{\tau,\bar\tau,\eta}$ gives the $S$-family of bordisms $(\rC^\pm_{\ell, (\tau',\bar\tau',\eta')\cdot (\tau,\bar\tau,\eta)},\Phi)$, and forgetting $i_0$ gives a family that is isomorphic to $(\rC_{\ell,\tau',\bar\tau',\eta'},\Phi\circ \Tran_{\tau,\bar\tau,\eta}^{-1})$ with isomorphism specified by the left-translation map $\Tran_{\tau,\bar\tau,\eta}$. These three $S$-families of bordisms gotten from~\eqref{eq:triplepath} coincide with the two projections out of the fibered product and the value of the functor~$\c$. The only subtle aspect in commutativity of the diagram~\eqref{eq:comp2iso} is again the choice of collars in the definitions of these families of bordisms, which in turn depends on the choices in the construction of the functor $\sC^\pm_M$. Since any pair of choices lead to isomorphic bordisms, we obtain canonical 2-commutativity data in~\eqref{eq:comp2iso}.
\ep

\begin{rmk}
The picture corresponding to the commutativity of~\eqref{eq:comp2iso} is
\beq\label{pic:pic1}
&&\begin{tikzpicture}[baseline=(basepoint)];
	\begin{pgfonlayer}{nodelayer}
		\node [style=none] (0) at (1.5, 1.75) {};
		\node [style=none] (1) at (1, 1.75) {};
		\node [style=none] (2) at (1.5, 0.25) {};
		\node [style=none] (3) at (1, 0.25) {};
		\node [style=none] (4) at (1.5, -1.25) {};
		\node [style=none] (5) at (1, -1.25) {};
		\node [style=none] (6) at (-2, 2) {};
		\node [style=none] (7) at (-2.5, 2) {};
		\node [style=none] (8) at (-2, 0.5) {};
		\node [style=none] (9) at (-2.5, 0.5) {};
		\node [style=none] (10) at (-2, 0) {};
		\node [style=none] (11) at (-2.5, 0) {};
		\node [style=none] (12) at (-2, -1.5) {};
		\node [style=none] (13) at (-2.5, -1.5) {};
		\node [style=none] (14) at (-2.5, 0.5) {};
		\node [style=none] (15) at (-2.5, 0.5) {};
				\node [style=none] (21) at (-.5, .25) {$\xrightarrow{\rm glue/compose}$};
		\node [style=none] (20) at (3.75,.75) {$(\rC^\pm_{\ell,\tau+\tau',\bar \tau+\bar \tau+\eta\eta',\eta+\eta'},\Phi)$};
\node [style=none] (19) at (-3.5,-.5) {$(\rC^\pm_{\ell,\tau,\bar \tau,\eta},\Phi)$};
\node [style=none] (18) at (-4.25,1.5) {$(\rC^\pm_{\ell,\tau',\bar \tau',\eta'},\Phi\circ \Tran_{\tau,\bar\tau,\eta}^{-1})$};
	\end{pgfonlayer}
	\begin{pgfonlayer}{edgelayer}
		\draw [bend left=270, looseness=0.75,thick] (0.center) to (1.center);
		\draw [bend left=270,thick,dashed] (2.center) to (3.center);
		\draw [bend left=270,thick,dashed] (4.center) to (5.center);
		\draw [thick] (0.center) to (4.center);
		\draw [thick] (1.center) to (5.center);
		\draw [bend left=90,thick] (0.center) to (1.center);
		\draw [bend left=90, looseness=0.75,thick] (2.center) to (3.center);
		\draw [bend left=90, looseness=0.75,thick] (4.center) to (5.center);
		\draw [bend left=270, looseness=0.75,thick] (6.center) to (7.center);
		\draw [bend left=270,thick,dashed] (8.center) to (9.center);
		\draw [bend left=90, looseness=0.75,thick] (6.center) to (7.center);
		\draw [bend left=90, looseness=0.75,thick] (8.center) to (9.center);
		\draw [thick] (8.center) to (6.center);
		\draw [thick] (9.center) to (7.center);
		\draw [bend left=270, looseness=0.75,thick] (10.center) to (11.center);
		\draw [bend left=270,thick,dashed] (12.center) to (13.center);
		\draw [bend left=90, looseness=0.75,thick] (10.center) to (11.center);
		\draw [bend left=90, looseness=0.75,thick] (12.center) to (13.center);
		\draw [thick] (12.center) to (10.center);
		\draw [thick] (13.center) to (11.center);
	\end{pgfonlayer}
	\path (0,0) coordinate (basepoint);
\end{tikzpicture}\label{eq:sgluecylinderpic}
\eeq
i.e., gluing of supercylinders with composition depicted vertically, where $(\ell,\tau,\bar\tau,\eta),(\ell,\tau',\bar\tau',\eta')\in ((\R_{>0}\times \HH^{2|1}_\pm)/\Z)(S)$; compare with the left picture in~\eqref{eq:gluecylinderpic}.
\end{rmk}

\begin{notation} Let $\sC^\pm_0$ denote the composition $\u\circ \sS^\pm_M$, 
$$
\sC^\pm_0\colon \R_{>0}\times \Map(\R^{2|1}/\Z,M)\to \Mor(2|1\EBord(M)).
$$
For an $S$-point of the source above, let $(\rC_{\ell,0}^\pm,\Phi)$ denote the associated $S$-family of (identity) $2|1$-Euclidean bordisms. 
\end{notation}

\begin{prop}\label{prop:adjunctionbordisms}
There is a 2-commutative diagrams of stacks
\beq
&&\resizebox{.95\textwidth}{!}{$\begin{tikzpicture}[baseline=(basepoint)];
\node (A) at (0,0) {$\Map(\R^{2|1}/\Z,M)\times (\R_{>0}\times \HH^{2|1}_\pm\times \HH^{2|1}_\pm\times \HH^{2|1}_\mp)/\Z$};
\node (B) at (11,0) {$\Mor(2|1\EBord(M))^{[2]}$};
\node (C) at (0,-1.5) {$\Map(\R^{2|1}/\Z,M)\times (\R_{>0}\times \HH^{2|1}_\pm)/\Z$};
\node (D) at (11,-1.5) {$\Mor(2|1\EBord(M))$};
\node (E) at (6.25,-.75) {$\twocommute$};
\draw[->] (A) to node [above] {$\sL^\pm_M\circ p_2 \times (\sC^\pm_M\circ p_1 \coprod \sC^\mp_M\circ p_3)$} (B);
\draw[->] (A) to (C);
\draw[->] (C) to node [below] {$\sL^\pm_M$} (D);
\draw[->] (B) to node [right] {$\c$} (D);
\path (0,-.75) coordinate (basepoint);
\end{tikzpicture}$} \label{eq:adjunctioncyl}
\eeq
where the top horizontal arrow is determined by the families $\sL^\pm_M$, $\sC^\pm_M$ and $\sC^\mp_M$ and the maps
\beq
p_i\colon \Map(\R^{2|1}/\Z,M)\times (\R_{>0}\times \HH^{2|1}_\pm\times \HH^{2|1}_\pm\times \HH^{2|1}_\mp)/\Z\to \Map(\R^{2|1}/\Z,M)\times (\R_{>0}\times \HH^{2|1}_\pm)/\Z\nonumber\\
p_1(\Phi,\ell,(\tau,\bar\tau,\eta),(\tau',\bar\tau',\eta'),(\tau'',\bar\tau'',\eta''))= (\Phi,\ell,\tau,\bar\tau,\eta), \\
 p_2(\Phi,\ell,(\tau,\bar\tau,\eta),(\tau',\bar\tau',\eta'),(\tau'',\bar\tau'',\eta''))= (\Phi\circ \Tran_{\tau,\bar\tau,\eta}^{-1},\tau',\bar\tau',\eta') \label{eq:targetofadj}\\
p_3\colon \Map(\R^{2|1}/\Z,M)\times (\R_{>0}\times \HH^{2|1}_\pm\times \HH^{2|1}_\pm\times \HH^{2|1}_\mp)/\Z\to \Map(\R^{2|1}/\Z,M)\times (\R_{>0}\times \HH^{2|1}_\mp)/\Z\nonumber\\
p_3(\Phi,\ell,(\tau,\bar\tau,\eta),(\tau',\bar\tau',\eta'),(\tau'',\bar\tau'',\eta''))= (\Phi\circ \Tran_{\tau,\bar\tau,\eta}^{-1}\circ L_{\tau',\bar\tau',\eta'}^{-1}, \ell,\tau'',\bar\tau'',\eta'')\nonumber
\eeq
and the left vertical arrow is determined by multiplication in $\E^{2|1}$,
\beq\label{eq:sfm2}
&&(\Phi,\ell,(\tau,\bar\tau,\eta),(\tau',\bar\tau',\eta'),(\tau'',\bar\tau'',\eta''))\mapsto (\Phi,\ell,(-\tau'',-\bar\tau'',-\eta'')\cdot (\tau',\bar\tau',\eta')\cdot(\tau,\bar\tau,\eta))
\eeq
\end{prop}
\bp
From the definitions of source and target data, the top map in~\eqref{eq:adjunctioncyl} does indeed yield an object in the fibered product that in terms of $S$-families of bordisms is
$$
\resizebox{\textwidth}{!}{$(\Phi,\ell,(\tau,\bar\tau,\eta),(\tau',\bar\tau',\eta'),(\tau'',\bar\tau'',\eta''))\mapsto \Big((\rL^\pm_{\ell,\tau',\bar\tau',\theta},\Phi\circ \Tran_{\tau,\bar\tau,\eta}^{-1}),\big((\rC^\pm_{\ell,\tau,\bar\tau,\eta},\Phi)\coprod (C^\mp_{\ell,\tau'',\bar\tau'',\eta''},\Phi\circ \Tran_{\tau,\bar\tau,\eta}^{-1}\circ \Tran_{\tau',\bar\tau',\eta'}^{-1})\big)\Big).$}
$$
The 2-commutativity data again comes from unpacking the definition of the functor $\c$. Similar to before, consider the diagram
\beq
\resizebox{.9\textwidth}{!}{$
\rS_\ell\coprod \rS_\ell\coprod \rS_\ell\coprod \rS_\ell \xrightarrow{i_0 \coprod i_{\tau,\bar\tau,\eta}\coprod i_{(\tau',\bar\tau',\eta')\cdot(\tau,\bar\tau,\eta)}\coprod i_{(\tau'',\bar\tau'',\eta'')^{-1}\cdot (\tau',\bar\tau',\eta')\cdot (\tau,\bar\tau,\eta)}} (S\times \R^{2|1})/\ell\Z \xrightarrow{\Phi} M.$}\nonumber
\eeq
In this case, forgetting any pair of factors in the coproduct gives the input data to an $S$-family of bordisms. For example, forgetting the middle two factors corresponds to the $S$-family of bordisms $(\rL^\pm_{\ell,(\tau'',\bar\tau'',\eta'')^{-1}\cdot (\tau',\bar\tau',\eta')\cdot (\tau,\bar\tau,\eta)},\Phi)$, which is the claimed composition.  The bordisms in~\eqref{eq:targetofadj} can similarly be extracted by forgetting the appropriate pair of inclusions and applying a translation action to $(S\times \R^{2|1})/\ell\Z$. 
\ep

\begin{rmk} The picture corresponding to the commutativity of~\eqref{eq:adjunctioncyl} is 
\beq\label{pic:pic2}
&&\resizebox{.95\textwidth}{!}{$
\begin{tikzpicture}[baseline=(basepoint)];
	\begin{pgfonlayer}{nodelayer}
		\node [style=none] (0) at (2.5, 0.5) {};
		\node [style=none] (1) at (2, 0.5) {};
		\node [style=none] (2) at (0.5, 0.5) {};
		\node [style=none] (3) at (0, 0.5) {};
		\node [style=none] (4) at (0.5, 0) {};
		\node [style=none] (5) at (0, 0) {};
		\node [style=none] (6) at (0.5, -1.5) {};
		\node [style=none] (7) at (0, -1.5) {};
		\node [style=none] (8) at (2.5, 0) {};
		\node [style=none] (9) at (2, 0) {};
		\node [style=none] (10) at (2.5, -1.5) {};
		\node [style=none] (11) at (2, -1.5) {};
		\node [style=none] (12) at (6.5, 0) {};
		\node [style=none] (13) at (6, 0) {};
		\node [style=none] (14) at (4.5, 0) {};
		\node [style=none] (15) at (4, 0) {};
		\node [style=none] (18) at (4.5, -1.25) {};
		\node [style=none] (19) at (4, -1.25) {};
		\node [style=none] (20) at (6.5, -1.25) {};
		\node [style=none] (21) at (6, -1.25) {};
		\node [style=none] (22) at (3.25, 0) {$\xrightarrow{\rm compose}$};
		\node [style=none] (23) at (4, -2) {$(\rC^\mp_{\ell,\tau'',\bar\tau'',\theta''},\Phi\circ L_{\tau,\bar\tau,\eta}^{-1}\circ L_{\tau',\bar\tau',\eta'}^{-1})$};
		\node [style=none] (24) at (0, -2) {$(\rC^\pm_{\ell,\tau,\bar\tau,\theta},\Phi)$};
		\node [style=none] (25) at (-1.75, 1.25) {$(\rL^\pm_{\ell,\tau',\bar\tau',\theta'},\Phi\circ L_{\tau,\bar\tau,\eta}^{-1})$};
		\node [style=none] (25) at (9, 1) {$(\rL^\pm_{\ell,(\tau'',\bar\tau'',\theta'')^{-1}\cdot (\tau',\bar\tau',\theta')\cdot (\tau,\bar\tau,\theta)},\Phi)$};
	\end{pgfonlayer}
	\begin{pgfonlayer}{edgelayer}
		\draw [bend right=105, looseness=1.25,thick] (0.center) to (1.center);
		\draw [bend left=90, looseness=0.75,thick] (0.center) to (1.center);
		\draw [bend right=90, looseness=1.50,thick] (0.center) to (3.center);
		\draw [bend right=90, looseness=1.50,thick] (1.center) to (2.center);
		\draw [bend left=90,thick] (2.center) to (3.center);
		\draw [bend right=90,thick] (2.center) to (3.center);
		\draw [bend left=270,thick] (4.center) to (5.center);
		\draw [bend left=285,thick,dashed] (6.center) to (7.center);
		\draw [bend left=90, looseness=0.75,thick,dashed] (4.center) to (5.center);
		\draw [bend left=90,thick] (6.center) to (7.center);
		\draw [thick] (6.center) to (4.center);
		\draw [thick] (7.center) to (5.center);
		\draw [bend left=270,thick] (8.center) to (9.center);
		\draw [bend left=270,thick,dashed] (10.center) to (11.center);
		\draw [bend left=90, looseness=0.75,thick,dashed] (8.center) to (9.center);
		\draw [bend left=75, looseness=0.75,thick] (10.center) to (11.center);
		\draw [thick] (10.center) to (8.center);
		\draw [thick] (11.center) to (9.center);
		\draw [bend right=105, looseness=1.25,thick] (12.center) to (13.center);
		\draw [bend left=90, looseness=0.75,thick,dashed] (12.center) to (13.center);
		\draw [bend right=90, looseness=1.50,thick] (12.center) to (15.center);
		\draw [bend right=90, looseness=1.50,thick] (13.center) to (14.center);
		\draw [bend left=90,thick,dashed] (14.center) to (15.center);
		\draw [bend right=90,thick] (14.center) to (15.center);
		\draw [bend left=270,thick,dashed] (18.center) to (19.center);
		\draw [bend left=75, looseness=0.75,thick] (18.center) to (19.center);
		\draw [bend left=270,thick,dashed] (20.center) to (21.center);
		\draw [bend left=75, looseness=0.75,thick] (20.center) to (21.center);
		\draw [thick] (15.center) to (19.center);
		\draw [thick] (14.center) to (18.center);
		\draw [thick] (13.center) to (21.center);
		\draw [thick] (12.center) to (20.center);
	\end{pgfonlayer}
	\path (0,.25) coordinate (basepoint);
\end{tikzpicture}
$}
\eeq
where the direction of the bordism is oriented vertically and $(\ell,\tau,\bar\tau,\eta),(\ell,\tau',\bar\tau',\eta')\in ((\R_{>0}\times \HH^{2|1}_\pm)/\Z)(S)$ and $(\ell,\tau'',\bar\tau'',\eta'')\in ((\R_{>0}\times \HH^{2|1}_\mp)/\Z)(S)$. 
\end{rmk}

The following shows that the restriction of the family of bordisms $\sC^\pm_M$ along the inclusion $\iHH_\pm^{2|1}\subset \HH^{2|1}$ factors as a composition involving $\sL^\pm_M$ and $\sR^\pm_M$. The proof uses similar techniques as above. 

\begin{prop}\label{prop:duality}
There is a 2-commutative diagrams of stacks
\beq
&&\resizebox{.95\textwidth}{!}{$
\begin{tikzpicture}[baseline=(basepoint)];
\node (A) at (0,0) {$\Map(\R^{2|1}/\Z,M)\times (\R_{>0}\times \HH^{2|1}_\pm\times \iHH^{2|1}_\pm)/\Z$};
\node (B) at (12,0) {$\Mor(2|1\EBord(M))^{[2]}$};
\node (C) at (0,-1.5) {$\Map(\R^{2|1}/\Z,M)\times (\R_{>0}\times \HH^{2|1}_\pm)/\Z$};
\node (D) at (12,-1.5) {$\Mor(2|1\EBord(M))$};
\node (E) at (6.25,-.75) {$\twocommute$};
\draw[->] (A) to node [above] {$(\sL_M^\pm \circ p_1\coprod \sC_0^\pm \circ q_2)\times (\sC_0^\pm \circ q_1\coprod \sR^\pm_M\circ q_2)$} (B);
\draw[->] (A) to node [left] {$\m$} (C);
\draw[->] (C) to node [below] {$\sC^\pm_M$} (D);
\draw[->] (B) to node [right] {$\c$} (D);
\path (0,-.75) coordinate (basepoint);
\end{tikzpicture}$} \label{eq:thickdiagram}
\eeq
where the top horizontal arrow is determined by the families $\sL^\pm_M$, $\sR^\pm_M$, $\sC^\pm_0$, and the maps 
\beq
p_i\colon \Map(\R^{2|1}/\Z,M)\times (\R_{>0}\times \HH^{2|1}_\pm\times \iHH^{2|1}_\pm)/\Z \to \times \Map(\R^{2|1}/\Z,M)\times (\R_{>0}\times \HH^{2|1}_\pm)/\Z\nonumber\\
p_1(\Phi,\ell,(\tau,\bar\tau,\eta),(\tau',\bar\tau',\eta'))=(\Phi,\ell,\tau,\bar\tau,\eta),\\
 p_2(\Phi,\ell,(\tau,\bar\tau,\eta),(\tau',\bar\tau',\eta'))=(\tau',\bar\tau',\eta',\Phi\circ \Tran_{\tau,\bar\tau,\eta}^{-1}), \label{eq:targetofdual1}\\
q_i\colon \colon  \Map(\R^{2|1}/\Z,M)\times (\R_{>0}\times \HH^{2|1}_\pm\times \iHH^{2|1}_\pm)/\Z\to \R_{>0}\times \Map(\R^{2|1}/\Z,M)\nonumber\\
q_1(\Phi,\ell,(\tau,\bar\tau,\eta),(\tau',\bar\tau',\eta'))=(\ell,\Phi), \\ 
q_2(\Phi,\ell,(\tau,\bar\tau,\eta),(\tau',\bar\tau',\eta'))=(\ell,\Phi\circ \Tran_{\tau,\bar\tau,\eta}^{-1} \circ \Tran_{\tau',\bar\tau',\eta'}^{-1})
\eeq
and $\m$ is~\eqref{eq:sfm}. 
\end{prop}

\begin{rmk} The picture corresponding to the commutativity of~\eqref{eq:thickdiagram} is 
\beq\label{pic:pic3}
&&\begin{tikzpicture}[baseline=(basepoint)];
	\begin{pgfonlayer}{nodelayer}
		\node [style=none] (0) at (-1.5, -0.25) {};
		\node [style=none] (2) at (-3, 0.25) {};
		\node [style=none] (3) at (-3.5, 0.25) {};
		\node [style=none] (4) at (-2, -0.25) {};
		\node [style=none] (5) at (-4.5, 0.25) {};
		\node [style=none] (6) at (-5, 0.25) {};
		\node [style=none] (13) at (-3, -0.25) {};
		\node [style=none] (14) at (-3.5, -0.25) {};
		\node [style=none] (15) at (0.5, 1.5) {};
		\node [style=none] (16) at (0, 1.5) {};
		\node [style=none] (17) at (0.5, 0) {};
		\node [style=none] (18) at (0, 0) {};
		\node [style=none] (19) at (0.5, -1.5) {};
		\node [style=none] (20) at (0, -1.5) {};
		\node [style=none] (21) at (2.55,.5) {$(\rC^\pm_{\ell,(\tau,\bar\tau,\eta)\cdot(\tau',\bar\tau',\eta')},\Phi)$};
		\node [style=none] (22) at (-6, 1) {$(\rL^\pm_{\ell,\tau,\bar\tau,\eta},\Phi)$};
		\node [style=none] (23) at (-2.5, -1.5) {$(\rR^\pm_{\ell,\tau',\bar\tau',\eta'},\Phi\circ\Tran_{\tau,\bar\tau,\eta}^{-1})$};
		\node [style=none] (24) at (-.75, 0) {$\xrightarrow{\rm compose}$};
	\end{pgfonlayer}
	\begin{pgfonlayer}{edgelayer}
		\draw [bend right=90, looseness=1.25,thick] (0.center) to (4.center);
		\draw [bend right=105, looseness=1.25,thick,dashed] (2.center) to (3.center);
		\draw [bend left=90, looseness=0.75,thick] (2.center) to (3.center);
		\draw [bend left=90,thick] (0.center) to (4.center);
		\draw [bend right=90, looseness=1.50,thick] (2.center) to (6.center);
		\draw [bend right=90, looseness=1.50,thick] (3.center) to (5.center);
		\draw [bend left=90,thick] (5.center) to (6.center);
		\draw [bend right=90,thick,dashed] (5.center) to (6.center);
		\draw [bend right=90, looseness=1.25,thick] (13.center) to (14.center);
		\draw [bend left=90,thick] (13.center) to (14.center);
		\draw [bend right=90, looseness=1.75,thick] (13.center) to (4.center);
		\draw [bend right=90, looseness=1.50,thick] (14.center) to (0.center);
		\draw [bend left=270, looseness=0.75,thick] (15.center) to (16.center);
		\draw [bend left=270,thick,dashed] (17.center) to (18.center);
		\draw [bend left=270,thick,dashed] (19.center) to (20.center);
		\draw[thick] (15.center) to (19.center);
		\draw[thick] (16.center) to (20.center);
		\draw [bend left=90,thick] (15.center) to (16.center);
		\draw [bend left=90, looseness=0.75,thick] (17.center) to (18.center);
		\draw [bend left=90, looseness=0.75,thick] (19.center) to (20.center);
	\end{pgfonlayer}
		\path (0,0) coordinate (basepoint);
\end{tikzpicture}
\eeq
where the identity bordisms have been omitted from the picture, and  $(\ell,\tau,\bar\tau,\eta),(\ell,\tau',\bar\tau',\eta')\in ((\R_{>0}\times \HH^{2|1}_\pm)/\Z)(S)$.
\end{rmk}

\begin{prop}\label{prop:tracerelations}
There is 2-commutative diagrams of stacks
\beq
&&\resizebox{.95\textwidth}{!}{$\begin{tikzpicture}[baseline=(basepoint)];
\node (A) at (0,0) {$(\Map(\R^{2|1}/\Z,M)\times (\R_{>0}\times \HH^{2|1}_-\times \iHH^{2|1}_+)/\Z)^\Z$};
\node (B) at (9,0) {$\Mor(2|1\EBord(M))^{[2]}$};
\node (C) at (0,-1.5) {$\mathcal{L}^{2|1}(M)$};
\node (D) at (9,-1.5) {$\Mor(2|1\EBord(M))$};
\node (E) at (3.25,-.75) {$\twocommute$};
\draw[->] (A) to node [above] {$\sL^\mp_M\circ p_1\times \sR^\pm_M\circ p_2$} (B);
\draw[->] (A) to (C);
\draw[->] (C) to (D);
\draw[->] (B) to node [right] {$\c$} (D);
\path (0,-.75) coordinate (basepoint);
\end{tikzpicture}$} \label{eq:tracerelations}
\eeq
where the $\Z$-fixed subspace has $S$-points
\beq\label{eq:periodic1}
&&\Phi\circ \Tran_{\tau,\bar\tau,\eta}^{-1}=\Phi\circ \Tran_{\tau',\bar\tau',\eta}^{-1} \iff \Phi=\Phi\circ (L_{(-\tau',-\bar\tau',-\eta')\cdot (\tau,\bar\tau,\eta)}^{-1})^{n},\quad n\in \Z
\eeq
the top vertical arrow is determined by the families $\sL^\pm_M$ and $\sR^\mp_M$ and the maps
\beq\label{eq:targetofdual}
&&\begin{array}{c} p_1(\Phi,\ell,(\tau,\bar\tau,\eta),(\tau',\bar\tau',\eta))=(\Phi,\ell,\tau,\bar\tau,\eta), \\ p_2(\Phi,\ell,(\tau,\bar\tau,\eta),(\tau',\bar\tau',\eta'))=(\Phi,\ell,\tau',\bar\tau',\eta')\end{array}
\eeq
and the left vertical arrow is determined by
\beq\label{eq:sfm3}
&&(\Phi,\ell,(\tau,\bar\tau,\eta),(\tau',\bar\tau',\eta'))\mapsto (\Phi,(\ell,\ell,0),(\tau',\bar\tau',\eta')^{\pm 1}\cdot(\tau,\bar\tau,\eta)^{\mp 1})\in  \mathcal{L}^{2|1}(M)(S)
\eeq
where $((\ell,\ell,0),(\tau',\bar\tau',\eta')^{\pm 1}\cdot(\tau,\bar\tau,\eta)^{\mp 1})\in s\Lat(S)$, and the map $\Phi$ is well-defined on the corresponding super torus because of condition~\eqref{eq:periodic1}. 
\end{prop}
\bp
The compatibility $\Phi\circ \Tran_{\tau,\bar\tau,\eta}^{-1}=\Phi\circ \Tran_{\tau',\bar\tau',\eta'}^{-1}$ in~\eqref{eq:periodic1} guarantees that the target of $(\rR^\pm_{\ell,\tau',\bar\tau',\eta'},\Phi)$ agrees with the source of~$(\rL^\mp_{\ell,\tau,\bar\tau,\eta},\Phi)$ so that the bordisms are indeed composable. By the definitions of $(\rR^\pm_{\ell,\tau',\bar\tau',\eta'},\Phi)$ and~$(\rL^\mp_{\ell,\tau,\bar\tau,\eta},\Phi)$, the composition is necessarily a bordism from the empty set to itself. Unpacking the definition of the functor $\c$ and examining the periodicity condition~\eqref{eq:periodic1}, one finds that this composition is the Euclidean supertorus in~$M$ associated to the lattice $((\ell,\ell,0),(\tau',\bar\tau',\eta')^{\pm 1}\cdot(\tau,\bar\tau,\eta)^{\mp 1})$ and map to~$M$ determined by~$\Phi$. 
\ep
\begin{rmk}  The picture corresponding to the commutativity of~\eqref{eq:tracerelations} is 
\beq\label{pic:pic4}
\begin{tikzpicture}[baseline=(basepoint)];
	\begin{pgfonlayer}{nodelayer}
		\node [style=none] (1) at (3.5, 0) {};
		\node [style=none] (2) at (3, 0) {};
		\node [style=none] (4) at (1.5, 0) {};
		\node [style=none] (5) at (1, 0) {};
		\node [style=none] (6) at (-0.5, 0.25) {};
		\node [style=none] (7) at (-1, 0.25) {};
		\node [style=none] (8) at (-2.5, 0.25) {};
		\node [style=none] (9) at (-3, 0.25) {};
		\node [style=none] (10) at (-0.5, -0.25) {};
		\node [style=none] (11) at (-1, -0.25) {};
		\node [style=none] (12) at (-2.5, -0.25) {};
		\node [style=none] (13) at (-3, -0.25) {};
		\node [style=none] (14) at (0.25,0) {$\xrightarrow{\rm compose}$};
		\node [style=none] (15) at (-4,-.75) {$(\rR^\pm_{\ell,\tau',\bar \tau',\eta'},\Phi)$};
				\node [style=none] (15) at (-4,.75) {$(\sL^\mp_{\ell,\tau,\bar \tau,\eta},\Phi)$};
								\node [style=none] (15) at (4,-1.5) {$(\rT_{(\ell,\ell,0), (\tau',\bar\tau',\eta')^{\pm 1}\cdot (\tau,\bar\tau,\eta)^{\mp 1}},\Phi)$};
	\end{pgfonlayer}
	\begin{pgfonlayer}{edgelayer}
		\draw [bend right=105, looseness=1.25,thick] (1.center) to (2.center);
		\draw [bend left=90, looseness=0.75,thick,dashed] (1.center) to (2.center);
		\draw [bend right=90, looseness=1.50,thick] (1.center) to (5.center);
		\draw [bend right=90, looseness=1.50,thick] (2.center) to (4.center);
		\draw [bend left=90,thick,dashed] (4.center) to (5.center);
		\draw [bend right=90,thick] (4.center) to (5.center);
		\draw [bend right=90, looseness=1.50,thick] (5.center) to (1.center);
		\draw [bend right=90, looseness=1.50,thick] (4.center) to (2.center);
		\draw [bend right=105, looseness=1.25,thick] (6.center) to (7.center);
		\draw [bend left=90, looseness=0.75,thick] (6.center) to (7.center);
		\draw [bend right=90, looseness=1.50,thick] (6.center) to (9.center);
		\draw [bend right=90, looseness=1.50,thick] (7.center) to (8.center);
		\draw [bend left=90,thick] (8.center) to (9.center);
		\draw [bend right=90,thick] (8.center) to (9.center);
		\draw [bend right=90, looseness=1.25,thick] (10.center) to (11.center);
		\draw [bend left=90,thick,dashed] (10.center) to (11.center);
		\draw [bend right=90, looseness=1.25,thick] (12.center) to (13.center);
		\draw [bend left=90,thick,dashed] (12.center) to (13.center);
		\draw [bend right=90, looseness=1.50,thick] (12.center) to (11.center);
		\draw [bend right=90, looseness=1.50,thick] (13.center) to (10.center);
	\end{pgfonlayer}
		\path (0,0) coordinate (basepoint);
\end{tikzpicture}
\eeq
where $(\ell,\tau,\bar\tau,\eta)\in ((\R_{>0}\times \HH^{2|1}_\mp)/\Z)(S)$ and $(\ell,\tau',\bar\tau',\eta'))\in ((\R_{>0}\times \HH^{2|1}_\pm)/\Z)(S)$.

\end{rmk}

\subsection{Additional relations from isometries of bordisms}

Define the map
\beq
\inv\colon \Map(\R^{2|1}/\Z,M)\times (\R_{>0}\times \HH^{2|1}_\pm)/\Z&\to& \Map(\R^{2|1}/\Z,M)\times (\R_{>0}\times \HH^{2|1}_\mp)/\Z,\nonumber\\
(\Phi,\ell,\tau,\bar\tau,\eta)&\mapsto& (\Phi\circ \Tran_{\tau,\bar\tau,\eta}^{-1},\ell,-\tau,-\bar\tau,-\eta).\label{eq:2morphsigma}
\eeq

\begin{lem} \label{lem:invflip}
There are equalities of maps of sheaves
$$
\act_\pm=\proj_\mp\circ \inv,\qquad \proj_\pm=\act_\mp \circ \inv.
$$. 
\end{lem}
\bp
This follows from~\eqref{eq:2morphsigma} and the definitions of $\proj_\pm$ and $\act_\pm$,
$$
\act_\pm(\Phi,\ell,\tau,\bar\tau,\eta)=(\Phi\circ \Tran_{\tau,\bar\tau,\eta}^{-1},\ell),\qquad \proj_\pm(\Phi,\ell,\tau,\bar\tau,\eta)=(\Phi,\ell)
$$
determined by restriction of the map~\eqref{eq:actproj}. 
\ep

\begin{prop}\label{prop:symmetry}
Let $\sL^\pm_M\circ\sigma$ denote the functor $\sL^\pm_M$ post-composed with the isomorphism in $\Mor(2|1\EBord(M))$ that exchanges the order of components of the target, i.e.,
$$
\t(\rL^\pm_{\ell,\tau,\bar\tau,\eta}\circ\sigma,\Phi)=(\rS_\ell^\mp,\Phi\circ \Tran_{\tau,\bar\tau,\eta}^{-1})\coprod (\rS_\ell^\pm,\Phi).
$$
There is a 2-commutative triangle
\beq
&&\begin{tikzpicture}[baseline=(basepoint)];
\node (A) at (0,0) {$\Map(\R^{2|1}/\Z,M)\times (\R_{>0}\times \HH^{2|1}_\pm)/\Z$};
\node (B) at (7,-.75) {$\Mor(2|1\EBord(M))$};
\node (C) at (0,-1.5) {$\Map(\R^{2|1}/\Z,M)\times (\R_{>0}\times \HH^{2|1}_\mp)/\Z$};
\node (E) at (3.25,-.75) {$\twocommute$};
\draw[->] (A) to node [above] {$\sL^\pm_M\circ \sigma$} (B);
\draw[->] (A) to node [left] {$\inv$} (C);
\draw[->] (C) to node [below] {$\sL^\mp_M$} (B);
\path (0,-.75) coordinate (basepoint);
\end{tikzpicture} \label{eq:symmetrybordism}
\eeq
for the left vertical arrow~\eqref{eq:2morphsigma}.
\end{prop} 

\bp
Starting with the data that determines $(\rL_{\ell,\tau,\bar\tau,\eta}^\pm,\Phi)$, 
\beq\nonumber
(S\times \R^{1|1})/\ell\Z\coprod (S\times \R^{1|1})/\ell\Z \xhookrightarrow{i_0 \coprod i_{\tau,\bar\tau,\eta}} (S\times \R^{2|1})/\ell\Z \xrightarrow{\Phi} M.
\eeq
exchanging the factors of the coproduct and acting on $(S\times \R^{2|1})/\ell\Z$ by $\Tran_{\tau,\bar\tau,\eta}^{-1}=\Tran_{-\tau,-\bar\tau,-\eta}$ provides the data determining $(\rL^\mp_{\ell,-\tau,-\bar\tau,-\eta},\Phi\circ \Tran_{\tau,\bar\tau,\eta}^{-1})$.
\ep

Below we will use (degenerate) notation $\sqrt{\Fl}$ to denote the maps of sheaves
\beq
\Map(\R^{2|1}/\Z,M)\times \R_{>0}&\xrightarrow{\sqrt{\Fl}} &\Map(\R^{2|1}/\Z,M)\times \R_{>0} \label{eq:Flonsheaves}\\
(\Phi,\ell)&\mapsto& (\Phi\circ \sqrt{\fl}{}^{-1},\ell).\nonumber\\
\Map(\R^{2|1}/\Z,M)\times (\R_{>0}\times \HH^{2|1}_\pm)/\Z&\xrightarrow{\sqrt{\Fl}} &\Map(\R^{2|1}/\Z,M)\times (\R_{>0}\times \HH^{2|1}_\pm)/\Z\label{eq:Flonsheaves2}\\
(\Phi,\ell,\tau,\bar\tau,\eta)&\mapsto& (\Phi\circ \sqrt{\fl}{}^{-1},\ell,-\tau,-\bar\tau,-i\eta).\nonumber\\
\Map(\R^{2|1}/\Z,M)\times (\R_{>0}\times \iHH^{2|1}_\pm)/\Z&\xrightarrow{\sqrt{\Fl}} &\Map(\R^{2|1}/\Z,M)\times (\R_{>0}\times \iHH^{2|1}_\pm)/\Z\label{eq:Flonsheaves3}\\
(\Phi,\ell,\tau,\bar\tau,\eta)&\mapsto& (\Phi\circ \sqrt{\fl}{}^{-1},\ell,-\tau,-\bar\tau,-i\eta).\nonumber
\eeq
where the particular map will be clear in context. All of these maps are induced by the isometry Euclidean supercylinders $\sqrt{\fl}$ defined in~\eqref{sec:sqrtfldef}. 

\begin{lem}\label{lem:orientation11restrict}
There are 2-commutative triangles
\beq
\begin{tikzpicture}[baseline=(basepoint)];
\node (A) at (0,0) {$\Map(\R^{2|1}/\Z,M)\times \R_{>0}$};
\node (B) at (7,-.75) {$\Ob(2|1\EBord(M))$};
\node (C) at (0,-1.5) {$\Map(\R^{2|1}/\Z,M)\times \R_{>0}$};
\node (E) at (3.25,-.75) {$\twocommute$};
\draw[->] (A) to node [above] {$\sS^\pm_M$} (B);
\draw[->] (A) to node [left] {$\sqrt{\Fl}$} (C);
\draw[->] (C) to node [below] {$\sS^\mp_M$} (B);
\path (0,-.75) coordinate (basepoint);
\end{tikzpicture} \label{eq:sqrtfl0}
\eeq
\beq
\begin{tikzpicture}[baseline=(basepoint)];
\node (A) at (0,0) {$\Map(\R^{2|1}/\Z,M)\times (\R_{>0}\times \HH^{2|1}_\pm)/\Z$};
\node (B) at (7,-.75) {$\Mor(2|1\EBord(M))$};
\node (C) at (0,-1.5) {$\Map(\R^{2|1}/\Z,M)\times (\R_{>0}\times \HH^{2|1}_\mp)/\Z$};
\node (E) at (3.25,-.75) {$\twocommute$};
\draw[->] (A) to node [above] {$\sC^\pm_M,\sL^\pm_M$} (B);
\draw[->] (A) to node [left] {$\sqrt{\Fl}$} (C);
\draw[->] (C) to node [below] {$\sC^\mp_M,\sL^\mp_M$} (B);
\path (0,-.75) coordinate (basepoint);
\end{tikzpicture} \label{eq:sqrtfl}
\eeq
\beq
\begin{tikzpicture}[baseline=(basepoint)];
\node (A) at (0,0) {$\Map(\R^{2|1}/\Z,M)\times (\R_{>0}\times \iHH^{2|1}_\pm)/\Z$};
\node (B) at (7,-.75) {$\Mor(2|1\EBord(M)).$};
\node (C) at (0,-1.5) {$\Map(\R^{2|1}/\Z,M)\times (\R_{>0}\times \iHH^{2|1}_\mp)/\Z$};
\node (E) at (3.25,-.75) {$\twocommute$};
\draw[->] (A) to node [above] {$\sR^\pm_M$} (B);
\draw[->] (A) to node [left] {$\sqrt{\Fl}$} (C);
\draw[->] (C) to node [below] {$\sR^\mp_M$} (B);
\path (0,-.75) coordinate (basepoint);
\end{tikzpicture} \nonumber
\eeq
\end{lem}
\bp
The 2-commuting data is given by the isometries in $\Ob(2|1\EBord(M))$ and $\Mor(2|1\EBord(M))$ determined by $\sqrt{\fl}$ defined in~\eqref{sec:sqrtfldef}. Specifically, the isometries of $S$-families
\beq
\sqrt{\Fl}(\rS_\ell^\pm,\Phi)&\simeq&(\rS_\ell^\mp,\Phi\circ \sqrt{\fl}{}^{-1})\nonumber\\
\sqrt{\Fl}(\rC^\pm_{\ell,\tau,\bar\tau,\eta},\Phi)&\simeq&(C^\mp_{\ell,-\tau,-\bar\tau,\sqrt{-1}\eta},\Phi\circ \sqrt{\fl}{}^{-1}) \nonumber\\
\sqrt{\Fl}(\rL^\pm_{\ell,\tau,\bar\tau,\eta},\Phi)&\simeq&(L^\mp_{\ell,-\tau,-\bar\tau,\sqrt{-1}\eta},\Phi\circ \sqrt{\fl}{}^{-1}) \nonumber\\
\sqrt{\Fl}(\rR^\pm_{\ell,\tau,\bar\tau,\eta},\Phi)&\simeq&(R^\mp_{\ell,-\tau,-\bar\tau,\sqrt{-1}\eta},\Phi\circ \sqrt{\fl}{}^{-1}). \nonumber
\eeq
Indeed, the isometry $\sqrt{\Fl}$ evaluated on any of the above families is determined by the specialization of diagram~\eqref{diag:sqrtflonobject} with $(\tau_\inn,\bar\tau_\inn,\eta_\inn)=0$. 
\ep

\subsection{Restriction of additional structures}

Analogously to~\eqref{eq:Flonsheaves}-\eqref{eq:Flonsheaves3}, we will use the notation for maps of sheaves
\beq
\dagger \colon \Map(\R^{2|1}/\Z,M)\times \R_{>0}&\to& \overline{\Map(\R^{2|1}/\Z,M)}\times \overline{\R_{>0}},\nonumber \\
(\Phi\colon (S\times \R^{2|1})/\ell \Z \to M)&\mapsto &(\overline{S}\times \overline{\R}^{2|1})/\overline{\ell}\Z \xrightarrow{\id_{\overline{S}}\times \overline{\orr}^{-1}} (\overline{S}\times \overline{\R}^{2|1})/\overline{\ell}\Z \xrightarrow{\overline{\Phi}} \overline{M}\simeq M).\label{eq:daggersheaf1}\\
\dagger \colon \Map(\R^{2|1}/\Z,M)\times (\R_{>0}\times \HH^{2|1}_\pm)/\Z &\to& \overline{\Map(\R^{2|1}/\Z,M)}\times (\overline{\R_{>0}}\times \overline{\HH}^{2|1}_\pm)/\Z \nonumber \\
(\Phi,\ell,(\tau,\bar\tau,\eta))&\mapsto &(\overline{\Phi}\circ \overline{\orr}^{-1},\overline{\ell},\overline{(\bar\tau,\tau,\eta)})\label{eq:daggersheaf2}\\
\dagger \colon \Map(\R^{2|1}/\Z,M)\times (\R_{>0}\times \iHH^{2|1}_\pm)/\Z &\to& \overline{\Map(\R^{2|1}/\Z,M)}\times (\overline{\R_{>0}}\times \overline{\iHH}^{2|1}_\pm)/\Z \nonumber \\
(\Phi,\ell,(\tau,\bar\tau,\eta))&\mapsto &(\overline{\Phi}\circ \overline{\orr}^{-1},\overline{\ell},\overline{(\bar\tau,\tau,\eta)})\label{eq:daggersheaf3}
\eeq
that cover complex conjugation on supermanifolds, where in this case all the above maps are determined by the map~\eqref{eq:daggerSpts} between families of Eulidean supercylinders. The reflection structure on $2|1\EBord(M)$ is constructed in Lemma~\ref{lem:extrastrEB}. 

\begin{prop}\label{prop:restrictRR}
The restriction of the reflection structure~\eqref{eq:21EBreflection} on $2|1\EBord(M)$ yields 2-commuting squares 
\beq
&&\begin{tikzpicture}[baseline=(basepoint)];
\node (A) at (0,0) {$\Map(\R^{2|1}/\Z,M)\times \R_{>0}$};
\node (B) at (5,0) {$\Ob(2|1\EBord(M))$};
\node (C) at (0,-1.5) {$\overline{\Map(\R^{2|1}/\Z,M)}\times \overline{\R_{>0}}$};
\node (D) at (5,-1.5) {$\Ob(2|1\EBord(M))$};
\draw[->] (A) to node [above] {$\sS^\pm_M$} (B);
\draw[->] (A) to node [left] {$\dagger$} (C);
\draw[->] (C) to node [below] {$\overline{\sS}^\mp_M$}  (D);
\draw[->] (B) to node [right] {$\dagger_0$} (D);
\node (E) at (2.5,-.75) {$\twocommute$};
\path (0,-.75) coordinate (basepoint);
\end{tikzpicture}\eeq
\beq
&&\begin{tikzpicture}[baseline=(basepoint)];
\node (A) at (0,0) {$\Map(\R^{2|1}/\Z,M)\times (\R_{>0}\times \HH^{2|1}_\pm)/\Z$};
\node (B) at (6,0) {$\Mor(2|1\EBord(M))$};
\node (C) at (0,-1.5) {$\overline{\Map(\R^{2|1}/\Z,M)}\times (\overline{\R_{>0}}\times \overline{\HH}^{2|1}_\pm)/\Z$};
\node (D) at (6,-1.5) {$\Mor(2|1\EBord(M))$};
\draw[->] (A) to node [above] {$\sC^\pm_M$} (B);
\draw[->] (A) to node [left] {$\dagger$} (C);
\draw[->] (C) to node [below] {$\overline{\sC}^\mp_M$}  (D);
\draw[->] (B) to node [right] {$\dagger_1$} (D);
\node (E) at (1.5,-.75) {$\twocommute$};
\path (0,-.75) coordinate (basepoint);
\end{tikzpicture}\eeq
\beq
&&\begin{tikzpicture}[baseline=(basepoint)];
\node (A) at (0,0) {$\Map(\R^{2|1}/\Z,M)\times (\R_{>0}\times \HH^{2|1}_\pm)/\Z$};
\node (B) at (6,0) {$\Mor(2|1\EBord(M))$};
\node (C) at (0,-1.5) {$\overline{\Map(\R^{2|1}/\Z,M)}\times (\overline{\R_{>0}}\times \overline{\HH}^{2|1}_\pm)/\Z$};
\node (D) at (6,-1.5) {$\Mor(2|1\EBord(M))$};
\draw[->] (A) to node [above] {$\sL^\pm_M$} (B);
\draw[->] (A) to node [left] {$\dagger$} (C);
\draw[->] (C) to node [below] {$\overline{\sL}^\mp_M$}  (D);
\draw[->] (B) to node [right] {$\dagger_1$} (D);
\node (E) at (1.5,-.75) {$\twocommute$};
\path (0,-.75) coordinate (basepoint);
\end{tikzpicture}\eeq
\beq
&&\begin{tikzpicture}[baseline=(basepoint)];
\node (A) at (0,0) {$\Map(\R^{2|1}/\Z,M)\times (\R_{>0}\times \iHH^{2|1}_\pm)/\Z$};
\node (B) at (6,0) {$\Mor(2|1\EBord(M))$};
\node (C) at (0,-1.5) {$\overline{\Map(\R^{2|1}/\Z,M)}\times (\overline{\R_{>0}}\times \overline{\iHH}^{2|1}_\pm)/\Z$};
\node (D) at (6,-1.5) {$\Mor(2|1\EBord(M)).$};
\draw[->] (A) to node [above] {$\sR^\pm_M$} (B);
\draw[->] (A) to node [left] {$\dagger$} (C);
\draw[->] (C) to node [below] {$\overline{\sR}^\mp_M$}  (D);
\draw[->] (B) to node [right] {$\dagger_1$} (D);
\node (E) at (1.5,-.75) {$\twocommute$};
\path (0,-.75) coordinate (basepoint);
\end{tikzpicture}\nonumber
\eeq

\end{prop}
\bp 
The 2-commuting data is given by the isometries in $\Ob(2|1\EBord(M))$ and $\Mor(2|1\EBord(M))$ 
\beq\label{eq:someisometriesconj}
\begin{array}{ccc}
\dagger(\rS_\ell^\pm,\Phi)&\simeq&(\rS_{\overline{\ell}}^\mp,\overline{\Phi\circ \orr^{-1}})\\
\dagger(\rC^\pm_{\ell,\tau,\bar\tau,\eta},\Phi)&\simeq&(C^\mp_{(\overline{\ell},\overline{(\bar\tau,\tau,\eta)})} ,\overline{\Phi\circ \orr^{-1}}) \\
\dagger(\rL^\pm_{\ell,\tau,\bar\tau,\eta},\Phi)&\simeq&(L^\mp_{\overline{\ell},\overline({\bar\tau,\tau,\eta})},\overline{\Phi\circ \orr^{-1}}) \\
\dagger(\rR^\pm_{\ell,\tau,\bar\tau,\eta},\Phi)&\simeq&(R^\mp_{\ell,\overline{(\bar\tau,\tau,\eta)}},\overline{\Phi\circ \orr^{-1}} ) \end{array}
\eeq
from specialization of the diagram~\eqref{diag:orRonobject}
\beq\nonumber
\begin{tikzpicture}[baseline=(basepoint)];
\node (BB) at (2.5,1.5) {$\rC_{\overline{\ell}}$};
\node (A) at (-1,0) {$\rS_{\overline{\ell}}$};
\node (AA) at (-1,1.5) {$\rS_{\overline{\ell}}$};
\node (B) at (2.5,0) {$\rC_{\overline{\ell}}$};
\node (C) at (6,0) {$\rS_{\overline{\ell}}$};
\node (CC) at (6,1.5) {$\rS_{\overline{\ell}}$};
\node (D) at (2.5,-1.25) {$M$};
\draw[->,right hook-latex] (A) to node [below] {$\overline{i_{0}}$} (B);
\draw[->,left hook-latex] (C) to node [below] {$\overline{i_{\tau,\bar\tau,\eta}}$} (B);
\draw[->,right hook-latex] (AA) to node [above=2pt] {$i_{\overline{0}}$} (BB);
\draw[->,left hook-latex] (CC) to node [above=2pt] {$i_{\overline{{\bar\tau,\tau,\eta}}}$} (BB);
\draw[->] (B) to node [right] {$\overline{\Phi}$} (D);
\draw[->] (B) to node [left] {$\overline{\orr}$} (BB);
\draw[->] (A) to node [left] {$\overline{\orr}$} (AA);
\draw[->] (C) to node [left] {$\overline{\orr}$} (CC);
\path (0,0) coordinate (basepoint);
\end{tikzpicture}
\eeq
where (in an abuse of notation) $\overline{\Phi}\colon \rC_\ell \to M$ is the map conjugate to $\Phi$ post-composed with the canonical real structure $\overline{M}\simeq M$ on the ordinary manifold $M$. To complete the argument, we note that $\orr$ exchanges the subspaces $\iHH^{2|1}_\pm\subset \R^{2|1}$, and hence the collar data for a bordism is exchanged under~\eqref{eq:someisometriesconj}. 
\ep

\begin{prop} 
Relative to the values in Proposition~\ref{prop:restrictRR}, the reflection structure $\dagger$ sends the isometries in Lemma~\ref{lem:orientation11restrict} to their inverses. 
\end{prop} 
\bp
This follows from Lemma~\ref{lem:commuteuptoflip}: the reflection functor commutes with the isometry $\sqrt{\fl}$ up to the spin flip, and $\sqrt{\fl}\circ \fl=\sqrt{\fl}{}^{-1}$. 
\ep

The renormalization group action on $2|1\EBord(M)$ is constructed in Lemma~\ref{lem:extrastrEB0}.

\begin{prop}\label{prop:RG11restrict}
For $\mu\in \R_{>0}$, the renormalization group action on $2|1\EBord(M)$ determines automorphisms of the families of bordisms $\sS^\pm_M$, $\sC^\pm_M$, $\sL^\pm_M$ and $\rR^\pm_M$ given by maps on $S$-points, 
\beq
\RG_\mu(\rS_\ell^\pm,\Phi)&\simeq&(\rS_{\mu^2\ell}^\pm,\Phi\circ \rg_\mu^{-1})\nonumber\\
\RG_\mu(\rC^\pm_{\ell,\tau,\bar\tau,\eta},\Phi)&\simeq&(\rC^\pm_{\mu^2\ell,\mu^2\tau,\mu^2\bar\tau,\mu\theta},\Phi\circ \rg_\mu^{-1}) \nonumber\\
\RG_\mu(\rL^\pm_{\ell,\tau,\bar\tau,\eta},\Phi)&\simeq&(\rL^\pm_{\mu^2\ell,\mu^2\tau,\mu^2\bar\tau,\mu\theta},\Phi\circ \rg_\mu^{-1}) \nonumber\\
\RG_\mu(\rR^\pm_{\ell,\tau,\bar\tau,\eta},\Phi)&\simeq&(\rR^\pm_{\mu^2\ell,\mu^2\tau,\mu^2\bar\tau,\mu\theta},\Phi\circ \rg_\mu^{-1}) \nonumber
\eeq
where $\rg_\mu\colon (S\times \R^{2|1})/\ell\Z \to (S\times \R^{2|1})/\mu^2\ell\Z$ is determined by the map~\eqref{Eq:21RGdefn}.
\end{prop}

\bp
The argument is analogous to the proof of Proposition~\ref{prop:restrictRR}, using a specialization of the diagram~\eqref{diag:rgonobject}. We note that $\rg_\mu$ preserves the subspaces $\iHH_\pm^{2|1}\subset \R^{2|1}$ and so the partition of the collar is preserved.
\ep

\subsection{Small bordisms and nearly constant supercylinders}\label{sec:smallbordisms}

\begin{defn}\label{defn:smallbordisms}
Consider the families of bordisms associated to the maps of stacks
\beq\label{eq:LM0}
\begin{array}{c} 
{\sf s}^\pm_M\colon \Map(\R^{0|1},M)\times \R_{>0}\to \Ob(2|1\EBord(M))\\
{\sf c}^\pm_M\colon \Map(\R^{0|1},M)\times (\R_{>0}\times \HH^{2|1}_\pm)/\Z\to \Mor(2|1\EBord(M))\\ 
{\sf l}^\pm_M\colon \Map(\R^{0|1},M)\times (\R_{>0}\times \HH^{2|1}_\pm)/\Z \to \Mor(2|1\EBord(M))\\ 
{\sf r}^\pm_M\colon \Map(\R^{0|1},M)\times (\R_{>0}\times \iHH^{2|1}_\pm)/\Z \to \Mor(2|1\EBord(M))
\end{array}
\eeq
by restriction of~\eqref{eq:sSM},~\eqref{eq:CM} and~\eqref{eq:LM} to the (finite-dimensional) subspaces
\beq\label{eq:constantpathsubspace}
\Map(\R^{0|1},M)\times \R_{>0} &\hookrightarrow& \Map(\R^{2|1}/\Z,M)\times \R_{>0},\\
(\Phi_0,\ell)&\mapsto& (\Phi\colon S\times \R^{2|1}/\ell\Z\xrightarrow{p}S\times \R^{0|1}\xrightarrow{\Phi_0} M,\ell).\nonumber
\eeq
\beq\nonumber
\Map(\R^{0|1},M)\times (\R_{>0}\times \HH^{2|1}_\pm)/\Z &\hookrightarrow& \Map(\R^{2|1}/\Z,M)\times (\R_{>0}\times \HH^{2|1}_\pm)/\Z,\\
(\Phi_0,\ell,\tau,\bar\tau,\eta)&\mapsto& (\Phi\colon S\times \R^{2|1}/\ell\Z\xrightarrow{p}S\times \R^{0|1}\xrightarrow{\Phi_0} M,\ell,\tau,\bar\tau,\eta).\nonumber
\eeq
\beq\nonumber
\Map(\R^{0|1},M)\times (\R_{>0}\times \iHH^{2|1}_\pm)/\Z &\hookrightarrow& \Map(\R^{2|1}/\Z,M)\times (\R_{>0}\times \iHH^{2|1}_\pm)/\Z,\\
(\Phi_0,\ell,\tau,\bar\tau,\eta)&\mapsto& (\Phi\colon S\times \R^{2|1}/\ell\Z\xrightarrow{p}S\times \R^{0|1}\xrightarrow{\Phi_0} M,\ell,\tau,\bar\tau,\eta).\nonumber
\eeq

We refer to the restricted families~\eqref{eq:LM0} as \emph{small bordisms} in $M$.
\end{defn}

\begin{rmk}
Morally, composition and disjoint union of small bordisms generates a subcategory $2|1\ebord(M)\subset 2|1\EBord(M)$ of small bordisms in $M$, and $2|1\eft^n(M)$ in Theorem~\ref{thm1} considers degree~$n$ field theories defined on these subcategories. However, formulating the notion of generating a subcategory in a category internal to symmetric monoidal stacks gets quite technical. 
\end{rmk}

In particular, the families~\eqref{eq:LM0} and Lemma~\ref{lem:orientation11restrict} determines maps of stacks
\beq\label{eq:supercyltoEB}
&&\Ob(\Ann^{2|1}(M))\to \Ob(2|1\EBord(M)),\quad \Mor(\Ann^{2|1}(M))\to \Mor(2|1\EBord(M))
\eeq
for the object and morphism stacks in Definition~\ref{defn:sC}.

\begin{prop}\label{prop:FTrestrict1}
The maps of stacks~\eqref{eq:supercyltoEB} are part of the data of an internal functor~\eqref{eq:maininclude} from the internal category of nearly constant supercylinders in~$M$ to the $2|1$-Euclidean bordism category over $M$. Under this functor, the restriction of the reflection $\dagger$ and RG-action $\RG_\mu$ on $2|1\EBord(M)$ (see Lemmas~\ref{lem:extrastrEB0} and~\ref{lem:extrastrEB}) are compatible with the structures on $\Ann^{2|1}(M)$ in Lemma~\ref{lem:structureonsP}.
\end{prop}

\bp
By Definition~\ref{defn:internalfunctor} of an internal functor, in addition to the maps of stacks~\eqref{eq:supercyltoEB} one requires compatibility data for composition and units. This is given by restriction of the 2-morphisms in Proposition~\ref{prop:stI}. The 2-morphisms in this proposition are canonical (coming from intersecting collars of bordisms), implying that the coherence property for these 2-isomorphisms is automatic. This completes the construction of the internal functor~\eqref{eq:maininclude}.

The compatibility of the functor~\eqref{eq:maininclude} with the additional structures follows from Propositions~\ref{prop:restrictRR} and~\ref{prop:RG11restrict}: the isomorphisms in these propositions restricted along~\eqref{eq:LM0} supply 2-commutativity data, 
\beq
\label{eq:RGOrRR}
&&\begin{tikzpicture}[baseline=(basepoint)];
\node (A) at (0,0) {$\Ann^{2|1}(M)$};
\node (B) at (3,0) {$2|1\EBord(M)$};
\node (C) at (0,-1.5) {$\Ann^{2|1}(M)$};
\node (D) at (3,-1.5) {$2|1\EBord(M)$};
\draw[->] (A) to (B);
\draw[->] (A) to node [left] {$\RG_\mu$} (C);
\draw[->] (C) to  (D);
\draw[->] (B) to node [right] {$\RG_\mu$} (D);
\node (E) at (1.5,-.75) {$\twocommute$};
\path (0,-.75) coordinate (basepoint);
\end{tikzpicture}\quad 
\begin{tikzpicture}[baseline=(basepoint)];
\node (A) at (0,0) {$\Ann^{2|1}(M)$};
\node (B) at (3,0) {$2|1\EBord(M)$};
\node (C) at (0,-1.5) {$\Ann^{2|1}(M)$};
\node (D) at (3,-1.5) {$2|1\EBord(M)$};
\draw[->] (A) to (B);
\draw[->] (A) to node [left] {$\dagger$} (C);
\draw[->] (C) to  (D);
\draw[->] (B) to node [right] {$\dagger$} (D);
\node (E) at (1.5,-.75) {$\twocommute$};
\path (0,-.75) coordinate (basepoint);
\end{tikzpicture}
\eeq
witnessing the claimed compatibility. 
\ep

\section{$2|1$-Euclidean field theories of degree~$n$}\label{sec:STmess}
In this section we compute the values of twisted $2|1$-Euclidean field theories over~$M$ on the bordisms constructed in the previous section. Following Remark~\ref{rmk:ultimate}, we expect this to characterize the ``half" of the data of a twisted $2|1$-Euclidean field theory, i.e., the Ramond sector values. 
We then specialize to the degree~$n$ twist relevant to conjecture~\eqref{eq:conjecture}. A complete construction of the degree~$n$ twist  has yet to appear; it is expected to be the unique supersymmetric extension of the 2-dimensional chiral free fermions constructed in~\cite[\S5-6]{ST11}, see Hypothesis~\ref{hyp:chiral}. Below we show that such a supersymmetric extension exists and is unique on the Ramond sector bordisms from the previous section. This partial construction suffices for our purposes. In Remark~\ref{rmk:generalproof}, we indicate a full construction of~$\Fer_n$ that depends upon a generators and relations presentation of $2|1\EBord(\pt)$ suggested in~\cite{ST11}. 

\subsection{The values of twisted field theories}
A \emph{twist} for the $2|1$-Euclidean bordism category over~$M$ is an internal functor $\twist\colon 2|1\EBord(M)\to \TA$ valued in the Morita category. A \emph{twisted $2|1$-Euclidean field theory over $M$} is an internal natural transformation~$E$
\beq
&&\begin{tikzpicture}[baseline=(basepoint)];
\node (A) at (0,0) {$2|1\EBord(M)$};
\node (B) at (5,0) {$\TA.$};
\node (C) at (2.5,0) {$E \Downarrow$};
\draw[->,bend left=15] (A) to node [above] {$\one$} (B);
\draw[->,bend right=15] (A) to node [below] {$\twist$} (B);
\path (0,0) coordinate (basepoint);
\end{tikzpicture}\label{eq:gentwistedEFT}
\eeq
We refer to~\S\ref{sec:STappen} for a review of these definitions from~\cite{ST11}. In the following series of propositions, we describe the values of a twisted $2|1$-Euclidean field theory on the families of bordisms constructed in the previous section. 

\begin{prop}\label{prop:valueoftwist}
A twist $\twist\colon 2|1\EBord(M)\to \TA$ determines data
\begin{enumerate}
\item[i.] a bundle of algebras $\twist(\sS^\pm_M)$ over $\Map(\R^{2|1}/\Z,M)\times \R_{>0}$;
\item[ii.] bundles of bimodules 
$$
_{\act_\pm^*\twist(\sS^\pm_M)}\twist(\sC^\pm_M)_{\proj_\pm^*\twist(\sS^\pm_M)}, \quad _\C\twist(\sL^\pm_M)_{\proj_\pm^*\twist(\sS^\pm_M)\otimes \act_\pm^*\twist(\sS^\mp_M)}, \quad
$$
over $\Map(\R^{2|1}/\Z,M)\times (\R_{>0}\times \HH^{2|1}_\pm)/\Z$ and 
$$
_{\proj_\pm^*\twist(\sS^\mp_M)\otimes \act_\pm^*\twist(\sS^\pm_M)}\twist(\sR^\pm_M)_\C
$$
over $\Map(\R^{2|1}/\Z,M)\times (\R_{>0}\times \iHH^{2|1}_\pm)/\Z$,
\item[iii.] a $\E^{2|1}\rtimes \Spin(2)\times \SL_2(\Z)$-equivariant vector bundle $\twist(\mathcal{L}^{2|1}(M))$ over the sheaf~$\mathcal{L}^{2|1}(M)$ defined in~\eqref{Eq:defnofsuperdoubleloops}. 
\end{enumerate}
\end{prop}
\bp
By Definition~\ref{defn:internalfunctor} of an internal functor, the twist determines maps of stacks
$$
\Ob(2|1\EBord(M))\to \Ob(\TA),\qquad \Mor(2|1\EBord(M))\to \Mor(\TA)
$$
that are compatible with the source and target functors. Evaluating these on the indicated families of objects of the source stacks recovers the statement in the proposition. 
\ep

\begin{prop}\label{prop:dataoftwistedFT} Using the notation from Proposition~\ref{prop:valueoftwist}, a $\twist$-twisted field theory~$E$ determines data
\begin{enumerate}
\item[i.] bundles of left modules $_{\twist(\sS^\pm_M)}E(\sS^\pm_M)$ over $\Map(\R^{2|1}/\Z,M)\times \R_{>0}$,
\item[ii.] maps of bimodule bundles over $\Map(\R^{2|1}/\Z,M)\times (\R_{>0}\times \HH^{2|1}_\pm)/\Z$ 
\beq
&&E(\sC^\pm_M)\colon \twist(\sC^\pm_M)\otimes_{\proj_\pm^*\twist(\sS^\pm_M)} \proj_\pm^*E(\sS^\pm_M)\to \act_\pm^*E(\sS^\pm_M),\nonumber\\
&& E(\sL^\pm_M)\colon \twist(\sL^\pm_M)\otimes_{\proj_\pm^*\twist(\sS^\pm_M)\otimes \act_\pm^*\twist(\sS^\mp_M)} \proj_\pm^*E(\sS^\pm_M)\otimes \act_\pm^*E( \sS^\mp_M)\to \underline{\C} ,\nonumber
\eeq
and maps of bimodule bundles over $\Map(\R^{2|1}/\Z,M)\times (\R_{>0}\times \iHH^{2|1}_\pm)/\Z$ 
\beq
&&E(\sR^\pm_M)\colon \twist(\sR^\pm_M)\to \proj_\pm^*E(\sS^\mp_M) \otimes \act_\pm^*E(\sS^\pm_M), \nonumber
 \eeq
for the projection and action maps $\proj^\pm$ and $\act^\pm$ defined in~\eqref{eq:actproj},  and
\item[iii.] a $\E^{2|1}\rtimes \Spin(2)\times \SL_2(\Z)$-invariant map of equivariant vector bundles 
$$
E(\mathcal{L}^{2|1}(M))\colon\twist(\mathcal{L}^{2|1}(M)) \to \underline{\C}\nonumber
$$
 over the sheaf $\mathcal{L}^{2|1}(M)$, taking the trivial equivariant structure on~$\underline{\C}$.
\end{enumerate}
\end{prop}
\bp
By Definition~\ref{defn:internalnattrans} of an internal natural transformation, a twisted field theory supplies the data of a map of stacks $\eta\colon \Ob(2|1\EBord(M))\to \Mor(\TA)$, and hence evaluating on the families $E(\sS^\pm_M)$ determines bundles of bimodules over $\Map(\R^{2|1}/\Z,M)\times \R_{>0}$
\beq\nonumber
\begin{array}{c}
\nonumber E(\sS^\pm_M) \in \Mor(\TA), \quad {\rm with} \quad \t(E(\sS^\pm_M))= \twist(\sS^\pm_M), \ \s(E(\sS^\pm_M))= \underline{\C}.
\end{array}
\eeq
The right action by $\underline{\C}$ is no additional data, so $E(\sS^\pm_M)$ are simply bundles of left modules. An internal natural transformation further specifies a compatibility 2-morphism~$\rho$ in~\eqref{diag:nu} for every family in $\Mor(2|1\EBord(M))$; this yields the data in (ii) and (iii) by the argument from~\cite[page 52]{ST11}. For example, $\rho$ evaluated on the family $\sC^\pm_M$ is a map of left modules
$$
E(\sC^\pm_M)\colon \twist(\sC^\pm_M)\otimes_{\proj_\pm^*\twist(\sC^\pm_M)} \proj_\pm^*E(\sS^\pm_M)\to \act_\pm^*E(\sS^\pm_M)
$$
which is the first map of bimodules listed above. 
\ep

\begin{rmk}
By Remarks~\ref{rmk:germs1} and~\ref{rmk:germs2}, the algebra and bimodule bundles in Proposition~\ref{prop:valueoftwist} and~\ref{prop:dataoftwistedFT} satisfy an additional property: they are isomorphic to the pullback of a bundle along a quotient that identifies shrinking collars. For example, the algebra bundle $\twist(\sS^\pm_M)$ pulls back along a quotient map $\R_{>0}\times \Map(\R^{2|1}/\Z,M)\to \R_{>0}\times \Map(\R^{2|1}/\Z,M)/{\sim}$ where the target is the sheaf of Euclidean superloops in $M$ with the 1-sided germ of a supercylinder. 
\end{rmk}

In the following two propositions, we record the data and property that a twist and a twisted field theory assigns to the compositions of bordisms computed in~\S\ref{sec:compositionEucbord}. 

\begin{prop}
A twist $\twist$ determines data
\begin{enumerate}
\item[i.] a map of $\twist(\sS^\pm_M)$-bimodules,
\beq
&&p_2^* \twist(\sC^\pm_M)\otimes_{p_2^*\twist(\sS^\pm_M)} p_1^* \twist(\sC^\pm_M)\to  \m^*\twist(\sC^\pm_M),\label{eq:mapofbimodules1}\\ 
&&{\rm over}\quad  \Map(\R^{2|1}/\Z,M)\times (\R_{>0}\times \HH^{2|1}_\pm\times \HH^{2|1}_\pm)/\Z \nonumber
\eeq
where $p_1,p_2$ are the maps in~\eqref{eq:composablesuperpath} and $\m$ is~\eqref{eq:sfm};
\item[ii.] a map of right $\twist(\sS^\pm_M)\otimes \twist(\sS^\mp_M)$-module bundles
\beq\label{eq:twistadjunction}
&&p_2^* \twist(\sL^\pm_M)\otimes_{p_2^*\twist(\sS^\pm_M)\otimes p_3^*\twist(\sS^\mp_M)} \twist (p_1^* \twist(\sC^\pm_M)\otimes p_3^*\twist(\sC^\mp_M))\to \m^*\twist(\sL^\pm_M),\label{eq:mapofbimodules2} \\ 
&&{\rm over}\quad  \Map(\R^{2|1}/\Z,M)\times (\R_{>0}\times \HH^{2|1}_\pm\times \HH^{2|1}_\pm\times \HH^{2|1}_\mp)/\Z\nonumber
\eeq
where $p_1,p_2,p_3$ are the maps in~\eqref{eq:targetofadj} and $\m$ is~\eqref{eq:sfm2};
\item[iii.] an isomorphism of vector bundles
\beq\label{Eq:twisttrace}
&&p_1^*\twist(\sL^\pm_M)\otimes_{p_1^*\twist(\sS^\pm_M)\otimes p_2^*\twist(\sS^\mp_M)} p_2^*\twist(\sR^\mp(M))\to \m^*\twist(\mathcal{L}^{2|1}(M)),\label{eq:mapofbimodules3}\\ &&{\rm over } \quad (\Map(\R^{2|1}/\Z,M)\times (\R_{>0}\times \HH^{2|1}_\pm\times \iHH^{2|1}_\mp)/\Z)^\Z\nonumber 
\eeq
for $p_1,p_2$ the maps in~\eqref{eq:targetofdual} and $\m$ is~\eqref{eq:sfm3}.
\end{enumerate}
The isomorphisms above associated to composition are subject to further coherence conditions, as dictated by the definition of an internal functor. 
\end{prop}
\bp
By Definition~\ref{defn:internalfunctor}, an internal functor assigns a map of bimodule bundles to any composition of bordisms. Then (i) comes from the composition in Proposition~\ref{prop:stI}, (ii) comes from the composition in Proposition~\ref{prop:adjunctionbordisms}, and (iii) comes from the composition in Proposition~\ref{prop:tracerelations}.
\ep
\begin{prop} \label{prop:maincompositionstatement}
A $\twist$-twisted field theory satisfies the conditions
\begin{enumerate}
\item[i.] the maps of bimodule bundles are equal 
\beq\label{eq:comppathbimodule}
&&p_1^* E(\sC^+_M)\circ p_2^* E(\sC^+_M)= \m^*E(\sC^+_M),\\ 
&&{\rm over}\quad  \Map(\R^{2|1}/\Z,M)\times (\R_{>0}\times \HH^{2|1}_\pm\times \HH^{2|1}_\pm)/\Z\nonumber
\eeq
relative to the map of bundles of bimodules~\eqref{eq:mapofbimodules1}, 
where $p_1,p_2$ are the maps~\eqref{eq:composablesuperpath} and $\m$ is~\eqref{eq:sfm};
\item[ii.] the maps of right module bundles are equal 
\beq\label{eq:bimoduleadjunction}
&&p_2^* E(\sL^+_M)\circ (p_1^* E(\sC^+_M)\otimes p_3^*E(\sC^-_M))= \m^*E(\sL^+_M), \\ 
&&{\rm over} \quad \Map(\R^{2|1}/\Z,M)\times (\R_{>0}\times \HH^{2|1}_\pm\times \HH^{2|1}_\pm\times \HH^{2|1}_\mp)/\Z\nonumber 
\eeq
relative to the map of right module bundles~\eqref{eq:mapofbimodules2}, where $p_1,p_2,p_3$ are the maps in~\eqref{eq:targetofadj} and $\m$ is~\eqref{eq:sfm2}; 
\item[iii.] the maps of vector bundles are equal 
\beq\label{Eq:bimoduletrace}
&&p_1^*E(\sL^\mp_M)\circ p_2^*E(\sR^\pm_M)=\m^*E(\mathcal{L}^{2|1}(M)),\\
&&{\rm over } \quad (\Map(\R^{2|1}/\Z,M)\times (\R_{>0}\times \HH^{2|1}_\mp\times \iHH^{2|1}_\pm)/\Z)^\Z\nonumber 
\eeq
relative to~\eqref{eq:mapofbimodules3}, for $p_1,p_2$ the projections to the factors in~\eqref{eq:targetofdual} and $\m$ is~\eqref{eq:sfm3}. 
\end{enumerate}
\end{prop}
\bp
By Definition~\ref{defn:internalnattrans}, an internal natural transformation requires that the data in Proposition~\ref{prop:dataoftwistedFT} satisfy compatibility properties for composition in the internal category $2|1\EBord(M)$. Specifically,~\eqref{eq:comppathbimodule} follows from compatibility with composition for the families $\sC^\pm_M$ as computed in Proposition~\ref{prop:stI}. The equality~\eqref{eq:bimoduleadjunction} follows from compatibility for composition of the families $\sL^\pm_M$ and $\sC^\pm_M\coprod \sC^\mp_M$ computed in Proposition~\ref{prop:adjunctionbordisms}. The equality~\eqref{Eq:bimoduletrace} follows from compatibility with the composition of $\sL^\mp_M$ and $\sR^\pm_M$ computed in Proposition~\ref{prop:tracerelations}. 
\ep

\begin{lem} \label{lem:isometryandtwist}
A twist determines the data of an isomorphism of right $\act_\pm^*\twist(\sS^\mp_M)\otimes \proj_\pm^*\twist(\sS^\pm_M)$-module bundles
\beq\label{eq:twistsymemtry}
\inv^*\twist(\sL^\mp_M)\simeq \twist(\sL^\pm_M\circ \sigma)\quad {\rm over} \quad \Map(\R^{2|1}/\Z,M)\times (\R_{>0}\times \HH^{2|1}_\pm)/\Z
\eeq
where $\inv$ is the map of sheaves~\eqref{eq:2morphsigma}. Relative to this compatibility data,  a twisted field theory is required to satisfy the compatibility property
\beq\label{eq:bimodulesymemtry}
\inv^*E(\sL^\mp_M)=E(\sL^\pm_M\circ \sigma)\quad {\rm over} \quad \Map(\R^{2|1}/\Z,M)\times (\R_{>0}\times \HH^{2|1}_\pm)/\Z.
\eeq
\end{lem}

\bp
The 2-commuting data in Proposition~\ref{prop:symmetry} is an isomorphism between maps to the stack $\Mor(2|1\EBord(M))$. A twist $\twist\colon 2|1\EBord(M)\to \TA$ therefore assigns to this isomorphism between maps to the stack and an isomorphism~\eqref{eq:twistsymemtry} between the corresponding bundles of modules. Arguing similarly, evaluating a twisted field theory imposes the property~\eqref{eq:bimodulesymemtry}.
\ep

As before, let $\sqrt{\Fl}$ denote the map of sheaves~\eqref{eq:Flonsheaves}-\eqref{eq:Flonsheaves3} generating a $\Z/4$-action. 

\begin{lem} \label{lem:Flstuff}
A twist determines an isomorphism of algebra bundles $\alpha_\pm\colon \twist(\sS^\pm_M)\to \sqrt{\Fl}^*\twist(\sS^\mp_M)$ over $\Map(\R^{2|1}/\Z,M)\times \R_{>0}$ with the property that $\alpha_\pm\circ\alpha_\mp$ is the parity involution. Relative to $\alpha_\pm$, there are isomorphisms of modules 
\beq\label{eq:valueofsqrtflontwist}
&&\twist(\sC^\pm_M)\simeq \sqrt{\Fl}^*\twist(\sC^\mp_M), \quad  \twist(\sL^\pm_M)\simeq \sqrt{\Fl}^*\twist(\sL^\mp_M), \quad \twist(\sR^\pm_M)\simeq \sqrt{\Fl}^*\twist(\sR^\mp_M)
\eeq
over $\Map(\R^{2|1}/\Z,M)\times (\R_{>0}\times \HH^{2|1}_\pm)/\Z$ and $\Map(\R^{2|1}/\Z,M)\times (\R_{>0}\times \iHH^{2|1}_\pm)/\Z$. A twisted field theory supplies the data of an isomorphism
\beq\label{eq:sqrtFlonsS}
_{\twist(\sS^\pm_M)}E(\sS^\pm_M)\xrightarrow{\sim} _{\sqrt{\Fl}^*\twist(\sS^\pm)} \sqrt{\Fl}^*E(\sS^\mp_M)
\eeq
that squares to the parity involution and with compatibility properties in the form of commutative diagrams of modules
\beq
\begin{tikzpicture}[baseline=(basepoint)];
\node (A) at (0,0) {$ \twist(\sC^\pm_M)\otimes_{\proj_\pm^*\twist(\sC^\pm_M)} \proj^*_\pm E(\sS^\pm_M)$};
\node (B) at (6,0) {$\act^*_\pm E(\sS^\pm_M)$};
\node (C) at (0,-1.5) {$\sqrt{\Fl}^*( \twist(\sC^\mp_M)\otimes_{\proj_\mp^*\twist(\sC^\mp_M)} \proj^*_\mp E(\sS^\mp_M))$};
\node (D) at (6,-1.5) {$ \sqrt{\Fl}^*\act^*_\mp E(\sS^\mp_M)$};
\draw[->] (A) to node [above] {$E(\sC^\pm_M)$} (B);
\draw[->] (A) to node [left] {$\simeq$} (C);
\draw[->] (C) to node [above] {$\sqrt{\Fl}^*E(\sC^\mp_M)$} (D);
\draw[->] (B) to node [right] {$\simeq$} (D);
\path (0,-.75) coordinate (basepoint);
\end{tikzpicture}\nonumber\\
\begin{tikzpicture}[baseline=(basepoint)];
\node (A) at (0,0) {$\twist(\sL^\pm_M)\otimes_{\proj_\pm^*\twist(\sS^\pm_M)\otimes \act_\pm^*\twist(\sS^\mp_M)} \proj_\pm^*E(\sS^\pm_M)\otimes \act_\pm^*E( \sS^\mp_M)$};
\node (B) at (8,0) {$\underline{\C}$};
\node (C) at (0,-1.5) {$\sqrt{\Fl}^*(\twist(\sL^\mp_M)\otimes_{\proj_\mp^*\twist(\sS^\mp_M)\otimes \act_\mp^*\twist(\sS^\pm_M)} \proj_\mp^*E(\sS^\mp_M)\otimes \act_\mp^*E( \sS^\pm_M))$};
\node (D) at (8,-1.5) {$\sqrt{\Fl}^*\underline{\C}$};
\draw[->] (A) to node [above] {$E(\sL^\pm_M)$} (B);
\draw[->] (A) to node [left] {$\simeq$} (C);
\draw[->] (C) to node [above] {$\sqrt{\Fl}^*E(\sL^\mp_M)$} (D);
\draw[->] (B) to node [right] {$\simeq$} (D);
\path (0,-.75) coordinate (basepoint);
\end{tikzpicture}\nonumber\\
\begin{tikzpicture}[baseline=(basepoint)];
\node (A) at (0,0) {$ \twist(\sR^\pm_M)$};
\node (B) at (6,0) {$\proj_\pm^*E(\sS^\mp_M) \otimes \act_\pm^*E(\sS^\pm_M)$};
\node (C) at (0,-1.5) {$\sqrt{\Fl}^*E(\sR^\mp_M)$};
\node (D) at (6,-1.5) {$ \sqrt{\Fl}^*(\proj_\mp^*E(\sS^\pm_M)\otimes \act_\mp^*E(\sS^\mp))$};
\draw[->] (A) to node [above] {$E(\sR^\pm_M)$} (B);
\draw[->] (A) to node [left] {$\simeq$} (C);
\draw[->] (C) to node [above] {$\sqrt{\Fl}^*E(\sR^\mp_M)$} (D);
\draw[->] (B) to node [right] {$\simeq$} (D);
\path (0,-.75) coordinate (basepoint);
\end{tikzpicture}\nonumber
\eeq
over $\Map(\R^{2|1}/\Z,M)\times (\R_{>0}\times \HH^{2|1}_\pm)/\Z$ and $\Map(\R^{2|1}/\Z,M)\times (\R_{>0}\times \iHH^{2|1}_\pm)/\Z$, where we again have suppressed the algebra isomorphisms $\alpha_\pm$.
\end{lem}

\bp
The isometries in Lemma~\ref{lem:orientation11restrict} give isomorphisms between maps to $\Ob(2|1\EBord(M))$ and $\Mor(2|1\EBord(M))$. When evaluating a twist on the family $\sS^\pm_M$, we obtain isomorphisms of algebra bundles~$\alpha_\pm=\twist(\sqrt{\Fl})$. To deduce the property that $\alpha_\pm\circ \alpha_\mp$ is the parity involution uses that a twisted field theory is a \emph{flip preserving} internal natural transformation (see~\cite[\S2.6]{ST11}) where the spin involution~\eqref{eq:flip} of super Euclidean families furnishes the flip $\Fl=\sqrt{\Fl}\circ \sqrt{\Fl}$. Similarly, the isomorphisms~\eqref{eq:valueofsqrtflontwist} come from applying the internal functor $\twist$ to the commutative diagrams in Lemma~\ref{lem:orientation11restrict}. Turning to a twisted field theory, $\sqrt{\Fl}$ provides the isomorphism of modules~\eqref{eq:sqrtFlonsS}. The remaining commutative diagrams follow from the diagram in Lemma~\ref{lem:orientation11restrict} and the definition of an internal natural transformation. 
\ep

\subsection{Reflection structures}\label{sec:reflrealEFT}
Next we unpack the data of a reflection structure from~\S\ref{sec:RPgeneral} for twists and twisted field theories. We recall the notation~\eqref{eq:conjsheaf} for the conjugation functor applied to a sheaf, and the maps of sheaves~\eqref{eq:daggersheaf1}-\eqref{eq:daggersheaf3} used to describe the restrictions of the reflection structure on $2|1\EBord(M)$ in Proposition~\ref{prop:restrictRR}.

\begin{lem}\label{lem:RPtwist}
A reflection structure for a twist $\twist\colon 2|1\EBord(M)\to \TA$ determines invertible $\twist(\sS^\pm_M)$-$\dagger^*\overline{\twist(\sS^\mp_M))}$-bimodules $\mathscr{R}^\pm\to \Map(\R^{2|1}/\Z,M)\times\R_{>0}$ with isomorphisms to the identity bimodule
\beq\label{eq:involutiondata}
\mathscr{R}^\pm\otimes_{\dagger^*\overline{\twist(\sS^\mp_M))}} \dagger^*\overline{\mathscr{R}}^\mp\simeq \id_{\twist(\sS^\pm_M)}.
\eeq
A reflection structure for $\twist$ further specifies isomorphisms of bimodules
\beq
&&\twist(\sC^\pm_M)\otimes_{\act_\pm^*\twist(\sS^\pm_M)} \act_\pm^*\mathscr{R}^\pm \xrightarrow{\sim} \proj_\pm^*\mathscr{R}^\pm \otimes_{\dagger^*\overline{\proj}_\mp^*\overline{\twist(\sS^\mp_M)}} \dagger^*\overline{\twist(\sC^\mp_M)},\nonumber
\eeq
\beq
&&\twist(\sL^\pm_M)\otimes_{\proj_\pm^*\twist(\sS^\pm_M)\otimes \act_\pm^*\twist(\sS^\mp_M)}(\proj_\pm^*\mathscr{R}^\pm\otimes \act_\pm^*\mathscr{R}^\mp) \xrightarrow{\sim} \dagger^*\overline{\twist(\sL^-_M)},\nonumber
\eeq
over $\Map(\R^{2|1}/\Z,M)\times (\R_{>0}\times \HH^{2|1}_\pm)/\Z$, and
\beq
&&\twist(\sR^+_M)\xrightarrow{\sim}(\proj_\pm^*\mathscr{R}^\mp \otimes \act_\pm^*\mathscr{R}^\pm)\otimes_{\dagger^*(\overline{\proj}_\mp^*\overline{\twist(\sS^\pm_M)}\otimes \overline{\act}_\mp^*\overline{\twist(\sS^\mp_M)})}  \dagger^*\overline{\twist(\sR^-_M)},
\eeq
over $\Map(\R^{2|1}/\Z,M)\times (\R_{>0}\times \iHH^{2|1}_\pm)/\Z$ satisfying an involutive property.
\end{lem}

\bp
From Definition~\ref{defn:RP0}, a reflection structure for the twist $\twist$ determines an invertible internal natural isomorphism, i.e., a map of stacks
$$
\eta\colon \Ob(2|1\EBord(\pt))\to \Mor(\TA)^\times
$$
valued in invertible bimodules together with coherence data $\rho$ satisfying properties. Applying $\eta$ to the family $\sS^\pm_M$, we obtain the bundle of invertible bimodules,~$\mathscr{R}^\pm$. Restricting $\rho$ to the families $\sC^\pm_M$, $\sL^\pm_M$ and $\sR^\pm_M$, a reflection structure provides the claimed isomorphisms of bimodules. A reflection structure also specifies involutive data (see Remark~\ref{rmk:involutivereflection}) which when evaluated on $\sS^\pm_M$ yields the isomorphism~\eqref{eq:involutiondata} and imposes involutive properties on the isomorphisms of bimodules.
\ep


\begin{lem} \label{lem:RPFT}
Fix a reflection structure for a twist $\twist\colon 2|1\EBord(M)\to \TA$. Then a reflection structure for a $\twist$-twisted field theory is an isomorphism of bundles of left $\twist(\sS^\pm_M)$-modules over $\Map(\R^{2|1}/\Z,M)\times\R_{>0}$
%
%
\beq
E(\sS^\pm_M) \stackrel{\sim}{\to} \mathscr{R}^\pm \otimes_{\dagger^*\overline{\twist(\sS^\mp_M)}}  \dagger^*\overline{E(\sS^\mp_M))}.\label{eq:RPdataonspt}
\eeq
These isomorphisms are required to satisfy compatibility conditions for the families $\sL^\pm(M)$, $\sR^\pm(M)$, and~$\sC^\pm(M)$, namely commutative squares, 
\beq
\begin{tikzpicture}[baseline=(basepoint)];
\node (A) at (0,0) {$\twist(\sC^\pm_M)\otimes_{\proj_\pm^*\twist(\sS^\pm_M)} \proj_\pm^*E(\sS^\pm_M)$};
\node (B) at (8,0) {$\act_\pm^*E(\sS^\pm_M)$};
\node (C) at (0,-1.5) {$\dagger^*\Big(\overline{\twist(\sC^\mp_M)} \otimes_{\overline{\proj}_\mp^*\overline{\twist(\sS^\mp_M)}} \overline{\proj}_\mp^*\overline{E(\sS^\mp_M)}\Big)$};
\node (D) at (8,-1.5) {$\dagger^*\overline{\act}_\mp^* \overline{E(\sS^\mp_M)}$};
\draw[->] (A) to node [above] {$E(\sC^\pm_M)$} (B);
\draw[->] (A) to node [left] {$\simeq$} (C);
\draw[->] (C) to node [below] {$\dagger^*\overline{E(\sC^\mp_M)}$} (D);
\draw[->] (B) to node [right] {$\simeq$} (D);
\path (0,-.75) coordinate (basepoint);
\end{tikzpicture}\label{eq:selfadjointdiagram}\\
\resizebox{.95\textwidth}{!}{$
\begin{tikzpicture}[baseline=(basepoint)];
\node (A) at (0,0) {$\twist(\sL^\pm_M)\otimes_{\proj_\pm^*\twist(\sS^\pm_M)\otimes \act_\pm^*\twist(\sS^\mp_M)} \proj_\pm^*E(\sS^\pm_M)\otimes\act_\pm^*E(\sS^\mp_M)$};
\node (B) at (9,0) {$\underline{\C}$};
\node (C) at (0,-1.5) {$\dagger^*(\overline{\twist(\sL^\mp_M)} \otimes_{\overline{\proj}_\mp^*\overline{\twist(\sS^\mp_M)}\otimes \overline{\act}_\mp^*\overline{\twist(\sS^\pm_M)}} \overline{\proj}_\mp^*\overline{E(\sS^\mp_M)}\otimes \overline{\act}_\mp^*\overline{E(\sS^\pm_M)}$};
\node (D) at (9,-1.5) {$\dagger^* \underline{\overline{\C}}$};
\draw[->] (A) to node [above] {$E(\sL^\pm_M)$} (B);
\draw[->] (A) to node [left] {$\simeq$} (C);
\draw[->] (C) to node [below] {$\dagger^*\overline{E(\sL^\mp_M)}$} (D);
\draw[->] (B) to node [right] {$\simeq$} (D);
\path (0,-.75) coordinate (basepoint);
\end{tikzpicture}$}\label{eq:RPdataonL0m}
\eeq
over $\Map(\R^{2|1}/\Z,M)\times (\R_{>0}\times \HH^{2|1}_\pm)/\Z$, and commutative squares over $\Map(\R^{2|1}/\Z,M)\times (\R_{>0}\times \iHH^{2|1}_\pm)/\Z$
\beq\resizebox{.95\textwidth}{!}{$
\begin{tikzpicture}[baseline=(basepoint)];
\node (A) at (0,0) {$\underline{\C}$};
\node (B) at (9,0) {$\twist(\sR^\pm_M)\otimes_{\proj_\pm^*\twist(\sS^\mp_M)\otimes \act_\pm^*\twist(\sS^\pm_M)} \proj_\pm^*E(\sS^\mp_M)\otimes\act_\pm^*E(\sS^\pm_M)$};
\node (C) at (0,-1.5) {$\dagger^* \underline{\overline{\C}}$};
\node (D) at (9,-1.5) {$\dagger^*(\overline{\twist(\sR^\mp_M)} \otimes_{\overline{\proj}_\mp^*\overline{\twist(\sS^\pm_M)}\otimes \overline{\act}_\mp^*\overline{\twist(\sS^\mp_M)}} \overline{\proj}_\mp^*\overline{E(\sS^\pm_M)}\otimes \overline{\act}_\mp^*\overline{E(\sS^\mp_M)}$};
\draw[->] (A) to node [above] {$E(\sR^\pm_M)$} (B);
\draw[->] (A) to node [left] {$\simeq$} (C);
\draw[->] (C) to node [below] {$\dagger^*\overline{E(\sR^\mp_M)}$} (D);
\draw[->] (B) to node [right] {$\simeq$} (D);
\path (0,-.75) coordinate (basepoint);
\end{tikzpicture}$}\nonumber
\eeq
where the vertical arrows use data from Lemma~\ref{lem:RPtwist} and~\eqref{eq:RPdataonspt}. In the above, we have suppressed the tensor products with the bimodules $\mathscr{R}^\pm$ that mediate between $\twist(\sS^\pm)$-modules and $\dagger^*\overline{\twist(\sS^\mp)}$-modules to avoid (further) cluttering the notation. 
\end{lem}
\bp
By Definition~\ref{defn:RP0}, reflection positivity data for a twisted field theory determines an isomorphism between the maps of stacks $\Ob(2|1\EBord(M))\to \Mor(\TA)$. Evaluating this isomorphism on the families $\sS^\pm_M$ yields the data~\eqref{eq:RPdataonspt}. 

For the data~\eqref{eq:RPdataonspt} to determine an isomorphism between internal natural transformations, there are compatibility requirements for each object of $\Mor(2|1\EBord(M))$. Evaluating on the families $\sC^\pm_M$, $\sL^\pm_M$ and $\sR^\pm_M$ recovers the claimed commutative squares as this compatibility data. 
\ep

\subsection{The degree~$n$ supersymmetric twist and its restrictions}\label{sec:susyextunique}

In~\S\ref{sec:appendfer} we review Stolz and Teichner's construction of the chiral free fermion twist $\Fer_n\colon 2\EBord(M)\to \TA$ as a functor out of the 2-dimensional Euclidean bordism category (without supersymmetry). We recall the functor $\mathcal{S}\colon 2\EBord(M)\to 2|1\EBord(M)$ from Lemma~\ref{lem:superfication} that sends a Euclidean bordism to a super Euclidean bordism. 

\begin{hyp} \label{hyp:chiral}
The $n$-chiral free fermions constructed in~\cite[\S5.3]{ST11} have a \emph{supersymmetric extension} denoted $s\Fer_n$
\beq
&&\begin{tikzpicture}[baseline=(basepoint)];
\node (A) at (0,0) {$2\EBord$};
\node (B) at (5,0) {$\TA$};
\node (C) at (2,-.5) {$\varphi\ \Downarrow$};
\node (D) at (0,-1.5) {$2|1\EBord(\pt)$};
\draw[->] (A) to node [left] {$\mathcal{S}$} (D);
\draw[->] (A) to node [above] {$\Fer_n$} (B);
\draw[->,dashed] (D) to node [below=5pt] {$s\Fer_n$} (B);
\path (0,-.75) coordinate (basepoint);
\end{tikzpicture}\label{eq:twistcont}
\eeq
that is unique up to an internal natural isomorphism $\varphi$. The reflection structure for the chiral fermions (see~\eqref{eq:FerRP}) has a supersymmetric extension that is unique up to isomorphism of internal natural transformations.
\end{hyp}

\begin{rmk}\label{Eq:chiraltosusy}
Hypothesis~\ref{hyp:chiral} is part of a larger expectation: 2-dimensional chiral conformal field theories should determine $\mathcal{N}=(0,1)$ supersymmetric ones, e.g., see~\cite[Page~12]{TopMoon} or \cite[\S1]{LinPei}. This connection is the bridge between Segal's early work~\cite{Segal_Elliptic,SegalCFT} and the Stolz--Teichner program; see~ \cite[page~50-51]{ST04} and~\cite[\S5.2]{ST11}. The essential idea is that a chiral theory has $\bar L_0=0$ (using the notation from~\S\ref{sec:motivateQFT}), and hence a supersymmetric extension is uniquely determined by taking $\bar G_0=0$. This manipulation in super Lie algebras exponentiates to the semigroup of cylinders. In particular, this argument applies to the chiral free fermions, see Proposition~\ref{prop:chiral1}.
%
\end{rmk}

\begin{notation}
In this subsection we will use the notation $s\Fer_n$ to denote the supersymmetric extension~\eqref{eq:twistcont} to avoid confusion with the (nonsupersymmetric) twist $\Fer_n$; in later sections we will revert to the simpler notation $\Fer_n\colon 2|1\EBord\to \TA$. 
\end{notation}

Our next goal is to verify Hypothesis~\ref{hyp:chiral} in the Ramond sector, i.e., for the $2|1$-Euclidean bordisms relevant to Theorem~\ref{thm1}. We begin with a definition.

\begin{defn} A supersemigroup representation $\varphi\colon \HH_\pm\to \End(V)$ is \emph{holomorphic} if~$\partial_{\bar \tau}\varphi=0$. Equivalently, for an $S$-point $(\tau,\bar\tau)\in \HH_\pm(S)$, we have $\varphi(\tau,\bar\tau)=\varphi(\tau)\colon S\to \End(V)$ is independent of $\bar\tau$. 
\end{defn}

The following has been used previously by Stolz and Teichner to lift a chiral field theory to a super Euclidean one, e.g., see \cite[Proof of Theorem 1.15, second part]{ST11}.

\begin{lem}\label{lem:holoext}
A holomorphic representation of the semigroup $\HH_\pm$ has a unique extension to a representation of the supersemigroup $\HH_\pm^{2|1}$ along the canonical inclusion $\HH_\pm\subset \HH_\pm^{2|1}$. 
\end{lem}

\bp
For $p$ the projection, consider the composition
\beq\label{eq:susyextension}
\tilde{\varphi}\colon \HH_\pm^{2|1}\xrightarrow{p} \HH_\pm\xrightarrow{\varphi} \End(V). 
\eeq
This composition is a super semigroup representation
$$
\tilde{\varphi}(\tau,\bar\tau,\eta)\cdot \tilde{\varphi}(\tau',\bar\tau',\eta')=\varphi(\tau)\cdot \tilde{\varphi}(\tau')=\varphi(\tau+\tau')=\tilde{\varphi}(\tau+\tau',\bar\tau+\bar\tau'+\eta\eta',\eta+\eta').
$$
We caution that the projection $p$ is \emph{not} a homomorphism, and the semigroup homomorphism property for $\tilde{\varphi}$ critically uses holomorphicity of $\varphi$. Since $p$ is a 1-sided inverse to the canonical inclusion $\HH_\pm\subset \HH_\pm^{2|1}$, the composition~\eqref{eq:susyextension} is indeed the unique extension of $\varphi$. 
\ep

We recall the Lie category $\Ann=\Ann(\pt)$ of Euclidean cylinders from Definition~\ref{defn:Path}. 


\begin{prop} \label{prop:chiral1}
Consider the restriction of $\Fer_n$ to the subcategory $\Ann$ of Euclidean cylinders. There is an extension to the category of Euclidean supercylinders~$\Ann^{2|1}$,
\beq\label{eq:sAnnextension}
&&\begin{tikzpicture}[baseline=(basepoint)];
\node (A) at (0,0) {$\Ann$};
\node (B) at (4,0) {$\TA$};
\node (C) at (0,-1.5) {$\Ann^{2|1}$};
\draw[->] (A) to node [above] {$\Fer_n$} (B);
\draw[->] (A) to node [left] {$\mathcal{S}$} (C);
\draw[->,dashed] (C) to node [below=5pt] {$s\Fer_n$} (B);
\node (E) at (1.4,-.5) {$\Downarrow$};
\path (0,-.75) coordinate (basepoint);
\end{tikzpicture}
\eeq
along $\mathcal{S}$, Stolz and Teichner's superfication functor~\eqref{eq:superfication}. Furthermore, the reflection structure~\eqref{eq:FerRP} on $\Fer_n$ has a unique extension to $s\Fer_n$ using the reflection structure~\eqref{eq:formdag} on the internal category $\Ann^{2|1}$.
\end{prop} 

\bp
It suffices to prove the lemma for $n=-1$. The restriction of~$\Fer_{-1}$ to $\Ann$ is characterized by the first family of bimodules in~\eqref{eq:degtwisbimod} and the anti-involution~\eqref{eq:Ferinvol}. The superfication functor~\eqref{eq:superfication} on Euclidean cylinders is given by the internal functor
$$
\resizebox{\textwidth}{!}{$
\Ann=\left(\begin{array}{c} ((\R_{>0}\times \HH_+)/\Z\coprod (\R_{>0}\times \HH_-)/\Z)\sq \Z/4 \\ \downarrow \downarrow \\ (\R_{>0}\coprod \R_{>0})\sq \Z/4 \end{array}\right) \xhookrightarrow{\mathcal{S}} \left(\begin{array}{c} ((\R_{>0}\times \HH^{2|1}_+)/\Z\coprod (\R_{>0}\times \HH^{2|1}_-)/\Z)\sq \Z/4\\ \downarrow\downarrow \\ (\R_{>0}\coprod \R_{>0})\sq \Z/4 \end{array}\right)=\Ann^{2|1}.$}
$$
On objects, $\mathcal{S}$ is the identity map, and on morphisms $\mathcal{S}$ is determined by the canonical inclusion of the reduced manifold. We define the extension $s\Fer_{-1}$ on $\Ob(\Ann^{2|1})$ to be equal to the value~\eqref{eq:Ferdefn} of $\Fer_{-1}$ on $\Ob(\Ann)$, i.e., the bundle of algebras $\Fer(\rS_\ell^\pm)$ where the generator of $\Z/4$ acts by the anti-involution~\eqref{eq:Ferinvol}. The value of $\Fer_{-1}$ on morphisms in $\Ann$ comes from an $\R_{>0}$-family of holomorphic $\HH$-actions on~$\Fer(\rS_\ell^\pm)$. By Lemma~\ref{lem:holoext}, these actions have a unique extension to $\HH_\pm^{2|1}$-actions, with formula given by
\beq\label{eq:sFerext}
\tilde{\varphi}_\ell^\pm(\tau,\bar\tau,\eta)=\varphi_\ell^\pm(\tau),
\eeq
where $\varphi_\ell^\pm(\tau)$ is determined by~\eqref{eq:holoHaction}.
With this fixed, $s\Fer_{-1}$ is defined in the same way as $\Fer_{-1}$: the value of $s\Fer_{-1}$ on $\Mor(\Ann^{2|1})$ is the family of bimodules gotten by viewing $\Fer(\rS_\ell^\pm)$ as a bimodule over itself with left-action twisted by the family of algebra automorphisms~\eqref{eq:sFerext}. The composition and unitor data are canonically determined by the fact that $\tilde\varphi^\pm_\ell(0)=\id$ and $\tilde\varphi^\pm_\ell(\tau)\cdot \tilde\varphi^\pm_\ell(\tau')=\tilde\varphi^\pm_\ell(\tau+\tau')$. The $\Z/4$-action on the resulting family of bimodules again comes from~\eqref{eq:Ferinvol}. This completes the construction of an internal functor $s\Fer_{-1}\colon \Ann^{2|1}\to \TA$, and (by inspection) the diagram~\eqref{eq:sAnnextension} strictly commutes. 

%

The reflection structure for $\Fer_{-1}$ is determined by a $*$-superalgebra structure~\eqref{eq:FerRP} on the algebra $\Fer(\rS_\ell^\pm)$ which then determines involutions on $\Fer(\rS_\ell^\pm)$ regarded as a bimodule. Since the objects in $\Ann$ and $\Ann^{2|1}$ are the same, this reflection structure uniquely extends to $s\Fer_{-1}$ using the same $*$-superalgebra structure. 
%
\ep

\begin{cor}\label{cor:FervaluesonRL}
Using~\eqref{eq:sFerext}, the values~\eqref{eq:degtwisbimod} of the internal functor $\Fer_{-1}\colon 2\EBord\to \TA$ have unique extensions 
to supersymmetric families, 
\beq
s\Fer_n(\rL_{\ell,\tau,\bar\tau,\eta}^\pm)=(\Fer_\ell^\mp)_{\Fer_\ell^\mp\otimes ^\tau\Fer_\ell^\mp}\qquad s\Fer_n(\rR_{\ell,\tau,\bar\tau,\eta}^\pm)= {}_{\Fer^\pm_\ell \otimes {}^\tau\Fer_\ell^\mp} \Fer_\ell^\mp\nonumber
\eeq
and these values are compatible with composition in $2|1\EBord$. Furthermore, the reflection structure~\eqref{eq:FerRP} has a unique extension to the above families.
\end{cor}

\begin{defn}\label{defn:Clifftwistfornormy}
Define the \emph{Clifford twist} as the internal functor $\Cl_n\colon \Ann\to \TA$ or $\Cl_n\colon \Ann^{2|1}\to \TA$ whose value on objects is the coproduct of trivial algebra bundles
$$
\underline{\cCl}_{+n}\coprod \underline{\cCl}_{-n} \to (\R_{>0}\coprod \R_{>0})\sq \Z/4=\Ob(\Ann)=\Ob(\Ann^{2|1})
$$
where the generator of $\Z/4$ acts by the algebra isomorphism
\beq\label{eq:Cliffinv0}
\alpha_\mp\colon  \cCl_{\mp 1}\to \cCl_{\pm 1},\quad v\mapsto  \pm iv, \quad v\in \C.
\eeq
The value of the Clifford twist on morphisms is the coproduct of trivial bimodule bundles gotten from viewing $\cCl_{\pm n}$ as a bimodule over itself, with $\Z/4$-equivariant structure again determined by~\eqref{eq:Cliffinv0}. Unit and composition data are determined by the canonical isomorphisms of bimodules $\cCl_{\pm n}\otimes_{\cCl_{\pm n}}\cCl_{\pm n}\simeq \cCl_{\pm n}$. 
\end{defn}

We recall the equivalence \cite[\S6]{ST11} between the Clifford twist and the restriction of the $\Fer_n$ to $\Ann$, see~\eqref{eq:ClFeriso} below. 

\begin{cor}\label{cor:ClFer}
There is an internal natural isomorphism between the supersymmetric extension $s\Fer_n$ and the Clifford twist, 
\beq
&&\begin{tikzpicture}[baseline=(basepoint)];
\node (A) at (0,0) {$\Ann^{2|1}$};
\node (B) at (5,0) {$\TA$};
\node (C) at (2.5,0) {$s\mathcal{F} \Downarrow$};
\draw[->,bend left=10] (A) to node [above] {$\Cl_n$} (B);
\draw[->,bend right=10] (A) to node [below] {$s\Fer_n$} (B);
\path (0,0) coordinate (basepoint);
\end{tikzpicture}\label{eq:sClFeriso}
\eeq
and $s\mathcal{F}$ restricts to the equivalence~\eqref{eq:ClFeriso} along $\mathcal{S}\colon \Ann\hookrightarrow \Ann^{2|1}$. 
Similarly, the values on the bordisms $\rL_{\ell,\tau,\bar\tau,\eta}^\pm$ and $\rR_{\ell,\tau,\bar\tau,\eta}^\pm$ in Corollary~\ref{cor:FervaluesonRL} are isomorphic (in $\Mor(\TA)$) with the trivial bundles of bimodules,
\beq
&&s\Fer_n(\rL_{\ell,\tau,\bar\tau,\eta}^\pm)\simeq(\cCl_{\pm n})_{\cCl_{\pm n}\otimes \cCl_{\mp n}}\qquad s\Fer_n(\rR_{\ell,\tau,\bar\tau,\eta}^\pm)\simeq {}_{\cCl_{\mp n}\otimes \cCl_{\pm n}} (\cCl_{\pm n}).\label{eq:CliffisoonLR}
\eeq
Under the above isomorphisms, the reflection structure on $s\Fer_n$ is compatible with the $*$-superalgebra structure~\eqref{eq:CliffRP} on~$\cCl_{\pm n}$.

\end{cor}
\bp
It again suffices to consider the case $n=-1$. The $\R_{>0}$-family of bimodules~\eqref{eq:Morita} implementing the equivalence~\eqref{eq:ClFeriso} have a holomorphic $\HH_\pm$-action~\eqref{eq:moritaHact}; by Lemma~\ref{lem:holoext} there is a unique extension to an action by~$\HH^{2|1}_\pm$. The $\Z/4$-equivariant structure on this bimodule also uniquely extends to a $\Z/4$-equivariant structure over $\Ob(\Ann^{2|1})$. This constructs the equivalence~\eqref{eq:sClFeriso} with the claimed properties. The same bimodule with $\HH_\pm^{2|1}$-action also constructs the isomorphisms~\eqref{eq:CliffisoonLR}. The reflection structure in all cases comes from the $*$-superalgebra structure on Clifford algebras, and so is compatible with the bimodule~\eqref{eq:Morita} implementing the isomorphisms. 
\ep

\begin{prop} \label{prop:chiral2}
The value of $\Fer_n$ over Euclidean tori admits an extension
\beq
&&\begin{tikzpicture}[baseline=(basepoint)];
\node (A) at (0,0) {$\Lat\sq (\Spin(2)\times \SL_2(\Z))$};
\node (B) at (6,0) {$\Mor(\TA)$};
\node (C) at (0,-1.5) {$s\Lat \sq (\Spin(2)\times \SL_2(\Z))$};
\draw[->] (A) to node [above] {$\Fer^{\otimes n}$} (B);
\draw[->] (A) to node [left] {$\mathcal{S}$} (C);
\draw[->,dashed] (C) to (B);
\node (E) at (1.7,-.7) {$\Downarrow$};
\path (0,-.75) coordinate (basepoint);
\end{tikzpicture}
\eeq
that is unique up to isomorphism. Equivalently, by~\eqref{eq:PfvalofFer}, the Pfaffian line bundle for the periodic-periodic spin structure has a unique extension to $2|1$-Euclidean tori. Furthermore, a section over $s\Lat \sq (\Spin(2)\times \SL_2(\Z))$ descends to a section over $s\Lat \sq (\Spin(2)\rtimes \E^{2|1} \times \SL_2(\Z))$ if and only if its restriction to $\Lat\sq (\Spin(2)\times \SL_2(\Z))$ is a holomorphic section of the Pfaffian. 
\end{prop} 

\bp 
Using the description of $\sLat$ from Lemma~\ref{lem:4component}, pulling back along the inclusion $\Lat\subset s\Lat$ induces an equivalence of categories between 
line bundles on $\Lat$ and line bundles on $s\Lat$. Hence, there is a unique isomorphism class of line bundle given by the trivial line. It therefore remains to consider extensions of equivariant structures for the trivial line on~$\Lat$. The presheaf $\GL_1(\C)$ of automorphisms of the trivial line bundle is characterized by
$$
\GL_1(\C)(S)=C^\infty(S)^{\ev,\times}.
$$
We have the isomorphism
$$
C^\infty(\sLat\times \Spin(2)\times \SL_2(\Z))^{\ev,\times}\simeq C^\infty(\Lat\times \Spin(2)\times \SL_2(\Z))^{\ev,\times}, 
$$
i.e., automorphism of line bundles over $s\Lat$ are determined by their restriction along $\Lat\hookrightarrow \sLat$. Hence, any $\Spin(2)\times \SL_2(\Z)$-equivariant structure for a line bundle on $s\Lat$ is determined by its restriction to $\Lat$. The remaining claim about descent of sections to the unpointed stacks is implied by Proposition~\ref{prop:maintorusprop} with $M=\pt$. 
\ep

\begin{rmk} \label{rmk:generalproof}
Assuming a generators and relations presentation of $2|1\EBord$ as described in~\cite[Proof of Theorem 1.15, second part]{ST11}, one can verify Hypothesis~\ref{hyp:chiral} in the framework of~\cite{ST11}. Propositions~\ref{prop:chiral1} and~\ref{prop:chiral2} establish the extension in the Ramond sector. By Remark~\ref{rmk:Beyond}, the Neveu--Schwarz sector supercylinders assemble into a moduli space that is isomorphic to the (non-supersymmetric) moduli of Euclidean cylinders. Hence, the supersymmetric extension on this sector is trivially unique. 
\end{rmk}

\subsection{The values of a degree~$n$ field theory}

\begin{cor} \label{cor:degreenvalues}
Assuming Hypothesis~\ref{hyp:chiral}, a degree~$n$ field theory $E\colon \one\Rightarrow \Fer_n$ as in~\eqref{eq:twistedEFT} determines 
\begin{itemize}
\item[i.] bundles of $\cCl_{\pm n}$-modules $E(\sS^\pm_M)$ over $\Map(\R^{2|1}/\Z,M)\times \R_{>0}$;
\item[ii.] maps of bundles of left $\cCl_{\pm n}$-modules over $\Map(\R^{2|1}/\Z,M)\times (\R_{>0}\times \HH^{2|1}_\pm)/\Z$ 
\beq
&&E(\sC^\pm_M)\colon \proj_\pm^*E(\sS^\pm_M)\to \act_\pm^*E(\sS^\pm_M),\nonumber
\eeq
compatible with composition~\eqref{eq:comppathbimodule};
\item[iii.] maps of vector bundles over $\Map(\R^{2|1}/\Z,M)\times (\R_{>0}\times \HH^{2|1}_\pm)/\Z$ 
\beq\label{eq:Clnliearpairing}
&& E(\sL^\pm_M)\colon \cCl_{\pm n} \otimes_{\cCl_{\pm n} \otimes \cCl_{\mp n}} \Big(\proj_\pm^*E(\sS^\pm_M)\otimes \act_\pm^*E( \sS^\mp_M)\Big)\to \underline{\C} ,
\eeq
satisfying the adjunction relation~\eqref{eq:bimoduleadjunction} and the symmetry relation~\eqref{eq:bimodulesymemtry}; 

\item[iv.] maps of bundles of left $\cCl_{\mp n}\otimes \cCl_{\pm n}$-modules over $\Map(\R^{2|1}/\Z,M)\times (\R_{>0}\times \iHH^{2|1}_\pm)/\Z$ 
\beq
&&E(\sR^\pm_M)\colon \cCl_{\mp n} \to \proj_\pm^*E(\sS^\mp_M) \otimes \act_\pm^*E(\sS^\pm_M), \nonumber
 \eeq
\item[v.] an equality of $\E^{2|1}\rtimes \Spin(2)\times \SL_2(\Z)$-invariant functions 
\beq\label{eq:equalityofsecttrivi}
p_1^*E(\sL^\mp_M)\circ p_2^*E(\sR^\pm_M)=\m^*E(\mathcal{L}^{2|1}(M))\in C^\infty((\R_{>0}\times \HH^{2|1}_\mp\times \iHH^{2|1}_\pm)/\Z)^\Z)
\eeq
relative to the isomorphism of trivial line bundles in~\eqref{eq:maptoPf}, where the left hand side of~\eqref{eq:equalityofsecttrivi} is an appropriately normalized Clifford supertrace, and the right hand side is the section $E(\mathcal{L}^{2|1}(M))$ pulled back to the standard trivialization of the Pfaffian line bundle on $\R_{>0}\times \HHz$ and then pulled back along $\m$ in~\eqref{eq:sfm3}. 
\end{itemize}
\end{cor}

\bp
Propositions~\ref{prop:dataoftwistedFT} and~\ref{prop:maincompositionstatement} describe the general data of a twist and twisted field theory. Proposition~\ref{prop:chiral1} and Corollary~\ref{cor:FervaluesonRL} gives values of the twist $\Fer_n$ on supercylinders, while Proposition~\ref{prop:chiral2} provides the values on supertori. This yields versions of the above data for modules over the algebra~\eqref{eq:Ferdefn}. Corollaries~\ref{cor:ClFer} and~\ref{prop:chiral2} and Lemma~\ref{lem:Tracenormalize1} allow one to rephrase these data in terms of modules over the Clifford algebras and the Clifford supertrace.
\ep

\begin{cor}\label{cor:RPFT}
Assuming Hypothesis~\ref{hyp:chiral}, a reflection structure on a degree~$n$ field theory $E\colon \one\to \Fer_n$ determines maps of bundles relative to the maps of sheaves~\eqref{eq:daggersheaf1} and~\eqref{eq:daggersheaf2}:
%
%
\begin{itemize}
\item[i.] an isomorphism of bundles of $\cCl_{\pm n}$-modules over $\Map(\R^{2|1}/\Z,M)\times \R_{>0}$ 
\beq
E(\sS^\pm_M) \stackrel{\sim}{\to} \dagger^*\overline{E(\sS^\mp_M)},\label{eq:RPdataonspt3}
\eeq
relative to the $*$-superalgebra structure on $\cCl_{\pm n}$;
\item[ii.] a commutative square of bundles of $\cCl_{\pm n}$-modules over $\Map(\R^{2|1}/\Z,M)\times (\R_{>0}\times \HH^{2|1}_\pm)/\Z$ 
\beq
&&\begin{tikzpicture}[baseline=(basepoint)];
\node (A) at (0,0) {$\proj_\pm^*E(\sS^\pm_M)$};
\node (B) at (6,0) {$\act_\pm^*E(\sS^\pm_M)$};
\node (C) at (0,-1.5) {$\dagger^* \overline{\proj}_\mp^*\overline{E(\sS^\mp_M)}$};
\node (D) at (6,-1.5) {$\dagger^* \overline{\act}_\mp^*\overline{E(\sS^\mp_M)}$};
\draw[->] (A) to node [above] {$E(\sC^\pm_M)$} (B);
\draw[->] (A) to node [left] {$\simeq$} (C);
\draw[->] (C) to node [below] {$\dagger^*\overline{E(\sC^\mp_M)}$} (D);
\draw[->] (B) to node [right] {$\simeq$} (D);
\path (0,-.75) coordinate (basepoint);
\end{tikzpicture}\label{eq:selfadjointdiagram2}
\eeq
\item[iii.] a commutative square of vector bundles over $\Map(\R^{2|1}/\Z,M)\times (\R_{>0}\times \HH^{2|1}_\pm)/\Z$ 
\beq
&&\begin{tikzpicture}[baseline=(basepoint)];
\node (A) at (0,0) {$\cCl_{\pm n} \otimes_{\cCl_{\pm n} \otimes \cCl_{\mp n}} \big( \proj_\pm^*E(\sS^\pm_M)\otimes\act_\pm^*E(\sS^\mp_M)\big)$};
\node (B) at (8,0) {$\underline{\C}$};
\node (C) at (0,-1.5) {$\dagger^*\Big(\overline{\cCl}_{\mp n} \otimes_{\overline{\cCl}_{\mp n} \otimes \overline{\cCl}_{\pm n}} \big(\overline{\proj}_\mp^*\overline{E(\sS^\mp_M)}\otimes \overline{\act}_\mp^*\overline{E(\sS^\pm_M)}\big)\Big)$};
\node (D) at (8,-1.5) {$\dagger^* \underline{\overline{\C}}$};
\draw[->] (A) to node [above] {$E(\sL^\pm_M)$} (B);
\draw[->] (A) to node [left] {$\simeq$} (C);
\draw[->] (C) to node [below] {$\dagger^*\overline{E(\sL^\mp_M)}$} (D);
\draw[->] (B) to node [right] {$\simeq$} (D);
\path (0,-.75) coordinate (basepoint);
\end{tikzpicture}\label{eq:RPdataonL0m2}
\eeq
over $\Map(\R^{2|1}/\Z,M)\times (\R_{>0}\times \HH^{2|1}_\pm)/\Z$.
\end{itemize}
\end{cor}
\bp
The proof is similar to that of the previous corollary, but as a specialization of Lemma~\ref{lem:RPFT} to the twist~$\Fer_n$, then translated into a statement for Clifford modules using Corollary~\ref{cor:ClFer}. 
\ep

\begin{cor}\label{cor:RPFT2}
Assuming Hypothesis~\ref{hyp:chiral}, a twisted field theory $E$ determines maps of bundles relative to isomorphisms of sheaves~\eqref{eq:Flonsheaves}-\eqref{eq:Flonsheaves3}:
%
%
\begin{itemize}
\item[i.] an isomorphism of bundles of $\cCl_{\pm n}$-modules over $\Map(\R^{2|1}/\Z,M)\times \R_{>0}$ 
\beq
E(\sS^\pm_M) \stackrel{\sim}{\to} \sqrt{\Fl}^*E(\sS^\mp_M),\label{eq:Fldataonspt3}
\eeq
relative to the anti-automorphism $\alpha\colon \cCl_{\pm n} \xrightarrow{\sim} \cCl_{\pm n}^\op=\cCl_{\mp n}$ defined in~\eqref{eq:Cliffinv0};
\item[ii.] a commutative square of bundles of $\cCl_{\pm n}$-modules over $\Map(\R^{2|1}/\Z,M)\times (\R_{>0}\times \HH^{2|1}_\pm)/\Z$ 
\beq
&&\begin{tikzpicture}[baseline=(basepoint)];
\node (A) at (0,0) {$\proj_\pm^*E(\sS^\pm_M)$};
\node (B) at (6,0) {$\act_\pm^*E(\sS^\pm_M)$};
\node (C) at (0,-1.5) {$\sqrt{\Fl}^* \proj_\mp^* E(\sS^\mp_M)$};
\node (D) at (6,-1.5) {$\sqrt{\Fl}^* \act_\mp^*E(\sS^\mp_M),$};
\draw[->] (A) to node [above] {$E(\sC^\pm_M)$} (B);
\draw[->] (A) to node [left] {$\simeq$} (C);
\draw[->] (C) to node [below] {$\sqrt{\Fl}^*E(\sC^\mp_M)$} (D);
\draw[->] (B) to node [right] {$\simeq$} (D);
\path (0,-.75) coordinate (basepoint);
\end{tikzpicture}\label{eq:Fldata2}
\eeq
\item[iii.] a commutative square of vector bundles over $\Map(\R^{2|1}/\Z,M)\times (\R_{>0}\times \HH^{2|1}_\pm)/\Z$ 
\beq
&&\begin{tikzpicture}[baseline=(basepoint)];
\node (A) at (0,0) {$\cCl_{\pm n} \otimes_{\cCl_{\pm n} \otimes \cCl_{\mp n}} \big( \proj_\pm^*E(\sS^\pm_M)\otimes\act_\pm^*E(\sS^\mp_M)\big)$};
\node (B) at (8,0) {$\underline{\C}$};
\node (C) at (0,-1.5) {$\sqrt{\Fl}^*\Big(\cCl_{\mp n} \otimes_{\cCl_{\mp n} \otimes \cCl_{\pm n}} \big(\proj_\mp^*E(\sS^\mp_M)\otimes \act_\mp^*E(\sS^\pm_M)\big)\Big)$};
\node (D) at (8,-1.5) {$\sqrt{\Fl}^* \underline{\overline{\C}}$};
\draw[->] (A) to node [above] {$E(\sL^\pm_M)$} (B);
\draw[->] (A) to node [left] {$\simeq$} (C);
\draw[->] (C) to node [below] {$\sqrt{\Fl}^*E(\sL^\mp_M)$} (D);
\draw[->] (B) to node [right] {$\simeq$} (D);
\path (0,-.75) coordinate (basepoint);
\end{tikzpicture}\label{eq:FldataonL0m2}
\eeq
over $\Map(\R^{2|1}/\Z,M)\times (\R_{>0}\times \HH^{2|1}_\pm)/\Z$.
\end{itemize}
\end{cor}
\bp
Using Lemma~\ref{lem:Flstuff}, the proof is completely analogous to the previous. 
\ep

The structures in Corollaries~\ref{cor:degreenvalues} and~\ref{cor:RPFT2} can be restricted along the inclusions 
\beq\label{eq:restrictalongthis}
&&M\times \R_{>0} \subset \Pi TM\times \R_{>0} \simeq \Map(\R^{0|1},M)\times \R_{>0} \subset \Map(\R^{2|1}/\Z,M)\times \R_{>0} 
\eeq
yielding bundles of $\cCl_{\pm n}$-modules 
\beq\label{eq:firstEE}
\HBR_\pm:=E(\sS^\pm_M)|_{M\times \R_{>0}}.
\eeq
In particular, the restriction of~\eqref{eq:Clnliearpairing} and~\eqref{eq:RPdataonspt3} combine to give a map
$$
\cCl_{\pm n}\otimes_{\overline{\cCl}_{\pm n}\otimes \cCl_{\pm n}} \Big(\overline{\HBR}_\pm \otimes \HBR_\pm\Big)\simeq \cCl_{\mp n}\otimes_{\cCl_{\mp n}\otimes \cCl_{\pm n}} \Big(\HBR_\mp \otimes \HBR_\pm\Big) \to \underline{\C}
$$
which determines a pairing for which the Clifford actions are adjoint,
\beq\label{eq:Lzeropairing3}\label{eq:Cliffadjoint} 
&&\langle-,-\rangle \colon \overline{\HBR}_+ \otimes_{\cCl_n} \HBR_+ \to \underline{\C},\qquad \langle \omega^*\cdot x,y\rangle =(-1)^{|\omega||x|}\langle x,\omega\cdot y\rangle, \ \omega\in \cCl_n.
\eeq

\begin{lem}\label{lem:hermitian}
The pairing~\eqref{eq:Lzeropairing3} is hermitian in the $\Z/2$-graded sense
\beq
\langle x,y\rangle =(-1)^{|x||y|}\overline{\langle y,x\rangle}.\label{eq:sesquilinear}
\eeq
\end{lem}
\bp
The argument uses similar ideas to the first part of \cite[Corollary 6.25]{HST}, though the signs here are different as we are working in the oriented setting (versus the unoriented framework of \cite{HST}). 
Consider the diagram 
\beq\label{eq:hermitianargument}
\begin{tikzpicture}[baseline=(basepoint)];
\node (AA) at (0,1.5) {$\overline{\cCl}_{\pm n}\otimes_{\overline{\cCl}_{\pm n}\otimes\cCl_{\pm n}} (\dagger^*\overline{\proj}_\pm^*\overline{E(\sS^\pm_M)}\otimes \act_\mp^*E(\sS^\pm_M))$};
\node (A) at (0,0) {$\cCl_{\mp n}\otimes_{\cCl_{\mp n}\otimes \cCl_{\pm n}}(\proj_\mp^*E(\sS^\mp_M)\otimes \act_\mp^*E(\sS^\pm_M))$};
\node (B) at (8,0) {$\underline{\C}$};
\node (C) at (0,-1.5) {$\dagger^*\Big(\overline{\cCl}_{\pm n}\otimes_{\overline{\cCl}_{\pm n}\otimes \overline{\cCl}_{\mp n}} (\overline{\proj}_\pm^*\overline{E(\sS^\pm_M)}\otimes \overline{\act}_\pm^*\overline{E(\sS^\mp_M)})\big)$};
\node (D) at (8,-1.5) {$\dagger^*\overline{\underline{\C}}$};
\node (DD) at (4,-3.5) {$\dagger^*\big(\overline{\cCl}_{\mp n}\otimes_{\overline{\cCl}_{\mp n}\otimes \overline{\cCl}_{\pm n}} (\overline{\act}_\pm^*\overline{E(\sS^\mp_M)}\otimes \overline{\proj}_\pm^*\overline{E(\sS^\pm_M)})$};
\draw[->] (AA) to (A);
\draw[->] (A) to node [above] {$E(\sL^\mp_M)$} (B);
\draw[->] (A) to (C);
\draw[->] (C) to node [below] {$\overline{E(\sL^\pm_M)}$} (D);
\draw[->] (B) to node [right] {$\simeq$} (D);
\draw[->] (C) to node [below] {$\sigma$} (DD);
\draw[->] (DD) to node [right=10pt] {$\overline{\inv}^*\overline{E(\sL^\mp_M)}$} (D);
\draw[->, bend left=10] (AA) to node [above=8pt] {restricts to $\langle-,-\rangle$} (B);
\path (0,0) coordinate (basepoint);
\end{tikzpicture}
\eeq
where the $*$-superalgebra structure on $\cCl_{\pm n}$ is used implicitly throughout. The indicated arrow in the upper right restricts to the pairing $\langle-,-\rangle$ by definition. 
Commutativity of the lower triangle follows from~\eqref{eq:bimodulesymemtry}.  Commutativity of the inner square is the condition ~\eqref{eq:RPdataonL0m2} on reflection positivity data.
The composition of arrows starting from the top left and using the lower triangle is the right hand side of~\eqref{eq:sesquilinear}, and hence the restriction of this diagram along~\eqref{eq:restrictalongthis} implies the claimed formula, where the sign in~\eqref{eq:sesquilinear} comes from the symmetry isomorphism~$\sigma$ for the (graded) tensor product in the lower triangle. 
\ep

See Definition~\ref{defn:RP0} for the notion of a reflection structure on a twisted field theory, and below we take the reflection structure~\eqref{eq:21EBreflection} on $2|1\EBord(M)$ and the reflection structure for the degree~$n$ twist determined by~\eqref{eq:FerRP}. We refer to~\eqref{eq:positivitypairing} for the signs involved in the positivity property for a $\Z/2$-graded hermitian pairing.

\begin{defn}\label{defn:RPFT} Assuming Hypothesis~\ref{hyp:chiral}, a $2|1$-Euclidean field theory of degree~$n$ with reflection structure is \emph{reflection positive} if the pairing~\eqref{eq:Lzeropairing3} is a (graded) positive-definite hermitian pairing on the super vector bundle~$\HBR_+$ over $M\times \R_{>0}$. 
\end{defn}


\begin{rmk}
For a reflection positive field theory,~\eqref{eq:Cliffadjoint} implies that the $\cCl_n$-action on $\HBR_+$ is automatically a $*$-representation relative to the hermitian inner product.
\end{rmk}

\begin{defn}\label{defn:STEFT}
Assuming Hypothesis~\ref{hyp:chiral},  for each $n\in \Z$ define the groupoid $2|1\EFT^n(M)$ as having objects degree~$n$, reflection positive $2|1$-Euclidean field theories, and morphisms isomorphisms between degree~$n$ field theories compatible with the reflection structure. \end{defn}

\section{Extracting geometric data from a degree~$n$ field theory}\label{sec:datafromEFT}

The main goal of this section is to prove Proposition~\ref{prop:mainprop}, extracting a superconnection with additional structure and property from an object of $2|1\EFT^n(M)$. The mechanism that produces this geometric data is the restriction of a field theory to Euclidean supercylinders,
\beq
\label{eq:twistedrestriction2}
&&\begin{tikzpicture}[baseline=(basepoint)];
\node (AA) at (-3,0) {$\Ann^{2|1}(M)$};
\node (A) at (0,0) {$2|1\EBord(M)$};
\node (B) at (4,0) {$\TA.$};
\node (C) at (2.1,0) {$\Downarrow$};
\draw[->] (AA) to (A);
\draw[->,bend left=10] (A) to node [above] {$\one$} (B);
\draw[->,bend right=10] (A) to node [below] {$\Fer_n$} (B);
\path (0,0) coordinate (basepoint);
\end{tikzpicture}
\eeq 
With the goal of making the geometric constructions in this section as flexible as possible, the main constructions are independent from the definitions in~\S\ref{sec:iota} and~\S\ref{sec:STmess}. This is achieved by phrasing all definitions in terms of a twisted representation of $\Ann^{2|1}(M)$ as in Definition~\ref{defn:degreenAnn} below, omitting reference to the intermediary degree~$n$ field theory in~\eqref{eq:twistedrestriction2}. At the end of each subsection we provide statements linking the representation theory of $\Ann^{2|1}(M)$ with degree~$n$ field theories under~\eqref{eq:twistedrestriction2}. 

\begin{defn}\label{defn:degreenAnn}
A \emph{degree~$n$ representation of $\Ann^{2|1}(M)$} is an internal natural transformation
\beq\label{Eq:Annrep}
&&\begin{tikzpicture}[baseline=(basepoint)];
\node (B) at (0,0) {$\Ann^{2|1}(M)$};
\node (C) at (5,0) {$\TA$};
\node (D) at (2.5,0) {$\Downarrow$};
\draw[->,bend left=12] (B) to node [above] {$\one$} (C);
\draw[->,bend right=12] (B) to node [below] {$\Cl_n$}  (C);
\path (0,-.05) coordinate (basepoint);
\end{tikzpicture}
\eeq
where $\Cl_n\colon \Ann^{2|1}(M)\to \TA$ is the pull back the Clifford twist (Definition~\ref{defn:Clifftwistfornormy}) along the functor $\Ann^{2|1}(M)\to \Ann^{2|1}(\pt)=\Ann^{2|1}$ induced by the canonical map $M\to \pt$. 
Degree~$n$ representations form a groupoid whose objects are internal natural transformations~\eqref{Eq:Annrep} and whose morphisms are isomorphisms between internal natural transformations. 
\end{defn}

\subsection{Clifford linear superconnections from degree~$n$ representations of $\tAnn^{2|1}(M)$}

Consider the restriction of a degree~$n$ representation~\eqref{Eq:Annrep} along the $\Z/4$-cover~\eqref{eq:sZ4descent}; we will also refer to the resulting twisted representation of $\tAnn^{2|1}(M)$ as a degree~$n$ representation. We caution that the notation~\eqref{eq:firstEE} in the previous section differs slightly from~$\HBR_\pm$ used in the following proposition and throughout this section.

\begin{prop}\label{prop:superconn}
A degree~$n$ representation of $\tAnn^{2|1}(M)$ determines vector bundles $\HBR_+\to M$ and $\HBR_-\to M$ with fiberwise $S^1\times \cCl_{n}$- and $S^1\times \cCl_{-n}$-actions, respectively. Furthermore, the bundles~$\HBR_\pm$ are endowed with the structure of a smooth 1-parameter family of superconnections $\A_\ell^\pm$ for $\ell\in \R_{>0}$ that are $\cCl_{\pm n}$-linear and $S^1$-invariant. 
\end{prop}

\bp
From Definition~\ref{defn:internalnattrans}, part of the data of an internal natural transformation $E$ in~\eqref{Eq:Annrep} is a map of stacks
$$
\Ob(\tAnn^{2|1}(M))\to \Mor(\TA)
$$
with specified source and target data. This is equivalent to bundles of left Clifford modules
\beq\label{eq:tildeHBR}
\cCl_{+n}\acts \widetilde\HBR_+\to \Map(\R^{0|1},M)\times \R_{>0},\quad \cCl_{-n}\acts \widetilde\HBR_-\to \Map(\R^{0|1},M)\times \R_{>0}
\eeq
over each component of $\Ob(\tAnn^{2|1}(M))$. The next piece of data in an internal natural transformation is an isomorphism between maps of stacks $\Mor(\tAnn^{2|1}(M))\to \Mor(\TA)$ (see~\eqref{diag:nu}), which is given by bundle maps
\beq
&&\cCl_{\pm n}\otimes_{\cCl_{\pm n}} \proj_\pm^*\Gamma(\widetilde\HBR_\pm)\to \act_\pm\Gamma(\widetilde\HBR_\pm)\label{eq:defrho}\\
&& \quad {\rm over} \quad \Map(\R^{0|1},M)\times (\R_{>0}\times \HH^{2|1}_\pm)/\Z.\nonumber
\eeq
The maps~\eqref{eq:defrho} are then further required to satisfy a condition with respect to composition in~$\tAnn^{2|1}(M)$. To unpack this in terms of the differential geometry of $M$, we first observe the isomorphism
\beq\label{eq:sectionsare}
&&\Gamma(\widetilde\HBR_\pm)\simeq \Omega^\bullet(M;C^\infty(\R_{>0}))\otimes_{C^\infty(M)} \Gamma(M,\HBR_\pm),\qquad \HBR_\pm=\widetilde\HBR_\pm|_{M\times\{1\}}
\eeq
of modules over $C^\infty(\Map(\R^{0|1},M)\times \R_{>0})\simeq \Omega^\bullet(M;C^\infty(\R_{>0}))$. This uses the fact that (up to isomorphism) a vector bundle over $\Map(\R^{0|1},M)$ pulls back from $M$ (see Lemma~\ref{lem:sVect}) and a vector bundle over $M\times \R_{>0}$ is isomorphic to one that pulls back along the projection $M\times \R_{>0}\to M$. 
Pulling back $\proj_\pm^*\widetilde{\HBR}_\pm$ and $\act_\pm^*\widetilde{\HBR}_\pm$ along $p$ gives the isomorphisms 
\beq
\Gamma(\proj_\pm^*\widetilde{\HBR}_\pm)&\simeq&\Omega^\bullet(M;C^\infty((\R_{>0}\times \HH^{2|1}_\pm)/\Z))\otimes_{\Omega^\bullet(M;C^\infty(\R_{>0}))} \Gamma(\widetilde{\HBR}_\pm)\nonumber\\
&\simeq&C^\infty((\R_{>0}\times  \HH^{2|1}_\pm)/\Z))\otimes \Omega^\bullet(M)\otimes_{C^\infty(M)} \Gamma(\HBR_\pm)\nonumber\\
&\simeq&C^\infty((\R_{>0}\times \HH^{2|1}_\pm)/\Z))\otimes \Omega^\bullet(M;\HBR_\pm)\label{eq:easydescript1}\\
\Gamma(\act_\pm^*\widetilde{\HBR}_\pm)&\simeq&\Omega^\bullet(M;C^\infty((\R_{>0}\times \HH^{2|1}_\pm)/\Z))\otimes_{\Omega^\bullet(M;C^\infty(\R_{>0}))} \Gamma(\widetilde{\HBR}_\pm)\nonumber\\
&\simeq&C^\infty((\R_{>0}\times \HH^{2|1}_\pm)/\Z))\otimes \Omega^\bullet(M)\otimes_{C^\infty(M)} \Gamma(\HBR_\pm)\nonumber\\
&\simeq&C^\infty((\R_{>0}\times \HH^{2|1}_\pm)/\Z))\otimes \Omega^\bullet(M;\HBR_\pm)\label{eq:easydescript2}
\eeq
using the projective tensor product throughout and the isomorphism~\eqref{eq:sectionsare}. The map of vector bundles~\eqref{eq:defrho} is a map between the $\Omega^\bullet(M;C^\infty((\R_{>0}\times \HH^{2|1}_\pm)/\Z))$-modules~\eqref{eq:easydescript1} and~\eqref{eq:easydescript2} for the module structure determined by the algebra maps $\proj_\pm^*$ and $\act_\pm^*$ determined by~\eqref{eq:actproj}. This map is $C^\infty((\HH^{2|1}_\pm\times \R_{>0})/\Z))$-linear map, affording the description of~\eqref{eq:defrho} in terms of a function valued in Clifford-linear endomorphisms,
\beq\label{eq:defrhotil}
&&\rho_\pm \in C^\infty((\R_{>0}\times \HH^{2|1}_\pm)/\Z;\End_{\cCl_{\pm n}}(\Omega^\bullet(M;\HBR_\pm))).
\eeq
Using the descriptions of $\proj_\pm^*$ and $\act_\pm^*$ from Lemma~\ref{lem:stmap}, $\rho_\pm$ must satisfy the property 
\beq \label{eq:Leibniz}
&&\resizebox{.93\textwidth}{!}{$
\rho_\pm(\ell,\tau,\bar\tau,\eta) (\alpha v)=(\alpha-\eta d\alpha)\rho_\pm(\ell,\tau,\bar\tau,\eta)(v),  \quad \begin{array}{l} \alpha\in \Omega^\bullet(M;C^\infty(\R_{>0}\times \HH_\pm)/\Z), \\ v\in C^\infty((\R_{>0}\times \HH^{2|1}_\pm)/\Z))\otimes \Omega^\bullet(M;\HBR_\pm). \end{array}$}
\eeq
Furthermore, Lemma~\ref{lem:concat} shows that $\rho_+$ and $\rho_-$ determine an $\R_{>0}$-family of super semigroup representations,
$$
\rho_\pm(\ell,\tau,\bar\tau,\eta)\circ \rho_\pm(\ell,\tau',\bar\tau',\eta')=\rho_\pm(\ell,\tau+\tau',\bar\tau+\bar\tau'+\eta\eta',\eta+\eta')
$$ 
since composition in $\tAnn^{2|1}(M)$ corresponds to multiplication in $\HH^{2|1}_\pm$. 

To further analyze the geometry of~\eqref{eq:defrhotil}, Taylor expand $\rho_\pm$ in the odd variable $-\eta\in C^\infty((\R_{>0}\times \HH^{2|1}_\pm)/\Z)$ to obtain even and odd endomorphisms $\aa_\pm(\ell,\tau,\bar\tau)$ and $\b_\pm(\ell,\tau,\bar\tau)$
\beq\label{eq:AB}
&&\rho_\pm(\ell,\tau,\bar\tau,\eta)=\aa_\pm(\ell,\tau,\bar\tau)- \eta \b_\pm(\ell,\tau,\bar\tau),  \qquad \begin{array}{l} \aa_\pm(\ell,\tau,\bar\tau)\in \End_{\cCl_{\pm n}}(\Omega^\bullet(M;\HBR_\pm))^\ev,\\ \b_\pm(\ell,\tau,\bar\tau)\in \End_{\cCl_{\pm n}}(\Omega^\bullet(M;\HBR_\pm))^\odd.\end{array}
\eeq
Restricting to $\R/\ell\Z\subset \HH_\pm/\ell\Z\subset (\R_{>0}\times \HH_\pm^{2|1})/\Z$, for each value of $\ell$ we obtain an action, $\R/\ell\Z\to \End_{\cCl_{\pm n}}(\Omega^\bullet(M;\HBR_\pm))$. By~\eqref{eq:Leibniz}, the $\R/\ell\Z$-action is linear over $\Omega^\bullet(M;C^\infty(\R_{>0}))$. Hence, the $\R/\ell\Z$-action is through Clifford-linear bundle endomorphisms
\beq\label{eq:circlact}
&&\R/\ell\Z\to \Gamma(\End_{\cCl_{\pm n}}(\HBR_\pm)), \\ 
&&t\mapsto \aa_\pm(\ell,t,t),\ \ \aa_\pm(\ell,s,s)\circ \aa_\pm(\ell,t,t)=\aa_\pm(\ell,s+t,s+t). \nonumber
\eeq
These actions necessarily vary smoothly with $\ell\in \R_{>0}$. So reparameterizing by $t\mapsto \ell t$ at each fiber $\ell\in \R_{>0}$, and choosing the isomorphism~\eqref{eq:sectionsare} to be $S^1$-equivariant, one obtains a single $S^1=\R/\Z$-action that commutes with the Clifford action.

Next consider the endomorphism of $\Omega^\bullet(M;\HBR_\pm)$
\beq
&&\A_\ell^\pm:=\b_\pm(\ell,0,0)\in \End_{\cCl_{\pm n}}(\Omega^\bullet(M;\HBR_\pm))^\odd. \label{eq:superconnconstrr}
\eeq
Since $\b_\pm(\ell,\tau,\bar\tau)$ is odd, so is $\A_\ell^\pm$. Furthermore, because $\b_\pm(\ell,\tau,\bar\tau)$ is $\cCl_{\pm n}$-linear, so is $\A_\ell^\pm$. Equation~\eqref{eq:Leibniz} shows that $\A_\ell^\pm$ is an endomorphism of $\Omega^\bullet(M;\HBR_\pm)$ for any $\ell\in \R_{>0}$, and in this way we view $\A_\ell^\pm$ as an $\R_{>0}$-family of operators. Equation~\eqref{eq:Leibniz} also shows that $\A_\ell^\pm$ commutes with the $\R/\ell\Z$-action~\eqref{eq:circlact}, and hence commutes with the $S^1$-action on $\HBR_\pm$. Finally,~\eqref{eq:Leibniz} applied to~\eqref{eq:AB} yields the pair of equations
$$
\aa_\pm(\ell,\tau,\bar\tau)(\alpha v)=\alpha \aa_\pm(\ell,\tau,\bar\tau)(v),\qquad \b_\pm(\ell,\tau,\bar\tau)(\alpha v)=d\alpha \aa_\pm(\ell,\tau,\bar\tau)(v)+(-1)^{|\alpha|}\alpha \b_\pm(\ell,\tau,\bar\tau)(v).
$$
Using that $\aa_\pm(\ell,0,0)=\id$, the value of the second equation at $(\tau,\bar\tau)=(0,0)$ is the Leibniz rule for~$\A_\ell^\pm$, and so we conclude that $\A_\ell^+$ and $\A_\ell^-$ are superconnections. This completes the proof of the proposition. 
\ep

\begin{lem} \label{prop:EFT1} Assuming Hypothesis~\ref{hyp:chiral}, the restriction of the degree~$n$ twist along~\eqref{eq:twistedrestriction2} is equivalent to the Clifford twist $\Cl_n$. Under this equivalence we have the identification
\beq\label{eq:HBRFT}
\widetilde{\HBR}_\pm\simeq E({\sf s}_M^\pm)
\eeq
where $\widetilde{\HBR}_\pm$ are the vector bundles~\eqref{eq:tildeHBR}, and ${\sf s}_M^\pm$ is the family of small bordisms defined in~\eqref{eq:LM0}.
\end{lem}
\bp
This is immediate from Definition~\ref{defn:Clifftwistfornormy}, Corollary~\ref{cor:ClFer}, and Corollary~\ref{cor:degreenvalues}. 
\ep

\begin{rmk} 
In light of Corollary~\ref{cor:ClFer}, an equivalent groupoid to the one in Definition~\ref{defn:degreenAnn} involves modules over the tensor powers of the free fermion algebras~\eqref{eq:Ferdefn}. The invertible bimodule~\eqref{eq:Morita} furnishes the equivalence between these two points of view. So without passing through the equivalence with the Clifford twist, the value~\eqref{eq:HBRFT} becomes
\beq\label{eq:Fermodextend}
E({\sf s}_M^\pm)\simeq ({}_{\Fer_\ell^\pm}(\cCl_{\pm 1} \bigotimes_{m>0} \Lambda^\bullet \C_{\ell,\mp m}))^{\otimes -n} \otimes_{\cCl_{\pm n}} \widetilde{\HBR}_\pm.
\eeq
 With the eventual goal of extracting $\KO(\!(q)\!)$-classes from $ \widetilde{\HBR}_\pm$, it is more convenient in this paper to work with Clifford modules as in~\eqref{eq:HBRFT}.
\end{rmk} 

\subsection{Flip preserving representations}
Next we consider descent data for a degree~$n$ representation of $\tAnn^{2|1}(M)$ along the $\Z/4$-cover~\eqref{eq:sZ4descent}. For degree~$n$ representations that restrict from field theories, this descent data must satisfy a condition related to the spin statistics theorem described in~\cite[Definition 6.44]{HST} and~\cite[Remark 6.45]{HST}.

We start with a general observation. The category of supermanifolds has a $\Z/2$-action that acts by the parity involution on ring of functions of a supermanifold; maps between supermanifolds are automatically equivariant for this $\Z/2$-action. By naturality, the objects and morphsims of any super Lie category also inherit canonical $\Z/2$-actions for which the source, target, unit and composition maps are $\Z/2$-equivariant. 
In the case of the super Lie category $\tAnn^{2|1}(M)$, these canonical $\Z/2$-actions have the geometric interpretation as the spin flip on super Euclidean manifolds (see Lemmas~\ref{lem:flip1} and~\ref{lem:structureonsP}). In particular, the $\Z/4$-cover~\eqref{eq:sZ4descent} factors as the composition
\beq
\tAnn^{2|1}(M)&=&\left(\begin{array}{c} (\R_{>0}\times \HH_+\coprod \HH_-)/\Z\times \Map(\R^{0|1},M) \\\s \downarrow \downarrow\t \\ (\R_{>0}\coprod \R_{>0})\times \Map(\R^{0|1},M) \end{array}\right) \nonumber\\ 
&\to& \left(\begin{array}{c} (\R_{>0}\times \HH_+\coprod \HH_-)/\Z\times \Map(\R^{0|1},M)\sq \Z/2 \\ \s\downarrow \downarrow\t \\ (\R_{>0}\coprod \R_{>0})\times \Map(\R^{0|1},M)\sq \Z/2 \end{array}\right) \label{eq:flippreservegroupoid}\\ 
&\to& \left(\begin{array}{c} (\R_{>0}\times \HH_+\coprod \HH_-)/\Z\times \Map(\R^{0|1},M)\sq \Z/4 \\ \s\downarrow \downarrow\t \\ (\R_{>0}\coprod \R_{>0})\times \Map(\R^{0|1},M)\sq \Z/4\end{array}\right)=\Ann^{2|1}(M) \nonumber
\eeq
where in the middle line the $\Z/2$-action is the canonical one that we identify with the spin flip $\Fl$; in the last line the generator of $\Z/4$ acts by the square root $\sqrt{\Fl}$. 


\begin{defn}
A degree~$n$ representation of $\Ann^{2|1}(M)$ is \emph{flip preserving} if on restriction to~\eqref{eq:flippreservegroupoid} the generator of $\Z/2$ acts by the parity involution on the Clifford algebras $\cCl_{\pm n}$ and by the parity involution on Clifford modules $\widetilde{\HBR}_\pm$ in the notation from Proposition~\ref{prop:superconn}. 
\end{defn}

\begin{lem} \label{lem:Z4}
A flip-preserving degree~$n$ representation of $\Ann^{2|1}(M)$ is the data of isomorphisms 
\beq\label{eq:Z4gen}
\alpha_\pm\colon \cCl_{\pm n}\to \cCl_{\mp n},\qquad 
{\sf f}_\pm \colon \widetilde{\HBR}_\pm\xrightarrow{\sim} \sqrt{\Fl}^*\widetilde{\HBR}_\mp.
\eeq
such that $\alpha_\pm\circ \alpha_\mp$ and ${\sf f}_\pm\circ {\sf f}_\mp$ are parity involutions. Furthermore, the isomorphisms~${\sf f}_\pm$ are $S^1$-equivariant relative to the inversion homomorphism $t\mapsto -t$ on $S^1$, and the superconnections in Proposition~\ref{prop:superconn} satisfy 
 \beq\label{eq:easyadsuperconn}
\A_\ell^\pm=(-i)^{\deg+1}\A_\ell^\mp
\eeq
relative to the isomorphisms ${\sf f}_\pm$. 
\end{lem}

\bp
Equivariance data for the generator $\sqrt{\Fl}$ of the $\Z/4$-action is equivalent to the isomorphisms~\eqref{eq:Z4gen}. The flip preserving condition requires $\alpha_\pm$ and ${\sf f}_\pm$ to square to the respective parity involutions. From the definition of an internal natural transformation, we have the condition on the super semigroup representations~\eqref{eq:defrhotil}
$$
\rho_\pm(\ell,\tau,\bar\tau,\eta)=\sqrt{\Fl}^*\rho_\mp(-\tau,-\bar\tau,\sqrt{-1}\eta)\in C^\infty(\HH^{2|1}_\pm,\End_{\cCl_{\pm n}}(\Gamma(\widetilde{\HBR}_\pm))).
$$
Using the notation from~\eqref{eq:AB} and the description of $\sqrt{\Fl}^*$ from~\eqref{Eq:sqflipfun} we obtain
$$
\aa_\pm(\ell,\tau,\bar\tau)-\eta \b_\pm(\ell,\tau,\bar\tau)=(-i)^{\deg}(\aa_\mp(-\tau,-\bar\tau)- i\eta \b_\mp(-\tau,-\bar\tau))
$$
where the equality implicitly uses the isomorphism~\eqref{eq:Z4gen}. From this we deduce that the isomorphism~\eqref{eq:Z4gen} is $\R$-equivariant for the actions~\eqref{eq:circlact} relative to the automorphism $t\mapsto -t$ of $\R$, and that the superconnections obey~\eqref{eq:easyadsuperconn}.
\ep

The upshot of Lemma~\ref{lem:Z4} is that for a flip-preserving degree~$n$ representation of~$\Ann^{2|1}(M)$, the vector bundles $\HBR_+$ and $\HBR_-$ in Proposition~\ref{prop:superconn} are determined by a single $S^1$-equivariant bundle with $\R_{>0}$-family of $S^1\times \cCl_n$-linear superconnections
\beq\label{eq:singleVA}
\HBR:=\HBR_+\to M,\qquad \A_\ell:=\A_\ell^+.
\eeq
There are two possible choices of $\alpha_\pm$ in~\eqref{eq:Z4gen} that differ by the parity automorphism; without loss of generality we assume $\alpha_\pm$ is the standard anti-automorphism~\eqref{eq:Cliffinv0}.

\begin{lem}
Assuming Hypothesis~\ref{hyp:chiral}, a degree~$n$ field theory~\eqref{eq:twistedrestriction2} determines a flip preserving degree~$n$ representation of $\Ann^{2|1}(M)$. 
\end{lem}
\bp
This follows from Lemma~\ref{lem:orientation11restrict} and the fact that twisted field theories are defined as flip preserving functors, see \cite[\S2.6]{ST11} and~\S\ref{sec:backFT} below. 
\ep

\subsection{Degree~$n$ representations with a pairing}
The restriction of the map~\eqref{eq:2morphsigma} along the inclusion $\Map(\R^{0|1},M)\subset \Map(\R^{2|1}/\Z,M)$ gives the map
\beq\label{eq:inv2}
&&\inv\colon \Map(\R^{0|1},M)\times (\R_{>0}\times \HH^{2|1}_\pm)/\Z \to \Map(\R^{0|1},M)\times (\R_{>0}\times \HH^{2|1}_\mp)/\Z,
\eeq
that we will also denote $\inv$. 

\begin{lem} 
The pullback of~\eqref{eq:defrho} along $\inv$ gives a map of vector bundles 
$$
 \act_\mp^*\widetilde{\HBR}_\pm \to \proj_\mp^*\widetilde{\HBR}_\pm
$$
over $\Map(\R^{0|1},M)\times (\R_{>0}\times \HH^{2|1}_+)/\Z$ commuting with the $\cCl_{\pm n}$-action. 
\end{lem}
\bp This follows from Lemma~\ref{lem:invflip}. 
\ep

\begin{defn} \label{defn:pairing}
A \emph{pairing} on a degree~$n$ representation of $\Ann^{2|1}(M)$ is a vector bundle map 
\beq\label{eq:sPpairing}
\pair\colon \cCl_n\otimes_{\cCl_{-n}\otimes \cCl_n} (\widetilde{\HBR}_-\otimes \widetilde{\HBR}_+)\to \underline{\C}
\eeq
over $\Map(\R^{0|1},M)\times \R_{>0}$ with the property that the diagram of vector bundles over $\Map(\R^{0|1},M)\times (\R_{>0}\times \HH^{2|1}_\pm)/\Z$ commutes 
\beq
&&\begin{tikzpicture}[baseline=(basepoint)];
\node (A) at (0,0) {$\cCl_n\otimes_{\cCl_{-n}\otimes \cCl_n}(\act_+^*\widetilde{\HBR}_-\otimes \proj_+^*\widetilde{\HBR}_+)$};
\node (C) at (0,-1.5) {$\proj_+^*(\cCl_n\otimes_{\cCl_{-n}\otimes \cCl_n}(\widetilde{\HBR}_-\otimes \widetilde{\HBR}_+))$};
\node (B) at (7,0) {$\act^*_+(\cCl_n\otimes_{\cCl_{-n}\otimes \cCl_n}(\widetilde{\HBR}_-\otimes \widetilde{\HBR}_+))$};
\node (D) at (7,-1.5) {$\underline{\C}.$};
\draw[->] (A) to node [above] {$\id\otimes \rho_+$} (B);
\draw[->] (A) to node [left] {$\inv^*\rho_-\otimes \id$} (C);
\draw[->] (C) to node [below] {$\proj_+^*\pair$} (D);
\draw[->] (B) to node [right] {$\act_+^*\pair$} (D);
\path (0,-.5) coordinate (basepoint);
\end{tikzpicture}\label{eq:adjunction}
\eeq
\end{defn}

\begin{prop} \label{eq:proppairing}
A pairing on a degree~$n$ representation of $\Ann^{2|1}(M)$ determines an $\R_{>0}$-family of $S^1\times \cCl_{\pm n}$-invariant pairings $\pair_\ell$ on $\HBR_\pm$ with the property that the 1-parameter family of superconnections $\A_\ell^+,\A_\ell^-$ from Proposition~\ref{prop:superconn} satisfy the adjunction formula,
\beq\label{eq:initialadjunction}
\pair_\ell(\A_\ell^-v,w)+(-1)^{|v|}\pair_\ell( v,\A_\ell^+w)=d\pair_\ell( v,w). 
\eeq
\end{prop}
\bp
Using the descriptions of $\proj^*_\pm \widetilde{\HBR}_\pm$ and $\act^*_\pm \widetilde{\HBR}_\pm$ in terms of modules over $\Omega^\bullet(M;C^\infty((\R_{>0}\times \HH_\pm^{2|1})/\Z))$ (see Lemma~\ref{lem:stmap}), commutativity of~\eqref{eq:adjunction} implies 
\beq
(\proj_+^*\pair_\ell)((\inv^*\rho_-)(\ell,\tau,\bar\tau,\eta)v,w)&=&\pair_\ell(\rho_-(\ell,-\tau,-\bar\tau,-\eta)v,w)\nonumber\\
&=&\pair_\ell( (\aa_-(\ell,-\tau,-\bar\tau)+\eta \b_-(\ell,-\tau,-\bar\tau))v,w)\nonumber\\
&=&(\act^*_+\pair_\ell)( v,\rho_+(\ell,\tau,\bar\tau,\eta)w)\nonumber\\
&=&\pair_\ell( v,\rho_+(\ell,\tau,\bar\tau,\eta)w)+\eta d\pair_\ell( v,\rho_+(\ell,\tau,\bar\tau,\eta)w)\nonumber\\
&=&\pair_\ell( v,(\aa_+(\ell,\tau,\bar\tau)-\eta \b_+(\ell,\tau,\bar\tau))w)\nonumber\\
&&+\eta d\pair_\ell ( v,(\aa_+(\ell,\tau,\bar\tau)-\eta \b_+(\ell,\tau,\bar\tau))w), \nonumber
\eeq
where the equalities above are between elements of $\Omega^\bullet(M;C^\infty((\R_{>0}\times \HH_\pm^{2|1})/\Z))$. Setting $\eta=0$ and $(\tau,\bar\tau)=(t,t)$, the pairing is seen to be $\R$-invariant for $\R\subset \HH_\pm$; this implies the claimed $S^1$-invariance. 
Collecting the terms involving $\eta$, setting $(\tau,\bar\tau)=(0,0)$, and using the definition~\eqref{eq:superconnconstrr} of the superconnections $\A_\ell^+$ and $\A_\ell^-$, equation~\eqref{eq:initialadjunction} follows. 
\ep

\begin{lem}\label{prop:EFT2} Assuming Hypothesis~\ref{hyp:chiral}, the restriction of a degree~$n$ field theory along~\eqref{eq:twistedrestriction2} equips the resulting degree~$n$ representation of $\Ann^{2|1}(M)$ with a pairing. 
\end{lem}

\bp The datum~\eqref{eq:sPpairing} follows from restricting $E(\sL^-_M)$ in Corollary~\ref{cor:degreenvalues} along~\eqref{eq:constantpathsubspace}. The property~\eqref{eq:adjunction} follows from~\eqref{eq:bimoduleadjunction} in Proposition~\ref{prop:maincompositionstatement}. 
\ep

\subsection{Reflection and real structures}

The functor $\dagger$ from Lemma~\ref{lem:structureonsP} determines a reflection structure on $\Ann^{2|1}(M)$ in the sense of Definition~\ref{defn:interreflection}. 

\begin{lem} A reflection structure for the Clifford twist of $\Ann^{2|1}(M)$ in the sense of Definition~\ref{defn:RPfancyrep} is determined by a $*$-superalgebra structure $*\colon \cCl_{\pm n}\to \overline{\cCl}_{\pm n}^\op=\overline{\cCl}_{\mp n}$ with the property
\beq\label{eq:reflectandflcompat}
\overline{\alpha}_\pm \circ (-)^*=(-1)^\F\circ (-)^*\circ \alpha_\mp,
\eeq
for parity involution $(-1)^\F$ and the anti-automorphism $\alpha_\pm$ of $\cCl_{\pm n}$, see~\eqref{eq:Cliffinv0}. In particular, a reflection structure for the Clifford twist determines a real structure $\cCl_{\pm n}\to \overline{\cCl}_{\pm n}$ on the Clifford algebras~$\cCl_{\pm n}$. 
\end{lem}
\bp
Unpacking Definition~\ref{defn:RPfancyrep}, a reflection structure for the twist $\Cl_n\colon \Ann^{2|1}(M)\to \TA$ is a $\Z/4$-equivariant isomorphism of trivial bundles of algebras 
$$
\cCl_{\pm n}\coprod \cCl_{\mp n} \simeq \dagger^*(\overline{\cCl}_{\mp n}\coprod \overline{\cCl}_{\pm n})
$$
over $\Ob(\tAnn^{2|1}(M))\simeq \Map(\R^{0|1},M)\times\R_{>0}\coprod \Map(\R^{0|1},M)\times \R_{>0}$ with additional involutive data. It suffices to consider the case $M=\pt$. An isomorphism of algebras 
\beq\label{eq:whatisstar}
\cCl_{\pm n}\simeq \overline{\cCl}_{\mp n}=\overline{\cCl}_{\pm n}^\op
\eeq
determines such an invertible bimodule, where the involutive data for the bimodule is inherited from an involutive property for~\eqref{eq:whatisstar}. Indeed, the above allows us to view $\cCl_n$ as a $\overline{\cCl}_{\pm n}^\op$-$\cCl_n$ bimodule. An involutive algebra isomorphism~\eqref{eq:whatisstar} is precisely the data of a $*$-superalgebra structure. By Lemma~\ref{lem:commuteuptoflip}, the $\Z/4$-equivariance property for this reflection structure is the compatibility~\eqref{eq:reflectandflcompat} between the $*$-superalgebra structure and the anti-involution~\eqref{eq:Cliffinv0}. Finally, we observe that the $*$-superalgebra structure and anti-involution combine to define a real structure $\cCl_{\pm n}\stackrel{\sim}{\to}\overline{\cCl}_{\pm n}$.
\ep

\begin{defn} Define a $*$-superalgebra structure (see Definition~\ref{defn:starsuper}) on $\cCl_n$ given by the $\C$-antilinear extension of the map on generators,
\beq\label{eq:superstarCl}
f_j^*=-if_j,\qquad e_j^*=ie_j.
\eeq
using the notation of generators for the Clifford algebra from~\eqref{eq:Clifford}. 
\end{defn}

\begin{rmk} Under the translation between $*$-superalgebras and $*$-structures on a $\Z/2$-graded algebras (e.g., see~\cite[page~91]{strings1}), \eqref{eq:superstarCl} corresponds to the standard $*$-structure~\eqref{eq:starCl} on the Clifford algebras. We observe that the standard anti-involution~\eqref{eq:Cliffinv0} and the $*$-superalgebra structure~\eqref{eq:superstarCl} satisfy the compatibility property~\eqref{eq:reflectandflcompat}. These structures combine to define the standard real structure $\cCl_{\pm n}\stackrel{\sim}{\to}\overline{\cCl}_{\pm n}$ determining the real Clifford algebras $\Cl_{\pm n}\subset \cCl_{\pm n}$. 
\end{rmk}

With the chosen reflection structure for the Clifford twist given by~\eqref{eq:superstarCl}, a reflection structure on a degree~$n$ representation $E$ of $\Ann^{2|1}(M)$ specifies isomorphisms of Clifford module bundle over $\Map(\R^{0|1},M)\times\R_{>0}$ (see Definition~\ref{defn:RPfancyrep})
\beq\label{eq:rsisos}
\widetilde{\HBR}_\mp \stackrel{\sim}{\to} \dagger^*\overline{\widetilde{\HBR}_\pm},
\eeq 
relative to the $*$-superalgebra structure~\eqref{eq:superstarCl} on $\cCl_{\pm n}$.

\begin{defn}\label{defn:rssP}
Suppose we are given a reflection structure on a degree~$n$ representation~$E$ of $\Ann^{2|1}(M)$ with a pairing~$\pair$. Then $E$ is \emph{reflection positive} if the pairing
$$
\langle-,-\rangle\colon \dagger^*\overline{\widetilde{\HBR}_+}\otimes_{\cCl_n} \widetilde{\HBR}_+\simeq \cCl_n\otimes_{\overline{\cCl}_n \otimes \cCl_n}(\dagger^*\overline{\widetilde{\HBR}_+}\otimes \widetilde{\HBR}_+)\simeq 
\cCl_n\otimes_{\cCl_{-n} \otimes \cCl_n}(\widetilde{\HBR}_-\otimes \widetilde{\HBR}_+)\xrightarrow{\pair}\underline{\C}
$$
 restricts along $M\times\R_{>0}\subset \Map(\R^{0|1},M)\times \R_{>0}$ to a positive hermitian pairing (in the $\Z/2$-graded sense), where the second isomorphism above uses~\eqref{eq:rsisos}. 
\end{defn}


\begin{prop} \label{prop:selfadjoint} 
For a degree~$n$ reflection positive representation $E$ of $\Ann^{2|1}(M)$, the associated 1-parameter family of $\cCl_{n}$-linear superconnections $\A_\ell$ is self-adjoint for each $\ell\in \R_{>0}$, the $S^1$-action is unitary, and the $\cCl_\pm$-action is self-adjoint with respect to the pairing $\langle-,-\rangle$ from Definition~\ref{defn:rssP}. 
\end{prop}

\bp
Relative to the isomorphisms~\eqref{eq:rsisos}, the reflection structure on $E$ requires that the super semigroup representations~\eqref{eq:defrhotil} satisfy
$$
\rho_\pm(\ell,\tau,\bar\tau,\eta)=(\dagger^*\overline{\rho}_\mp)(\ell,\tau,\bar\tau,\eta)=i^\deg\overline{\rho}_\mp(\ell,\bar\tau,\tau,\eta)
$$
Taylor expanding as in~\eqref{eq:AB}, we find that the superconnection $\A_\ell^-$ is sent to $-i^{\deg+1}\overline{\A}_\ell^+$ under the isomorphism~\eqref{eq:rsisos}. Then by the adjunction formula~\eqref{eq:initialadjunction} and the definition of the hermitian pairing, we obtain the desired self-adjointness property for $\A_\ell^+$. The unitary property for the $S^1$-action follows similarly. Self-adjointness of the Clifford action follows from the definition of the reflection structure for the Clifford twist using the $*$-superalgebra structure~\eqref{eq:superstarCl} and the definition of the pairing $\langle-,-\rangle$. 
\ep

\begin{prop} \label{prop:realA} 
A reflection positive degree~$n$ representation $E$ of $\Ann^{2|1}(M)$ equips the vector bundle $\HBR=\HBR_+$ with a real structure such that the pairings in Definition~\ref{defn:rssP} comes from ordinary (i.e., not $\Z/2$-graded) real pairings on the underlying real bundle of $\Cl_n$-modules and the 1-parameter family of self-adjoint superconnections $\A_\ell=\A_\ell^+$ is a real family of superconnections.
\end{prop}
\bp
The isomorphism~\eqref{eq:Z4gen} and the reflection structure~\eqref{eq:rsisos} combine to give an isomorphism
\beq\label{eq:Realstructure}
\widetilde{\HBR}_\pm\xrightarrow\sim \dagger^*(\overline{\sqrt{\Fl}})^*\overline{\widetilde{\HBR}_\pm}
\eeq
that when restricted to $M\times \R_{>0}\subset \Map(\R^{0|1},M)\times \R_{>0}$ is an $\R_{>0}$-family of real structures~$\HBR_\pm\simeq \overline{\HBR}_\pm$, relative to the real structure on $\cCl_n$ specified by the $*$-structure~\eqref{eq:superstarCl} and anti-automorphism~\eqref{eq:Cliffinv}. Below, let $x\mapsto \overline{x}$ denote the image of a section under this family of real structures. 
We claim that the pairing satisfies the condition
\beq\label{eq:realcondition}
\langle x,y\rangle=(-1)^{|\overline{x}|} \overline{\langle \overline{x},\overline{y}\rangle}
\eeq
for $x,y\in \Gamma(\HBR_+)$. This follows from~\eqref{eq:compositionsofstructure2}: commuting $\dagger$ and $\sqrt{\Fl}$ introduces a sign depending on parity. The condition~\eqref{eq:realcondition} implies that the pairing is real for pairings between even sections, zero for pairings between even and odd sections, and purely imaginary for pairings between odd sections. The translation between a super Hilbert space and a $\Z/2$-graded Hilbert space is (compare~\cite[page~91]{strings1})
\beq\label{Eq:metric}
\langle v,w\rangle_0:=\left\{\begin{array}{ll} \langle v,w\rangle & v,w \in \Gamma(\HBR_+)^\ev \\ i^{-1}\langle v,w\rangle & v,w\in \Gamma(\HBR_+)^\odd \\ 0 & {\rm else} \end{array}\right.
\eeq
and so we see that the condition~\eqref{eq:realcondition} leads to a real pairing in the usual sense under this translation. Finally, the super semigroup representations~\eqref{eq:defrhotil} must pull back to themselves along the real structure~\eqref{eq:Realstructure} 
\beq\label{eq:realdataonC}
\rho_\pm(\ell,\tau,\bar\tau,\eta)=(\dagger^*(\overline{\sqrt{\Fl}})^* \overline{\rho}_\pm)(\ell,\tau,\bar\tau,\eta).
\eeq
In particular, the Taylor component $\b_+(\ell,0,0)$ (as in~\eqref{eq:AB}) pulls back to itself. This is precisely the condition that the superconnection is real. 
\ep

\begin{prop} \label{prop:EFT3}
Assuming Hypothesis~\ref{hyp:chiral}, the restriction of a reflection structure on a degree~$n$, $2|1$-Euclidean field theory along~\eqref{eq:twistedrestriction2} determines a reflection structure on a degree~$n$ representation of $\Ann^{2|1}(M)$. The restriction of a reflection positive degree~$n$, $2|1$-Euclidean field theory along~\eqref{eq:twistedrestriction2} determines a reflection positive degree~$n$ representation of~$\Ann^{2|1}(M)$. 
\end{prop}
\bp
Hypothesis~\ref{hyp:chiral} implies that the reflection structure for the degree~$n$ twist comes from the $*$-superalgebra structure on $\cCl_n$. If we then further examine the data in Proposition~\ref{lem:RPtwist} and Corollary~\ref{cor:degreenvalues}, the restriction of a reflection positive field theory along~\eqref{eq:constantpathsubspace} determines a reflection positive degree~$n$ representation for the reflection structure on $\Ann^{2|1}(M)$ given by $\dagger$. The symmetry of the pairing in Definition~\ref{defn:rssP} follows from Lemma~\ref{lem:hermitian}. Finally, we observe that the positivity condition in Definition~\ref{defn:RPFT} is equivalent to the positivity condition in Definition~\ref{defn:rssP}. 
\ep

\subsection{Trace class representations} \label{sec:traceclassreps}

We recall the super Euclidean double loop space $\mathcal{L}^{2|1}(M)$ from Definition~\ref{defn:doubleloop} with its action by $\E^{2|1}\rtimes\Spin(2)\times \SL_2(\Z)$ and the span~\eqref{eq:supercyltosuperloop} that regards a subspace of supercylinders in $M$ as supertori in $M$. Define the subspace of nearly constant super double loops,
\beq\label{eq:constsuperdouble}
&&\mathcal{L}_0^{2|1}(M):=\Map(\R^{0|1},M)\times s\Lat\subset \Map(\R^{2|1}/\Z^2,M)\times s\Lat=\mathcal{L}^{2|1}(M).
\eeq
The action by $\E^{2|1}\rtimes\Spin(2)\times \SL_2(\Z)$ preserves the inclusion~\eqref{eq:constsuperdouble}. The assignment $M\mapsto \mathcal{L}_0^{2|1}(M)\sq (\E^{2|1}\rtimes\Spin(2)\times \SL_2(\Z))$ is natural in the manifold $M$, and hence there is a canonical map to the stack of super Euclidean tori,
\beq\label{eq:sloopover}
&&\mathcal{L}_0^{2|1}(M)\sq (\E^{2|1}\rtimes \Spin(2)\times \SL_2(\Z)) \to s\Lat\sq (\E^{2|1}\rtimes \Spin(2)\times \SL_2(\Z) ).
\eeq
Over the target of the above map, we take the extensions of the Pfaffian line bundle specified by Proposition~\ref{prop:chiral2}. Explicitly, this line bundle is classified by the map of stacks
\beq\label{eq:superPf}
&&\sLat\sq \Spin(2)\times \SL_2(\Z)\to \pt\sq U(1), \quad \Spin(2)\times \SL_2(\Z)\to U(1),
\eeq
determined by the homomorphism specified in~\eqref{eq:Pfhomo}. 
We use the same notation $\Pf$ for the pullback of the line bundle classified by~\eqref{eq:superPf} along the map of stacks~\eqref{eq:sloopover}. 
To give an explicit description of the sections of $\Pf$ over $\mathcal{L}_0^{2|1}(M)$, we adopt notation 
$$
C^\infty_{j/2,k/2}(\Lat):=\{f\in C^\infty(\Lat) \mid (u,\bar u)\cdot f=u^j\bar u^k f,\ (u,\bar u)\in \Spin(2)\}\subset \C^\infty(\Lat)
$$
for functions on $\Lat$ that transform with weight $(j,k)$ for the $\Spin(2)\simeq U(1)$-action on $\Lat$. 

\begin{rmk}
The $\SL_2(\Z)$-invariant elements $C^\infty_{j/2,k/2}(\Lat)^{\SL_2(\Z)}$ determine global sections of $\Pf^{-j/2}\otimes \overline{\Pf}{}^{-k/2}$ where $\Pf$ and $\overline{\Pf}$ are the Pfaffian line bundle and its conjugate over the moduli stack of (non-super) Euclidean spin tori; see \S\ref{sec:back1}. 
\end{rmk}

\begin{proof}[Proof of Proposition~\ref{prop:maintorusprop}]
There is an injection
$$
\Gamma(\mathcal{L}_0^{2|1}(M)\sq (\E^{2|1}\rtimes \Spin(2)\times \SL_2(\Z));\Pf^{\otimes n})\hookrightarrow C^\infty(\Map(\R^{0|1},M)\times \sLat)
$$ 
that realizes sections as functions on $\Map(\R^{0|1},M)\times \sLat$ that transform appropriately for the action by $\E^{2|1}\rtimes \Spin(2)\times \SL_2(\Z)$. Specifically, sections are invariant under the $\E^{2|1}\times \SL_2(\Z)$-action and transform with weight for the the $\Spin(2)\simeq U(1)$-action. Using the identification,
$$
C^\infty(\Map(\R^{0|1},M)\times \sLat)\simeq C^\infty(\sLat;C^\infty(\Map(\R^{0|1},M)))\simeq C^\infty(\Lat;\Omega^\bullet(M)[\lambda_1,\lambda_2]/(\lambda_1\lambda_2))
$$
that Taylor expands in the odd variables $\lambda_1,\lambda_2\in C^\infty(\sLat)$, any function can be written as
\beq
\omega&=&(2\vol)^{\deg/2}(\omega_0+\lambda_1(2\vol)^{1/2}\omega_1+\lambda_2(2\vol)^{1/2}\omega_2),\label{eq:itsomega}\\
&& \omega_i\in C^\infty(\Lat;\Omega^\bullet(M))\simeq \Omega^\bullet(M;C^\infty(\Lat))\nonumber
\eeq
where $\vol:=\frac{\ell_1\bar\ell_2-\bar\ell_1\ell_2}{2i}\in C^\infty(\Lat)^{\SL_2(\Z)}$ is a real-valued nonvanishing function. By \cite[Lemma~3.19]{DBEChern}, the function $\omega$ is $\E^{2|1}$-invariant if and only if\footnote{Equation~\eqref{eq:initialE21invar} differs from \cite[Lemma~3.19]{DBEChern} by a factor of $1/2$ because the normalization for $\omega$ chosen above has an extra factor of $2^{\deg/2}$, which will be more convenient below.}
\beq\label{eq:initialE21invar}
\dR\omega_0=0,\quad \partial_{\bar \ell_i}\omega_0=\frac{1}{2}\dR\omega_i,\quad i=1,2.
\eeq
Next we consider $\Spin(2)$-equivariance: for $(u,\bar u)\in U(1)\simeq \Spin(2)$ we will use that
$$
(u,\bar u)\cdot\vol^{1/2}=u\bar u\vol^{1/2}, \ (u,\bar u)\cdot\lambda_i=\bar u \lambda_i, \ (u,\bar u)\cdot \alpha=\bar u^{-k} \alpha, \ \alpha\in \Omega^k(M)\subset \Omega^k(M;C^\infty(\Lat))
$$
where in the last calculation, we view a $k$-form as a differential form as an element of $\Omega^k(M;C^\infty(\Lat))$ via the inclusion of constant functions $\C\hookrightarrow C^\infty(\Lat)$. 
From the above, we deduce that for $\omega_i$ to be in the $n^{\rm th}$ weight space of of the $U(1)$-action, we require 
$$
\omega_0\in \bigoplus_{k+j=n} \Omega^k(M;C^\infty_{j/2,0}(\Lat)),\quad \omega_i\in \bigoplus_{k+j=n-1} \Omega^k(M;C^\infty_{j/2,1}(\Lat)).
$$
We further observe that $\omega_0$ is required to be $\SL_2(\Z)$-invariant, so by Lemma~\ref{lem:euctoribackground} and \eqref{eq:powerofell}, restriction along $\HHz\times \R_{>0}\hookrightarrow \Lat$ gives 
\beq
Z&:=&\ell^{(\deg-n)/2} (\omega_0|_{\HHz\times \R_{>0}})=(\ell^{(\deg-n)/2}(2\vol)^{-\deg/2} \omega)|_{\iHH\times \R_{>0}}\nonumber\\
&=&\ell^{-n/2}(2\im(\tau))^{-\deg/2}(\omega|_{\HHz\times \R_{>0}})\in \bigoplus_{j+k=n} \Omega^j(M;\mmf^{k,0})\label{eq:defnofZ}
\eeq
where in the second line we used that $\vol|_{\HHz\times \R_{>0}}=\im(\tau)\ell$. Unpacking the notation from~\S\ref{sec:backtot}, $Z$ is a differential form valued in $C^\infty(\HHz\times \R_{>0})$ with the transformation property
\beq\label{eq:modulartransformZ}
Z\mapsto (cz+d)^{(\deg-n)/2}Z,\qquad z=\tau/\ell
\eeq
for the $\MP_2(\Z)$-action on coefficients. We note the diffeomorphism
$$
\HHz\times\R_{>0}\xrightarrow{\sim} \HHz\times \R_{>0}, \qquad (\tau,\bar\tau,\ell)\mapsto (\tau/\ell,\bar\tau/\ell,\ell\im(\tau))
$$
giving coordinates on $\HHz\times \R_{>0}$ in terms of the conformal modulus~$(z,\bar z)=(\tau/\ell,\bar\tau/\ell)$ and total volume~$v=\ell\im(\tau)$ of the associated torus. By \cite[Lemma~3.21]{DBEChern}, the data of $\omega_1,\omega_2$ are equivalent to
\beq
Z_v&\in& \bigoplus_{j+k=n-1} \Omega^{j}(M;\mmf^{(k-2),-2})\label{eq:Zv}\\
Z_{\bar z}&\in& \bigoplus_{j+k=n-1} \Omega^{j}(M;\mmf^{k,-2} )\label{eq:Ztau}
\eeq
i.e., differential forms with values in $C^\infty(\HHz\times \R_{>0})$ with transformation properties 
$$
Z_{\bar z}\mapsto (cz+d)^{(\deg-n-1)/2}(c\bar z+d)Z_{\bar z},\qquad Z_v\mapsto (cz+n)^{(\deg-n-3)/2}(c\bar z+d)Z_v
$$
for the $\MP_2(\Z)$-action on coefficients. By \cite[Lemma~3.21]{DBEChern}, property~\eqref{eq:initialE21invar} is equivalent to 
$$
\partial_{\bar z} Z=\dR Z_{\bar z},\qquad \partial_v Z=\dR Z. 
$$
This completes the proof. 
\ep


The element $\Gamma_n\in \cCl_n$ is defined in~\eqref{eq:Gamma} and determines the Clifford supertrace~\eqref{eq:Cliffordsupertrace}.

\begin{defn}\label{defn:tracestru}
A \emph{trace class} degree~$n$ representation of $\Ann^{2|1}(M)$ is the data of a degree~$n$ representation $E$ of $\Ann^{2|1}(M)$ and a section $(Z,Z_\tau,Z_v)\in \Gamma(\mathcal{L}^{2|1}_0(M);\Pf^{\otimes n})$. We require these data to satisfy the properties: 
\begin{enumerate}
\item[(i)] the (ordinary) trace exists 
\beq\label{eq:ordtrace}
\Tr(\Gamma_n\circ \rho_+|_{M\times (\R_{>0}\times \iHH_+)/\Z})\in C^\infty(M;C^\infty((\R_{>0}\times \iHH_+)/\Z))
\eeq
and 
\item[(ii)] the supertrace has the compatibility condition
\beq\label{eq:traceclass}\label{eq:constpartitionfun2}
&&\resizebox{.9\textwidth}{!}{$Z=\phi(\ell,z)^{-n}(2\im(z))^{-\deg/2}\Tr_{\cCl_n}(\rho_+|_{\Map(\R^{0|1},M)\times(\R_{>0}\times \iHH_+)/\Z})\in \Omega^\bullet(M;(C^\infty(\R_{>0}\times \iHH_+)/\Z))$}
\eeq
\end{enumerate}
where $\Tr_{\cCl_n}$ denotes the Clifford supertrace and 
$$
\phi(\ell,\tau)=\prod_{n>0} (1-e^{2\pi i \tau/\ell})\in C^\infty(\R_{>0}\times \iHH_+)
$$
is closely related to the Dedekind $\eta$-function. 
\end{defn}

\begin{rmk}
The condition~\eqref{eq:constpartitionfun2} is equivalent to requiring $Z$ to be equal to the Clifford supertrace of the super semigroup representation in $\Fer^{\pm n}_\ell$-modules, see~\eqref{eq:Fermodextend}. 
\end{rmk}

\begin{prop} \label{prop:trace}
Assuming Hypothesis~\ref{hyp:chiral}, the restriction of a degree~$n$ field theory along~\eqref{eq:twistedrestriction2} determines a trace class degree~$n$ representation of $\Ann^{2|1}(M)$.
\end{prop} 

\begin{proof}
The argument follows \cite[Lemma~3.20]{ST11}. In the terminology of \cite{STTraces}, the bordisms in the image of the composition
\beq\label{eq:thickfamily}
\Map(\R^{0|1},M)\times(\R_{>0}\times \iHH_+)/\Z\hookrightarrow \Mor(\Ann^{2|1}(M))\to \Mor(2|1\EBord(M))
\eeq
are \emph{thick} morphisms in the bordism category: this follows directly from the factorization in Proposition~\ref{prop:duality}. The categorical trace of a thick morphism is an endomorphism of the monoidal unit. Proposition~\ref{prop:tracerelations} shows that this categorical trace is the image of the composition
\beq\label{eq:imageofthick}
&&\Map(\R^{0|1},M)\times(\R_{>0}\times \iHH_+)/\Z\to \mathcal{L}_0^{2|1}(M)\to \Mor(2|1\EBord(M)).
\eeq
Applying a degree~$n$ field theory to these families, Corollary~\ref{cor:degreenvalues} shows that the categorical trace of the family~\eqref{eq:thickfamily} is sent to the composition 
\beq\label{eq:itsCliffodsupertr}
&&\cCl_n\otimes_{\cCl_n\otimes \cCl_n} \cCl_n\xrightarrow{\id_{\cCl_n}\otimes {\sf r}_M^+}\cCl_n\otimes_{\cCl_n\otimes \cCl_n} (\widetilde{\HBR}_+\otimes \widetilde{\HBR}_-)\xrightarrow{\sigma}\cCl_n\otimes_{\cCl_n\otimes \cCl_n} (\widetilde{\HBR}_-\otimes \widetilde{\HBR}_+) \xrightarrow{\mathring{\sf l}_M^-} \underline{\C},
\eeq
for the families of small bordisms ${\sf r}_M^+$ and ${\sf l}_M^-$ from~\eqref{eq:LM0}, where $\mathring{\sf l}_M^-$ denotes the restriction of ${\sf l}_M^-$ along $0 \in \HH_-^{2|1}$. 
The element $\Gamma_n$ provides an identification $\cCl_n\otimes_{\cCl_n\otimes \cCl_n} \cCl_n\simeq \underline{\C}$ with the trivial line bundle, and the composition~\eqref{eq:itsCliffodsupertr} becomes the Clifford supertrace. Using that nuclear Fr\'echet spaces have the approximation property~\cite[page~109]{Schaefer}, we find that the categorical trace above is independent of the choice of thickener~\cite[Theorem 1.7]{STTraces}. Finally, Lemma~\ref{lem:Tracenormalize1} translates between the Clifford trace and the trace associated to $\Fer_n$, giving the property~\eqref{eq:constpartitionfun2}. The argument for the trace property~\eqref{eq:ordtrace} is similar, but one modifies the outgoing boundary of the supercylinder in the image of~\eqref{eq:imageofthick} by the spin flip automorphism. This changes the supertrace to an ordinary trace. 

To recover the claimed compatibility formula~\eqref{eq:constpartitionfun2}, the above argument gives
\beq\label{eq:constpartitionfun3}
&&\omega|_{\HHz\times \R_{>0}}=\ell^{n/2}\phi(\ell,\tau)^{-n}\Tr_{\cCl_n}(\rho_+|_{\Map(\R^{0|1},M)\times(\R_{>0}\times \iHH_+)/\Z})
\eeq
for $\omega$ the function in~\eqref{eq:itsomega} determining the section of $\Pf^{-\otimes n}$. The prefactor before the Clifford supertrace in~\eqref{eq:constpartitionfun3} is the $-n$th power of the character of the invertible bimodule~\eqref{eq:Morita} as a function on $\R_{>0}\times \iHH_+$, see~\eqref{eq:traceformula}. Then by~\eqref{eq:defnofZ}, we obtain~\eqref{eq:constpartitionfun2}. 
%
%
%
\ep

\section{Energy cutoffs and the cocycle map}\label{sec:cocyclemap}

In this section we use the notation~\eqref{eq:singleVA} for the super vector bundle $\HBR$ with 1-parameter family of superconnections $\A_\ell$ associated with a degree~$n$ representation of $\Ann^{2|1}(M)$.  For $\ell\in \R_{>0}$, we also set $\rho_\ell(\tau,\bar\tau,\eta):=\rho_+(\ell,\tau,\bar\tau,\eta)$, i.e.,  $\rho_\ell$ is the restriction of $\rho_+$ of~\eqref{eq:defrhotil} to a fixed~$\ell\in \R_{>0}$. 

\subsection{Extracting a sequence of superconnections} \label{sec:sequenceofsuperconn}

\begin{lem} \label{lem:weightspaces}
A degree~$n$, reflection positive representation of $\Ann^{2|1}(M)$ (in the sense of Definition~\ref{defn:rssP}) determines an orthogonal direct sum of metrized super vector bundles 
\beq\label{eq:weightspaces}
\HBR \simeq \bigoplus_{k\in \Z} \HBR_k,\qquad \HBR_k\to M,
\eeq
with respect to which the action of $\HH^{2|1}_+/\ell\Z$ is given by 
\beq\label{eq:weightspacesupersemi}
\rho_\ell(\tau,\bar \tau,\eta)=\bigoplus_{k\in \Z} e^{2\pi i k\tau/\ell} \rho_{\ell,k}(\im(\tau),\theta),
\eeq
where for each $(\ell,k)\in \R_{>0}\times \Z$, $\rho_{\ell,k}\colon \R^{1|1}_{\ge 0}\to \End(\Omega^\bullet(M;\HBR_{k,\ell}))$ is a $\cCl_{-n}$-linear super semigroup representation. If the representation is also trace class, then $\HBR_k=\{0\}$ for $k\ll 0$. 
\end{lem}
\bp
The exact sequence~\eqref{eq:itsexact} of super Lie groups gives exact sequences of super semigroups for each $\ell\in \R_{>0}$, 
\beq\label{eq:exactseq}
\begin{array}{ccccccccc}
1&\to &\R/\ell\Z&\to& \HH_+^{2|1}/\ell\Z&\to& \R^{1|1}_{\ge 0}&\to& 1,\\
&&t&\mapsto & (t,t,0), \ (\tau,\bar \tau,\eta) &\mapsto & (2\im(\tau),\eta).
\end{array}
\eeq
By Proposition~\ref{prop:selfadjoint} the action of $\R/\ell\Z$ on $\HBR$ is a unitary, and hence we obtain~\eqref{eq:weightspaces} as an orthogonal direct sum of weight spaces. For $p_k\colon \HBR\to \HBR_k$ the orthogonal projection onto a weight space, define the super semigroup representation 
\beq\nonumber
\rho_{\ell,k}\colon \HH^{2|1}_+/\ell\Z&\to& \End(\Omega^\bullet(M;\HBR_{k,\ell}))\\ 
\rho_{\ell,k}(\tau,\bar \tau,\eta)&:=& e^{-2\pi i k\tau/\ell}\otimes (p_k\circ \rho_\ell (\tau,\bar \tau,\eta)\circ p_k)\label{eq:weightspacerep}
\eeq
where $e^{-2\pi i k\tau/\ell}\colon \HH^{2|1}_+/\ell\Z\to \End(\C)$ a 1-dimensional representation. The subgroup $\R/\ell\Z\subset \HH^{2|1}_+/\ell\Z$ acts trivially, and so the exact sequence~\eqref{eq:exactseq} shows that~\eqref{eq:weightspacerep}  is determined by a representation of the quotient. This gives the desription~\eqref{eq:weightspacesupersemi}. 

Finally, the trace class condition implying $\HBR_k=\{0\}$ for sufficiently negative $k$ follows from a standard argument, e.g., see \cite[page~43]{ST11} or \cite[Proof of Theorem~3.3.14]{ST04}. In brief, $\rho_{\ell,k}$ needs to have trace zero for $k\ll 0$ if the trace~\eqref{eq:ordtrace} is to converge. Continuity of this trace in the parameter $(\tau,\bar \tau,0)$ implies that the identity operator $\rho_{\ell,k}(0,0,0)=\id_{\HBR_{\ell,k}}$ is the zero operator for $k\ll 0$, and hence $\HBR_k=\{0\}$ for $k\ll 0$. 
\ep

\begin{lem}
The 1-parameter family of superconnections from Proposition~\ref{prop:superconn} decomposes with respect to~\eqref{eq:weightspaces} as an infinite direct sum of self-adjoint, real Clifford linear superconnections,
\beq
\A_\ell=\bigoplus_{k\in \Z} \A_{\ell,k},\qquad \A_{\ell,k}\colon \Omega^\bullet(M;\HBR_k)\to \Omega^\bullet(M;\HBR_k).\label{eq:sequenceofsuperconns}
\eeq
\end{lem}
\bp
The infinite direct sum of superconnections comes from Taylor expanding the right hand side of~\eqref{eq:weightspacesupersemi} in the odd variable $\eta$ and applying the argument in the proof of Proposition~\ref{prop:superconn}, where self-adjointness follows from Proposition~\ref{prop:selfadjoint} and the real structure is from Proposition~\ref{prop:realA}. 
\ep

\begin{defn}\label{defn:smoothindexbundleEFT} A reflection positive degree~$n$ representation of $\Ann^{2|1}(M)$ \emph{admits smooth index bundles} if for all $(\ell,k)\in \R_{>0}\times \Z$, the superconnections~\eqref{eq:sequenceofsuperconns} admit smooth index bundles in the sense of Definition~\ref{defn:cutoffs}. \end{defn}

\begin{prop} \label{prop:Ktate}
A reflection positive degree~$n$ representation $E$ of $\Ann^{2|1}(M)$ that admits smooth index bundles determines a class
\beq\label{eq:Eindex}
[E]\in \KO^n(M)(\!(q)\!).
\eeq
\end{prop}
\bp
We claim it suffices to consider the case that $\mathcal{E}_k=\{0\}$ for $k<0$. To see this, multiply the representation~\eqref{eq:weightspacesupersemi} by $q^J$ for a sufficiently large power so that the resulting representation $q^JE$ satisfies $\mathcal{E}_k=\{0\}$ for $k<0$. The argument below then constructs an index $[q^JE]\in \KO^n(M)[\![q]\!]$. Finally, define~\eqref{eq:Eindex} as $[E]:=q^{-J}[q^JE]\in \KO^n(M)(\!(q)\!).$ Such a global power $J$ may not exist if~$M$ is not compact, but by choosing an exhaustion function on~$M$ and constructing compatible classes on compact subsets we obtain the same reduction. 

Below, let $\A_k:=\A_{1,k}$, i.e., the $k$th 1-parameter family of superconnections~\eqref{eq:sequenceofsuperconns} specialized to $\ell=1$. We apply the index bundle construction to the sequence of superconnections $\A_k$ analogously to the proof Lemma~\ref{lem:KOqmap}. For each $N$, define an open cover of $X$ indexed by $\vec{\lambda}=(\lambda_0,\dots,\lambda_N)\in \R_{>0}^N$
\beq\nonumber
U^{\vec{\lambda}}&:=&\{x\in X\mid \lambda_k\notin{\rm Spec}((\A_{k}^{[0]})^2_x)\}\subset X \ {\rm for} \ k\le N\}.
\eeq
Over each component of this cover one has finite-rank vector bundles $\HBR_k^{<\lambda_k}\to U^{\vec{\lambda}}$, and on overlaps $U^{\vec{\lambda}}\bigcap U^{\vec{\lambda'}}$ there are odd endomorphisms for $k\le N$,
$$
f_k^{\lambda_k\lambda_k'}\in \End(U^{\vec{\lambda}}\bigcap U^{\vec{\lambda'}},(\HBR_k^{<\lambda})^\perp),\qquad \lambda_k<\lambda_k'
$$
defined on the orthogonal complement to the inclusion $\HBR_k^{<\lambda}\hookrightarrow \HBR_k^{<\lambda'}$. Hence, for each $N$ Proposition~\ref{prop:appenindexbundle} gives a map 
\beq\label{eq:maptoholim2}
\left[\{U^{\vec\lambda},\bigoplus_{k=0}^N \HBR_k^{<\lambda_k},\bigoplus_{k=0}^Nf_k^{\lambda_k\lambda_k'}\}_{\vec\lambda}\right]\colon X\to (\KO\wedge \{1,q,\dots q^N\}_+ )^{n}
\eeq
determined by a topological functor (see~\eqref{eq:constructtoKO}) from the \v{C}ech category of $\{U^{\vec{\lambda}}\}$ to a topological category whose classifying space represents the $n$th space in the $N$th spectrum in the homotopy colimit diagram~\eqref{eq:KTateascolim}. A refinement of an open cover determines a functor between \v{C}ech categories inducing an equivalence; so by \cite[Proposition~2.1]{Segalclassifying} the maps constructed above come equipped with canonical homotopy compatibilities for increasing~$N$. Therefore, the maps~\eqref{eq:maptoholim2} are part of the data of  a map into  the homotopy colimit diagram~\eqref{eq:KTateascolim} and hence defines a map into the spectrum~$\KO[\![q]\!]$.
\ep

\subsection{Modular invariance and the Pontryagin character}

\begin{defn}\label{defn:infgen}
A self-adjoint representation of $\Ann^{2|1}(M)$ is \emph{infinitesimally generated} if $\rho_\ell(\tau,\bar \tau,\eta)$ is determined by the formula
\beq\label{eq:infgenerator}
\rho_\ell(\tau,\bar \tau,\eta)=\bigoplus_{k\in \Z} e^{2\pi i k \tau/\ell} \rho_{\ell,k}(2\im(\tau),\eta)=\bigoplus_{k\in \Z} e^{2\pi i k \tau/\ell} e^{-2\im(\tau) \A_{\ell,k}^2+\eta \A_{\ell,k}}
\eeq
for the sequence of self-adjoint, Clifford linear superconnections,
$$
\A_{\ell,k}:=p_k\circ \A_\ell\circ p_k
$$
where $\A_\ell$ is the 1-parameter family of self-adjoint superconnections on $\HBR$ from~\eqref{eq:singleVA} and $p_k\colon \HBR\to \HBR_k$ is the orthogonal projection.
\end{defn}

\begin{rmk}\label{rmk:ODEfail}
Assuming existence and uniqueness of solutions to ordinary differential equations, a self-adjoint representation of $\Ann^{2|1}(M)$ is automatically infinitesimally generated, e.g., if the vector bundles $\HBR_k$ are finite-dimensional for all~$k$. However, existence and uniqueness can fail in general topological vector spaces, see \cite[Example~5.6.1]{Hamilton} for counterexamples in Fr\'echet spaces. From the physics point of view, one expects field theories to be determined by infinitesimal data. This is essential if one is to make sense out of physical quantities like energy as arising from eigenvalues of operators, which is a prerequisite for discussing energy cutoffs in effective field theory. 
To force existence and uniqueness, one could modify the functional analytic framework for field theories, e.g., to tame Fr\'echet spaces and tame maps, or some variant of Hilbert bundles and continuous maps. A more Wilsonian approach would be to consider representations in Ind-finite vector bundles with compatibilities via the energy filtration. Such conditions quickly become technical, and the correct analytical context for extending down to points in conjecture~\eqref{eq:conjecture} is not yet clear. As such, we prefer to leave these analytical choices flexible, whence Definition~\ref{defn:infgen}. 
\end{rmk}

\begin{prop} For a sequence of superconnections extracted from an infinitesimally generated trace class degree~$n$ representation of $\Ann^{2|1}(M)$, the normalized Clifford supertrace
$$
Z=\phi(\ell,\tau)^{-n}\sum q^k \sTr_{\cCl_n}(e^{-(\A_{\ell,k}(2\im(\tau)))^2})\in \Omega^\bullet(M;C^\infty(\iHH\times \R_{>0}))
$$
transforms as~\eqref{eq:modulartransformZ} for the $\MP_2(\Z)$-action on $\iHH_+\times \R_{>0}$, i.e., transforms as a section of~$\Pf^{-\otimes n}$. In the above, $\A_{\ell,k}(2\im(\tau))$ is the rescaling~\eqref{eq:rescale} of the  superconnection $\A_{\ell,k}$ by~$2\im(\tau)$. 

\end{prop}
\bp
Using the definition of a trace class representation, the equality~\eqref{eq:constpartitionfun2} from Proposition~\ref{defn:tracestru} yields
\beq
Z&=&(2\im(\tau))^{-\deg/2}\phi(\ell,\tau)^{-n}\sTr_{\cCl_n}(\rho_+(\tau,\bar \tau,0))\nonumber\\
&=&\phi(\ell,\tau)^{-n}\sum q^k (2\im(\tau))^{-\deg/2}  \sTr_{\cCl_n}(e^{-2\im(\tau)\A_{\ell,k}^2}) \nonumber\\
&=&\phi(\ell,\tau)^{-n}\sum q^k  \sTr_{\cCl_n}(e^{-(\A_{\ell,k}(2\im(\tau))^2}), \label{eq:thisishow}
\eeq
where $Z$ is modular by~\eqref{eq:defnofZ}. This completes the proof. 
\ep

\subsection{Homotopy compatibility of the index bundle and the partition function}

\begin{lem}\label{lem:1sthomotopyforcoccyle}
There is a uniquely determined coboundary between the Pontryagin character of the cocycle constructed in Proposition~\ref{prop:Ktate} and the differential form 
\beq\label{eq:1sthomotopyforcoccyle}
\phi(1,\tau)^{-n}\sum  q^k \sTr_{\Cl_n}(e^{-\A_k(2\im(\tau))^2})\in \Omega^\bullet(M;C^\infty(\R_{>0})(\!(q)\!))=Z|_{\ell=1}
\eeq
where above we identify the rescaling paramter $2\im(\tau)\in C^\infty(\R_{>0})\subset C^\infty(\R_{>0})(\!(q)\!)$ with the standard coordinate. 
\end{lem}
\bp
The existence of index bundles for $\A_k$ implies that for each $(\lambda,k)$ there are uniquely determined coboundaries (see Lemma~\ref{lem:anotherdamnhomotopy})
\beq\label{eq:somecoboundaries2}
&&d\omega_k^\lambda=\sTr_{\Cl_n}(e^{-\A_k(2\im(\tau))^2})|_{U_\lambda}-\Ch(\nabla^{<\lambda_k})\in \Omega(U_{\lambda_k};\C[u^{\pm 2}]).
\eeq
The Pontryagin character of the cocycle constructed in Proposition~\ref{prop:Ktate} is a \v{C}ech--de~Rham cocycle for each $N$ (compare~\eqref{eq:limitofCherns}) with compatibilities between different $N$ coming from refinements of the cover. For each $N$, restricting the coboundaries~\eqref{eq:somecoboundaries2} along a refinement of covers gives a coboundary in the \v{C}ech--de~Rham complex between the Pontryagin character of the index bundle and~\eqref{eq:1sthomotopyforcoccyle}.
\ep

\begin{lem}\label{lem:2ndhomotopyforcoccyle}
Given an infinitesimally generated, self-adjoint representation of $\Ann^{2|1}(M)$, there are Chern--Simons form associated to the $\ell$-parameter family of superconnections, giving
$$
\sum q^k  \sTr_{\cCl_n}(e^{-(\A_{\ell,k}(2\im(\tau))^2})=\sum q^k  \sTr_{\cCl_n}(e^{-(\A_{1,k}(2\im(\tau))^2})+d \sum q^k\zeta_k \in \Omega^\bullet(M;C^\infty(\R_{>0})(\!(q)\!))
$$
\end{lem} 
\bp
The supertrace $\sTr(e^{-\A_{\ell,k}(2\im(\tau))^2})$ varies smoothly with $\ell$. The lemma then follows from the standard construction of Chern--Simons forms applied to the 1-parameter family of superconnections~$\A_{\ell,k}$ for the parameter $\ell\in \R_{>0}$, 
\beq
\zeta_k(\ell)&=&\int_{1}^{\ell} \sTr_{\cCl_n}\left(\frac{d\A_{s,k}(2\im(\tau))}{ds} e^{-\A_{s,k}(2\im(\tau))^2}\right)ds\in \Omega^\bullet(M;C^\infty(\R_{>0}))\nonumber \\
d\zeta_k(\ell)&=&\sTr_{\cCl_n}(e^{-(\A_{\ell,k}(2\im(\tau)))^2})-\sTr_{\cCl_n}(e^{-(\A_{1,k}(2\im(\tau))^2})\in \Omega^\bullet(M;C^\infty(\R_{>0}))\nonumber
\eeq
where (as before), $2\im(\tau)\in C^\infty(\R_{>0})$ is identified with the standard coordinate function in the coefficient ring. 
\ep

\subsection{The groupoid $2|1\eft^n(M)$ and the proof of Theorem~\ref{thm1}}\label{sec:proofofthm1}

By Definition~\ref{defn:21eft} of the groupoid $2|1\eft^n(M)$, an object $E\in 2|1\eft^n(M)$ determines the geometric data in Proposition~\ref{prop:mainprop} satisfying the cutoff property in Definition~\ref{defn:admitscutoffsEFT}.

\begin{proof}[Proof of Theorem~\ref{thm1}]
Proposition~\ref{prop:Ktate} gives a cocycle $\Ind (E)$ representing a class $[\Ind (E)] \in \KO^{-n}(M)(\!(q)\!)$ depending on choices of cutoffs. Choose a $\MP_2(\Z)$-equivariant, $\partial_{\bar z}$-preimage~$T$ of~$Z_{\bar z}$ so that 
\beq\label{eq:TheMFstuff}
&&Z-dT\in \Omega^\bullet(M;\mathcal{O}(\HHz)\otimes C^\infty(\R)))^{\MP_2(\Z)},\quad (Z-dT)|_{\vol=1}\in \Omega^\bullet(M;\mathcal{O}(\HHz))^{\MP_2(\Z)}.
\eeq
The restriction to $\vol=1$ is the map on coefficients induced by restricting to the subspace $\HHz\subset \HHz\times \R_{>0}\subset \Lat$ of lattices with volume~1. In light of~\eqref{eq:thisishow}, we may choose $T$ so that the $q$-expansion of $(Z-dT)|_{\vol=1}$ has finitely many negative powers of~$q$ over any subset of~$M$ contained in a compact subspace.

Lemmas~\ref{lem:1sthomotopyforcoccyle} and~\ref{lem:2ndhomotopyforcoccyle} give a coboundary mediating between a cocycle representative of the Pontryagin character of $E$ and $Z \in \Omega^\bullet(M;C^\infty(\HHz\times \R_{>0}))^{\MP_2(\Z)}$. Together with $T$ we obtain a coboundary between the Pontryagin character of the index bundle and $(Z-dT)|_{\vol=1}$. Explicitly, for each $N$ and $\vec\lambda$, we have an equality of elements in $\Omega^\bullet(U^{\vec\lambda};\C(\!(q)\!)/q^{N+1})$
$$
((Z-dT)|_{\vol=1})|_{U^{\vec\lambda}}=\sum^N q^k \sTr(e^{-(\nabla^{<\lambda_k})^2})-\phi(\tau)^{-n}d(\sum^N q^k(\omega_k^\lambda+\zeta_k|_{U^{\vec\lambda}})+T|_{U^{\vec\lambda}})|_{\vol =1}
$$
between the restriction of the globally defined form and the component of the \v{C}ech--de~Rham cocycle. Above, the $\omega_k^\lambda$ come from Lemma~\ref{lem:1sthomotopyforcoccyle} and the $\zeta_k$ are from Lemma~\ref{lem:2ndhomotopyforcoccyle}. The cocycles $\Ind (E)$ and $Z-dT$ together with the above homotopy above determines a class $[E]\in \KO_{\MF}^n(M)$. 

A priori this class depends on the choice of $\partial_{\bar z}$ preimage~$T$ and the choice of cutoffs $\{\vec\lambda\}$ for each $N$ determining the locally finite cover $\{U^{\vec\lambda}\}$. The independence from choices of cutoffs follows from general index bundle considerations, see Proposition~\ref{prop:appenindexbundle}. On the other hand, any pair of choices for $T$ differ by an element of $\Omega^\bullet(M;\MF)$. This changes the representative $Z-dT$ of the class $[Z]\in \H^\bullet(M;\MF)$ by a coboundary in $\Omega^\bullet(M;\MF)$ while simultaneously modifying the homotopy compatibility data~$T$. This determines a homotopic map to $n$th space in the spectrum representing $\KO_\MF$, and hence does not change the class~$[E]\in \KO_{\MF}^n(M)$.
\ep

\subsection{A field theory representative of the $\KO_\MF$-Euler class} \label{sec:EulercocycleKOMF}

We have the following recipe for constructing objects in $2|1\eft^n(M)$. 

\begin{prop} \label{prop:howtoconstruct}
Consider data
\begin{enumerate}
\item[i.] Fr\'echet vector bundles $\HBR_k\to M$ with fiberwise $\cCl_n$-action and $\HBR_k=\{0\}$ for~$k\ll 0$;
\item[ii.] smooth, 1-parameter families of superconnections $\A_{k,\ell}$ on $\HBR_k$ for~$\ell\in \R_{>0}$ that are $\cCl_n$-linear;
\item[iii.] hermitian metrics on each $\HBR_k$ relative to which the $\cCl_n$-action is self-adjoint and the superconnection $\A_{k,\ell}$ is self-adjoint;
\item[iv.] real structures on $\HBR_k$ and $\cCl_n$ compatible with the Clifford action and superconnection; and
\item[v.] differential forms 
$$
(Z,Z_{\bar z},Z_v)\in \Omega^\bullet(M;C^\infty(\HHz\times \R_{>0})[u^{\pm 2}]),\qquad |u^2|=-4
$$
where $Z$ has total degree~$n$, $Z_{\bar z},Z_v$ have total degree $n-1$, they satisfy the $\MP_2(\Z)$-invariance conditions~\eqref{eq:defnofZ}, \eqref{eq:Zv}, and \eqref{eq:Ztau}, and 
$$
\dR Z=0,\quad \partial_{\bar z}Z=\dR Z_{\bar z}, \quad \partial_v Z=\dR Z_v
$$
where $v\in C^\infty(\R_{>0})$ is the standard real coordinate, and $(z,\bar z)\in C^\infty(\HHz)$ are the standard complex coordinates. 
\end{enumerate}
Next suppose that the data (i)-(v) satisfy the properties 
\begin{enumerate}
\item[i.] The endormorphism $\bigoplus q^k e^{-(\A_{\ell,k}(2\im(\tau))^2}$ of $\Omega^\bullet(M;\HBR)=\Omega^\bullet(M;\bigoplus \HBR_k)$ has a trace (i.e., is nuclear) and there an equality of differential forms
\beq\label{Eq:constructTC}
Z= \phi(\ell,\tau)^{-n}\sum q^k  \sTr_{\cCl_n}(e^{-(\A_{\ell,k}(2\im(\tau))^2}),\label{eq:thisishow2}
\eeq
and
\item[ii.] the superconnections $\A_{k,\ell}$ admit index bundles in the sense of Definition~\ref{defn:cutoffs}. 
\end{enumerate}
Then the above determines an object $E\in 2|1\eft^n(M)$. 
\end{prop} 

\bp
The sequence of superconnections determines an infinitesimally generated degree~$n$ representation of $\Ann^{2|1}(M)$ using~\eqref{eq:infgenerator}. The hermitian inner product, self-adjointness, and real structures guarantee that this representation is reflection positive. The differential forms $(Z,Z_{\bar z},Z_v)$ and property~\eqref{Eq:constructTC} imply that the representation is also trace class. Finally, the existence of index bundles imply the cutoff condition. Hence, by Definition~\ref{defn:21eft} we obtain an object~$E\in 2|1\eft^n(M)$. 
\ep

Let $V\to M$ be a real vector bundle with geometric string structure. In~\S\ref{sec:KOMFEuler} we constructed the $\KO_\MF$-Euler class of $V$. We promote this to a construction of an object $\Eu(V)\in 2|1\eft^n(M)$ whose image under the cocycle map~\eqref{eq:cocycle} is the $\KO_\MF$-Euler class of~$V$. Using~\eqref{eq:infgenerator}, we construct an infinitesimally generated degree~$n$ representation of $\Ann^{2|1}(M)$ from 1-parameter family of superconnections determined by
\beq\label{eq:Eulervb2}
&&\HBR=\bS_V\otimes \bigotimes_{k\ge 0} \Lambda_{q^k} V,\qquad \A_k=\nabla_k+2\pi \ell H,\qquad \nabla_k=\nabla^{\bS_V\otimes \bigotimes_{k\ge 0} \Lambda_{q^k} V}|_{\HBR_k}
\eeq
where $H\in \Omega^3(M)\hookrightarrow \Omega^3(M;\End(\HBR_k))$ is part of the data of the geometric string structure. The $\Cl_n$-action on $\HBR$ is inherited from the action on $\bS_V$, see~\eqref{eq:spinorsdefn2}. This determines data (i)-(iv) in Proposition~\ref{prop:howtoconstruct}.

From~\eqref{eq:Eulervb2} we compute
\beq
Z&=&\phi(z)^{-{\rm dim}(V)} \sum_{k\ge 0} q^k \sTr_{\cCl_n}(e^{-\nabla_k^2+\frac{2\pi \ell}{2\im(\tau)}dH})=\Ch(\nabla^{\bS_V\otimes \bigotimes_{k\ge 0} \Lambda_{q^k} V}) e^{\frac{p_1}{\im(z)}}\label{eq:Eulervalontorus}\\
&=&\Pf(F)\Witt(\nabla^V)^{-1} e^{\frac{\pi p_1}{\im(z)}}\nonumber\\
&=&\Pf(F)\Witt^*(\nabla^V)^{-1}\in \Omega^\bullet(M;C^\infty(\HH\times \R_{>0})[u^{\pm 1}])^{\MP_2(\Z)}\nonumber
\eeq
using~\eqref{eq:nonholoWit},~\eqref{eq:ChernofEuler},~\eqref{eq:thisishow2}, where $z=\tau /\ell$. From this we find
\beq\label{eq:Eulerform2}
&&\resizebox{.92\textwidth}{!}{$\partial_v Z=0,\quad \partial_{\bar z} Z=\frac{i p_1/2}{\im(z)^2}\Pf(F)\Witt(\nabla^V)^{-1} e^{\frac{p_1}{\im(z)}}=d\left(\frac{i H/2}{\im(z)^2}\Pf(F)\Witt(\nabla^V)^{-1} e^{\frac{p_1}{\im(z)}}\right)$}
\eeq
which supplies the data (v) in Proposition~\ref{prop:howtoconstruct},
$$
Z_v=0,\qquad Z_{\bar z}=\frac{i H/2}{\im(\tau)^2}\Pf(F)\Witt(\nabla^V)^{-1} e^{\frac{p_1}{\im(\tau)}}
$$
and these data satisfy property (i). 


\begin{defn}
Define the object $\Eu(V)\in 2|1\eft^n(M)$ from the data~\eqref{eq:Eulervb2} and~\eqref{eq:Eulerform2}. 
\end{defn}


\begin{proof}[Proof of Theorem~\ref{thm2}] 
To verify the construction of the object $\Eu(V)$, the discussion above provides the data (i)-(v) in Proposition~\ref{prop:howtoconstruct} as well as property (i); it remains to verify property (ii). But this follows from the fact that the vector bundles $\HBR_k$ are finite rank, and superconnection on finite rank bundles always admit index bundles by~\cite[Theorem~9.7]{BGV}. Using the translation between the modular and the holomorphic Witten class~\eqref{eq:nonholoWitdelta}, we observe that $\Eu(V)\in 2|1\eft^{{\rm dim}(V)}(M)$ refines the $\KO_\MF$-Euler using the description of this class from~\S\ref{sec:KOMFEuler} and the construction of the cocycle map~\eqref{eq:cocycle}. 
\ep

\begin{rmk} \label{rmk:comparetoSTEuler}
To compare the object $\Eu(V)$ above with Stolz and Teichner's proposed $\TMF$-Euler class from~\cite[Theorem 1.0.3]{ST04}, one must extend the $S^1\times \cCl_n$-action on $\HBR$ to a $S^1\rtimes \Cl(L\R^n)$-action via~\eqref{eq:Fermodextend}, i.e., tensor $\HBR$ with the invertible bimodule~\eqref{eq:Morita}. The resulting bundle of $S^1\rtimes \cCl(L\R^n)$-modules on $M$ recovers Stolz and Teichner's stringor bundle restricted to the constant loops. The extension to nearly constant super loops encodes a version of the stringor connection. 
%
%
\end{rmk}

\subsection{A partial construction of the supersymmetric $\sigma$-model}\label{sec:partialsigmamodel}

For a family of Riemannian string manifolds $\pi\colon X\to M$ with $n$-dimensional fibers, we use Proposition~\ref{prop:howtoconstruct} to provide the data of an object $\sigma(X)\in 2|1\eft^{-n}(M)$. When $M=\pt$ we verify that this data satisfies the required properties to define such an object. 

To begin, use~\eqref{eq:infgenerator} to define the infinitesimally generated degree~$n$ representation of $\Ann^{2|1}(M)$ from 1-parameter family of superconnections determined by
\beq\label{eq:sigmamodelsuperconn}
&&\HBR:=\bigoplus_{k\ge 0} q^k \pi_*(\bS\otimes W_k(T(X/M)))\to M,\qquad \A_{\ell,k}:=\B(W_k(T(X/M)))(\ell)+2\pi \ell^2 H
\eeq
where $\B(W_k(T(X/M)))(\ell)$ is the $\ell$-rescaled Bismut superconnection for the Dirac operator $\slashed{D}\otimes \nabla^{W_k(T(X/M))}$, and $\A_{\ell,k}$ modifies this by the datum $H\in \Omega^3(X)$ coming from a geometric string structure. This determines data (i)-(iv) in Proposition~\ref{prop:howtoconstruct}.

We regard the normalized trace of the associated degree~$n$ representation of $\Ann^{2|1}(M)$ as a formal power series, 
\beq\label{Eq:formaltraces}
&&Z=\phi(z)^n\sum_{k\ge 0}q^k \sTr_{\cCl_{-n}} (e^{-\A_{\ell,k}(2\im(\tau))^2})=\phi(z)^n\sum_{k\ge 0}q^k \sTr_{\cCl_{-n}} (e^{-\A_{k}(2\vol)^2})
\eeq
using that $2\im(z)\ell=2\vol$.
\begin{lem}
The formal sum of traces~\eqref{Eq:formaltraces} has the limit
$$
\lim_{\vol \to 0}Z=\int_{X/M} \Witt^*(\nabla^{X/M})
$$
and there are uniquely determined forms $\alpha_k(t)$ satisfying 
\beq\label{eq:Zstuffhomotopy}
Z=\int_{X/M} \Witt^*(\nabla^{X/M})+d\left(\sum_{k \ge 0} q^k \int_0^{2\vol} \alpha_k(t)dt\right).
\eeq
\end{lem}
\bp
This follows from Bismut's local index theorem~\cite[Theorems~10.23 and~10.32]{BGV}, where the $\alpha_k(t)$ are the Chern--Simons forms associated with rescaling superconnections. Specifically, we require the version of the local index theorem for vector bundles with superconnection~\cite[Theorems 5.33 and 5.41]{Kahle}. Using this, the above formula follows from~\eqref{eq:nonholoWitdelta} and~\eqref{eq:BismutstheoremforWG}.
%
\ep

\begin{defn} Define the differential forms in $\Omega^\bullet(M;C^\infty(\R_{>0})(\!(q)\!))$,
\beq
Z_v&=&\sum_{k \ge 0} q^k \alpha_k(2\vol)  \label{eq:failforms}\\ 
Z_{\bar z}&=&\int_{X/M}\frac{-i H/2}{\im(z)^2}\Witt^*(\nabla^{X/M})+\partial_{\bar z} \left(\sum_{k \ge 0} q^k \int_0^{2\vol} \alpha_k(t)dt\right),\label{eq:failforms2}
\eeq
where $H$ comes from the data of a geometric string structure. 
\end{defn}

\begin{lem}
The forms~\eqref{eq:failforms} and~\eqref{eq:failforms2} satisfy
$$
\partial_vZ=dZ_v,\quad \partial_{\bar z}Z=dZ_{\bar z}. 
$$
\end{lem}
\bp
This follows directly from~\eqref{eq:Zstuffhomotopy}. 
\ep

\begin{prop}\label{prop:almostsigma}
Given a family $\pi\colon X\to M$ of Riemannian string manifolds, the data~\eqref{eq:sigmamodelsuperconn},~\eqref{Eq:formaltraces}, \eqref{eq:failforms} and \eqref{eq:failforms2} determine an object $\sigma(X)\in 2|1\eft^{-n}(M)$ if and only if $(Z,Z_v,Z_{\bar z})$ satisfy the $\MP_2(\Z)$-equivariance properties~\eqref{eq:defnofZ}, \eqref{eq:Zv}, and \eqref{eq:Ztau}. 
\end{prop}
\bp
This follows from Proposition~\ref{prop:howtoconstruct}: the required data are supplied as indicated, property (ii) follows from standard results for Bismut superconnections (see Proposition~\ref{prop:Bismutcutoffs}), and so it remains to verify property (i). 
\ep

\begin{proof}[Proof of Theorem~\ref{thm3}]
The superconnections over $M=\pt$ are the twisted Dirac operators comprising the Dirac--Ramond operator (see Definition~\ref{defn:DiracRamond}) and we write $\slashed{D}_k=\A_k$. Specializing the construction above to $M=\pt$, we find that $Z$ is independent of $\vol$ by the McKean--Singer theorem~\cite{McKeanSinger}: the supertrace
$$
\sum q^k \sTr_{\cCl_{-n}}(e^{-2\im(\tau)\ell\slashed{D}_k})=\sum q^k \sTr_{\cCl_{-n}}(e^{-2\vol\slashed{D}_k})
$$
is independent of the parameter $\vol$ and equal to Chern character of the formal sum of indicies. Hence, $Z_v=0$ in this case and the second term in~\eqref{eq:failforms2} also vanishes. The $\MP_2(\Z)$-equivariance properties for $Z$ and $Z_{\bar z}$ are clear, and so by Proposition~\ref{prop:almostsigma} we obtain a class $\sigma(X)\in 2|1\eft^{-n}(\pt)$. Agreement with the $\KO_\MF$-valued analytic index is clear from construction. 
\ep


\appendix 

\section{Pfaffian lines and modular forms}\label{sec:backtot}
Below we review the description of modular forms as sections of tensor powers of a Pfaffian line bundle. This fact is classical, e.g., see~\cite{Quillen_det},~\cite[\S4]{Freed_Det} and~\cite[\S5]{ST11}. However, we require explicit formulas that do not seem to exist in the literature.

\subsection{Conventions for modular forms}\label{Sec:MF} 

Throughout, $\tau\in \HHz\subset \C$ is the standard coordinate on the upper half plane.

\begin{defn}\label{defn:modularform}
For $k\in \Z$, define the vector space of \emph{weight $k$ weak modular forms} as
$$
\MF_k^\omega:=\left\{f\in\mathcal{O}(\HHz)\mid f\left(\frac{a\tau+b}{c\tau+d}\right)=(c\tau+d)^kf(\tau), \ \left[\begin{smallmatrix} a & b \\ c & d\end{smallmatrix} \right]\in \SL_2(\Z)\right\}.
$$ 
The \emph{ring of weak modular forms} is $\MF^\omega:=\bigoplus_{k\in \Z} \MF_k^\omega$ with multiplication inherited from multiplication of functions on $\HHz$. 
A weak modular form is a \emph{weakly holomorphic} if $f(\tau)$ is meromorphic as $\tau\to i\infty$. Let $\MF=\bigoplus_k \MF_k$ denote the ring of weakly holomorphic modular forms. 
\end{defn}

When regarding modular forms and weak modular forms as graded rings, the graded components are defined as 
\beq\label{eq:MFgraded}
\MF^k:=\left\{\begin{array}{ll} \MF_{-k/2} &  k\ {\rm even} \\ 0 & l \ {\rm odd} \end{array}\right.\qquad (\MF^\omega)^k:=\left\{\begin{array}{ll} \MF_{-k/2}^\omega &  k\ {\rm even} \\ 0 & l \ {\rm odd} \end{array}\right.
\eeq
which results in a graded-commutative product structure. 

\begin{rmk}
An equivalent description of the graded ring $\MF^\omega$ is the $\SL_2(\Z)$-invariants
\beq\label{eq:MFassubspace}
&&\MF^\omega\simeq (\mathcal{O}(\HHz)[u^{\pm 1}])^{\SL_2(\Z)},\qquad |u|=-2
\eeq
for the action by fractional linear transformations on $\HHz$ and $u\mapsto u/(c\tau+d)$. 
\end{rmk}

A modular form can equivalently be expressed as a function of $q=e^{2\pi i \tau}$ on the unit punctured disk in $\C$. Taylor expanding defines an injective ring homomorphism
\beq
q\hbox{-}{\rm expand}\colon \MF\hookrightarrow \C(\!(q)\!).\label{eq:qexpand}
\eeq
Let $\MF^\Z_k\subset \MF_k$ denote the subgroup of \emph{integral modular forms} of weight~$k$, defined as having $q$-expansions in integral Laurent series, $\Z(\!(q)\!)\subset \C(\!(q)\!)$.

\subsection{The Pfaffian line over Euclidean spin tori}\label{sec:back1}

We recall the description of a Euclidean spin torus from Remark~\ref{rmk:chiralD} in terms of a complex manifold $\C/\Lambda$ with a chosen square root of its canonical bundle. In this description, the chiral Dirac operator is the $\bar\partial$-operator twisted by the canonical bundle. 

The chiral Dirac operator for the even spin structure is invertible (see~\S\ref{sec:dumbEuc} for spin structures on the torus), and hence the Pfaffian line bundle is canonically trivialized by its $\SL_2(\Z)\times \Spin(2)$-invariant Pfaffian section.\footnote{This trivialization does not respect the Quillen metric: the norm of the Pfaffian section over tori with even spin structure is essentially given by the three Weber modular functions with appropriate normalizations; e.g., this can be seen by specializing \cite[Proposition 4.10]{Freed_Det} to real points in the Jacobian.} 
The chiral Dirac operator for the odd spin structure has a 1-dimensional kernel, and so the Pfaffian line bundle is the $\SL_2(\Z)\times \Spin(2)$-equivariant line bundle over $\Lat$ whose fiber at $(\ell_1,\ell_2)\in \Lat$ is the dual space of harmonic spinors. Since $\Spin(2)\simeq U(1)$ acts on harmonic spinors through the weight~1 representation, the Pfaffian line $\Pf\to \Lat\sq \SL_2(\Z)\times \Spin(2)$ is classified by the map of stacks
\beq\label{eq:Pfhomo}
&&\Lat \sq \SL_2(\Z)\times \Spin(2)\to \pt\sq U(1),\ \  \SL_2(\Z)\times \Spin(2)\to \Spin(2)\simeq U(1)\xrightarrow{(-)^{-1}} U(1)
\eeq
determined by the projection homomorphism, the identification $\Spin(2)\simeq U(1)$, and inversion on $U(1)$. The inversion comes from the fact that the Pfaffian line is \emph{dual} to harmonic spinors. Tensor powers of the Pfaffian are then determined by the homomorphism
\beq\label{eq:canonicalline}
&&\SL_2(\Z)\times \Spin(2)\to \Spin(2)\simeq U(1)\xrightarrow{k} U(1)\iff \Pf^{\otimes -k}\to \Lat\sq \SL_2(\Z)\times \Spin(2)
\eeq
where the first arrow is the projection and the second arrow is the $k$-fold covering map,~$z\mapsto z^{k}$ for $k\in \Z$. Global sections of $\Pf^{\otimes -k}$ are therefore functions on $\Lat$ that transform with weight $k$ under the $U(1)\simeq \Spin(2)$-action and are invariant under the $\SL_2(\Z)$-action. 

There is a line bundle $\overline{\Pf}^{\otimes -k}\to \Lat\sq \SL_2(\Z)\times \Spin(2)$ defined analogously to~\eqref{eq:canonicalline}, but using the conjugate character $z\mapsto \overline{z}^k$ of $U(1)$. The (positive) square root of the volume of a lattice determines a function 
$$
\vol^{-1/2}\in C^\infty(\Lat),\qquad \vol(\ell_1,\ell_2)=\frac{\ell_1\bar\ell_2-\bar\ell_1\ell_2}{2i}
$$ 
that descends to a nonvanishing global section of $\Pf\otimes \overline{\Pf}$, i.e., a trivialization. Hence there is an isomorphism of line bundles
$$
\Pf\simeq \overline{\Pf}^\vee= \overline{\Pf}^{\otimes -1}
$$
over $\Lat \sq \SL_2(\Z)\times \Spin(2)$ determined by $\vol^{-1/2}$.

 For future reference, we observe the functors
$$
\Lat\sq \Spin(2)\times \SL_2(\Z)\xrightarrow{\rm forget} \Lat\sq \SO(2)\times \SL_2(\Z)\xrightarrow{\vol} \R_{>0} 
$$
where the first arrow forgets the odd spin structure (or equivalently, is determined by the double cover $\Spin(2)\to \SO(2)$), and the second takes the volume of a Euclidean torus. The fiber at $v\in \R_{>0}$ determines a $\Spin(2)\times \SL_2(\Z)$-invariant subspace of $\Lat$ of lattices whose associated Euclidean spin tori have volume~$v$.

\subsection{Sections of the Pfaffian line and modular forms}\label{sec:back2}
Next we explain how~$\Pf^{\otimes k}$ and its sections are related to modular forms. 
Here it is useful to work with a different presentation of the stack of Euclidean tori with periodic-periodic spin structure. Let $\MP_2(\Z)$ denote the metaplectic group, the nontrivial double cover of~$\SL_2(\Z)$. 
This group has the explicit description 
$$
\MP_2(\Z)\simeq \left\{\left[\begin{array}{cc} a & b \\ c & d \end{array}\right]\in \SL_2(\Z),\epsilon\in \mathcal{O}(\HHz)^\times \mid \epsilon(z)^2=cz+d\right\}
$$
where multiplication is determined by multiplication in $\SL_2(\Z)$, multiplication in $\mathcal{O}(\HHz)^\times$, and the action of $\SL_2(\Z)$ on $\mathcal{O}(\HHz)^\times$ by fractional linear transformations. Below we will use the notation $\sqrt{cz+d}$ for the function $\epsilon(z)$ when describing an element of~$\MP_2(\Z)$. Define the quotient stack $\HHz\times\R_{>0}\sq \MP_2(\Z)$ for the action
$$
\left(\left[\begin{array}{cc} a & b \\ c & d \end{array}\right],\sqrt{cz+d}\right)\cdot (\tau,\ell)=\left(\frac{a\tau+b}{c\tau+d},\ell|c\tau+d|\right).
$$
where we observe that this action factors through the homomorphism $\MP_2(\Z)\to \SL_2(\Z)$. The quotient stack $\HHz\times\R_{>0}\sq \MP_2(\Z)$ supports line bundles indexed by $k\in \Z$ classified by functors $\HHz\sq \MP_2(\Z)\to \pt\sq \C^\times$ determined by
\beq\label{eq:sqmodularform}
\HHz\times \MP_2(\Z)\to \C^\times ,\qquad \left(\tau,\left[\begin{array}{cc} a & b \\ c & d \end{array}\right],\sqrt{cz+d}\right)\to (c\tau+d)^{k/2}. 
\eeq
Similarly, the maps
\beq\label{eq:sqmodularformcong}
\HHz\times \MP_2(\Z)\to \C^\times ,\qquad \left(\tau,\left[\begin{array}{cc} a & b \\ c & d \end{array}\right],\sqrt{cz+d}\right)\to (c\overline{\tau}+d)^{l/2}. 
\eeq
determine conjugate line bundles indexed by $l\in \Z$. 

\begin{defn} 
For $(k,l)\in \Z\times \Z$, let $\mmf^{k,l}$ denote the global sections of the tensor product of the line bundles~\eqref{eq:sqmodularform} and~\eqref{eq:sqmodularformcong}. Hence, we have
\beq\label{Eq:mmf}
&&\resizebox{.93\textwidth}{!}{$
\mmf^{k,l}=\left\{f\in C^\infty(\HHz\times \R_{>0})\mid f\left(\frac{a\tau+b}{c\tau+d},\ell|c\tau+d|\right)=(c\tau+d)^{-k/2}(c\overline{\tau}+d)^{-l/2}f(\tau,\ell)\right\}$}
\eeq
\end{defn}

\begin{lem} \label{lem:euctoribackground}
There is an equivalence of stacks
\beq\label{eq:appenequivtori}
\HHz\times\R_{>0}\sq \MP_2(\Z)\stackrel{\sim}{\to} \Lat\sq \SL_2(\Z)\times \Spin(2)
\eeq
and the pullback of the lines $\Pf^{\otimes -k}$ from~\eqref{eq:canonicalline} are isomorphic to the lines~\eqref{eq:sqmodularform}. Similarly, the pullback of~$\overline{\Pf}^{\otimes -k}$ are isomorphic to the lines~\eqref{eq:sqmodularformcong}. In particular, we obtain the isomorphism on global sections,
$$
\mmf^{k,l}\simeq \Gamma(\Lat\sq (\Spin(2)\times \SL_2(\Z));\Pf^{\otimes k}\otimes \overline{\Pf}^{\otimes l}).
$$
\end{lem}
\bp
Consider the functor between Lie groupoids defined by the maps
\beq
\HHz\times \R_{>0}\to \Lat,&& (\tau,\ell)\mapsto (\ell\tau,\ell)\nonumber
\eeq
and
\beq
\HHz\times \R_{>0}\times \MP_2(\Z)&\to& \SL_2(\Z)\times \Spin(2),\nonumber\\ 
\left(\tau,\ell,\left[\begin{array}{cc} a & b \\ c & d \end{array}\right],\sqrt{c\tau+d}\right)&\mapsto& \left(\left[\begin{array}{cc} a & b \\ c & d \end{array}\right],\frac{(c\tau+d)^{1/2}}{|c\tau+d|^{1/2}}\right).\nonumber
\eeq
A short calculation shows that these maps fit together to determine a fully faithful and essentially surjective functor~\eqref{eq:appenequivtori}, yielding the claimed equivalence of stacks. 

To verify the isomorphism of line bundles, we view sections as smooth functions on $\HHz\times \R_{>0}$ with transformation properties. Then multiplication by $\ell^{k/2}$ provides the desired isomorphism between sections; indeed, for $f(\tau,\ell)\in C^\infty(\HHz\times \R_{>0})$, the transformation property for $\ell^{k/2}f(\tau,\ell)$ to be a section of the pullback of $\Pf^{\otimes k}$ along~\eqref{eq:appenequivtori} is the equality (using~\eqref{eq:canonicalline})
\beq\label{eq:powerofell}
&&\begin{array}{lll}
\left(\left[\begin{array}{cc} a & b \\ c & d \end{array}\right],\sqrt{c\tau+d}\right)\cdot (\ell^{k/2}f(\tau,\ell))&=&(\ell|c\tau+d|)^{k/2}f\left(\frac{a\tau+b}{c\tau+d},\ell|c\tau+d|\right)\\
&=&\frac{|c\tau+d|^{k/2}}{(c\tau+d)^{k/2}}\ell^{k/2}f(\tau,\ell), 
\end{array}
\eeq
which is equivalent to $f(\tau,\ell)$ transforming with weight $-k/2$, via~\eqref{eq:sqmodularform}. The argument for $\overline{\Pf}{}^{\otimes k}$ is analogous. 
\ep

The moduli stack of elliptic curves is presented by the complex analytic quotient stack $\HHz\sq \SL_2(\Z)$ for the $\SL_2(\Z)$-action by fractional linear transformations. Define a line bundle $\Mf_k \to \HHz\sq \SL_2(\Z)$ by the map
\beq\label{eq:sqmodularform1}
\SL_2(\Z)\to \mathcal{O}(\HHz)^\times ,\qquad \left[\begin{array}{cc} a & b \\ c & d \end{array}\right]\to (c\tau+d)^{k}. 
\eeq
Holomorphic sections of $\Mf_k$ are weight $k$ weak modular forms, whence the notation. 
The moduli stack of elliptic curves with periodic-periodic spin structure is presented by the complex analytic quotient stack $\HHz\sq \MP_2(\Z)$, where the action on $\HHz$ factors through the homomorphism $\MP_2(\Z)\to \SL_2(\Z)$. There is an evident functor $\HHz\sq\MP_2(\Z)\to \HHz\sq \SL_2(\Z)$; geometrically, this forgets the spin structure on the elliptic curve. The pullback of the line bundle~\eqref{eq:sqmodularform1} along this forgetful functor admits a square root denoted $\Mf_{k/2}$ determined by the character~\eqref{eq:sqmodularform}. Taking underlying smooth stacks, there is a functor from Euclidean spin tori to elliptic curves given by the composition, 
\beq\label{eq:forgettoell}
&&f\colon \Lat\sq \Spin(2)\times \SL_2(\Z) \simeq \HHz\times \R_{>0}\sq \MP_2(\Z)\to \HHz\sq \MP_2(\Z)\to \HHz\sq \SL_2(\Z)
\eeq
where the first arrow is induced by projection. Then we obtain the following corollary to Lemma~\ref{lem:euctoribackground}. 

\begin{cor}\label{cor:modularcompare}
The pullback of the line bundle determined by~\eqref{eq:sqmodularform1} along the composition~\eqref{eq:forgettoell} is isomorphic to the line determined by~\eqref{eq:canonicalline}:
$$
f^*\Mf_k\simeq \Pf^{\otimes 2k}.
$$ 
In particular, weakly holomorphic modular forms of weight $k$ pull back to sections of~$\Pf^{\otimes 2k}$. More generally, smooth sections of $\Mf_{k/2}\otimes \overline{\Mf}_{l/2}$ over $\HHz\sq \MP_2(\Z)$ pull back to smooth sections of $\Pf^{\otimes k}\otimes \overline{\Pf}^{\otimes l}$. 
\end{cor}

\subsection{Geometry of the Dedekind $\eta$-function and 2nd Eisenstein series}\label{sec:etageometry}

Recall the Dedekind $\eta$-function,
\beq\label{eq:etadefn}
\eta(q):=q^{1/24}\prod_{n>0}(1-q^n),\qquad q=e^{2\pi i z}.
\eeq
The following lemma helps explain certain normalizations in the physics literature. Concretely, if a function on $\HHz$ transforms by a $U(1)$-phase under the $\MP_2(\Z)$-action, multiplying by a power of~$\eta$ gives a function that transforms as a modular form. Rephrasing in the geometry of stacks yields the following. 

\begin{lem}\label{lem:Dedekind}
The Dedekind $\eta$-function determines an isomorphism between the line bundle over $\HHz\sq \MP_2(\Z)$ determined by~\eqref{eq:sqmodularform} and the line bundle determined by the character
\beq\label{eq:mu24}
\mu \colon  \MP_2(\Z)\to \MP_2(\Z)/[\MP_2(\Z),\MP_2(\Z)]\simeq \Z/24\subset U(1)
\eeq
given by embedding the abelianization of $\MP_2(\Z)$ as 24th roots of unity. 
\end{lem}
\bp
The transformation formula for the $\eta$-function under the action by $\MP_2(\Z)$ is
$$
\eta(\gamma z)=(\mu(\gamma)^{-1}\sqrt{(cz+d)})\eta(z),\qquad \gamma\in \MP_2(\Z)
$$
and hence $\eta$ is a global section of the tensor product of the line bundle determined by~\eqref{eq:sqmodularform} and the dual of the line bundle determined by~\eqref{eq:mu24}. In view of~\eqref{eq:etadefn}, $\eta(z)$ is a nonvanishing function and therefore determines a nonvanishing section of this tensor product of line bundles over $\HHz\sq \MP_2(\Z)$. This proves the lemma.
\ep
%

There are two common forms of the 2nd Eisenstein series
\beq
E_2(z)&=&\sum_{n\ne 0} \frac{1}{n^2}+\sum_{m\ne 0}\sum_{n\in \Z} \frac{1}{(mz+n)^2},\nonumber\\
 E_2^*(z,\bar z)&=&E_2(z)-\frac{\pi}{\im(z)}=E_2(z)-\frac{2\pi i}{z-\bar z}.\label{Eq:2ndEisen}
\eeq
We recall~\cite[page 19]{ZagierMF} that $E_2(z)$ is holomorphic but not modular, whereas $E_2^*(z,\bar z)$ is not holomorphic but is weight~2,
$$
E_2^*\left(\frac{az+b}{cz+d}\right)=(cz+d)^2E_2^*(z). 
$$
This can be rephrased as follows. 

\begin{lem}\label{eq:2ndEisen}
The nonholomorphic 2nd Eisenstein series $E_2^*$ defines a smooth section of the $4$th tensor power of the line bundle determined by~\eqref{eq:sqmodularform}. This smooth section satisfies
\beq\label{eq:E2isntholo}
\partial_{\bar z}E_2^*=-\frac{2\pi i}{(z-\bar z)^2}.
\eeq
\end{lem}

\section{Super background}

\subsection{Real structures and super hermitian pairings}\label{sec:superherm}
The category of \emph{super vector spaces} is the category of $\Z/2$-graded vector spaces over $\C$ equipped with the graded tensor product. The adjective ``super" refers to the sign in the braiding isomorphism for this tensor product
$$
\sigma\colon V\otimes W\stackrel{\sim}{\to} W\otimes V \qquad v\otimes w\mapsto (-1)^{|v||w|} w\otimes v \quad v\in V, \ w\in W
$$
where $v$ and $w$ are homogeneous elements of degree $|v|,|w|\in \Z/2$. 
Let $V^\ev\oplus V^\odd$ denote the direct sum decomposition of $V$ into its even and odd subspaces. The \emph{grading involution} is the linear map $(-1)^{\sf F}\colon V \to V$ that acts by $+1$ on $V^\ev$ and~$-1$ on $V^\odd$. 


There is an involution on the category of super vector spaces that on objects reverses complex structures, $V\mapsto \overline{V}$, i.e., $\overline{V}$ is gotten from $V$ by precomposing the action of $\C$ with complex conjugation. 

\begin{defn}
A \emph{real structure} on a super vector space $V$ is an isomorphism $r\colon V\stackrel{\sim}{\to} \overline{V}$ of super vector spaces such that $\overline{r}\circ r=\id_V$. The \emph{real subspace} associated with a real structure is the fixed point set $V_\R\subset V$ for $r$. 
\end{defn}


\begin{defn}\label{defn:sesquilinear} A \emph{hermitian form} on a super vector space $V$ is the data of a map $\langle-,-\rangle\colon \overline{V}\otimes V\to \C$ making the diagram commute, 
\beq
&&\begin{tikzpicture}[baseline=(basepoint)];
\node (A) at (0,0) {$\overline{V}\otimes V$};
\node (B) at (4,0) {$\C$};
\node (C) at (0,-1.25) {$V\otimes \overline{V}$};
\node (D) at (4,-1.25) {$\overline{\C}$};
\draw[->] (A) to node [above] {$\langle-,-\rangle$} (B);
\draw[->] (A) to node [left] {$\sigma$} (C);
\draw[->] (C) to node [below] {$\overline{\langle-,-\rangle}$} (D);
\draw[->] (B) to node [right] {$(\overline{\phantom{A}})$} (D);
\path (0,-.75) coordinate (basepoint);
\end{tikzpicture}
\eeq
or equivalently, satisfying the formula $\langle x,y\rangle =(-1)^{|x||y|}\overline{\langle y,x\rangle}$ for homogeneous elements. A hermitian form is \emph{positive} if 
\beq
\langle x,x\rangle>0, & x\ne 0 \ {\rm even} \quad {\rm and} \quad i^{-1}\langle x,x\rangle>0, & x\ne 0 \ {\rm odd.}\label{eq:positivitypairing}
\eeq
A \emph{real} hermitian form on a super vector space with real structure is a hermitian form that when restricted to $(V_\R)^\ev\subset V$ takes values in $\R\subset \C$ and when restricted to $(V_\R)^\odd\subset V$ takes values in $i^{-1}\R\subset \C$. 
\end{defn}

A positive hermitian form in the above sense determines a positive pairing $(-,-)$ in the usual sense on the underlying (ungraded) vector space $V$ by the formula
\beq
(x,y)=\left\{\begin{array}{ll} 0 & |x|\ne |y| \\ \langle x,y\rangle & x,y \ {\rm even} \\  i\langle x,y\rangle & x,y \ {\rm odd.} \end{array}\right.\label{eq:ordinarypairing}
\eeq
When $V$ is endowed with a real structure, this translation sends a real hermitian form to a real inner product on $V_\R$.

\begin{defn}
Let $V$ be a super vector space equipped with a positive hermitian form~$\langle-,-\rangle$. The \emph{super adjoint} of an operator $T\colon V\to V$ is characterized by
\beq
\langle x,T y\rangle =(-1)^{|x||A|}\langle T^*x,y\rangle. \label{eq:superadjoint}
\eeq
\end{defn}

In the translation to the ordinary pairing~\eqref{eq:ordinarypairing}, we have
\beq
&&T^*=\left\{\begin{array}{ll} T^\dagger & T \ {\rm even} \\ iT^\dagger & T\ {\rm odd}\end{array}\right.\label{eq:superadjointtrans}
\eeq
where $(-)^\dagger$ denotes the usual adjoint with respect to the positive pairing $(-,-)$.

\subsection{Super algebras and their modules} \label{sec:supertrace}
A \emph{superalgebra} is an algebra object in the category of super (i.e., $\Z/2$-graded) vector spaces equipped with the $\Z/2$-graded tensor product. Super algebras are the objects of a monoidal category, where the tensor product of superalgebras~$A$ and $B$ is the super vector space $A\otimes B$ with multiplication
$$
(a\otimes b)\cdot (a'\otimes b')=(-1)^{|b||a'|} aa'\otimes bb'. 
$$
The \emph{opposite} of a superalgebra $A$ is a superalgebra $A^\op$ with the same underlying super vector space and multiplication 
\beq
a\cdot_{\op} b:=(-1)^{|a||b|}b\cdot a\label{eq:opposite}
\eeq
where $b\cdot a$ is the multiplication in $A$. 

\begin{defn}\label{defn:starsuper}
A \emph{$*$-superalgebra} is a superalgebra together with the data of a homomorphism $(-)^*\colon A\to \overline{A}{}^\op$ that squares to the identity on $A$, where the opposite algebra is taken in the graded sense.
\end{defn}

Given a superalgebra $A$, let ${}_A\Mod$ and $\Mod_A$ denote the categories of left and right $A$-modules, respectively, whose objects are super vector spaces with a (graded) $A$-action. These categories have symmetric monoidal structures from direct sum of $A$-modules. Similarly, let ${}_A\Mod_B$ denote the category of $A-B$-bimodules. There are canonical equivalences of categories 
\beq
\Mod_{A^\op\otimes B^\op}\simeq {}_{A}\Mod_{B^{\op}} \simeq {}_{A\otimes B}\Mod\label{eq:leftright}
\eeq 
where a left $A\otimes B$-bimodule $W$ is sent to the $A-B^\op$-bimodule with the same underlying super vector space and action
$$
a\cdot w\cdot b:=(-1)^{|w||b|} (a\otimes b) \cdot w
$$
with similar formulas defining a right $A^\op\otimes B^\op$-module. 
We often use notation like ${}_AW$ for an object of ${}_A\Mod$, and ${}_AW_B$ for an object of ${}_A\Mod_B$. 

The \emph{supertrace} of an endomorphism $T\colon V\to V$ is
\beq
\sTr(T):= {\sf Tr}((-1)^{\sf F}\circ T)\label{eq:supetrace}
\eeq
where ${\sf Tr}$ is the ordinary trace. 
The \emph{super commutator} of elements in a superalgebra is
\beq
[a,b]=ab-(-1)^{|a||b|}ba\qquad a,b\in A.\label{eq:supercommutator}
\eeq
A \emph{supertrace} on a superalgebra is a map of super vector spaces 
\beq
\sTr_A\colon A\to \C,\qquad \sTr_A([a,b])=0\label{eq:superalgsupertrace}
\eeq
that vanishes on super commutators, and hence is determined by a linear map~$A/[A,A]\to \C$. 

A supertrace on a superalgebra $A$ induces a supertrace on the category of (finite-dimensional) $A$-modules. For right $A$-modules, this supertrace comes from the composition
\beq
\End_A(W)\stackrel{\sim}{\leftarrow} W\otimes_A W^\vee \to A/[A,A]\stackrel{\sTr_A}{\to}  \C \label{eq:HattoriStallings}
\eeq
where $W^\vee:={\rm Hom}_A(W,A)$ is the dual (left) $A$-module, and the map to $A/[A,A]$ is the Hattori--Stallings trace induced by the evaluation pairing,
$$
{\rm eval}\colon W\otimes_\C W^\vee \to A. 
$$
One obtains traces for bimodules (and left modules) from traces on $A\otimes B^\op$ using~\eqref{eq:leftright}. 



\subsection{(Real) supermanifolds}
Stolz and Teichner's supersymmetric field theories are couched in the language of supermanifolds whose structure sheaves are defined over~$\C$; this is called the category of $cs$-supermanifolds in~\cite{DM}. The majority of what we require is covered in the concise introduction~\cite[\S4.1]{ST11}, but we establish a little notation below. 

By definition, any supermanifold is locally isomorphic to the supermanifold~$\R^{n|m}$ characterized by the global sections of its structure sheaf,~$C^\infty(\R^{n|m})\simeq C^\infty(\R^n;\C)\otimes_\C \Lambda^\bullet \C^m$. We emphasize that this is a sheaf of super (i.e., $\Z/2$-graded) algebras over $\C$. 

Let $\Pi$ denote the functor
\beq\label{Eq:Pifunctor}
\Pi\colon \Vect_{/\Mfld} \to \SMfld,\qquad (E\to N)\mapsto \Pi E
\eeq
that sends a complex vector bundle $E\to N$ over an ordinary manifold $N$ to the supermanifold $\Pi E$ whose structure sheaf is the sheaf of sections of $\Lambda^\bullet E^\vee$, viewed as a sheaf of $\Z/2$-graded algebras. For a supermanifold $S$, let $S_{\rm red}$ denote its reduced manifold; this is an ordinary manifold obtained (as locally ringed space) by the quotient of the structure sheaf of~$S$ by its nilpotent ideal.  Every supermanifold can be written as $S\simeq \Pi E$, where $E\to S_{\rm red}$ is a complex vector bundle~\cite{Batchelor}. In particular, the supermanifolds $\R^{n|m}$ arise as the trivial rank $m$-bundle over~$\R^n$, i.e., $E=\C^n\times \R^n\to \R^n=(\R^{n|m})_{\rm red}$. 

One frequently performs computations in supermanifolds using the \emph{functor of points}. This formalism views a supermanifold as a representable presheaf.

\begin{ex} 
The representable presheaf associated with~$\R^{n|m}$ assigns to a supermanifold~$S$ the set 
\beq\label{eq:RnmSpts}
&&\R^{n|m}(S):=\{t_1,t_2,\dots,t_n\in C^\infty(S)^{\ev}, \theta_1,\theta_2,\dots,\theta_m\in C^\infty(S)^\odd \mid (t_i)_{\rm red}=\overline{(t_i)}_{\rm red}\}
\eeq
where $(t_i)_{\rm red}$ denotes the restriction of a function to the reduced manifold $S_{\rm red}\hookrightarrow S$, and $\overline{(t_i)}_{\rm red}$ is the conjugate of the complex-valued function $(t_i)_{\rm red}$ on the smooth manifold~$S_{\rm red}$. We use this functor of points description throughout the paper, typically with roman letters denoting even functions and greek letters denoting odd functions. 
\end{ex}
\begin{ex}\label{ex:R21points}
We record two equivalent descriptions of the functor of points of~$\R^{2|1}$
\beq
\R^{2|1}(S)&\simeq& \{x,y \in C^\infty(S)^\ev,\ \theta \in C^\infty(S)^\odd\mid (x)_{\rm red}=\overline{(x)}_{\rm red}, (y)_{\rm red}=\overline{(y)}_{\rm red}\}\label{eq:r211}\\
&\simeq& \{z,w \in C^\infty(S)^{\ev}, \theta\in C^\infty(S)^\odd\mid (z)_{\rm red}=\overline{(w)}_{\rm red}\}.\label{eq:r212}
\eeq
The isomorphism between~\eqref{eq:r211} and~\eqref{eq:r212} is $(x,y)\mapsto (x+iy,x-iy)=(z,w)$. Below we adopt the standard (though potentially misleading) notation $\overline{z}:=w$. 
\end{ex}

\begin{defn}
A \emph{super Lie group} is a group object in supermanifolds. 
\end{defn}

By definition, a super Lie group structure on a supermanifold~$G$ comes from endowing the sets $G(S)$ with the structure of a group natural in $S$. Most of the super Lie groups in this paper will be defined in this manner. We also note that every ordinary Lie group determines a super Lie group. 

\begin{ex}\label{ex:U1}
Consider the super Lie group determined by the ordinary Lie group $U(1)\subset \C$, the unit complex numbers. We use the functor of points description 
\beq
&&U(1)(S)=\{u,\bar u\in C^\infty(S)^{\ev}\mid (u)_{\rm red}=\overline{(\bar u)}_{\rm red}, \  u \bar u=1 \},\label{Eq:spin2}
\eeq
analogous to the notation for the $S$-points of $\R^{2|1}$ introduced after~\eqref{eq:r212}.
\end{ex}

\begin{defn}
An \emph{open cover} of a supermanifold $S$ is a collection $\{\varphi\colon U_i\to S\}$ where $\varphi_i$ is a local isomorphism for each $i$ and $\{\varphi_{\rm red}\colon (U_i)_{\rm red}\to S_{\rm red}\}$ is an ordinary open cover of the manifold~$S_{\rm red}$. 
\end{defn}

\begin{defn}\label{defn:realsmfld}
There is an involution on the category of supermanifolds, 
\beq\label{eq:conjugationfunctor}\label{eq:conjsmfld}
(\overline{\phantom{-}})\colon \SMfld\to \SMfld,\qquad S\mapsto \overline{S}
\eeq
sending a supermanifold to the same locally ringed space but whose sheaf has complex numbers acting through complex conjugation. A \emph{real structure} on a supermanifold $S$ is an isomorphism $r_S\colon S\stackrel{\sim}{\to} \bar S$ such that $\overline{r}\circ r=\id_S$. A map $f\colon S'\to S$ between supermanifolds with real structure is \emph{real} if $r_S\circ f=f\circ r_{S'}$. 
\end{defn}

Since~\eqref{eq:conjugationfunctor} is a functor, a super Lie group $G$ has a conjugate super Lie group $\overline{G}$ (which typically is not isomorphic to $G$). A \emph{real structure} on a super Lie group is a real structure on its underlying supermanifold such that all the structure maps are real. 


\begin{ex}
A real structure on a supermanifold presented as $S= \Pi E$ amounts to a choice of real structure on the complex vector bundle~$E\to S_{\rm red}$. In particular, real structures need not exist, and when they do exist they need not be unique. However, ordinary manifolds regarded as supermanifolds have a canonical real structure where the map $M\to \overline{M}$ comes from complex conjugation of $\C$-valued functions on $M$. Similarly, ordinary Lie groups regarded as super Lie groups have canonical real structures. 
\end{ex}

\begin{ex} \label{ex:real}
For a complex vector bundle $E\to N$, a real structure on $E$ (as a vector bundle) determines a real structure on $\Pi E$ (as a supermanifold) via complex conjugation
$$
C^\infty(\Pi E)\simeq \Gamma(N;\Lambda^\bullet E^\vee)\simeq \Gamma(N;(\Lambda E_\R^\vee)\otimes \C) \qquad (E_\R)\otimes \C\simeq E. 
$$
Indeed, complex conjugation is a $\C$-antilinear map $C^\infty(\Pi E)\to C^\infty(\Pi E)$, determining an isomorphism of supermanifolds $\Pi E\stackrel{\sim}{\to} \overline{\Pi E}$. 
\end{ex}

Finally, we require a technical lemma. Given a vector bundle $V\to M$, we obtain a vector bundle $p^*V\to \Map(\R^{0|1},M)$ by pulling back along the projection $p\colon \Map(\R^{0|1},M)\simeq \Pi TM\to M$. 


\begin{lem} \label{lem:sVect}For any vector bundle $\V\to \Map(\R^{0|1},M)$, there is an isomorphism $\V \simeq p^*V$ where $V\to M$ is a vector bundle on $M$.
\end{lem}

\bp
Let $i\colon M\hookrightarrow \Map(\R^{0|1},M)$ denote the inclusion of the reduced manifold. Given a vector bundle $\V\to \Map(\R^{0|1},M)$, our candidate for $V\to M$ in the statement of the lemma is $V:=i^*\V$. Indeed, the composition $i\circ p$ is homotopic to the identity on $\Map(\R^{0|1},M)$. For example, consider the homotopy $\R\times \Pi TM_\C\to \Pi TM_\C$ that rescales the odd fibers of $\Pi TM_\C\to M$ by a scalar $t\in \R$; in the functor of points this map is $(t,x,\psi)\mapsto (x,t\psi)$ for $t\in \R(S)$ and $(x,\psi)\in \Pi TM_\C(S)\simeq \Map(\R^{0|1},M)(S)$. Now we use that pullbacks along homotopic maps yield isomorphic vector bundles,
concluding that $p^*V=p^*i^*\V\simeq (i\circ p)^*\V$ is isomorphic to $\V$. 
\ep

\section{Stolz and Teichner's twisted geometric field theories}\label{sec:STappen}

We refer to the companion paper in this series~\cite[Appendix~B]{DBEEFT} for a more detailed overview of the geometric bordism categories constructed in~\cite{ST11}. Below we collect the main definitions and notation used in this paper, as well as the definition of a reflection structure for a twisted field theory (see~\S\ref{sec:RPgeneral}). The presentation below emphasizes that twisted field theories fit into the larger context of twisted representations of internal categories; this larger context includes previously studied examples in Lie groupoid geometry and twisted equivariant K-theory. 

\subsection{Categories internal to stacks} \label{sec:internalcategories}


\begin{defn}[{\cite[Definition~2.13]{ST11}}]\label{defn:internalcat}
A \emph{category ${\sf C}$ internal to stacks} is the data of stacks $\Ob({\sf C})$ and $\Mor({\sf C})$ with source, target, unit, and composition functors
\beq\label{eq:structuremapsinternal}
&&\s,\t\colon \Mor({\sf C})\to \Ob({\sf C})\quad \u\colon \Ob({\sf C})\to \Mor({\sf C}),\quad \c\colon \Mor({\sf C})\times_{\Ob({\sf C})} \Mor({\sf C})\to \Mor({\sf C})
\eeq
along with additional coherence data (e.g., an associator and unitors) satisfying compatibility properties. \emph{A category internal to symmetric monoidal stacks} and \emph{a category internal to symmetric monoidal stacks with involution} are defined analogously. 
\end{defn} 

Stacks form a $(2,1)$-category, and correspondingly we will refer to maps between stacks as 1-morphisms and isomorphisms between these 1-morphisms as 2-morphisms. The definitions below apply to categories internal to stacks, categories internal to symmetric monoidal stacks, or categories internal to symmetric monoidal stacks with involution. We will simply refer to an internal category when all three of these options are in play. 

\begin{rmk}
For the categorical foundations in~\cite{ST11} (building on~\cite{Pseudocategories}), the composition functor uses the strict fibered product of stacks. This creates problems when attempting to define the correct notion of equivalence of internal categories, as the composition functor does not play well with equivalence of stacks. However, there is an easy fix: if the source and target functors $\s$ and $\t$ are fibrations of stacks, this strict fibered product coincides with the usual (homotopy) fibered product of stacks. Geometric bordism categories have this fibration property; see Remark~\ref{rmk:fibration}. 
\end{rmk}

\begin{defn}[{\cite[Definition~2.18]{ST11}}]\label{defn:internalfunctor}
For internal categories ${\sf C}$, ${\sf D}$, an \emph{internal functor} $F=(F_0,F_1,\mu,\varepsilon)\colon {\sf C}\to {\sf D}$ consists of 1-morphisms $F_0\colon \Ob({\sf C})\to \Ob({\sf D})$, $F_1\colon \Mor({\sf C})\to \Mor({\sf D})$, and 2-morphisms $\mu$ and $\epsilon$,
\beq
&&\begin{tikzpicture}[baseline=(basepoint)];
\node (A) at (0,0) {$\Mor({\sf C})\times_{\Ob({\sf C})} \Mor({\sf C})$};
\node (B) at (4,0) {$\Mor({\sf C})$};
\node (C) at (0,-1.5) {$\Mor({\sf D})\times_{\Ob({\sf D})}\Mor({\sf D})$};
\node (D) at (4,-1.5) {$\Mor({\sf D})$};
\node (E) at (2,-.75) {$\mu \ \twocommute$};
\draw[->] (A) to node [above] {$\c$} (B);
\draw[->] (A) to node [left] {$F_1\times F_1 $} (C);
\draw[->] (C) to node [below] {$\c$} (D);
\draw[->] (B) to node [right] {$F_1$} (D);
\path (0,-.75) coordinate (basepoint);
\end{tikzpicture}
\qquad \begin{tikzpicture}[baseline=(basepoint)];
\node (A) at (0,0) {$\Ob({\sf C})$};
\node (B) at (4,0) {$\Ob({\sf D})$};
\node (C) at (0,-1.5) {$\Mor({\sf C})$};
\node (D) at (4,-1.5) {$\Mor({\sf D}).$};
\node (E) at (2,-.75) {$\varepsilon \ \twocommute$};
\draw[->] (A) to node [above] {$F_0$} (B);
\draw[->] (A) to node [left] {$\u$} (C);
\draw[->] (C) to node [below] {$F_1$} (D);
\draw[->] (B) to node [right] {$\u$} (D);
\path (0,-.75) coordinate (basepoint);
\end{tikzpicture}\label{eq:internalfunctordiag}
\eeq
These data are required to satisfy further properties involving the associator and unitors. 
\end{defn}

\begin{defn}[{\cite[Definition~2.19]{ST11}}]\label{defn:internalnattrans}
An \emph{internal natural transformation} $(\eta,\rho)\colon F\Rightarrow G$ between internal functors $F,G\colon {\sf C}\to {\sf D}$ is a 1-morphism $\eta\colon \Ob({\sf C})\to \Mor({\sf D})$ and 2-morphism~$\rho$
\beq
&&\begin{tikzpicture}[baseline=(basepoint)];
\node (A) at (0,0) {$\Mor({\sf C})$};
\node (B) at (6,0) {$\Mor({\sf D})\times_{\Ob({\sf D})}\Mor({\sf D})$};
\node (C) at (0,-1.5) {$\Mor({\sf D})\times_{\Ob({\sf D})}\Mor({\sf D})$};
\node (D) at (6,-1.5) {$\Mor({\sf D}).$};
\node (E) at (3,-.75) {$\rho \ \twocommute$};
\draw[->] (A) to node [above] {$\eta\circ \t \times F_1$} (B);
\draw[->] (A) to node [left] {$G_1\times \eta\circ \s $} (C);
\draw[->] (C) to node [below] {$\c$} (D);
\draw[->] (B) to node [right] {$\c$} (D);
\path (0,-.75) coordinate (basepoint);
\end{tikzpicture}\label{diag:nu}
\eeq
These data are required to satisfy additional coherence properties. Define an \emph{isomorphism between internal natural transformations} $(\eta,\rho)\Rightarrow (\eta',\rho')$ as a 2-morphism $\eta\Rightarrow \eta'$ along which $\rho'$ pulls back to~$\rho$. 
\end{defn}

\begin{ex}\label{defn:superLiecat}
A \emph{super Lie category} is a category internal to supermanifolds whose source and target maps are submersions. This is a version of Definition~\ref{defn:internalcat} in which $\Ob({\sf C})$ and $\Mor({\sf C})$ are supermanifolds, i.e., representable stacks. There can be no higher coherence data in a super Lie category, so the associator and unitor are both necessarily trivial. A \emph{smooth functor} $F\colon {\sf C}\to {\sf D}$ between super Lie categories is a pair of maps of supermanifolds $F_0\colon \Ob({\sf C})\to \Ob({\sf D})$, $F_1\colon \Mor({\sf C})\to \Mor({\sf D})$ satisfying the axioms for a functor. A \emph{smooth natural transformation} $\eta\colon F\Rightarrow G$ is a map of supermanifolds $\eta\colon \Ob({\sf C})\to \Mor({\sf D})$ satisfying the axioms for a natural transformation between functors. By the Yoneda lemma for stacks, a smooth functor determines an internal functor, and a smooth natural transformation determines an internal natural transformation. 
\end{ex}

\begin{ex}\label{ex:Liegroupoiod}
As a special case of the previous example, every (super) Lie groupoid determines an internal category. In particular, every super Lie group $G$ determines the super Lie groupoid $*\sq G$, which in turn provides a category internal to stacks. For the (discrete) Lie groupoid determined by an ordinary manifold $M$, a graded gerbe on~$M$ determines an invertible algebra bundle (possibly after replacing $M$ with an equivalent Lie groupoid, see~\cite[Example 1.74]{Freedalg}). Hence, gerbes provide examples of twists for the internal category associated to~$M$. We refer to~\cite{Stoffeltwist} for a related connection between $1|1$-dimensional field theories and twisted K-theory. 
\end{ex}

%

\subsection{Twisted representations of internal categories}

\begin{defn}[Sketch of {\cite[Definition~5.1]{ST11}}] \label{defn:Morita}
Define the \emph{internal Morita category}~$\TA$ as the following category internal to symmetric monoidal stacks with involution. The stack $\Ob(\TA)$ classifies locally trivial bundles of nuclear Fr\'echet superalgebras and isomorphisms of algebra bundles. The stack $\Mor(\TA)$ classifies locally trivial bundles of nuclear Fr\'echet bimodules over topological algebras. Composition in $\TA$ is the projective tensor product of bimodule bundles,
$$
\c \colon \Mor(\TA)\times_{\Ob(\TA)}\Mor(\TA)\to \Mor(\TA),\qquad \c(B,B')=B\otimes_A B'.
$$ 
Viewing an algebra bundle as a bimodule over itself constructs the unit map $\u\colon \Ob(\TA)\to \Mor(\TA)$. Endow $\Ob(\TA)$ and $\Mor(\TA)$ with monoidal structures from the projective tensor product of algebra and bimodule bundles over the structure sheaf of the base supermanifold. Involutions on the stacks $\Ob(\TA)$ and $\Mor(\TA)$ come from the parity involution on bundles, and all structure maps are $\Z/2$-equivariant for these involutions. 
\end{defn}

\begin{rmk} \label{rmk:differ}
Definition~\ref{defn:Morita} deviates from the one in~\cite{ST11} in two ways. First, we require sheaves of topological vector spaces to be locally free, i.e., vector bundles. Stolz and Teichner's motivation to work with more general sheaves is explained in Remark~\ref{rmk:concordance}. Second, we specialize from general locally convex topological vector spaces to nuclear Fr\'echet spaces. We make this choice in order to guarantee that certain formal categorical structures in the functorial definition of quantum field theory agree with standard structures in physics: Fr\'echet spaces have the \emph{approximation property}~\cite[page~109]{Schaefer}, leading to a well-behaved theory for traces of nuclear operators, see~\cite[\S4.5]{STTraces}. It also happens that all the desired examples in this paper result in vector bundles of nuclear Fr\'echet spaces. We refer to~\cite[Appendix~2]{costbook} for a succinct overview of nuclear Fr\'echet spaces. 
\end{rmk}

\begin{rmk}
One can forget an involution or forget a monoidal structure on a stack. Hence Definition~\ref{defn:Morita} produces three internal categories: a category internal to symmetric monoidal stacks with involution, a category internal to symmetric monoidal stacks, and a category internal to stacks. We will use the notation $\TA$ to denote any one of these three options, where the particular choice should be clear from context. 
\end{rmk}


\begin{defn}
A \emph{twist} for an internal category ${\sf C}$ is an internal functor $\twist\colon {\sf C}\to \TA$. The \emph{trivial twist} is the functor
\beq\label{ex:one}
\one \colon {\sf C}\to \TA
\eeq
that is constant to the monoidal unit of $\TA$, e.g., all objects $Y\in \Ob({\sf C})(S)$ are sent to the trivial algebra bundle on~$S$, all objects $\Sigma\in \Mor({\sf C})(S)$ are sent to the trivial bimodule bundle over $S$, and all the morphisms in $\Ob({\sf C})(S)$ and $\Mor({\sf C})(S)$ to the trivial automorphism of the trivial bundle. 
\end{defn}

\begin{defn} \label{lem:fancyrep}
A \emph{$\twist$-twisted representation} of an internal category ${\sf C}$ is an internal natural transformation~$E$
\beq\label{eq:twistrepofinternal}
&&\begin{tikzpicture}[baseline=(basepoint)];
\node (A) at (0,0) {${\sf C}$};
\node (B) at (4,0) {$\TA.$};
\node (C) at (2.1,0) {$E \Downarrow$};
\draw[->,bend left=12] (A) to node [above] {$\one$} (B);
\draw[->,bend right=12] (A) to node [below] {$\twist$} (B);
\path (0,0) coordinate (basepoint);
\end{tikzpicture}
\eeq
Isomorphisms between twists are internal natural isomorphisms $\twist\Rightarrow\twist'$. Isomorphisms between twisted representations are isomorphisms between internal natural transformations~$E\simeq E'$. 
\end{defn}

\begin{ex}
The standard definition of a representation of a Lie groupoid~\cite[Definition~1.7.1]{Mackenzie} determines a twisted representation of the associated internal category (see Example~\ref{ex:Liegroupoiod}) for the trivial twist $\twist=\one$.  
\end{ex}

\begin{ex} 
When ${\sf C}$ is a Lie groupoid, Freed's invertible algebra bundles over ${\sf C}$ \cite[Definition 1.59]{Freedalg} determine twists for the associated internal category. Invertible algebra bundles determine twists for equivariant K-theory~\cite[Definition 1.78]{Freedalg}, see also Karoubi~\cite{KaroubiFrench} and Donovan--Karoubi \cite{DonovanKaroubi}. 
\end{ex} 

\begin{ex}\label{ex:easytwistedrep1}
For an ordinary manifold $M$ regarded as an internal category (see Example~\ref{ex:Liegroupoiod}), a twist is a bundle of algebras $A$ over $M$. A twisted representation is a bundle of $A$-modules $E\to M$.
\end{ex}

\begin{ex} \label{ex:easytwistedrep2}
For a super Lie group $G$ regarded as an internal category (see Example~\ref{ex:Liegroupoiod}), a twist is an algebra $A$ together with a bundle of $A-\C$ bimodules over $G$ with a map of bundles covering multiplication of $G$. When $A=\C$ this is the same data as a $\C^\times$-central extension of $G$. A twisted representation is an $A$-module~$E$ with (weak) $G$-action; when $A=\C$ this is the same data as a projective representation of~$G$. 
\end{ex}

\subsection{Reflection structures for twisted representations}\label{sec:RPgeneral}

The following notions of internal involution and equivariance data expand on~\cite[Definition 6.47]{HST}. In the following definition we emphasize that $\mathscr S$ is just a usual category, not an internal category. 

\begin{defn}[{\cite[Definition~B.1]{FreedHopkins}}]
For a category ${\mathscr S}$, an \emph{involution} is a functor $(\overline{\phantom{A}})\colon {\mathscr S}\to {\mathscr S}$ and a natural transformation $\varepsilon\colon \id_{\mathscr S}\Rightarrow (\overline{\phantom{A}})\circ (\overline{\phantom{A}})$ satisfying a compatibility condition for the triple composition of $(\overline{\phantom{A}})$ with itself. 
\end{defn}

\begin{ex}\label{ex:involutiononSMF}
Below we will consider two involutions on the category of supermanifolds. The first is the trivial involution, and the second is the \emph{conjugation functor} $(\overline{\phantom{A}})$ from Definition~\ref{defn:realsmfld} that reverses the complex structure on the structure sheaf of a supermanifold. In both of these cases, the involution squares to a functor that is equal to the identity, and so we take $\varepsilon$ to be the identity. 
\end{ex}



Below we view (symmetric monoidal) stacks (with involution) as a category fibered over $\SMfld$, and so one can ask for a map of stacks covering an involution of $\SMfld$. 

\begin{defn}\label{defn:interinvol}
Fix an involution $\SMfld\to \SMfld$ on the category of supermanifolds. For an internal category ${\sf C}$, an \emph{internal involution} is the data of an internal functor $\tau\colon {\sf C} \to {\sf C}$ covering the fixed involution of supermanifolds, an internal natural isomorphism $\eta \colon \id_{\sf C}\Rightarrow \tau\circ\tau$ covering the isomorphism data for the involution on $\SMfld$, and an isomorphism $\delta$ between internal natural isomorphisms $\tau\simeq \tau\circ \id_{\sf C} \Rightarrow \id_{\sf C}\circ \tau\simeq \tau$ in the diagram
\beq
\begin{tikzpicture}[baseline=(basepoint)];
\node (A) at (0,0) {${\sf C}$};
\node (B) at (3,0) {${\sf C}$};
\node (C) at (6,0) {${\sf C}$};
\node (D) at (9,0) {${\sf C}.$};
\node (E) at (3,.5) {$\eta \ \Downarrow$};
\node (E) at (6,-.5) {$\eta \ \Uparrow$};
\draw[->] (A) to node [above] {$\tau$} (B);
\draw[->,bend left] (A) to node [above] {$\id_{\sf C}$} (C);
\draw[->,bend right] (B) to node [below] {$\id_{\sf C}$} (D);
\draw[->] (B) to node [above] {$\tau$} (C);
\draw[->] (C) to node [above] {$\tau$} (D);
\path (0,-.75) coordinate (basepoint);
\end{tikzpicture}\nonumber
\eeq
The isomorphism $\delta$ satisfies a compatibility condition with $\eta$ for the quadruple composition of $\tau$ with itself. 
\end{defn}

\begin{defn}\label{defn:interreflection}
 A \emph{reflection} for an internal category ${\sf C}$ is an internal involution relative to the conjugation functor~\eqref{eq:conjsmfld} on $\SMfld$. 
\end{defn}

\begin{lem} \label{lem:TAreal}
The internal category $\TA$ admits a reflection from complex conjugation of algebras, bimodules and bimodule maps,
\beq
(\overline{\phantom{-}})\colon \TA\to \TA, \ \ (A\to S)\mapsto (\overline{A}\to \overline{S}),\label{eq:realTVS}
\eeq
where $(\overline{\phantom{-}})\circ (\overline{\phantom{-}})=\id$. 
\end{lem}



The following generalizes \cite[Definition~B.6]{FreedHopkins} to internal categories. 

\begin{defn} \label{defn:equivariancedata}
Let $({\sf C},\tau_{\sf C},\eta_{\sf C},\delta_{\sf C})$ and $({\sf D},\tau_{\sf D},\eta_{\sf D},\delta_{\sf D})$ be internal categories with internal involution, and let $F\colon {\sf C}\to {\sf D}$ be a functor. \emph{Equivariance data} for $F$ is an internal natural transformation~$\phi$
\beq
\begin{tikzpicture}[baseline=(basepoint)];
\node (A) at (0,0) {${\sf C}$};
\node (B) at (4,0) {${\sf D}$};
\node (C) at (0,-1.5) {${\sf C}$};
\node (D) at (4,-1.5) {${\sf D}$};
\node (E) at (2,-.75) {$\phi \ \twocommute$};
\draw[->] (A) to node [above] {$F$} (B);
\draw[->] (A) to node [left] {$\tau_{\sf C}$} (C);
\draw[->] (C) to node [below] {$F$} (D);
\draw[->] (B) to node [right] {$\tau_{\sf D}$} (D);
\path (0,-.75) coordinate (basepoint);
\end{tikzpicture}\label{eq:equivariantfunctor}
\eeq
%
and an isomorphism $\beta\colon \phi^2\circ \eta_{\sf C}\rightarrow \eta_{\sf D}$ between internal natural transformations,
\beq
\begin{tikzpicture}[baseline=(basepoint)];
\node (A) at (0,0) {${\sf C}$};
\node (B) at (4,0) {${\sf D}$};
\node (C) at (0,-1.5) {${\sf C}$};
\node (D) at (4,-1.5) {${\sf D}$};
\node (F) at (0,-3) {${\sf C}$};
\node (G) at (4,-3) {${\sf D}$};
\node (E) at (2,-.75) {$\phi \ \twocommute$};
\draw[->,bend right=50] (A) to node [left] {$\id_{\sf C}$} (F);
\draw[->,bend left=50] (B) to node [right] {$\id_{\sf D}$} (G);
\draw[->] (A) to node [above] {$F$} (B);
\draw[->] (A) to node [left] {$\tau_{\sf C}$} (C);
\draw[->] (C) to node [below] {$F$} (D);
\draw[->] (B) to node [right] {$\tau_{\sf D}$} (D);
\draw[->] (C) to node [left] {$\tau_{\sf C}$} (F);
\draw[->] (F) to node [below] {$F$} (G);
\draw[->] (D) to node [right] {$\tau_{\sf D}$} (G);
\node (H) at (2,-2.25) {$\phi \ \twocommute$};
\node (I) at (-.45,-1.5) {$\stackrel{\eta_{\sf C}}{\Rightarrow}$};
\node (I) at (4.45,-1.5) {$\stackrel{\eta_{\sf D}}{\Leftarrow}$};
\path (0,-1.5) coordinate (basepoint);
\end{tikzpicture}\label{eq:equivariantfunctor2}
\eeq
where $\beta$ satisfies a compatibility condition with $\delta_{\sf C}$ and $\delta_{\sf D}$. 
\end{defn}
\begin{rmk}\label{rmk:involutivereflection}
In the examples of interest in this paper, the coherence data $\eta$ and $\delta$ in Defintion~\ref{defn:interinvol} will be trivial. In such cases, the data in~\eqref{eq:equivariantfunctor2} becomes an isomorphism between $\phi^2$ and the identity on $F$, i.e., $\phi$ is an involutive 2-morphism. 
\end{rmk}

\begin{ex}
The trivial twist $\one\colon {\sf C}\to \TA$ has canonical equivariance data for any choice of reflection structure on the source. 
\end{ex}

Categories internal to stacks also have a notion of equivariance data for internal natural transformations, defined as follows. 

\begin{defn}  \label{defn:equivariancedata2}
Let  $({\sf C},\tau_{\sf C},\eta_{\sf C},\delta_{\sf C})$ and $({\sf D},\tau_{\sf D},\eta_{\sf D},\delta_{\sf D})$ be internal categories with internal involution, and let $F,G\colon {\sf C}\to {\sf D}$ be internal functors with equivariance data $(\phi_F,\beta_F)$ and $(\phi_G,\beta_G)$. \emph{Equivariance data} for an internal natural transformation $E\colon F\Rightarrow G$ is an isomorphism $\sigma\colon \phi_G\circ E\simeq E\circ \phi_F$  of natural transformations satisfying a compatibility condition with~$\beta_F$ and~$\beta_G$.
\end{defn}

The above combine to define reflection structures on twisted representations of internal categories. 

\begin{defn} \label{defn:RPfancyrep}
Given an internal category ${\sf C}$ with reflection structure $(\tau,\eta,\delta)$, a \emph{reflection structure} for a twist $\twist\colon {\sf C}\to \TA$ is equivariance data for $\twist$ (in the sense of Definition~\ref{defn:equivariancedata}) relative to the reflection structure on $\TA$ from Lemma~\ref{lem:TAreal}. Isomorphisms between twists with reflection structure are internal natural isomorphisms $\twist\Rightarrow\twist'$ with equivariance data (in the sense of Definition~\ref{defn:equivariancedata2}). 

Given a $\twist$-twisted representation $E$~\eqref{eq:twistrepofinternal} and a reflection structure for $\twist$, a \emph{reflection structure} for $E$ is equivariance data for the internal natural transformation~$E$ in the sense of Definition~\ref{defn:equivariancedata2}. Isomorphisms between twisted representations with reflection structures are isomorphisms between internal natural transformations $E\simeq E'$ that are compatible with the equivariance data~$\sigma$.
\end{defn}

\begin{ex}
Building on Example~\ref{ex:easytwistedrep1} of an ordinary manifold $M$ regarded as an internal category, a reflection structure for a twist is a real structure on a bundle of algebras~$A$ over~$M$. A reflection structure for a twisted representation $E$ is a real structure on~$E$ for which the~$A$-action is real. 
\end{ex}

\begin{ex} 
Building on Example~\ref{ex:easytwistedrep2} of a super Lie group $G$ regarded as an internal category, a reflection structure for a twist determined by $A$ is a real structure on $A$ and a compatible real structure on the bundle of $A-\C$ bimodules over~$G$. A reflection structure for a twisted representation $E$ is a real structure on $E$ compatible with the real structure on $A$ and the action by $G$. When $A=\C$, the real structure for the twist is the same data as a $\Z/2$-extension of $G$, and a real structure on $E$ is the same data as a real projective representation. 
\end{ex}


%


\subsection{Twisted geometric field theories and reflection structures}\label{sec:backFT}\label{sec:21EB}

See~\S\ref{sec:EBord} for a review of rigid geometries. 

\begin{defn}[{Sketch of \cite[Definitions~2.21, 2.46, end of 2.48, and 4.4]{ST11}}]\label{defn:21EB}
Given a rigid geometry $(G,\M)$ with involution, let $(G,\M)\Bord(M)$ denote Stolz and Teichner's $(G,\M)$-bordism category over $M$ as a category internal to symmetric monoidal stacks with involution. Hence, $(G,\M)\Bord(M)$ includes the data of symmetric monoidal stacks with involution
\beq\label{Eq:objectmor21bord}
\Ob((G,\M)\Bord(M)),\qquad \Mor((G,\M)\Bord(M))
\eeq
as well as source, target, unit and compositions functors between symmetric monoidal stacks with involution
\beq\label{Eq:comp21bord}
\begin{array}{c}\s,\t\colon \Mor((G,\M)\Bord(M))\to \Ob((G,\M)\Bord(M)),\\ \u\colon \Ob((G,\M)\Bord(M))\to \Mor((G,\M)\Bord(M)),\\
\c\colon\Mor((G,\M)\Bord(M))^{[2]}\to \Mor((G,\M)\Bord(M)).\end{array}
\eeq
There is associator and unitor data that satisfy further compatibility properties. 
\end{defn}

Above and throughout, we use the notation 
\beq\nonumber
\resizebox{\textwidth}{!}{$
\Mor((G,\M)\Bord(M))^{[k]}:=\underbrace{\Mor((G,\M)\Bord(M))\times_{\Ob((G,\M)\Bord(M))} \cdots \times_{\Ob((G,\M)\Bord(M))} \Mor((G,\M)\Bord(M))}_{k-{\rm fold\ fibered\ product}}
$}
\eeq
to denote the $k$-fold fibered product of $\Mor((G,\M)\Bord(M))$ over $\Ob((G,\M)\Bord(M))$ with respect to the maps $\s$ and $\t$; hence, $S$-points of $\Mor((G,\M)\Bord(M))^{[d]}$ are $k$-tuples of composable bordisms. 
\begin{rmk}\label{rmk:fibration}
In \cite[Appendix~B.3]{DBEEFT}, it is shown that the source and target functor~\eqref{Eq:comp21bord} in the construction of $(G,\M)\Bord(M)$ are fibrations of stacks, and hence the strict fibered product agrees with the weak fibered product of stacks. We use this fact throughout when considering compositions of bordisms. 
\end{rmk}

We sketch the construction of $(G,\M)\Bord(M)$ as an internal category. The stack $\Ob((G,\M)\Bord(M))$ has objects specified by data: (i) a $(G,\M)$-pair $(Y^c,Y)$ over $S$, (ii) a map~$\Phi\colon Y\to M$, and (iii) a partition of~$Y\setminus Y^c$:
\beq\label{Eq:somemoredata}
&&(Y^c\subset Y),\qquad Y\xrightarrow{\Phi} M,\qquad Y\setminus Y^c\simeq Y^-\coprod Y^+.
\eeq
The $S$-family $Y\to S$ is a \emph{collar} and the partition specifies an \emph{incoming collar} $Y^-$ and an \emph{outgoing collar} $Y^+$. The \emph{core} $Y^c\to S$ is required to be proper over~$S$. These objects admit an involution determined by the involution on the $(G,\M)$-geometry, and a monoidal structure from disjoint union of supermanifolds. Isomorphisms in $\Ob((G,\M)\Bord(M)$ are spans of $(G,\M)$-isometric embeddings $Y\hookleftarrow W \hookrightarrow Y'$ satisfying compatibility properties with $Y^c$ and $Y^\pm$, up to an equivalence relation that allows for shrinking of collars.

Similarly, an object of the stack $\Mor((G,\M)\Bord(M))$ is specified by data
\beq\label{eq:vaguedata}
&&Y^c_\inn\subset Y_\inn\hookrightarrow \Sigma\hookleftarrow Y_\out\supset Y^c_\out,\qquad \Phi\colon \Sigma\to M,
\eeq
where $(Y^c_{\inn/\out}, Y_{\inn/\out})$ are $(G,\M)$-pairs associated with the source and target of the bordism. There is an analogously defined core $\Sigma^c:=\Sigma\setminus (i_{\rm in}(Y_{\rm in}^-)\cup i_{\rm out}(Y_{\rm out}^+))$ that is required to be proper over~$S$. Isomorphisms in $\Mor((G,\M)\Bord(M))$ are again specified by spans of isometric embeddings $\Sigma\hookleftarrow W \hookrightarrow \Sigma'$ with compatibility properties subject to an equivalence relation that allows for shrinking of collars. We note that the complete definition replaces the inclusions in~\eqref{eq:vaguedata} with spans, which is essential for the fibration property described in Remark~\ref{rmk:fibration}. We will only require bordisms of the form~\eqref{eq:vaguedata} in this paper.

\begin{defn}[{\cite[Definition~5.2]{ST11}}]\label{defn:FTmain2}
A \emph{twisted geometric field theory} is a twisted representation of $(G,\M)\Bord(M)$, i.e., an internal natural transformation~$E$
\beq
&&\begin{tikzpicture}[baseline=(basepoint)];
\node (A) at (0,0) {$(G,\M)\Bord(M)$};
\node (B) at (5,0) {$\TA$};
\node (C) at (2.5,0) {$E \Downarrow$};
\draw[->,bend left=15] (A) to node [above] {$\one$} (B);
\draw[->,bend right=15] (A) to node [below] {$\twist$} (B);
\path (0,0) coordinate (basepoint);
\end{tikzpicture}\label{eq:twistedGFT}
\eeq
where $\twist$ is a \emph{twist functor} and $\one$ is~\eqref{ex:one}. Twisted field theories are the objects of a groupoid whose morphisms are isomorphisms between internal natural transformations. 
\end{defn}

Reflection structures for field theories are special cases of Definitions~\ref{defn:equivariancedata} and~\ref{defn:equivariancedata2}.

\begin{defn}\label{defn:RP0}
Given a reflection $\dagger$ of $(G,\M)\Bord(M)$, a \emph{reflection structure} for a twist and \emph{reflection structure} for a twisted field theory is given by Definition~\ref{defn:RPfancyrep}.
\end{defn}

\begin{lem} \label{ex:orientationGM}
A reflection structure on a geometric bordism category $(G,\M)\Bord(M)$ is determined by a map $\dagger\colon \M\to \overline{\M}$ and a homomorphism $\alpha\colon G\to \overline{G}$ satisfying 
\beq\label{eq:anotherdagger}
\overline{\dagger}\circ\dagger=\id_\M,\qquad \dagger|_{\M^c}=\overline{\M}^c\subset \overline{\M},\qquad \overline{\alpha}\circ \alpha=\id_G.
\eeq
\end{lem}

\begin{proof}
The statement follows from the argument in~\cite[Lemma 6.19]{HST}. Given an $S$-family of $(G,\M)$-supermanifolds $(Y\to S)$, the map~\eqref{eq:anotherdagger} gives a unique $(G,\M)$-stucture for the $\overline{S}$-family $(\overline{Y}\to \overline{S})$. This extends to an internal functor 
\beq
&&(G,\M)\Bord(M)\to (G,\M)\Bord(M),\qquad (S\leftarrow Y\to M)\mapsto (\overline{S}\leftarrow \overline{Y}\to \overline{M}\simeq M).\label{eq:realstru}
\eeq
The identification $\overline{M}\simeq M$ uses that $M$ is an ordinary manifold, and hence has a canonical real structure when regarded as a supermanifold. Via the forgetful functor $(G,\M)\Bord(M)\to \SMfld$ that sends an $S$-family of bordisms to $S$, the internal functor~\eqref{eq:realstru} covers the conjugation functor on supermanifolds. 
\ep

The assignment $M\mapsto (G,\M)\Bord(M)$ is natural in $M$: a smooth map $f\colon M'\to M$ induces an internal functor
\beq\label{eq:Ebordnat}
f_*\colon (G,\M)\Bord(M')\to (G,\M)\Bord(M)
\eeq
by precomposing the map $\Phi$ in~\eqref{eq:vaguedata} with $f$. This implies that twisted geometric field theories (with reflection structures) determine a contravariant functor from manifolds to groupoids. 

We will primarily be interested in the \emph{$2|1$-dimensional Euclidean bordism category over~$M$}, denoted  $2|1\EBord(M)$, defined as the geometric bordism category for the geometry $(\E^{2|(0,1)}\rtimes \Spin(2),\R^{2|1},\R^{1|1})$ from Definition~\ref{defn:sEuc} with involution from the spin flip, see Lemma~\ref{lem:flip1}. We describe some structures on the $2|1$-Euclidean bordism category and the affect of these structure on twisted field theories.

\begin{lem} \label{lem:extrastrEB0}
For $\mu\in \R_{>0}$, the internal category $2|1\EBord(M)$ admits an internal endofunctor $\RG_\mu\colon 2|1\EBord(M)\to 2|1\EBord(M)$ that dilates a super Euclidean bordism by $\mu$. Furthermore, for $\mu,\mu'\in \R_{>0}$ there are internal natural isomorphisms $\RG_\mu\circ \RG_{\mu'}\simeq \RG_{\mu\mu'}$.
\end{lem} 
\begin{proof}[Sketch of proof.]
The functors $\RG_\mu$ are inherited from the endofunctors on $2|1$-Euclidean supermanifolds given by Lemma~\ref{lem:RGgeneral}. This determines maps on the object and morphisms stacks of $2|1\EBord(M)$, and it remains to supply the 2-morphisms in~\eqref{eq:internalfunctordiag} and verify their compatibilities. 
\ep

\begin{defn}\label{defn:RG21} The \emph{renormalization group flow} is the (weak) $\R_{>0}$-action on twisted $2|1$-Euclidean field theories over $M$ gotten from pulling back internal natural transformations~\eqref{eq:twistedGFT} along the endofunctors of the $2|1$-Euclidean bordism category in Lemma~\ref{lem:extrastrEB0}. 
\end{defn}

\begin{lem} \label{lem:extrastrEB}
The internal category $2|1\EBord(M)$ admits a reflection  
\beq\label{eq:21EBreflection}
\dagger \colon 2|1\EBord(M)\to 2|1\EBord(M),\qquad \dagger\circ\dagger\simeq \id.
\eeq
\end{lem} 
\begin{proof}
The involution of $2|1\EBord(M)$ is inherited from the involution $2|1$-Euclidean supermanifolds in Definition~\ref{defn:dagger}.
\ep

\subsection{The degree~$n$ twist on the Euclidean bordism category}\label{sec:appendfer}

We review the \emph{free fermion twist} (also called the \emph{degree twist}) from \cite[\S5-6]{ST11}, which is a particular collection of $\otimes$-invertible internal internal functors
\beq\label{eq:chiraltwist}
\Fer_n \colon 2\EBord \to \TA,\qquad \Fer_n\otimes \Fer_n\simeq \Fer_{n+m}
\eeq
indexed by $n\in \Z$ together with the indicated internal natural isomorphisms. Above, $2\EBord=2\EBord(\pt)$ is the \emph{2-dimensional Euclidean bordism category} (over $M=\pt$), i.e., the geometric bordism category associated to the 2-dimensional Euclidean geometry defined in~\S\ref{sec:dumbEuc}. The morphisms in this bordism category are (roughly) flat 2-manifolds with spin structure and geodesic boundary. 

\begin{rmk}
Free fermions haven been studied extensively in the mathematical literature in a variety of languages, e.g., see~\cite{Wassermann,SegalCFT,ST04,DouglasHenriques}. 
\end{rmk}

We begin by reviewing some families of Euclidean bordisms from~\cite[\S3.2]{ST11}. There are maps of stacks~\cite[page~34]{ST11}
$$
\sS^\pm\colon \R_{>0}\to \Ob(2\EBord),\qquad \ell\mapsto \rS_\ell^\pm
$$
sending $\ell\in \R_{>0}$ to a circumference $\ell$ Euclidean circle with periodic spin structure, and choice of collar (indicated by the $\pm$) analogous to~\eqref{eq:collardataI}. In analogy to~\eqref{eq:CM} and~\eqref{eq:LM}, there are maps of stacks~\cite[page~34-35]{ST11}
\beq
\sC^\pm\colon (\R_{>0}\times \HH_\pm)/\Z\to \Mor(2\EBord), && (\ell,\tau)\mapsto \rC_{\ell,\tau}^\pm\nonumber \\
\sL^\pm\colon (\R_{>0}\times \HH_\pm)/\Z\to \Mor(2\EBord),&& (\ell,\tau)\mapsto \rL_{\ell,\tau}^\pm\label{eq:someEucbordisms} \\
\sR^\pm\colon (\R_{>0}\times \iHH_\pm)/\Z\to \Mor(2\EBord),&& (\ell,\tau)\mapsto \rR_{\ell,\tau}^\pm\nonumber 
\eeq
classifying families of bordisms whose cores are flat cylinders of diameter $\ell$ and height $\im(\tau)$. The distinction in the 3 functors~\eqref{eq:someEucbordisms} comes from the choice of source and target data. Composition among the families~\eqref{eq:someEucbordisms} is determined by addition in $\HH_\pm$; we refer to the pictures~\eqref{pic:pic1},~\eqref{pic:pic2}, and \eqref{pic:pic3} for depictions of some of these compositions. Formulas can be recovered from~\eqref{eq:comp2iso}, \eqref{eq:adjunctioncyl}, \eqref{eq:thickdiagram}, and \eqref{eq:tracerelations}. For example, concatenation of Euclidean cylinders yields~\cite[Proposition 3.5]{ST11}, see~\eqref{pic:pic1}
$$
\rC_{\ell,\tau'}^\pm\circ \rC_{\ell,\tau}^\pm\simeq \rC_{\ell,\tau'+\tau}^\pm.
$$

Following~\cite[\S5.3]{ST11}, we will describe the values of $\Fer_{-1}$ in~\eqref{eq:chiraltwist} on the bordisms~\eqref{eq:someEucbordisms} and their compositions. 
Define algebra bundles $\Fer(\sS^\pm)\to \R_{>0}$ whose fiber at~$\rS^\pm_\ell$ is~\cite[Equation~6.1]{ST11}
\beq\label{eq:Ferdefn}
\Fer_\ell^\pm :=\Fer(\rS_\ell^\pm)=\cCl_{\pm 1} \otimes \bigotimes_{m\in \N} \cCl(\C_{\ell,\pm m}\oplus \C_{\ell,\mp m}),
\eeq
for the following definitions of the above factors. For $m\in \R$, let $\C_{\ell,m}$ denote the 1-dimensional representation of~$\R$ with action by $e^{2\pi i mt/\ell}$ for $t\in \R$. The Clifford algebras $\cCl(\C_{\ell,m}\oplus \C_{\ell,-m})$ are determined by the $\R$-invariant (hyperbolic) pairing
\beq\label{eq:Cliffordpairing}
&&\omega\colon \C_{\ell,m}\oplus \C_{\ell,-m}\to \C,\qquad  \omega((v_m,v_{-m}),(w_m,w_{-m}))=v_mw_{-m}+v_{-m}w_m.
\eeq
Finally, the infinite tensor product in~\eqref{eq:Ferdefn} is the \emph{restricted} tensor product,
i.e., the closure of finite sums of tensor products $\otimes_m v_m$ where $v_m=1$ for all but finitely many $m$. Throughout we will use the canonical isomorphism
$$
(\Fer_\ell^\pm)^\op=\Fer_\ell^\mp
$$
whose map on underlying sets is the identity. 

Next construct an $\HH_\pm$-action on the algebras~\eqref{eq:Ferdefn}, following~\cite[page~59]{ST11}. Let $\tau\in \HH_\pm$ act by $e^{2\pi i m \tau/\ell }$ on $\C_{\ell, m}$. When applied diagonally to $\C_{\ell,m}\oplus \C_{\ell,-m}$, this $\HH_\pm$-action preserves the pairing~\eqref{eq:Cliffordpairing}, and so determines an automorphism 
\beq\label{eq:holoHaction}
\varphi_{\ell,m}^\pm(\tau)\colon \cCl(\C_{\ell,m}\oplus \C_{\ell,-m})\to \cCl(\C_{\ell,m}\oplus \C_{\ell,-m}),\qquad \tau \in \HH_\pm
\eeq
extending the $\R$-action on~\eqref{eq:Ferdefn}. Varying $m$, these automorphisms assemble into an $\HH_\pm$-family of automorphisms $\varphi^\pm_\ell(\tau) \colon \Fer_\ell^\pm\to \Fer_\ell^\pm$. Define bundles of bimodules over $\R_{>0}\times \HH_\pm$ whose fibers are denoted $_{^\tau\Fer_\ell^\pm}(\Fer_\ell^\pm)_{\Fer_\ell^\pm}$ and defined by taking $\Fer_\ell^\pm$ as a bimodule over itself where the left action is twisted by the automorphism $\varphi^\pm_\ell(\tau)$. Since the algebra  maps $\varphi_{\ell,m}(\tau+\ell)=\varphi_{\ell,m}(\tau)$ are equal, this bundle of bimodules descends to the quotients $(\R_{>0}\times \HH_\pm)/\Z$. The twist~\eqref{eq:chiraltwist} assigns the bimodule bundles over the families~\eqref{eq:someEucbordisms} whose fibers at~$(\ell,\tau)$ are
\beq
\Fer_{-1}(\sC^\pm)_{\ell,\tau}=\Fer_{-1}(\rC_{\ell,\tau}^\pm)&=&_{^\tau\Fer_\ell^\mp}(\Fer_\ell^\mp)_{\Fer_\ell^\mp}\nonumber\\
\Fer_{-1}(\sL^\pm)_{\ell,\tau}=\Fer_{-1}(\rL_{\ell,\tau}^\pm)&=&(\Fer_\ell^\mp)_{\Fer_\ell^\mp\otimes ^\tau\Fer_\ell^\pm}\label{eq:degtwisbimod}\\
\Fer_{-1}(\sR^\pm)_{\ell,\tau}=\Fer_{-1}(\rR_{\ell,\tau}^\pm)&=&_{\Fer^\pm_\ell \otimes {}^\tau\Fer_\ell^\mp} \Fer_\ell^\mp.\nonumber
\eeq
The value of~\eqref{eq:chiraltwist} on the moduli of Euclidean tori is the Pfaffian line bundle~\cite[page~56]{ST11}
\beq\label{eq:PfvalofFer}
\Fer_{-1}(\Lat\sq (\E^2\rtimes \Spin(2)\times \SL_2(\Z)))=\Pf,
\eeq
which we shall regard as an $\E^2\rtimes \Spin(2)\times \SL_2(\Z)$-equivariant structure on the trivial bundle on lattices $\Lat$; we refer to~\S\ref{sec:backtot} for a review of this line bundle. 

\begin{rmk}\label{rmk:antiperiodic}
The values on the antiperiodic versions of the Euclidean cylinders~\eqref{eq:someEucbordisms} follows a similar discussion as above, but starting with the bundle of algebras over $\R_{>0}\coprod  \R_{>0}$ whose fiber at $(\ell,\pm)$ is~\cite[Equation~6.1]{ST11}
\beq\label{eq:antiperFer}
 \bigotimes_{m\in \N_0+1/2} \cCl(\C_{\ell,\pm m}\oplus \C_{\ell,\mp m}).
\eeq
The value of $\Fer_{-1}$ on the other components~\eqref{eq:pointedEuctori} of the moduli stack of spin tori is again the Pfaffian line, which on these components is canonically trivialized by the Pfaffian section (the Dirac operator is invertible for the even spin structure). 
\end{rmk}

A priori, an internal functor~\eqref{eq:chiraltwist} requires additional data for compositions and units. Most of this is supplied from the fact the compositions of the bimodules~\eqref{eq:degtwisbimod} corresponds to composition of the algebra maps $\varphi^\pm_\ell(\tau)$, e.g., there are canonical bimodule maps 
\beq\label{eq:explicitcomposition}
&&\Fer_{-1}(\rC_{\ell,\tau'})\circ \Fer_{-1}(\rC_{\ell,\tau})= \Fer_{-1}(\rC_{\ell,\tau'})\otimes_{^\tau\Fer_\ell^\mp} \Fer_{-1}(\rC_{\ell,\tau})\to \Fer_{-1}(\rC_{\ell,\tau'+\tau})
\eeq
since $\varphi_\ell^\pm(\tau')\circ \varphi^\pm_\ell(\tau)=\varphi_\ell^\pm(\tau'+\tau).$ The remaining data to be specified corresponds to compositions that glue cylinders into tori~\cite[Equation~3.11]{ST11} (see the picture~\eqref{pic:pic4}). For such a composition the twist assigns the data of an isomorphism of line bundles over~$\R_{>0}\times \iHH_\pm$
\beq\label{eq:Fertrace}
 \Fer_{-1}(\sL^\mp)_{\ell,\tau}\circ \Fer_{-1}(\sR^\pm)_{\ell,\tau}\simeq \Fer^\pm_\ell \otimes_{\Fer^\pm_\ell \otimes {}^\tau\Fer_\ell^\mp} \Fer_\ell^\mp \stackrel{{\rm data}}{\simeq}\Pf|_{(\ell,\tau)}
\eeq
between the indicated tensor product and the restriction of the Pfaffian along the inclusion $\R_{>0}\times \iHH_+\hookrightarrow \Lat$. The above is essentially a choice of trace on the bimodules ${}^\tau\Fer_\ell^\pm$, which we characterize in Lemma~\ref{lem:Tracenormalize1} below. 

We obtain further structure on the algebras~\eqref{eq:Ferdefn} from evaluating the twist on isometries, e.g., the isometry $\rS_\ell^\pm\simeq\rS_\ell^\mp$ determined by the action of $\sqrt{-1}\in U(1)\simeq \Spin(2)$ (analogous to Lemma~\ref{lem:orientation11restrict}) gives an isomorphism 
\beq\label{eq:Ferinvol}
&&\cCl_{\mp 1} \otimes \bigotimes_{m\in \N} \cCl(\C_{\ell,\mp m}\oplus \C_{\ell,\pm m})\stackrel{{\rm data}}{\simeq} \cCl_{\pm 1} \otimes \bigotimes_{m\in \N} \cCl(\C_{\ell,\pm m}\oplus \C_{\ell,\mp m}),
\eeq
that squares to the parity involution. The map~\eqref{eq:Ferinvol} is unique up to sign, given by
\beq\label{eq:Cliffinv}
 \cCl_{\mp 1}\to \cCl_{\pm 1},\quad v\mapsto \pm iv, \quad v\in \C
\eeq
 $$
 \cCl(\C_{\ell, \mp m}\oplus \C_{\ell,\pm m})\to \cCl(\C_{\ell, \pm m}\oplus \C_{\ell,\mp m}),\qquad  (v,w)\mapsto (iw,iv), \quad (v,w)\in \C_m\oplus \C_{-m}).
$$
Similarly, a reflection structure on $\Fer_{-1}$ is determined by the conjugate-linear involutive isomorphism of algebras, 
\beq
 \cCl_{\mp 1}&\to& \overline{\cCl}_{\pm 1},\quad v\mapsto  \pm i\overline{v}, \quad v\in \C\nonumber\\
 \cCl(\C_{\ell, \mp m}\oplus \C_{\ell,\pm m})&\to& \overline{\cCl(\C_{\ell, \pm m}\oplus \C_{\ell,\mp m})}\simeq \cCl(\overline{\C}_{\ell, \mp m}\oplus \overline{\C}_{\ell,\pm m})\label{eq:FerRP}\\
 (v,w)&\mapsto& (\pm i\overline{w},\mp i\overline{v}), \quad (v,w)\in \C_{\mp m}\oplus \C_{\pm m})\nonumber
\eeq
i.e., the $\Z/2$-graded version of the $*$-structure~\eqref{eq:starCl} on the Clifford algebras. The above determines conjugate-linear maps between the bimodules~\eqref{eq:degtwisbimod}, and so determines the data gotten by restriction of a reflection structure, compare Lemma~\ref{lem:RPtwist}.

\begin{rmk}
The degree~$n$ twist is defined as a functor out of the Euclidean bordism category in~\cite[page~56]{ST11}, and then described as a functor out of the ``conformal Euclidean" bordism category in~\cite[\S6]{ST11}. We indicate the translation between these descriptions. To begin, the conformal Euclidean bordism category includes the renormalization group action as an isometry, i.e., the conformal Euclidean geometry is $(\E^2\rtimes \C^\times,\R^2,\R)$ under the identification $\Spin(2)\times \R_{>0}\simeq \C^\times$. Decomposing into Fourier modes gives an $\R$-equivariant algebra isomorphism between~\eqref{eq:Ferdefn} and
the Clifford algebras on the spaces of (even) sections of the spinor bundle over the periodic circle~\eqref{eq:Fernoncircle}.
This recovers the values of Stolz and Teichner's Euclidean twist on objects from~\cite[page~56]{ST11}. We remark that~\eqref{eq:antiperFer} is the version of~\eqref{eq:Fernoncircle} where the spinor bundle is for the anti-periodic spin structure. The algebra bundles~\eqref{eq:Ferdefn} are isomorphic to trivial bundles with fiber $\Fer_{-1}(\rS^\pm_1)$ (i.e., the fiber at $\ell=1$); on the negative circle $\rS^-_\ell$ this recovers the conformal Euclidean twist~\cite[Equation~6.1]{ST11}. The trivialization is $\HH_\pm$-equivariant relative to homomorphisms between semigroups at each fiber. This sketches the data of a 2-commutative triangle
\beq
&&\begin{tikzpicture}[baseline=(basepoint)];
\node (A) at (0,0) {$2\EBord$};
\node (B) at (5,-.75) {$\TA.$};
\node (C) at (2,-.75) {$\Downarrow$};
\node (D) at (0,-1.5) {$2\CEBord(\pt)$};
\draw[->] (A) to node [left] {forget} (D);
\draw[->] (A) to node [above] {$\Fer_{-1}$} (B);
\draw[->] (D) to node [below=5pt] {${\rm c}\Fer_{-1}$} (B);
\path (0,-.75) coordinate (basepoint);
\end{tikzpicture}\nonumber
\eeq
where $\Fer_{-1}$ is~\eqref{eq:chiraltwist} defined in~\cite[page~56]{ST11}, whereas ${\rm c}\Fer_{-1}$ is the conformal version of this functor from~\cite[\S6]{ST11}, and the 2-commuting data is essentially determined by the $\HH_\pm$-equivariant isomorphism of algebra bundles over $\R_{>0}\coprod \R_{>0}$. Said differently, the degree twist in dimension~2 is RG-invariant, which in this framework demands the additional (descent) data in the 2-commuting diagram above. 
\end{rmk}

We recall the internal $\Ann=\Ann(\pt)$ of Euclidean cylinders from Definition~\ref{defn:Path}. 

\begin{lem}[{\cite[Equation~6.2]{ST11}}]\label{lem:ClFeriso}
Consider the restriction of $\Fer_n$ to the subcategory~$\Ann\hookrightarrow 2\EBord$ of Euclidean cylinders. There is an internal natural isomorphism 
\beq
&&\begin{tikzpicture}[baseline=(basepoint)];
\node (A) at (0,0) {$\Ann$};
\node (B) at (5,0) {$\TA$};
\node (C) at (2.5,0) {$\mathcal{F} \Downarrow$};
\draw[->,bend left=10] (A) to node [above] {$\Cl_n$} (B);
\draw[->,bend right=10] (A) to node [below] {$\Fer_n$} (B);
\path (0,0) coordinate (basepoint);
\end{tikzpicture}\label{eq:ClFeriso}
\eeq
from the internal functor $\Cl_n$ in Definition~\ref{defn:Clifftwistfornormy}. Furthermore, the internal natural isomorphism $\mathcal{F}$ has a reflection structure relative to the reflection structures on the internal functors determined by the conjugate linear anti-involution~\eqref{eq:FerRP} of $\Fer_n$, and the $*$-superalgebra structure
\beq\label{eq:CliffRP}
 \cCl_{\pm n}&\to& \overline{\cCl}_{\mp n},\quad v\mapsto  \mp i\overline{v}, \quad v\in \C^n.
\eeq
\end{lem}
\bp
To construct an internal natural isomorphism~$\mathcal{F}$ in~\eqref{eq:ClFeriso} it suffices to consider the case $n=-1$. The basic ingredient is a bundle of invertible bimodules over $\R_{>0}\coprod \R_{>0}$ with fibers~\cite[Equation~6.2]{ST11},
\beq\label{eq:Morita}
{}_{\Fer_\ell^\mp}(\cCl_{\mp 1} \bigotimes_{m>0} \Lambda^\bullet \C_{\ell,\pm m})_{\cCl_{\mp 1}}
\eeq
where $\cCl_{\pm 1}$ is regarded as a bimodule over itself, and $\Lambda^\bullet \C_{\ell,\pm m}$ is a left $\cCl(\C_{\ell,\pm m}\oplus \C_{\ell,\mp m})$-module for the action determined by the Lagrangian subspace $\C_{\ell,\pm m}\subset \C_{\ell,\pm m}\oplus \C_{\ell,\mp m}$.\footnote{For example, see~\cite[Definition 2.2.4]{ST04} for Clifford modules constructed out of Lagrangian subspaces.} The bundle~\eqref{eq:Morita} provides part of the data of a map of stacks $\Ob(\Ann)\to \Mor(\TA)^\times$; the remaining data is a $\Z/4$-equivariant structure determined by the action of $\pm i$ on $\C_{\ell,\pm m}$ and the map $\cCl_{\pm 1}\to \cCl_{\mp 1}$ in~\eqref{eq:Cliffinv}. This furnishes the first part of the data of the internal natural isomorphism~\eqref{eq:ClFeriso}. The remaining data are $\Z/4$-equivariant isomorphisms of bimodules over $(\R_{>0}\times \HH_\pm)/\Z$
\beq\label{eq:Moritamapstuff}
&& _{^\tau\Fer_\ell^\mp}(\Fer_\ell^\mp) \otimes_{\Fer_\ell^\mp} (\cCl_{\mp 1} \bigotimes_{m>0} \Lambda^\bullet \C_{\ell,\pm m})_{\cCl_{\mp 1}}\to {}_{\Fer_\ell^\mp} (\cCl_{\mp 1} \bigotimes_{m>0} \Lambda^\bullet \C_{\ell,\pm m})\otimes_{\cCl_{\mp 1}} {\cCl_{\mp 1}}.
\eeq
Identifying the source with the bimodule~\eqref{eq:Morita} but with left action twisted by $\varphi^\pm_\ell(\tau)$, the map~\eqref{eq:Moritamapstuff} is given by
\beq\label{eq:moritaHact}
\C_{\ell,\pm m}\to \C_{\ell,\pm m}, \qquad v\mapsto e^{\pm 2\pi i \tau m/\ell}v. 
\eeq
This completes the construction of the natural isomorphism~\eqref{eq:ClFeriso}. In this case, the reflection structure on $\mathcal{F}$ is a property, corresponding to the compatibility of the anti-linear anti-involutions~\eqref{eq:FerRP} and~\eqref{eq:CliffRP}. 
\ep

Finally, we characterize the isomorphism~\eqref{eq:Fertrace} specified in Stolz and Teichner's construction; this map is implicit in~\cite[Lemma 2.3.14]{ST04} and~\cite[Proof of Theorem 1.16]{ST11}, though the presentation there is somewhat terse. For our purposes, one can take the following as a definition rather than a lemma. 
Up to a constant phase that can be absorbed into the trivialization of~$\Pf$, the following is also consistent with conventions in the physics literature, e.g., see~\cite[page~6]{GJF2}, as well as~\cite[Equation~3.3.10]{ST04}. 

\begin{lem}[{\cite[\S6]{ST11}}]\label{lem:Tracenormalize1}
There is a commutative triangle of line bundles over $\R_{>0}\times \iHH_+$
\beq
&&\begin{tikzpicture}[baseline=(basepoint)];
\node (A) at (0,0) {$\Fer^\mp_\ell \otimes_{\Fer^\mp_\ell \otimes {}^\tau\Fer_\ell^\pm} \Fer_\ell^\pm$};
\node (B) at (5,0) {$\underline{\C}\simeq \Pf|_{\R_{>0}\times \iHH_+}$};
\node (C) at (0,-1.5) {$\cCl_{\mp 1}\otimes_{\cCl_{\mp 1}\otimes \cCl_{\pm 1}} \cCl_{\pm 1}$};
\draw[->] (A) to node [above] {$\simeq$} (B);
\draw[->] (C) to node [left] {$\simeq$} (A);
\draw[->,dashed] (C) to (B);
\path (0,-.75) coordinate (basepoint);
\end{tikzpicture}\label{eq:maptoPf}
\eeq
where the upper isomorphism is the data~\eqref{eq:Fertrace} coming from the internal functor $\Fer$, and the trivialization of the Pfaffian over $\R_{>0}\times \iHH_+\hookrightarrow \Lat$ is the canonical one (from restricting the trivial bundle on $\Lat$). Then the dashed arrow is given by the map on sections
\beq
C^\infty(\R_{>0}\times \iHH_+; \cCl_{-1}/[\cCl_{-1},\cCl_{-1}]) &\to& C^\infty(\R_{>0}\times \iHH_+)\nonumber\\
 \Gamma&\mapsto& \ell^{1/2}\prod_{m=1}^\infty (1-e^{2 \pi i m \tau/\ell})\label{eq:traceformula}
\eeq
where $ \Gamma=\Gamma_1\in \cCl_1$ as defined in~\eqref{eq:Gamma}. 
\end{lem}
\begin{proof}[Sketch of proof.]
The invertible bimodule~\eqref{eq:Morita} supplies the vertical isomorphism in~\eqref{eq:maptoPf}. The upper isomorphism in~\eqref{eq:maptoPf} is a fiberwise trace on the bimodules ${}^\tau\Fer_\ell^\pm$. This is determined by a trace on the Clifford algebras modified by the character of the invertible bimodule~\eqref{eq:Morita}. Since the character of each factor $\Lambda^\bullet \C_{\ell,m}$ is $1-e^{2\pi i\tau m/\ell}$, the total character is the claimed infinite product. The factor of $\ell^{1/2}$ comes from the fact that the zero modes in~\eqref{eq:Fernoncircle} are the Clifford algebra associated with $\ell$ times the standard metric; see the description~\eqref{eq:Fernoncircle}. This Clifford algebra is isomorphic to $\cCl_{\mp 1}$, but the natural Clifford supertrace differs by a factor of $\ell^{1/2}$. Indeed, by~\cite[Proposition 3.2.1]{BGV} a Clifford supertrace is determined by a (metric) trivialization of $\Lambda^{\rm top}\R^n$, which introduces the factor of $\ell^{1/2}$.
%
%
\ep

\section{The families index in KO-theory}\label{sec:KOindex}

\subsection{Clifford algebras, spin structures, and Dirac operators} \label{sec:Cliffordandspin}
Define the (complex) Clifford algebras
\beq
\cCl_n:=\left\{\begin{array}{ll} \langle f_1,\dots f_n \mid [f_j,f_k]=-2\delta_{jk}\rangle & n \ge 0 \\ \langle e_1,\dots e_n \mid [e_j,e_k]=2\delta_{jk}\rangle & n<0\end{array}\right.
\label{eq:Clifford}
\eeq
where $\{e_1,\dots, e_n\}$ or $\{f_1,\dots,f_n\}$ denotes the standard basis of $\R^n$. These Clifford algebras have an evident real structure and also a $*$-structure determined by
\beq
f_j^\dagger=-f_j,\qquad e_j^\dagger=e_j.\label{eq:starCl}
\eeq
Define the element $\Gamma_n$ of the $n^{\rm th}$ Clifford algebra
\beq\label{eq:Gamma}
\Gamma_n:=\left\{\begin{array}{ll} 2^{-n/2}f_1f_2\dots f_n & n \ge 0 \\  2^{n/2}e_1e_2\dots e_n & n<0.\end{array}\right.
\eeq 
The \emph{Clifford supertrace} of a $\cCl_{-n}$-linear endomorphism $T\colon V\to V$ is (see~\S\ref{sec:supertrace})
\beq\label{eq:Cliffordsupertrace}
\sTr_{\cCl_n}(T):=\sTr(\Gamma_n\circ T).
\eeq

\begin{defn} A \emph{spin structure} for a rank $n$ real, oriented, metrized vector bundle $V\to M$ is a $\Spin(n)$-principal bundle $P_V\to M$ equipped with a double cover $P_V\to {\rm Fr}(V)$ of the oriented frame bundle of $V$ that is equivariant relative to $\Spin(n)\to \SO(n)$. 
\end{defn}

The collection of spin structures on $V$ forms a category (possibly empty) whose morphisms are isomorphisms of principal bundles compatible with the map to the frame bundle of $V$. A spin structure on a family $\pi\colon X\to M$ of Riemannian manifolds is a spin structure on the vertical tangent bundle $T(X/M)=\ker(d\pi)\subset TM$. We will require two types of bundles of Clifford modules built out of a spin structure. In the following, for a metrized vector bundle $V\to M$ the bundle $\Cl(V)\to M$ is the bundle of fiberwise Clifford algebras. 

\begin{defn}
Given a vector bundle with spin structure $V\to M$, the \emph{Euler bundle} is
\beq\label{eq:spinorsdefn2}
\bS_V=P_V \times_{\Spin(n)} \Cl_{-n},
\eeq
which is a bundle of invertible $\Cl(V)$-$\Cl_{-n}$-bimodules. Hence, one obtains a left $\Cl_n=\Cl_{-n}^\op$-module structure on $\bS_V$ that graded commutes with the $\Cl(V)$-action. 
For $\pi\colon X\to M$ a bundle of Riemannian manifolds with chosen spin structure on $T(X/M)$, the \emph{spinor bundle} is 
\beq\label{eq:spinorsdefn}
\bS_{X/M}=P_{T(X/M)} \times_{\Spin(n)} \Cl_{n},
\eeq
which is a bundle of invertible $\Cl(T(X/M))$-$\Cl_n$-bimodules. The right $\Cl_n$-action determines a left $\Cl_{-n}=\Cl_n^\op$-action that graded commutes with the left $\Cl(T(X/M))$-action. This allows for the construction of the $\Cl_{-n}$-linear Dirac operator supported on sections of~$\bS_{X/M}$. 
\end{defn}


\subsection{Clifford modules and KO-theory}
We overview a model for the KO-spectrum due to~\cite[\S1.3-1.4]{Cheung} and~\cite[Theorem~7.1]{HST}, see also~\cite[\S4]{DBEIndex}. Let $\mathcal{H}_n$ denote an infinite-dimensional $\Z/2$-graded separable Hilbert space with a self-adjoint action by~$\Cl_n$ with the property that every irreducible representation of $\Cl_n$ appears with infinite multiplicity in $\mathcal{H}_n$. Let $\mathcal{B}(\mathcal{H}_n)$ denote the algebra of $\Cl_n$-linear bounded operators on $\mathcal{H}_n$ equipped with the norm topology. 

\begin{defn} 
Let $\mathcal{M}_n$ be the category whose objects are finite-dimensional $\Cl_{n}$-submodules $W\subset \mathcal{H}_n$ and whose morphisms are
\beq
&&{\rm Mor}(W_0,W_1):=\left\{ \begin{array}{cl} \begin{array}{c} \Cl_{n}\otimes \Cl_{-1}-{\rm action\ on\ } W_0^\perp \\ {\rm extending\ the }\ \Cl_{n}{\rm-action}\end{array}  & {\rm if} \ W_0\subseteq W_1\subset \mathcal{H}_n \\
\null\\
\emptyset & {\rm else} \end{array}\right.
\eeq
where $W_0^\perp\subseteq W_1$ denotes the orthogonal complement of the subspace $W_0\subseteq W_1$. Composition for a nested inclusion $W_0\subseteq W_1\subseteq W_2$ uses the direct sum of Clifford modules to obtain a $\Cl_{n}\otimes \Cl_{-1}$-module structure on the orthogonal complement of $W_0\subseteq W_2$. Endow the objects and morphisms in this category with the subspace topology for the embeddings
\beq
{\rm Obj}(\mathcal{M}_n)\hookrightarrow \mathcal{B}(\mathcal{H}_n),\qquad {\rm Mor}(\mathcal{M}_n)\hookrightarrow \mathcal{B}(\mathcal{H}_n)^{\times 3} \label{eq:bddops}
\eeq
where the map identifies a submodule $W_0\subset \mathcal{H}_n$ with a finite rank projection operator, and $(W_0,W_1,e)\in {\rm Mor}(\mathcal{M}_n)$ determines a pair of finite rank projection operators and a bounded operator.
\end{defn} 

\begin{defn} The \emph{nerve} of a topological category $\mathcal{C}$ is the simplicial space whose $k$-simplices are $N_k\mathcal{C}=\mathcal{C}_k$ for $k=0,1$, and for $k>1$ are length $k$ chains of composable morphisms
\beq
N_k\mathcal{C}&=&\{x_0\stackrel{f_1}{\to}x_1\stackrel{f_2}{\to}x_2\stackrel{f_3}{\to}\cdots \stackrel{f_k}{\to}x_k\mid x_i\in \mathcal{C}_0, f_i\in \mathcal{C}_1\}\nonumber\\&=&\mathcal{C}_1\times_{\mathcal{C}_0}\cdots \times_{\mathcal{C}_0} \mathcal{C}_1,\nonumber 
\eeq
topologized as the fibered product.
The face maps in  $N_\bullet \mathcal{C}$ are determined by composing morphisms, and degeneracies are determined by the unit map in $\mathcal{C}$. 
\end{defn}


The \emph{geometric realization} of a simplicial space $Z_\bullet$ is 
\beq
| Z_\bullet | :=(\coprod_n Z_n\times \Delta^n)/{\sim}\qquad (f^*x,t)\sim (x,f_*t) \quad \forall f\colon [l]\to [k]\label{eq:fatrealizae}
\eeq
where $\Delta^n$ is the standard $n$-simplex, and $f\colon [l]\to [k]$ is an order-preserving map. Geometric realization is natural: a map of simplicial spaces $F\colon Z_\bullet\to Y_\bullet$ determines a map between realizations, $|F|\colon |Z_\bullet|\to |Y_\bullet|$. 

\begin{prop}[{\cite[\S1.3-1.4]{Cheung}, \cite[Theorem~7.1]{HST}, \cite[\S4]{DBEIndex}}]\label{prop:ABScat}
The classifying space ${\rm B} \mathcal{M}_n$ represents the functor $\KO^{n}(-)$. 
\end{prop}
Given a topological space  $Z$ with open cover $\{U_i\}_{i\in I}$, for a nonempty finite subset  $\sigma\subset I$ define $U_\sigma=\bigcap_{i\in \sigma} U_i$ as the intersection.

\begin{defn} \label{ex:opencover0}
For $Z$ a topological space and $\{U_i\}_{i\in I}$ an open cover, the \emph{\v{C}ech category}, denoted $\check{{\rm C}}(U_i)$, is the topological category whose objects are $\check{{\rm C}}(U_i)_0=\coprod_{\sigma} U_\sigma$ and whose morphisms are the space $\check{{\rm C}}(U_i)_1=\coprod_{\sigma\subseteq \tau} U_\tau$. The source and target maps are determined by the projection and inclusion $U_\tau\subseteq U_\sigma$, the unit is induced by the identity map $U_\sigma \to U_\sigma$, and composition comes from nested inclusions. 
\end{defn}

\begin{lem}[{\cite[Proposition~4.12]{Segalclassifying}, \cite[Theorem~2.1]{DuggerIsaksen}}]\label{lem:Cechrealize}
For $Z$ a topological space and $\{U_i\}_{i\in I}$ an open cover, the canonical map 
\beq
{\rm B} \check{C}(U_i)\stackrel{\sim}{\to} Z,\label{eq:Segalwe}
\eeq 
is a homotopy equivalence. 
 \end{lem}
 
 \begin{cor}\label{eq:cormaptoKO}
For $Z$ a topological space and $\{U_i\}_{i\in I}$ an open cover, a continuous functor 
\beq\label{eq:constructtoKO}
 \check{{\rm C}}(U_i)\to \mathcal{M}_n
\eeq
determines a map $Z\to B\mathcal{M}_n\simeq \KO^n$ into the $n$th space of $\KO$ spectrum.
\end{cor}
\bp This follows from Proposition~\ref{prop:ABScat} and Lemma~\ref{lem:Cechrealize}. 
\ep

\subsection{Super connections and and the families index bundle}\label{sec:cutoffconstrution}

\begin{defn}\label{defn:superconn}
Let $\HBR\to M$ be a real vector bundle with a fiberwise $\Cl_n$-action. A \emph{Clifford linear superconnection}~$\A$ on $\HBR$ is an odd, $\Cl_n$-linear map of Fr\'echet spaces satisfying a Leibniz rule
\beq
&&\A\colon \Omega^\bullet(M;\HBR)\to \Omega^\bullet(M;\HBR),\quad \A(\alpha\cdot v)=d\alpha \cdot v+(-1)^{|\alpha|} \alpha \cdot \A v\label{eq:ordsuperconn}
\eeq
where $v\in \Omega^\bullet(M;\HBR)$ and $\alpha\in \Omega^\bullet(M)$. 
\end{defn}

Using the $\Z$-grading on forms, a superconnection can be written as 
\beq
\A=\sum \A^{[k]},\qquad \A^{[k]}\colon \Omega^\bullet(M;\HBR)\to \Omega^{\bullet+k}(M;\HBR)\label{eq:superconnectionZ}
\eeq 
where $\A^{[1]}=\nabla$ is an ordinary connection on $\HBR$, and $\A^{[2k]}\in \Omega^{2k}(M;\End(\HBR)^{\odd})$, $\A^{[2k+1]}\in \Omega^{2k+1}(M;\End(\HBR)^{\ev})$ are endomorphism-valued forms. 

\begin{defn}
Suppose that $\langle-,-\rangle\colon\overline{\HBR}\otimes \HBR\to \underline{\C}$ is a positive (super) hermitian pairing. Then a superconnection $\A$ on $\HBR$ is \emph{self-adjoint} if 
\beq
\langle \A^{[1]}x,y\rangle+(-1)^{|x||y|}\langle x,\A^{[1]}y\rangle=d\langle x,y\rangle,\qquad (\A^{[k]})^*=i^{k+1}\overline{\A}{}^{[k]},\qquad k\ne 1,\label{eq:supersa}
\eeq
where $(\A^{[k]})^*$ denotes the super adjoint of $\A^{[k]}$. 
\end{defn}

\begin{rmk}
Translating from super adjoints to ordinary adjoints as in~\eqref{eq:superadjointtrans}, the above is equivalent to the standard definition of self-adjoint superconnection~\cite[page 117]{BGV}, i.e., the usual metric preserving condition for the ordinary connection $\A^{[1]}=\nabla$ and 
$$
(\A^{[k]})^\dagger=\left\{\begin{array}{ll} \A^{[k]} & k=1\mod 4\\ -\A^{[k]} & k=1\mod 4\\ -\A^{[k]} & k=2\mod 4\\ \A^{[k]} & k=0 \mod 4\end{array}\right. \qquad\qquad k\ne 1
$$
where $(\A^{[k]})^\dagger$ denotes the ungraded adjoint relative to an ordinary hermitian inner product.
\end{rmk}

For a superconnection $\A$, define 
\beq
&&\A(u):=u^{-1/2}\A^{[0]}+\A^{[1]}+u^{1/2}\A^{[2]}+\dots+ u^{j/2}\A^{[j]}+\dots\quad\nonumber
\eeq
for a formal parameter $u^{-1/2}$ of degree $+1$. 

\begin{defn}\label{defn:chernform}
The \emph{Chern form} of a $\Cl_n$-linear superconnection is the Clifford supertrace (which need not exist for infinite rank bundles)
\beq\label{eq:theChernchar}
&&\Ch(\A):=\sTr_{\cCl_n}(u^{-n/2} e^{-u\A(u)^2})=\sTr(\Gamma_n u^{-n/2} e^{-u\A(u)^2}) \in \Omega^\bullet(M;\C[u^{\pm 2}]))
\eeq
for the Clifford supertrace defined in~\eqref{eq:Cliffordsupertrace}. 
\end{defn}
The absence of fractional powers of~$u$ in $\Ch(\A)$ follows from the fact that the supertrace of an odd endomorphism is zero. Self-adjoint superconnections have Chern forms with only even powers of~$u$. The following is a standard algebraic argument, e.g., \cite[\S1.5]{BGV}.

\begin{lem}
When it exists, the Chern form of a $\Cl_n$-linear superconnection $\A$ is closed and has total degree $n$.
\end{lem}

Given a self-adjoint, $\Cl_{n}$-linear superconnection~$\A$ on a Fr\'echet vector bundle $\HBR\to M$, let $\mathcal{H}_m$ denote the Hilbert space completion of the fiber of $\mathcal{E}$ at $m\in M$. Define the subset $U_\lambda\subset M$ for $\lambda\in \R_{>0}$ as follows. A point $m\in M$ belongs to $U_\lambda$ if $\lambda$ is not an eigenvalue of the (non-negative) operator $(\A^{[0]})^2_m$ on $\mathcal{H}_m$:
\beq
U_\lambda=:\{m\in M\mid \lambda\notin {\rm Spec}((\A^{[0]})^2_m)\}\subset M.\label{eq:Ulambda2}
\eeq
Over each subset $U_\lambda\subset M$, define a subset $\mathcal{H}^{<\lambda}\subset \coprod_{m\in U_\lambda} \mathcal{H}_m$ as the direct sum of $\nu$-eigenspaces~$\mathcal{H}_{m}^\nu$ of $(\A^{[0]})^2_m$ with eigenvalue $\nu<\lambda$
\beq
 \mathcal{H}^{<\lambda}_m:=\bigoplus_{\nu<\lambda} \mathcal{H}_{m}^\nu \subset \mathcal{H}_m.\label{eq:bundleofClnmodules}
 \eeq
Letting $m$ vary, $\mathcal{H}^{<\lambda}$ denotes the $\Cl_{n}$-invariant subset of~$\coprod_{m\in U_\lambda} \mathcal{H}_m$ whose fibers are $ \mathcal{H}^{<\lambda}_m$. Define the maps (of sets) over $U_\lambda$
\beq\label{eq:fiberwiseincludeproject}
p^\lambda\colon \coprod_{m\in U_\lambda} \mathcal{H}_m\to \mathcal{H}^{<\lambda},\qquad i_\lambda\colon \mathcal{H}^{<\lambda}\to \coprod_{m\in U_\lambda} \mathcal{H}_m
\eeq
where $p^\lambda$ is the fiberwise orthogonal projection, and $i_\lambda$ is the fiberwise inclusion. 

\begin{defn} \label{defn:cutoffs}
Suppose that $\A$ is a self-adjoint $\Cl_n$-linear superconnection for which the differential form-valued Chern character~\eqref{eq:theChernchar} is well-defined. Then $\A$ \emph{admits smooth index bundles} if the following properties are satisfied:
\begin{enumerate}
\item the collection of subsets $\{U_\lambda\}_{\lambda\in \R_{>0}}$ defined in~\eqref{eq:Ulambda2} form a smooth open cover of~$M$;
\item the subsets $\mathcal{H}^{<\lambda}\subset \HBR|_{U_\lambda}\to U_\lambda$ are smooth vector subbundles;
\item the composition defines a smooth connection:
$$
\nabla^{<\lambda}:=p^\lambda\circ \A^{[1]}|_{U_\lambda}\circ p^\lambda\colon \Gamma(\mathcal{H}^{<\lambda})\to \Omega^1(U_\lambda;\mathcal{H}^{<\lambda}),
$$
\item the limit of Chern--Simons forms
\beq
\eta_{\lambda}:=\lim_{r\to \infty} \int_1^r \sTr_{\Cl_n}(\frac{d\A_{\lambda}(s)}{ds}\exp(\A_{\lambda}(s)))ds \in \Omega^\bullet(U_\lambda ;C^\infty(\R_{>0})[u^{\pm 1}])\label{eq:improperintegral}
\eeq
exists and satisfies
\beq\label{eq:satisfies}
d\eta_\lambda=\Ch(\A)|_{U_\lambda}-\Ch(\nabla^{<\lambda})
\eeq
where in~\eqref{eq:improperintegral}, $\A_{\lambda}(s)$ is the 1-parameter family of superconnections on $\HBR|_{U_\lambda} \oplus \Pi \mathcal{H}^{<\lambda}$,
\beq\label{eq:nullconcordFT}
&&\A_{\lambda}(s)=\A_{\lambda}(s)^{[0]}+\A^{[1]} \oplus \nabla^{<\lambda} +\sum_{i>1} s^{-i} \A^{[i]} \oplus 0_{\Pi \mathcal{H}^{<\lambda}},\quad \A_{\lambda}(s)^{[0]}=\left[\begin{array}{cc} \A^{[0]} & s\cdot i_\lambda \\ s \cdot p^\lambda & 0 \end{array}\right]
\eeq
where $i_\lambda$ and $p^\lambda$ are the maps~\eqref{eq:fiberwiseincludeproject}.
\end{enumerate}
\end{defn}

\begin{rmk}\label{rmk:indexvsdiff}
Conditions (1) and (2) are required to obtain an index of $\A$ valued in~$\KO^n(M)$. Conditions (3) and (4) are required to relate a \v{C}ech-de~Rham cocycle representing the Chern character of this index bundle to the Chern character of~$\A$. We also note that for a superconnection on a bundle over $M=\pt$, conditions (1)-(3) are automatic; the cutoff condition in this case is therefore convergence of the integral~\eqref{eq:nullconcordFT}, which in turn requires that the eigenvalues $\A^{[0]}$ grow sufficiently quickly. 
\end{rmk}

\begin{prop}\label{prop:Bismutcutoffs}
Any finite-rank self-adjoint $\Cl_n$-linear superconnection admits smooth index bundles. The Clifford linear Bismut superconnection adapted to a family of Clifford linear Dirac operators admits smooth index bundles. 
\end{prop}
\bp
For $\Cl_0=\R$, the finite rank case is~\cite[Theorem~9.7]{BGV} and for Bismut superconnections this follows from~\cite[Proposition~9.10 and Theorem~9.26]{BGV}, see also \cite[\S7.2]{FreedLott}. The straightforward extension of these results to Clifford linear superconnections is~\cite[\S8.2-8.3]{DBEIndex}.
\ep

For a superconnection admitting smooth index bundles, we review the index bundle construction from~\cite[\S8.2-8.3]{DBEIndex}. 

\begin{prop}[{\cite[\S8.2]{DBEIndex}}] \label{prop:appenindexbundle}
Suppose a $\Cl_n$-linear superconnection admits smooth index bundles. Then a choice of locally finite open cover $\{U_\lambda\}_{\lambda \in \Lambda}$ associated with a discrete subset $\Lambda\subset \R_{>0}$ determines a continuous functor, 
\beq\label{eq:KOindexspec}
 \check{{\rm C}}(U_\lambda)\to \mathcal{M}_n\quad \implies \quad \Ind_\Lambda(\A)\colon M\to \KO^n,
\eeq
whose realization is a map to the $n$th space in the $\KO$-spectrum. Given choices $\Lambda\subset \Lambda' \subset \R_{>0}$ related by a refinement of open covers, there is a uniquely specified homotopy between the associated maps~\eqref{eq:KOindexspec}.
\end{prop}
\bp
The $\Cl_n$-modules $\mathcal{H}^{<\lambda}\to U_\lambda$ determine classes $[\mathcal{H}^{<\lambda}]\in \KO^n(U_\lambda)$ for each $\lambda\in \R_{>0}$. On overlaps $U_\lambda\bigcap U_{\mu}$ for $\lambda<\mu$ there are canonical inclusions
\beq
g_{\lambda\mu}\colon \mathcal{H}^{<\lambda}\hookrightarrow \mathcal{H}^{<\mu}\label{eq:compatibilitydata}
\eeq
and the orthogonal complement to the image of $g_{\lambda\mu}$ has a $\Cl_{n}\otimes \Cl_{-1}$-action extending the $\Cl_n$-action, determined by the invertible odd endomorphisms
\beq\label{eq:compatibilitydata2}
e^{\lambda\mu}|_{\mathcal{H}_b^{\nu}}= \frac{1}{\sqrt{\nu}} \A^{[0]}|_{\mathcal{H}_b^{\nu}},\quad \mathcal{H}_b^{\nu}\subset (\mathcal{H}^{<\lambda})^\perp\subset \mathcal{H}^{<\mu}
\eeq
that define the action of the generator $e\in \Cl_{-1}$. The inclusions $g_{\lambda\mu}$ and odd endomorphisms $e^{\lambda\mu}$ satisfy a compatibility property on triple intersections. By the contractibility of the orthogonal group of a separable infinite dimensional Hilbert space~\cite{Kuiper}, the data~\eqref{eq:bundleofClnmodules},~\eqref{eq:compatibilitydata} and~\eqref{eq:compatibilitydata2} determine a continuous functor~\eqref{eq:constructtoKO} uniquely up to contractible choice. Hence by Corollary~\ref{eq:cormaptoKO}, we obtain a map to the $n$th space in the $\KO$-spectrum~\eqref{eq:KOindexspec}.
A refinement of open covers associated with $\Lambda\subset \Lambda'\subset \R_{>0}$, determines a continuous functor between the associated \v{C}ech categories that after realization is a homotopy equivalence. This leads to a pair of functors~\eqref{eq:constructtoKO} and restriction of the operators~\eqref{eq:compatibilitydata2} furnish a natural isomorphism between these functors. By \cite[Proposition~2.1]{Segalclassifying}, the resulting maps $\Ind_\Lambda(\A),\Ind_{\Lambda'}(\A)\colon M\to \KO^{-n}$ are homotopic. 
\ep

The Pontryagin character of the class associated with~\eqref{eq:KOindexspec} is determined by the local data $\Ch(\nabla^{<\lambda})$ with compatibilities on overlaps coming from Chern--Simons forms; more precisely we obtain a cocycle in the \v{C}ech--de~Rham complex denoted
\beq\label{eq:appencdR}
&&\Ch(\Ind_\Lambda(\A))=[U_\lambda,\Ch(\nabla^{<\lambda}),\CS(\nabla^{<\lambda},\nabla^{<\lambda'})]\in C^*(\{U_\lambda\},\Omega^\bullet).
\eeq

\begin{lem} \label{lem:anotherdamnhomotopy}
The differential forms~\eqref{eq:improperintegral} witness a coboundary in the \v{C}ech--de~Rham complex between~\eqref{eq:appencdR} and the (globally defined) cocycle $\Ch(\A)$.
\end{lem}
\begin{proof}[Proof sketch.] 
This follows from~\eqref{eq:satisfies} and the defining properties of Chern--Simons forms; we refer to \cite[\S8.3]{DBEIndex} for details. 
\ep

\bibliographystyle{amsalpha}
\bibliography{references}

\end{document}

%% file: sample.tikzstyles
\tikzstyle{new edge style 0}=[-, dashed=]